\documentclass[11pt]{article}
\usepackage{apacite}
\usepackage{natbib}
\usepackage{epsfig,lscape}
\usepackage{amssymb,amsfonts,amsmath,amsthm}
\usepackage{rotating}
\usepackage{bbm}
\usepackage{multibib}
\newcites{supp}{Supplement references}

\usepackage{threeparttable}
\usepackage{booktabs}
\usepackage{subcaption}
\usepackage{multirow}
\usepackage{graphicx}
\usepackage[T1]{fontenc}
\usepackage{helvet}
\usepackage{comment}
\usepackage{xcolor}
\usepackage{amsthm}
\usepackage{apacite}
\usepackage{natbib}
\usepackage{statmath}
\usepackage{amsmath}
\usepackage{amssymb}
\usepackage{algorithm}
\usepackage{algpseudocode}
\usepackage{appendix}
\usepackage{amsmath}
\usepackage{mathtools}
\usepackage{multicol}
\usepackage{algorithm,algpseudocode}
\usepackage{lipsum}
\usepackage{hyperref}

\algrenewcommand\algorithmicrequire{\textbf{Input:}}
\algrenewcommand\algorithmicensure{\textbf{Output:}}

\makeatletter

\newcommand{\LL}{\textcolor{red}}

\allowdisplaybreaks[4]

\newcommand\mytext[1]{\text{\scriptsize{#1}}}

\newtheorem*{condition*}{\conditionnumber}

\providecommand{\conditionnumber}{}

\newtheorem{thm}{Theorem}
\newtheorem{lem}{Lemma}
\newtheorem{prop}{Proposition}

\def\argmin{\mathrm{argmin}}
\def\min{\mathrm{min}}
\def\max{\mathrm{max}}
\def\diag{\mathrm{Diag}}
\def\supp{\mathrm{supp}}

\def\me{\mathrm e}

\def\c{\mathrm{c}}
\def\T{ {\mathrm{\scriptscriptstyle T}} }
\def\F{\mathrm{F}}
\def\V{\mathrm{V}}
\def\supp{\mathrm{supp}}
\def\exp{\mathrm{exp}}
\def\log{\mathrm{log}}

\def\Tr{\mathrm{Tr}}
\def\Obj{\mathrm{Obj}}

\newenvironment{prf}
{\noindent \textbf{Proof.}}{\hfill $\Box$ \vspace{.1in}}

\newtheorem{assumptionA}{Assumption}

\newtheorem{conditionB}{Condition}

\newtheorem{conditionC}{Condition}

\newtheorem{conditionE}{Condition}

\theoremstyle{definition}
\newtheorem{rem}{Remark}

\usepackage{thmtools}
\usepackage{thm-restate}
\allowdisplaybreaks

\title{Inference for misspecified high-dimensional model}
\begin{document}
\begin{titlepage}

\begin{center}
{\Large Debiased Prediction Inference with Non-sparse Loadings in Misspecified High-dimensional Regression Models}

\vspace{.1in} {\large Libin Liang\footnotemark[1] and Zhiqiang Tan\footnotemark[1]}

\vspace{.1in}
\today
\end{center}

\footnotetext[1]{Department of Statistics, Rutgers University. Address: 110 Frelinghuysen Road,
Piscataway, NJ 08854. E-mails: ll866@stat.rutgers.edu, ztan@stat.rutgers.edu.}

\paragraph{Abstract.}

High-dimensional regression models with regularized sparse estimation are widely applied.
For statistical inferences, debiased methods are available about single coefficients or predictions with sparse new covariate vectors (also called loadings),
in the presence of possible model misspecification. However, statistical inferences about predictions with non-sparse loadings
are studied only under the assumption of correctly specified models.
In this work, we develop debiased estimation and associated Wald confidence intervals for
predictions with general loadings, allowed to be non-sparse, from possibly misspecified high-dimensional regression models.
Our debiased estimator involves estimation of a debiasing vector,
which is the general loading left-multiplied by the non-centered precision matrix in the linear model (LM) setting or
the inverse Hessian of the objective function at the target coefficient vector in the generalized linear model (GLM) setting.
We propose suitable estimators of the precision matrix or the inverse Hessian respectively in the LM or GLM settings and, for the first time,
establish a root-n asymptotic expansion for the debiased prediction and justify associated Wald confidence intervals
under sparsity conditions on the precision matrix or the inverse Hessian which are comparable to the conjunction of sparsity conditions
required for inferences about all single coefficients in existing works.
We also provide numerical results which further demonstrate the validity of our proposed confidence intervals for predictions with general loadings
from possibly misspecified regression models.

%\vspace{-.2in}

\end{titlepage}

\addtocontents{toc}{\protect\setcounter{tocdepth}{0}}

\section{Introduction}

High-dimensional regression models are widely applied with regularized sparse estimation such as Lasso estimation \citep{hastie_book}.
But statistical inferences about coefficients or predictions are challenging, due to the fact that
the bias of regularized sparse estimation can be large compared to the standard deviation.
To address such problems, various methods and theory have been proposed. See Section~\ref{sec:related-work} for a discussion of related works.

% Hastie et al. book, https://www.amazon.com/Statistical-Learning-Monographs-Statistics-Probability/dp/0367738333

In this work, we aim to obtain valid confidence intervals of the prediction from a new covariate vector (also called loading), allowed to be non-sparse,
for possibly misspecified high-dimensional regression models. With model misspecification, the target coefficient vector at the population level is
defined as a minimizer of a population objective function. Specifically, let $Y_i\in\mathbb{R}$ be a response and $X_i\in\mathbb{R}^p$ be a covariate vector for $i=1,\ldots, n$. Assume that $(X_i,Y_i)$, $i=1,\ldots, n$, are independent and identically distributed (i.i.d.) as $(X_0, Y_0)$.
A correctly specified linear model (LM) indicates that
\begin{align*}
 Y_i= X_i^\T \beta^* + \varepsilon_i, \quad i=1,\ldots, n,  %\label{eq:LM-correct}
\end{align*}
where the noise $\varepsilon_i$ is mean-zero conditionally on $X_i$, and $\beta^*$ represents the true coefficient vector.
In a misspecified LM, the target (not true) coefficient vector associated with least-squares estimation is defined as
\begin{align}
    \bar\beta = \argmin_{\beta \in\mathbb{R}^p} \mathbb{E}[\frac{1}{2}(Y_0-X_0^\T \beta)^2]. \label{eq:target_beta_LM}
\end{align}
A correctly specified generalized linear model (GLM) (more precisely, a GLM with correct mean restriction) indicates that
\begin{align*}
 E(Y_i | X_i) = g^\prime(X^\T \beta^*), \quad i=1,\ldots,n,
\end{align*}
where $\beta^*$ represents the true coefficient vector and
$g^\prime(\cdot)$, the derivative of a convex function $g(\cdot)$, is an increasing function, called the inverse link function.
For logistic regression, $g(b) = \log(1+\me^b)$.
In a misspecified GLM, the target (not true) coefficient vector associated with convex M-estimation
%(e.g., maximum likelihood estimation with a canonical link)
is defined as
\begin{align}
    \bar\beta = \argmin_{\beta \in\mathbb{R}^p}\mathbb{E}[\phi(Y_0,X_0^\T\beta)], \label{eq:target_beta_GLM}
\end{align}
where $\phi(a,b)$ is a function convex and three-times differentiable in $b$.
For maximum likelihood estimation with a canonical link, $\phi(Y_0,X_0^\T\beta) = -Y_0 X_0^\T\beta + g(X_0^\T\beta)$.

For a general new covariate vector (or loading) $\xi\in\mathbb{R}^p$,
the target of point prediction is defined as $\xi^\T\bar\beta$. In addition, we denote $\Sigma = \mathbb{E}[X_0 X_0^\T]\in \mathbb{R}^{p\times p}$ and $\Omega^* = \Sigma^{-1}$, which, by some abuse of terminology, are the non-centered variance and precision matrices of $X_0$ respectively. For misspecified GLM, we further denote $H=\mathbb{E}[\phi_2(Y_0,X_0^\T\bar\beta)X_0X_0^\T]$ and $M^*= H^{-1}$, as the counterparts of $\Sigma$ and $\Omega^*$,
where $\phi_2(a,b)=\frac{\partial^2\phi(a,b)}{\partial^2 b}$.
The LM setting is a special case of GLM setting, with $\phi (a,b) = \frac{1}{2} (a-b)^2$.
Nevertheless, we discuss LM and GLM settings separately, because the former is simpler and allows stronger results to be obtained.
A complication in GLM setting is that $M^*$ (but not $\Omega^*$) depends on unknown $\bar\beta$
and hence estimation of $M^*$ is affected by that of $\bar\beta$. See Remarks \ref{rem:est_error} and \ref{rem:CLIME_weighted_nodewise_regression_M}
for relevant discussion.

We employ a debiased estimator for the target prediction, $\xi^\T\bar\beta$, and construct a Wald confidence interval in the misspecified LM and GLM settings.
Let $\hat{\beta}\in \mathbb{R}^p$ be
an initial regularized sparse estimator of $\bar\beta$, for example, a Lasso estimator computed from $\{(X_i,Y_i)\}_{i=1,\ldots,n}$. Our debiased estimator is of the form
\begin{align}
    \widehat{\xi^\T\beta} = \xi^\T\hat{\beta} + \tilde{u}^\T\frac{1}{n}\sum_{i=1}^n X_i(Y_i-X_i^\T\hat{\beta}) ,  \label{eq:debiased_lm}
\end{align}
in misspecified LM setting and is of the form
\begin{align}
    \widehat{\xi^\T\beta} = \xi^\T\hat{\beta} -  \tilde{v}^\T\frac{1}{n}\sum_{i=1}^n X_i\phi_1(Y_i,X_i^\T\hat{\beta}), \label{eq:debiased_glm}
\end{align}
in misspecified GLM setting, where $\phi_1(a,b)=\frac{\partial\phi(a,b)}{\partial b}$,
$\tilde{u}\in\mathbb{R}^p$ is an estimator of $u^*=\Omega^*\xi$ and
$\tilde{v}\in\mathbb{R}^p$ is an estimator of $v^*=M^* \xi$. We refer to $u^*$ or $v^*$ as a debiasing vector.
The estimator $\widehat{\xi^\T\beta}$ in \eqref{eq:debiased_lm} or \eqref{eq:debiased_glm} can be motivated from
the decomposition in (\ref{debiased_est_diff}) or (\ref{debiased_est_diff_glm}) later.
Moreover, our estimator of the asymptotic variance of $\widehat{\xi^\T\beta}$ is of the form
\begin{align}
    \hat{V} = \frac{1}{n} \sum_{i=1}^n (\tilde{u}^\T X_i)^2(Y_i-X_i^\T\hat{\beta})^2,\label{eq:debiased_lm_v}
\end{align}
in misspecified LM setting and is of the form
\begin{align}
    \hat{V} = \frac{1}{n} \sum_{i=1}^n (\tilde{v}^\T X_i)^2\phi^2_1(Y_i,X_i^\T\hat{\beta}),\label{eq:debiased_glm_v}
\end{align}
in misspecified GLM setting.
For $0< \alpha < 1$, in the LM and GLM settings, our $1-\alpha$ Wald confidence interval is constructed as
\begin{align}
   \left[\widehat{\xi^{\T}\beta} - z_{1-\frac{\alpha}{2}}\sqrt{\frac{\hat{V}}{n}}, \widehat{\xi^{\T}\beta} + z_{1-\frac{\alpha}{2}}\sqrt{\frac{\hat{V}}{n}}\right], \label{CI}
\end{align}
where $z_{1-\frac{\alpha}{2}}$ is $1-\frac{\alpha}{2}$ percentile of the standard normal distribution.

To justify our Wald confidence intervals in LM setting, we propose a proper choice of $\tilde{u}$ such that under suitable regularity and sparsity conditions, we have
\begin{subequations}
\label{result_LM}
    \begin{equation}
        \frac{\widehat{\xi^{\T}\beta} - \xi^\T\bar\beta}{\sqrt{V}}= \frac{1}{n}\sum_{i=1}^n  \frac{(\xi^{\T}\Omega^* X_i)(Y_i-X_i^\T\bar\beta)}{\sqrt{V}} + o_p( n^{-\frac{1}{2}}), \label{result_LM_1}
    \end{equation}
    \begin{equation}
        \frac{\hat{V}}{V} =1+o_p(1), \label{result_LM_2}
    \end{equation}
\end{subequations}
where $V=\mathbb{E}[( \xi^{\T} \Omega^* X_i)^2(Y_i-X_i^\T\bar\beta)^2]$ is the asymptotic variance (after rescaling).
Similarly, to justify our confidence intervals in GLM setting, we propose a proper choice of $\tilde{v}$ such that under regularity and sparsity conditions, we have
\begin{subequations}
    \begin{equation}
        \frac{\widehat{\xi^{\T}\beta} -\xi^\T\bar\beta}{\sqrt{V}} = -\frac{1}{n}\sum_{i=1}^n \frac{(\xi^{\T} M^* X_i)\phi_1(Y_i,X_i^\T\bar\beta)}{\sqrt{V}}
     + o_p(  n^{-\frac{1}{2}}), \label{result_GLM_1}
    \end{equation}
    \begin{equation}
        \frac{\hat{V}}{V}=1+o_P(1), \label{result_GLM_2}
    \end{equation}
    \label{result_GLM}
\end{subequations}
where $V=\mathbb{E}[(\xi^{\T}M^* X_i)^2\phi^2_1(Y_i,X_i^\T\bar\beta)]$ is the asymptotic variance (after rescaling).

Our estimators $\tilde{u}$ and $\tilde{v}$ for debiasing vectors $u^*$ and $v^*$ are defined by replacing $\Omega^*$ and $M^*$ with their estimators, i.e., $\tilde{u} = \tilde{\Omega}\xi$ and $\tilde{v}=\tilde{M}\xi$ where $\tilde{\Omega}\in\mathbb{R}^{p\times p}$ and $\tilde{M}\in\mathbb{R}^{p\times p}$ are estimators of $\Omega^*$ and $M^*$ respectively. To facilitate theoretical analysis, we require that $\tilde{\Omega}$ and $\tilde{M}$ are obtained from a new set of covariate vectors $X_i^\prime\in \mathbb{R}^p$ and (in case of $\tilde{M}$) responses $Y_i^\prime$ for $i=1,\ldots, m$. The new set $\{(X_i^\prime, Y_i^\prime)\}_{i=1,\ldots m}$ 
satisfy $(X_i^\prime, Y_i^\prime) \sim (X_0,Y_0)$ for each $i$ and are independent of the original samples $\{(X_i, Y_i)\}_{i=1,\ldots, n}$, from which
$\hat\beta$ is computed. Alternatively, $\hat{\beta}$ can also be computed from the combination of $\{(X_i,Y_i)\}_{i=1,\ldots,n}$ and $\{(X^\prime_i,Y^\prime_i)\}_{i=1,\ldots,m}$. In practical applications, we apply a cross-fitting technique to realize the independence of two samples.
The cross-fitting version of our method is discussed in Section \ref{sec:sparse_omega} and \ref{sec:sparse_M}.

To highlight main contributions, we present our method and theory in the preceding setup with two independent samples.
We denote $s_0=\Vert \bar\beta\Vert_0$. In the LM setting, we denote $s_{1j}=\Vert \Omega_{\cdot j}^*\Vert_0$ for $j=1,\ldots, p$ and $s_1 = \max_{j=1,\ldots, p}\{s_{1j}\}$. In the GLM setting, we denote $s_{1j}=\Vert M^*_{\cdot j}\Vert_0$ for $j=1,\ldots, p$ and $s_1 = \max_{j=1,\ldots, p}\{s_{1j}\}$.
In both settings, we further denote $S_\xi=\supp(\xi)$ (the support of $\xi$) and $s_{1\xi} = (s_{1j}: j\in S_\xi)^\T$ (a vector).
See Section~\ref{sec:notation} for more details of our notation.

In the LM setting, our estimator $\tilde{\Omega}$ can be defined by directly applying CLIME \citep{clime} on the sample $\{ X_i^\prime\}_{i=1,\ldots,m}$
as shown in Algorithm \ref{alg:omega_est_prop},
and the corresponding results are summarized in Proposition \ref{prop:lm} (deduced from Theorem \ref{thm_clime} later).
Our method can also be employed with an alternative estimator $\tilde\Omega$ defined by a two-stage algorithm
with the first stage being a variation of Lasso nodewise regression \citep{zhang2014confidence, van_de_Geer_2014},
and the corresponding results can be deduced from Theorem \ref{thm_two_stage} later.

\begin{prop}
\label{prop:lm}
     Suppose that $m$ and $n$ are within a constant factor of each other, Assumptions \ref{sigma_eigenvalue}, \ref{sub_gaussian_residual} and \ref{ass:sparse_omega} are satisfied,
     $ \Vert \hat{\beta}-\bar\beta\Vert_1  =O_p(s_0\sqrt{\frac{\log(p\vee n)}{n}})$ and $(\hat{\beta}-\bar\beta)^\T\hat{\Sigma}(\hat{\beta}-\bar\beta)=O_p(s_0\frac{\log(p\vee n)}{n})$ with $\hat{\Sigma}=\frac{1}{n}\sum_{i=1}^nX_iX_i^\T$, and $\tilde{\Omega}$ is defined by Algorithm \ref{alg:omega_est_prop} with $\rho_n=A_0\lambda_n$ and a sufficiently large constant $A_0$.
     Then the following results hold for $\widehat{\xi^\T\beta}$ in (\ref{eq:debiased_lm}) and $\hat{V}$ in (\ref{eq:debiased_lm_v}). \\
    \noindent (i)
    If further Assumption \ref{ass:sub_gaussian_entries} is satisfied, then
(\ref{result_LM}) holds provided $\max\{s_0, \Vert s_{1\xi}\Vert_\infty \frac{\Vert \xi\Vert_1}{\Vert \xi\Vert_2}, \min\{\sqrt{\Vert s_{1\xi}\Vert_1},\Vert s_{1\xi}\Vert_\infty\} \sqrt{ s_0}\}\tau_n=o(1)$. \\
\noindent (ii)
If further Assumption \ref{ass:sub_gaussian_LM} is satisfied,
then (\ref{result_LM}) holds provided $\max\{s_0, \min\{\sqrt{\Vert s_{1\xi}\Vert_1},\Vert s_{1\xi}\Vert_\infty\} \sqrt{ s_0} \}\tau_n=o(1)$.\\
Here $\lambda_n = \sqrt{\frac{\log(p\vee n)}{n}}$ and $\tau_{n} =\frac{\log(p\vee n)}{\sqrt{n}}$.
\end{prop}

\begin{algorithm}[t]
     \caption{CLIME estimator}
     \vspace{-0.8em}
     \begin{flushleft}
    \hspace*{\algorithmicindent}\textbf{Input: } $\{X^\prime_i\}_{i=1,\ldots, m}$ and $\rho_n >0$. \\
 \hspace*{\algorithmicindent}\textbf{Output: }$\tilde{\Omega}\in\mathbb{R}^{p\times p}$.
    \end{flushleft}
    \vspace{-1.2em}
    \begin{algorithmic}[1]
     % Input:
     % Output:
    \State Compute $\hat{\Sigma}^\prime=\frac{1}{m}\sum_{i=1 }^m X^\prime_iX_i^{\prime \T}$.
    \State Obtain $\tilde{G}\in \mathbb{R}^{p\times p}$ by solving
        \begin{align*}
            \tilde{G} = &\argmin_{\text{$G\in\mathbb{R}^{p\times p}$}}  \sum_{j=1}^p \Vert G_{\cdot j}\Vert_{1}\\
            &\text{s.t. } \Vert \hat{\Sigma}^\prime G - I_p\Vert_{\max} \leq \rho_n.
        \end{align*}
   \State Compute $\tilde{\Omega} \in\mathbb{R}^{p\times p}$ such that $\tilde{\Omega}_{ij}=\tilde{G}_{ij}\textbf{1}_{[|\tilde{G}_{ij}| \leq |\tilde{G}_{ji}|]} + \tilde{G}_{ji}\textbf{1}_{[|\tilde{G}_{ij}| > |\tilde{G}_{ji}|]}$.
    \end{algorithmic}
    \label{alg:omega_est_prop}
\end{algorithm}

\begin{algorithm}[t]
    \caption{Extended CLIME estimator}
    \vspace{-0.8em}
     \begin{flushleft}
    \hspace*{\algorithmicindent} \textbf{Input: }$(X^\prime_i, Y_i^\prime)_{i=1,\ldots,m}$ and $\rho_n >0$. \\
 \hspace*{\algorithmicindent} \textbf{Output: }$\tilde{M}\in\mathbb{R}^{p\times p}$.
  \end{flushleft}
   \vspace{-1em}
    \begin{algorithmic}[1]
    \State Compute $\hat{\beta}^\prime $ from $(X^\prime_i, Y^\prime_i)_{i=1,\ldots, m}$ similarly as $\hat\beta$ from $(X_i, Y_i)_{i=1,\ldots, n}$.
     \State Compute $\hat{H}^\prime=\frac{1}{m}\sum_{i=1}^m\phi_2(Y_i^\prime,X_i^{\prime \T} \hat{\beta}^\prime)X_i^\prime X_i^{\prime \T}$.
    \State Obtain $\tilde{G}^J\in \mathbb{R}^{p\times p}$ by solving
        \begin{align*}
            \tilde{G} = &\argmin_{\text{$G\in\mathbb{R}^{p\times p}$}}  \sum_{j=1}^p \Vert G_{\cdot j}\Vert_{1}\\
            &\text{s.t. } \Vert \hat{H}^\prime G - I_p\Vert_{\max} \leq \rho_n.
        \end{align*}
    \State Compute $\tilde{M}\in\mathbb{R}^{p\times p}$ as $(\tilde{M})_{ij}=\tilde{G}_{ij}\textbf{1}_{ [ |\tilde{G}_{ij}| \leq |\tilde{G}_{ji}|]} + \tilde{G}_{ji}\textbf{1}_{[|\tilde{G}_{ij}| > |\tilde{G}_{ji}|]}$.
    \end{algorithmic}
    \label{alg:est_M_prod}
\end{algorithm}

In the GLM setting, our estimator $\tilde{M}$ can be defined by an extended version of CLIME as shown in Algorithm \ref{alg:est_M_prod},
and the corresponding results are summarized in Proposition \ref{pro:glm_setting} (deduced from Theorem \ref{thm_clime_M} later).
Similarly as in the LM setting, an alternative estimator $\tilde M$ can be defined by a two-stage algorithm,
and the corresponding results can be deduced from Theorem \ref{thm_two_stage_M} later.

\begin{prop}
\label{pro:glm_setting} Suppose that $m$ and $n$ are within a constant factor of each other,
Assumptions \ref{ass:H_eigen_value}, \ref{ass:residulal_glm_lower},  \ref{bounded_higher_differentiate} and \ref{ass:sparse_M} are satisfied, $\Vert \hat{\beta}-\bar\beta\Vert_1=O_p(s_0\sqrt{\frac{\log(p\vee n)}{n}})$ and $(\hat{\beta}-\bar\beta)^\T\hat{\Sigma}(\hat{\beta}-\bar\beta)=O_p(s_0\frac{\log(p\vee n)}{n})$, $\Vert \hat{\beta}^\prime-\bar\beta\Vert_1 = O_p(s_0\sqrt{\frac{\log(p\vee n)}{n}})$ and  $(\hat{\beta}^\prime-\bar\beta)^\T \hat{\Sigma}^\prime (\hat{\beta}^\prime-\bar\beta)=O_p(s_0\frac{\log(p \vee n )}{n})$  for $\hat{\beta}^\prime$ in Algorithm \ref{alg:est_M_prod}  and $\hat{\Sigma}^\prime=\frac{1}{m}\sum_{i=1 }^m X^\prime_iX_i^{\prime \T}$,
and $\tilde{M}$ is defined by Algorithm \ref{alg:est_M_prod} with $\rho_n = \sqrt{s_0}A_0\lambda_n$ and a sufficiently large constant $A_0$.
Then the following results hold for $\widehat{\xi^\T\beta}$ in (\ref{eq:debiased_glm}) and $\hat{V}$ in (\ref{eq:debiased_glm_v}). \\
\noindent (i) If further Assumption \ref{ass:sub_gaussian_entries_GLM} is satisfied, then (\ref{result_GLM}) holds provided
$\max\{s_0\sqrt{\log( n)}, \Vert s_{1\xi}\Vert_\infty \sqrt{s_0}\frac{\Vert \xi\Vert_1}{\Vert \xi\Vert_2} ,\min\{ \sqrt{ \Vert s_{1\xi}\Vert_1 }  ,  \Vert s_{1\xi}\Vert_\infty \}s_0\}\tau_n=o(1)$. \\
\noindent (ii)
If further Assumption \ref{ass:sub_gaussian_GLM} is satisfied, then (\ref{result_GLM}) holds provided $\max\{s_0 \sqrt{\log( n) },\min\{ \sqrt{ \Vert s_{1\xi}\Vert_1}  ,  \Vert s_{1\xi}\Vert_\infty \}s_0 \}\tau_n = o(1)$.\\
Here we denote $\lambda_n = \sqrt{\frac{\log(p\vee n)}{n}}$ and $\tau_n=\frac{\log(p\vee n)}{\sqrt{n}}$.
\end{prop}

The primary difference between results (i) and (ii) in Proposition \ref{prop:lm} and \ref{pro:glm_setting} is as follows.
The sparsity condition in (i) depends on $\frac{\Vert \xi\Vert_1}{\Vert \xi\Vert_2}$
and hence may fail for a non-sparse loading $\xi$ with large $\frac{\Vert \xi\Vert_1}{\Vert \xi\Vert_2}$,
for example, $\xi$ being a vector with all entries of similar magnitudes and  $\frac{\Vert \xi\Vert_1}{\Vert \xi\Vert_2} \approx \sqrt{p}$.
On the other hand, provided with the stronger assumption that $X_0$ is a sub-gaussian vector,
the sparsity condition in (ii) is independent of $\frac{\Vert \xi\Vert_1}{\Vert \xi\Vert_2}$,
and hence can be satisfied for a non-sparse loading $\xi$ with large $\frac{\Vert \xi\Vert_1}{\Vert \xi\Vert_2}$, subject to
sufficiently small $s_0$ and $s_1$.
In fact, by taking only $\Vert s_{1\xi}\Vert_\infty$ in the term $\min\{ \sqrt{ \Vert s_{1\xi}\Vert_1}  ,  \Vert s_{1\xi}\Vert_\infty \}$ in the sparsity conditions and the fact that $\Vert s_{1\xi}\Vert_\infty\leq s_1$,
Proposition \ref{prop:lm} (ii) provides valid inference about $\xi^\T \bar\beta$ under the sparsity condition  $\max\{s_0, s_1\sqrt{ s_0} \}\tau_n=o(1)$,
and Proposition \ref{pro:glm_setting} (ii) provides valid inference about $\xi^\T \bar\beta$ under the sparsity condition $\max\{s_0 \sqrt{\log(n)},s_1s_0 \}\tau_n = o(1)$, for any general loading $\xi$.
This is the first time such results are established for valid prediction inference using regularized sparse estimation with possibly misspeicifed regression models. Our key technical innovation behind Propositions \ref{prop:lm} (ii) and \ref{pro:glm_setting} (ii) is to control the $L_2$ norm error of the estimator $\tilde u$ or $\tilde v$ for the debiasing vector $u^*$ or $v^*$ via the operator norm error of the associated matrix estimator $\tilde \Omega$ or $\tilde M$. See Sections \ref{sec:high-level-LM} and \ref{sec:high_level_ideas_M} for high-level explanations.

We further discuss implications of our results and compare with various existing results, in Section~\ref{sec:related-work} and Remarks \ref{rem:on_sparisity_condition}, \ref{rem:on_clime_nodewise}, \ref{rem:on_sparsity_conditoin_M}, and \ref{rem:CLIME_weighted_nodewise_regression_M}.

\subsection{Related works and comparison} \label{sec:related-work}

There has been an extensive literature on debiased methods using regularized sparse estimation in high-dimensional regression models.
We discuss several works directly related to our work.

\cite{zhang2014confidence} and \cite{van_de_Geer_2014} introduced debiased methods for inferring (e.g., constructing confidence intervals) about a single coefficient,
$\beta_j^*$, of $\bar\beta$ in high-dimensional LM and GLM respectively.
A generalization is discussed in \cite{Ning}.
Their debiased estimators for $\bar\beta_j$ are of form (\ref{eq:debiased_lm}) or (\ref{eq:debiased_glm}) with the new covariate vector (or loading) $\xi$ set to $e_j$, whose entries are 0 except 1 in the $j$th entry. Their methods are shown to work under suitable sparsity conditions on $\bar\beta$ and the $j$th column of $\Omega^*$ or $M^*$ (i.e., $u^*= \Omega^* e_j$ or $v^* = M^* e_j$) while allowing misspecification of the regression model.
However, the associated theory cannot be directly extend to inferring about prediction with a general, non-sparse loading $\xi$.

\cite{JMLR:v15:javanmard14a} proposed another debiased method for single $\bar\beta_j$ in high-dimensional LM, without requiring sparsity of the $j$th single column of $\Omega^*$. Their debiased estimator for $\bar\beta_j$ is also of form (\ref{eq:debiased_lm}) with $\xi=e_j$, but their estimator $\tilde{u}$ for $u^*$, the $j$th column of $\Omega^*$, is obtained by minimizing the empirical second moment of $\tilde{u}^\T X_i$, $i=1,\ldots, n$, under certain constraints.
However, their method is only shown to work when the regression model is correctly specified, and
their theory cannot be extended to prediction inference with a general loading $\xi$.

\cite{cai2020optimal} introduced a debiased method for prediction inference with a general loading $\xi$ in LM,
and \cite{guo2021inference} proposed an extension to logistic regression.
Their estimators for $\xi^\T\bar\beta$ are of form (\ref{eq:debiased_lm}).
The estimator $\tilde u$ of \cite{cai2020optimal} is obtained by incorporating extra constraints in the method of \cite{JMLR:v15:javanmard14a},
whereas that of \cite{guo2021inference} is defined by exploiting similar ideas in \cite{cai2020optimal}, but after re-weighting the second term in (\ref{eq:debiased_glm}), $\tilde{v}^\T\frac{1}{n}\sum_{i=1}^n X_i\phi_1(Y_i,X_i^\T\hat{\beta})$, as $\tilde{v}^\T\frac{1}{n}\sum_{i=1}^n X_i\phi_2(Y_i,X_i^\T\hat{\beta})^{-1}\phi_1(Y_i,X_i^\T\hat{\beta})$.
The methods of \cite{cai2020optimal} and \cite{guo2021inference} are justified without requiring sparsity of $\Omega^*$ or $M^*$,
but only when the regression model is correctly specified.
In contrast, by Proposition \ref{prop:lm} (ii) or Proposition \ref{pro:glm_setting} (ii), our method is justified for any general loading $\xi$ under suitable sparsity of all columns of $\Omega^*$ or $M^*$ in LM or GLM respectively (in addition to suitable sparsity of $\bar\beta$),
while allowing the regression model to be misspecified.

\cite{zhu2018linear} proposed an inference method for prediction with a general loading $\xi\in\mathbb{R}^p$ in the high-dimensional LM setting. Their
method starts from a transformation as follows:
\begin{align*}
   X_i^\T \bar\beta  = Z_i (\xi^\T\bar\beta) + \tilde{W}^\T_i\alpha^* ,
\end{align*}
where $Z_i=\xi^\T X_i/\Vert\xi\Vert_2^2 \in \mathbb{R}$, $\tilde{W}_i=U_a^\T W_i \in \mathbb{R}^{p-1}$, $\alpha^*=U_a^\T\bar\beta$, $W_i=Q_\xi X_i$, $Q_\xi=I_p-P_\xi$, $P_\xi = \xi\xi^\T/\Vert \xi\Vert_2^2$, and $U_a\in \mathbb{R}^{p \times (p-1)}$ is an orthogonal matrix such that $U_a^\T U_a=I_{p-1}$ and $Q_\xi = U_aU_a^\T$. Let $\gamma^*\in \mathbb{R}^{p-1}$ be the target coefficient in the linear regression of $Z_i$ on $\tilde{W}_i$.
\cite{zhu2018linear} proposed estimators, $\hat\alpha$ and $\hat\gamma$, for $\alpha^*$ and $\gamma^*$ and showed that under suitable conditions,
\begin{align*}
    S_n = \sqrt{n}\frac{\sum_{i=1}^n(Z_i-\tilde{W_i}\hat{\gamma})(Y_i - Z_i (\xi^\T\bar\beta) - \tilde{W}_i\hat{\alpha} )}{\Vert Z-\tilde{W}\hat{\gamma}\Vert_2 \Vert Y - Z (\xi^\T\bar\beta) - \tilde{W}\hat{\alpha}\Vert_2}
\end{align*}
converges to $\mathcal{N}(0,1)$ in distribution. Then $S_n$ can be applied for inference of $\xi^\T\bar\beta$.
In contrast with our work, their theory requires that LM is correctly specified and relies on directly the sparsity of $\gamma^*$ (which depends on $\xi$),
even though not on the sparsity of $\bar\beta$.

\cite{yliu} proposed an estimation method for prediction with a general loading $\xi\in\mathbb{R}^p$ in the high-dimensional LM setting.
Their method starts from a transformation related to that in \cite{zhu2018linear} but in a matrix form:
\begin{align*}
    X \bar\beta  = \tilde Z \psi^* + \Gamma \theta^* ,
\end{align*}
where $X=(X_1,\ldots,X_n)^\T \in \mathbb{R}^{n\times p}$, $\tilde Z = X \xi/\Vert \xi\Vert_2 \in \mathbb{R}^{n\times 1}$, $\psi^* =\xi^\T\bar\beta/\Vert \xi\Vert_2\in \mathbb{R}$ and $\theta^* = \Gamma^{-1 }X Q_\xi\bar\beta$ for a properly chosen matrix $\Gamma\in\mathbb{R}^{n\times n}$.
\cite{yliu} then employed Lasso or other penalized estimators, $\hat\psi$ and $\hat\theta$, for $\psi^*$ and $\theta^*$
in the transformed linear regression of $Y=(Y_1,\ldots,Y_n)^\T$ on $(\tilde Z, \Gamma)$, and derived
estimation error bounds for $\hat\psi$ provided that either $\bar\beta$ is sparse or  $(\lambda_1(X Q_\xi X^\T), \ldots, \lambda_n(X Q_\xi X^\T))^\T$ is sparse, where $\lambda_i(X Q_\xi X^\T)$ is the $i$th eigenvalue of $X Q_\xi X^\T$.
In contrast with our work, their theory requires that LM is correctly specified. In addition, their work does not provide inference such as confidence intervals for $\psi^*=\xi^\T \bar\beta/\Vert \xi\Vert_2$.

\subsection{Organization and notation} \label{sec:notation}

The rest of the paper is organized as follows. In Section \ref{sec:high_dimensional_linear_model}, we introduce the debiased estimator of $\xi^\T\bar\beta$ and present the theoretical properties in misspecified high-dimensional LM setting, depending on an estimator of $\Omega^*$ subject to certain convergence rates. We provide specific estimators for $\Omega^*$ satisfying the desired convergence rates and present theoretical results for the debiased estimators of $\xi^\T\bar\beta$ incorporating the proposed estimators of $\Omega^*$.
In Section \ref{sec:high_dimen_glm}, we discuss our method and theory in misspecified high-dimensional GLM setting.
In Section \ref{sec:outline_of_proof}, we give proof outlines of our main theoretical results.
In Section \ref{sec:numerical_study}, we present numerical studies to support our theoretical results.

We introduce the following notations used in the article. We denote
$\phi_k(a,b)=\frac{\partial^k\phi(a,b)}{\partial^k b}$ for $k=1,2,3$. For a vector $a=(a_1,\ldots,a_p)^\T\in\mathbb{R}^p$, we denote as $a_j$ the $j$th entry of $a$ and define $\Vert a\Vert_1=\sum_{j=1}^p |a_j|$, $\Vert a\Vert_2=\sqrt{\sum_{j=1}^p a_j^2}$, $\Vert a\Vert_\infty = \max_{j=1,\ldots,p}\{|a_j|\}$ and $\Vert a \Vert_0=|\supp(a)|$, where $\supp(a)$ is the support of $a$. For a matrix $A\in \mathbb{R}^{p\times q}$, we denote as $A_{\cdot j}$ the $j$th column of $A$ and define $\|A\|_\max = \max_{1\leq i\leq p,1\leq j\leq q}\{|A_{ij}|\}$, $\| A\|_{L_1} = \sum_{1\leq i\leq p, 1\leq j\leq q} |A_{ij}|$, $\| A \|_\F = \sqrt{\sum_{1\leq i\leq p, 1\leq j\leq q}A^2_{ij}}$, $\| A\|_{1}=\inf\{c\geq 0:\Vert Av\Vert_1 \leq c\Vert v\Vert_1 \text{ for all $v\in \mathbb{R}^q$}\}=\max_{1\leq j\leq q}\{\sum_{i=1}^p |A_{ij}|\}$, and $\Vert A\Vert_2$ is the $L_2$ operator norm of $A$,
i.e., $\Vert A\Vert_2 = \inf\{c\geq 0:\Vert Av\Vert_2 \leq c\Vert v\Vert_2 \text{ for all $v\in \mathbb{R}^q$}\}$. For $j=1,\ldots,p$, $e_j\in \mathbb{R}^p$ is a vector with 1 in the $j$th entry and 0 in other entries. For any subset $S\subset \{1,2,\ldots, p\}$ and matrix $M\in \mathbb{R}^{p\times p}$, $M_{S,S}\in \mathbb{R}^{|S|\times |S|}$ is a submatrix of $M$ containing rows and columns with indexes belonging to $S$.

Given two positive sequences $\{a_k\}_{k\geq 1}$ and $\{b_k\}_{k\geq 1}$, $a_k=o(b_k)$ if $\lim_{k\rightarrow \infty}\frac{a_k}{b_k}=0$ and $a_k\asymp b_k$ indicates that there exist constants $c>0$ and $K\geq 1$ such that $a_k \leq c b_k$ and $b_k\leq ca_k$ for all $k\geq K$. Given random variables sequence $\{X_k\}_{k\geq 1}$ and constants sequence $\{a_k\}_{k\geq 1}$, we denote $X_k=O_p(a_k)$ if there exists $C>0$ such that for any $\epsilon>0$, there exists $N>0$ such that $\mathbb{P}(|X_k|\leq C|a_k|)>1-\epsilon$ for any $k>N$ and we denote $X_k = o_p(a_k)$ if for any $\epsilon > 0$, $\lim_{k\to\infty}\mathbb{P}(|X_k|>\epsilon |a_k|)=0$. Given random variable sequence $\{X_k\}_{k\geq 1}$ and random variable $X$, we denote $X_k \rightarrow_d X$ if the random variables sequence $\{X_k\}_{k\geq 1}$ converges to $X$ in distribution.

\begin{comment}
\LL{For a matrix $A\in \mathbb{R}^{p\times p}$, $\diag(A)\in \mathbb{R}^p$ is a vector whose entries are diagonal entries of $A$. For a vector $a\in\mathbb{R}^p$, $\diag(a)\in \mathbb{R}^{p\times p}$ is a matrix whose diagonal entries are the entries in $a$ and off-diagonal entries are zero. $\mathcal{U}(a,b)$ is a uniform distribution from $[a,b]$. $e_j\in\mathbb{R}^p$ for $j=1,\ldots,p$ is a vector with one at the $j$-th position and zero everywhere else. We denote $X = (X_1,\ldots,x_n)^\T \in \mathbb{R}^{n\times p}$, $X_{d}$ for $d=1,\ldots, p$ to be the vector comprised of the
$d$ column of $X$ and $X_{-d}$ for $d=1,\ldots, p$ to be the matrices comprised of the columns of $X$ except the $d$ column.}
\end{comment}

\section{Debiased prediction and inference in misspecified LM settings}
\label{sec:high_dimensional_linear_model}

\subsection{High-level ideas} \label{sec:high-level-LM}

Before moving to the cross-fitting version of our method, we discuss main ideas behind
our method and theory as presented in the Introduction.
Suppose that
$\Vert \hat{\beta}-\bar\beta\Vert_1=O_p(s_0\sqrt{\frac{\log(p\vee n)}{n}})$, and $(\hat{\beta}-\bar\beta)^\T\hat{\Sigma} (\hat{\beta}-\bar\beta)=O_p(s_0\frac{\log(p\vee n)}{n})$,
and suitable regularity conditions are satisfied as in Proposition \ref{prop:lm}.
We focus on the justification of asymptotic expansion (\ref{result_LM_1}).
Consider the following decomposition of the difference between $\widehat{\xi^\T\beta}$ in (\ref{eq:debiased_lm}) and $\xi^\T\bar\beta$:
\begin{align}
    \widehat{\xi^\T\beta} - \xi^\T\bar\beta &= u^{* \T}\frac{1}{n}\sum_{i=1}^n X_i(Y_i-X_i^\T\bar\beta)  + (\hat{\Sigma} u^* -\xi)^\T(\bar\beta-\hat{\beta})  + (\tilde{u}-u^*)^\T\frac{1}{n} \sum_{i=1}^n X_i(Y_i-X_i^\T\bar\beta)\notag \\
    &- (\tilde{u}-u^*)^\T\hat{\Sigma}(\hat{\beta}-\bar\beta) , \label{debiased_est_diff}
\end{align}
where $\hat{\Sigma}=\frac{1}{n}\sum_{i=1}^nX_iX_i^\T$.
The first term on the right-hand side of (\ref{debiased_est_diff}) gives the leading term in (\ref{result_LM_1}), by the definition $u^*= \Omega^* \xi$.
Note that $u^{* \T} X_i(Y_i-X_i^\T\bar\beta)$, $i=1,\ldots,n$, are mean-zero i.i.d.~variables by the definition of $\bar\beta$ in \eqref{eq:target_beta_LM}.
By standard arguments,
\begin{align*}
   \frac{\sqrt{n}}{\sqrt{V}} (u^{* \T}\frac{1}{n}\sum_{i=1}^n X_i(Y_i-X_i^\T\bar\beta) ) \rightarrow_d \mathcal{N}(0, 1),
\end{align*}
and $V \asymp \Vert \xi\Vert^2_2$.
To obtain asymptotic expansion (\ref{result_LM_1}), it suffices to control the remaining terms on the right-hand side of (\ref{debiased_est_diff})
as $o_p(\frac{\Vert \xi\Vert_2}{\sqrt{n}})$.
We distinguish two different scenarios, mainly depending on different sub-gaussian assumptions on $X_0$, Assumption \ref{ass:sub_gaussian_entries} or \ref{ass:sub_gaussian_LM}.

\textbf{First scenario.}\;
Suppose that Assumption \ref{ass:sub_gaussian_entries} holds, i.e., (i) the entries of $X_0$ are sub-gaussian random variables and
(ii) $X_0^\T \Omega^*\xi$ is a sub-gaussian random variable with sub-gaussian norm being $O(\Vert \Omega^* \xi \Vert_2)$, which is $O(\Vert  \xi \Vert_2)$ provided Assumption \ref{sigma_eigenvalue}(ii).
The independence of $\tilde{\Omega}$ with $\{(X_i, Y_i)\}_{i=1,\ldots, n}$ is not required,
i.e., $\tilde\Omega$ can also be defined by Algorithm \ref{alg:omega_est_prop} with $\{(X^\prime_i, Y^\prime_i)\}_{i=1,\ldots, m}$ replaced by $\{(X_i, Y_i)\}_{i=1,\ldots, n}$.

The term $(\hat{\Sigma} u^* -\xi)^\T(\bar\beta-\hat{\beta})$ can be controlled by the H\"{o}lder inequality: $|(\hat{\Sigma} u^* -\xi)^\T(\bar\beta-\hat{\beta})|\leq \Vert \hat{\Sigma} u^* -\xi\Vert_\infty \Vert \bar\beta-\hat{\beta}\Vert_1$. The $j$th entry of $\hat{\Sigma} u^* -\xi$ can be expressed as $\frac{1}{n}\sum\limits_{i=1}^n(X_{ij}(X_i^\T \Omega^* \xi)-\mathbb{E}[X_{ij}(X_i^\T \Omega^* \xi)])$.
Under Assumption \ref{ass:sub_gaussian_entries} and \ref{sigma_eigenvalue}(ii), it can be shown that
 $\Vert \hat{\Sigma} u^* -\xi\Vert_\infty=O_p(\sqrt{\frac{\log(p \vee n)}{n}}\Vert \xi\Vert_2)$ (Lemma \ref{lem:bound_of_random_quan} (v)).
Together with the assumption on $\Vert \hat{\beta}- \bar\beta\Vert_1$, we have
\begin{align}
    (\hat{\Sigma} u^* -\xi)^\T(\bar\beta-\hat{\beta})=O_p(s_0\frac{\log(p\vee n)}{n}\Vert\xi\Vert_2) . \label{eq:sigma_u_beta}
\end{align}

The term  $(\tilde{u}-u^*)^\T\frac{1}{n} \sum\limits_{i=1}^n X_i(Y_i-X_i^\T\bar\beta)$ can also be controlled by the H\"{o}lder inequality:
\begin{align*}
 \left| (\tilde{u}-u^*)^\T\frac{1}{n} \sum\limits_{i=1}^n X_i(Y_i-X_i^\T\bar\beta) \right| \leq \Vert \tilde{u}-u^*\Vert_1\Vert \frac{1}{n} \sum\limits_{i=1}^n X_i(Y_i-X_i^\T\bar\beta)\Vert_\infty.
\end{align*}
By the definition of $\bar\beta$ in (\ref{eq:target_beta_LM}), we have $\mathbb{E}[X_{0j}(Y_0-X_0^\T\bar\beta)] = 0$ for $j=1,\ldots, p$.
Under Assumption \ref{ass:sub_gaussian_entries}(i) and \ref{sub_gaussian_residual}, it can be shown that $\Vert \frac{1}{n} \sum\limits_{i=1}^n X_i(Y_i-X_i^\T\bar\beta)\Vert_\infty=O_p(\sqrt{\frac{\log(p\vee n)}{n}})$ (Lemma \ref{lem:bound_of_random_quan} (iv)). Hence
\begin{align}
    (\tilde{u}-u^*)^\T\frac{1}{n} \sum\limits_{i=1}^n X_i(Y_i-X_i^\T\bar\beta) = O_p(\Vert \tilde{u}-u^*\Vert_1 \sqrt{\frac{\log(p\vee n)}{n}}) . \label{eq:u_Y_Xbeta_LM}
\end{align}

The term $(\tilde{u}-u^*)^\T\hat{\Sigma}(\bar\beta-\hat{\beta})$ can be controlled by the Cauchy--Schwartz inequality:
\begin{align}
    |(\tilde{u}-u^*)^\T\hat{\Sigma}(\hat{\beta}-\bar\beta)|&\leq  \sqrt{(\tilde{u}-u^*)^\T\hat{\Sigma} (\tilde{u}-u^*)} \sqrt{(\hat{\beta}-\bar\beta)^\T\hat{\Sigma}(\hat{\beta}-\bar\beta)} .\label{inq:cross_term_u_beta}
\end{align}
Under Assumption \ref{ass:sub_gaussian_entries} (i), it can be shown that $\Vert \hat{\Sigma}-\Sigma\Vert_{\max}=O_p(\sqrt{\frac{\log(p\vee n)}{n}})$ (Lemma \ref{lem:bound_of_random_quan} (i)).
Then $(\tilde{u}-u^*)^\T\hat{\Sigma} (\tilde{u}-u^*)=(\tilde{u}-u^*)^\T(\hat{\Sigma}-\Sigma) (\tilde{u}-u^*)+(\tilde{u}-u^*)^\T\Sigma (\tilde{u}-u^*) \leq  \Vert\tilde{u}-u^*\Vert_1^2 \Vert\hat{\Sigma}-\Sigma\Vert_{\max}+(\tilde{u}-u^*)^\T\Sigma (\tilde{u}-u^*)=O_p(\Vert\tilde{u}-u^*\Vert_1^2 \sqrt{\frac{\log(p\vee n)}{n}} + \Vert\tilde{u}-u^*\Vert_2^2 )$. Together with the assumption on $(\hat{\beta}-\bar\beta)^\T\hat{\Sigma} (\hat{\beta}-\bar\beta)$, we have
\begin{align}
    (\tilde{u}-u^*)^\T\hat{\Sigma}(\hat{\beta}-\bar\beta) =O_p(\sqrt{\Vert\tilde{u}-u^*\Vert_1^2 \sqrt{\frac{\log(p\vee n)}{n}} + \Vert\tilde{u}-u^*\Vert_2^2} \sqrt{\frac{s_0\log(p\vee n) }{n}}) . \label{eq:cross_term_lm}
\end{align}

From (\ref{eq:sigma_u_beta}), (\ref{eq:u_Y_Xbeta_LM}) and (\ref{eq:cross_term_lm}),
the non-leading terms on the right-hand side of (\ref{debiased_est_diff}) can be controlled through $\Vert \tilde{u}-u^*\Vert_{1}$ and $\Vert \tilde{u}-u^*\Vert_{2}$,
which in turn can be controlled by
$\Vert \tilde{\Omega}_{S_\xi,S_\xi}-\Omega^*_{S_\xi,S_\xi}\Vert_1\Vert\xi\Vert_1$ and $\Vert \tilde{\Omega}_{S_\xi,S_\xi}-\Omega^*_{S_\xi,S_\xi}\Vert_2\Vert\xi\Vert_2$ respectively.
Therefore, if $\frac{\Vert\xi\Vert_1}{\Vert\xi\Vert_2}$ is bounded,
then  asymptotic expansion (\ref{result_LM_1}) can be obtained from (\ref{debiased_est_diff}) by controlling
$\Vert \tilde{\Omega}_{S_\xi,S_\xi}-\Omega^*_{S_\xi,S_\xi}\Vert_1$ and $\Vert \tilde{\Omega}_{S_\xi,S_\xi}-\Omega^*_{S_\xi,S_\xi}\Vert_2$.

\vspace{.05in}
\textbf{Second scenario.}\;
With possibly large $\frac{\Vert\xi\Vert_1}{\Vert \xi\Vert_2}$, we develop a different approach
to control the terms $(\tilde{u}-u^*)^\T\frac{1}{n} \sum\limits_{i=1}^n X_i(Y_i-X_i^\T\bar\beta)$ and $(\tilde{u}-u^*)^\T\hat{\Sigma}(\bar\beta-\hat{\beta})$.
Suppose that
Assumption \ref{ass:sub_gaussian_LM} holds (i.e., $X_0$ is a sub-gaussian vector) and $\tilde{\Omega}$ is independent of $\{(X_i, Y_i)\}_{i=1,\ldots, n}$.

To control the term $(\tilde{u}-u^*)^\T\frac{1}{n} \sum\limits_{i=1}^n X_i(Y_i-X_i^\T\bar\beta)$, we note that
conditionally on $\tilde u$,
$\mathbb{E}[(\tilde{u}-u^*)^\T X_{i}(Y_i-X_i^\T\bar\beta)] = 0$ by the definition of $\bar\beta$ and the independence of $\tilde{\Omega}$ with $(X_i,Y_i)$ for $i=1,\ldots, n$,
and $(\tilde{u}-u^*)^\T X_{i}$ is a sub-gaussian variable with sub-gaussian norm being $O(\Vert \tilde{u}-u^* \Vert_2)$ by Assumption \ref{ass:sub_gaussian_LM}.
Then conditionally on $\tilde u$,
 $(\tilde{u}-u^*)^\T X_{i}(Y_i-X_i^\T\bar\beta)$ is a centered sub-exponential random variable with sub-exponential norm being $O(\Vert \tilde{u}-u^*\Vert_2)$
 provided that $Y_i-X_i^\T\bar\beta$ is sub-gaussian for $i=1,\ldots,n$.
By sub-exponential concentration properties, we have (Lemma \ref{lem:bound_of_random_quan} (x))
\begin{align}
    (\tilde{u}-u^*)^\T\frac{1}{n} X_i(Y_i-X_i^\T\bar\beta) = O_p(\Vert \tilde{u}-u^*\Vert_2 \sqrt{\frac{\log(p\vee n)}{n}}). \label{inq:u_y_xi_beta_vector_ass}
\end{align}
Compared to (\ref{eq:u_Y_Xbeta_LM}), the right-hand side of (\ref{inq:u_y_xi_beta_vector_ass}) depends on $\Vert \tilde{u}-u^*\Vert_2$, but not $\Vert \tilde{u}-u^*\Vert_1$.

To control the term $(\tilde{u}-u^*)^\T\hat{\Sigma}(\bar\beta-\hat{\beta})$,
it suffices, from (\ref{inq:cross_term_u_beta}), to control $(\tilde{u}-u^*)^\T\hat{\Sigma} (\tilde{u}-u^*)$, which can be expressed as $\frac{1}{n}\sum_{i=1}^n (( \tilde{u}-u^*)^\T X_i)^2$.
By Assumption \ref{ass:sub_gaussian_LM} and the independence of $\tilde{\Omega}$ with $X_i$, we have that $((\tilde{u}-u^*)^\T X_i)^2$ is a sub-exponential random variable with sub-exponential norm being $O(\Vert \tilde{u}-u^*\Vert^2_2)$, and hence $( \tilde{u}-u^*)^\T\hat{\Sigma} (\tilde{u}-u^*) =O_p(\Vert \tilde{u}-u^*\Vert_2^2)$ (Lemma \ref{lem:bound_of_random_quan} (xi)).
Together with (\ref{inq:cross_term_u_beta}) and the assumption on $(\hat{\beta}-\bar\beta)^\T\hat{\Sigma} (\hat{\beta}-\bar\beta)$, we have
\begin{align}
    (\tilde{u}-u^*)^\T\hat{\Sigma}(\hat{\beta}-\bar\beta)=O_p(\Vert\tilde{u}-u^*\Vert_2\sqrt{\frac{s_0\log(p\vee n) }{n}} ) . \label{eq:cross_u_sigma_beta}
\end{align}
Compared to (\ref{eq:cross_term_lm}), the right-hand side of (\ref{eq:cross_u_sigma_beta}) depends on $\Vert \tilde{u}-u^*\Vert_2$, but not $\Vert \tilde{u}-u^*\Vert_1$.

From (\ref{eq:sigma_u_beta}), (\ref{inq:u_y_xi_beta_vector_ass}) and (\ref{eq:cross_u_sigma_beta}),
the non-leading terms on the right-hand side of (\ref{debiased_est_diff}) can be controlled through $\Vert \tilde{u}-u^*\Vert_{2}$,
independently of $\Vert \tilde{u}-u^*\Vert_{1}$. Hence asymptotic expansion (\ref{result_LM_1}) can be obtained from (\ref{debiased_est_diff}) by controlling
$\Vert (\tilde{\Omega}-\Omega^*)_{S_\xi,S_\xi}\Vert_2$ regardless of the magnitude of $\frac{\Vert \xi\Vert_1}{\Vert\xi\Vert_2}$.

\vspace{.05in}
\textbf{Summary.}\;
To complete the justification of asymptotic expansion (\ref{result_LM_1}), the key is to find a suitable estimator $\tilde{\Omega}$ such that
the \textit{matrix norms} $\Vert \tilde{\Omega}_{S_\xi,S_\xi}-\Omega^*_{S_\xi,S_\xi}\Vert_1$ and $\Vert \tilde{\Omega}_{S_\xi,S_\xi}-\Omega^*_{S_\xi,S_\xi}\Vert_2$ are
$O_p( \sqrt{\frac{\log(p\vee n)}{n}})$, up to some sparsity pre-factors.
We take $\tilde \Omega$ to be the CLIME estimator in Algorithm \ref{alg:omega_est_prop} and show that
$\Vert \tilde{\Omega}_{S_\xi,S_\xi}-\Omega^*_{S_\xi,S_\xi}\Vert_1=O_p(\Vert s_{1\xi}\Vert_\infty\sqrt{\frac{\log(p\vee n)}{n}})$ and $\Vert \tilde{\Omega}_{S_\xi,S_\xi}-\Omega^*_{S_\xi,S_\xi}\Vert_2=O_p(\min\{\sqrt{\Vert s_{1\xi}\Vert_1}, \Vert s_{1\xi}\Vert_\infty\}\sqrt{\frac{\log(p\vee n)}{n}})$ (Lemma \ref{lem:error_rate_omega_clime}).
Then we obtain the sparsity conditions as stated in Proposition \ref{prop:lm} (i) and (ii),
corresponding to respectively the first and second scenarios above.

\begin{comment}
Because $|(\hat{\Sigma} u^* -\xi)^\T(\bar\beta-\hat{\beta})|\leq \Vert \hat{\Sigma} u^* -\xi\Vert_\infty \Vert \bar\beta-\hat{\beta}\Vert_1$ and $\Vert \hat{\Sigma} u^* -\xi\Vert_\infty$ can be controlled under a mild regular assumption from the concentration property, $|(\hat{\Sigma} u^* -\xi)^\T(\bar\beta-\hat{\beta})|$ can be easily controlled given $\Vert \bar\beta-\hat{\beta}\Vert_1$ is controlled. For establishing the asymptotic normality of $\widehat{\xi^\T\beta} - \xi^\T\bar\beta$, it is important to find $\tilde{u}$ such that $\tilde{u}$ is a good estimator of $u^*$. In the following discussion, we will construct the debiased estimator of $\xi^\T\bar\beta$ by giving $\tilde{u}=\hat{\Omega}\xi$ where $\hat{\Omega}$ is an oracle estimator of $\Omega^*$. We justify our debiased estimator assuming the oracle estimator satisfies certain conditions. We give or introduce algorithms to compute the estimators of $\Omega^*$ which can be incorporated into the proposed debiased estimator.
\end{comment}

\subsection{Debiased prediction with cross-fitting in LM settings}

\label{sec:sparse_omega}

We give our debiased estimator of $\xi^\T\bar\beta$ with cross-fitting in the misspecified LM setting.
Specifically, we divide the sample into two complementary subsets of $\{1,\ldots, n\}$, that is $J,J^\c \subset\{1,\ldots,n\}$ and $J\cup J^\c = \{1,\ldots,n\}$. We assume that $|J|>\rho n$  and $|J^\c|>\rho n$, where $0<\rho < 1$. We denote as $X_J\in\mathbb{R}^{|J|\times p}$ and $Y_J \in \mathbb{R}^{|J|}$ the covariates matrix and response vector with row index belonging to $J$ and as $X_{J^\c}\in\mathbb{R}^{|J^\c|\times p}$ and $Y_{J^\c} \in \mathbb{R}^{|J^\c|}$ the covariate matrix and response vector with row index belonging to $J^\c$.

Let $\hat{\beta}\in \mathbb{R}^p$ be
an initial regularized sparse estimator of $\bar\beta$, for example, a Lasso estimator computed from $\{(X_i,Y_i)\}_{i=1,\ldots,n}$.
Our debiased estimator of $\xi^\T\bar\beta$ is defined as
\begin{align}
    \widehat{\xi^{\T}\beta} = \xi^\T \hat{\beta} + \xi^\T\hat{\Omega}^{J^\c} \frac{1}{n}\sum_{i\in J}X_i(Y_i-X_i^\T\hat{\beta}) + \xi^\T\hat{\Omega}^{J} \frac{1}{n}\sum_{i\in J^\c}X_i(Y_i-X_i^\T\hat{\beta}) , \label{debiased_est_1}
\end{align}
where $\hat{\Omega}^J$ and $\hat{\Omega}^{J^\c}$ are estimators of $\Omega^*$ based on $X_J$ and $X_{J^\c}$ respectively, which are specified later in Section \ref{sec:omega_estimators}.
By construction, $\hat{\Omega}^{J^\c}$ is independent of $(X_J,Y_J)$ and $\hat{\Omega}^J$ is independent of $(X_{J^\c}, Y_{J^\c})$,
which is important for theoretical analysis in the second scenario as discussed in Section~\ref{sec:high-level-LM}.
To construct a Wald confidence interval for $\xi^\T\bar\beta$, our estimator of the asymptotic variance for $\widehat{\xi^{\T}\beta}$ is defined as
\begin{align}
    \hat{\V} = \frac{1}{n} \sum_{i\in J} (\xi^\T\hat{\Omega}^{J^\c}X_i)^2(Y_i-X_i^\T\hat{\beta})^2 + \frac{1}{n}\sum_{i\in J^\c}(\xi^\T\hat{\Omega}^{J}X_i)^2(Y_i-X_i^\T\hat{\beta})^2. \label{eq:definition_V_hat}
\end{align}

\subsection{Theoretical properties of debiased prediction in LM settings}
\label{sec:theorey_LM}

We give the theoretical properties of our debiased estimator of $\xi^\T\bar\beta$ and the variance estimator $\hat V$ introduced in Section \ref{sec:sparse_omega},
assuming certain convergence rates on $\hat\beta$, $\hat{\Omega}^J$ and $\hat{\Omega}^{J^\c}$.
Specifically, we assume that $\hat{\beta}$ satisfies the following condition, which is satisfied by Lasso estimation under suitable regularity conditions
\citep{lasso_and_dantzig, hastie_book}.

\begin{conditionB}
\label{beta_hat_condition}
With probability at least $1-\delta_1(n)$, the estimator $\hat{\beta}$ satisfies
\begin{align*}
    \Vert \hat{\beta}-\bar\beta\Vert_1 &\leq C_1 s_0\sqrt{\frac{\log(p\vee n)}{n}}, \notag \\
    \frac{1}{n}\sum_{i=1}^n(X_i^\T \hat{\beta} -X_i^\T \bar\beta)^2 &\leq C_1 s_0\frac{\log(p \vee n)}{n},
\end{align*}
where $\delta_1(n)$ satisfies $\delta_1(n)\rightarrow 0$ as $n\rightarrow \infty$ and $C_1$ is some positive constant.
\end{conditionB}

We also assume that the estimators $\hat{\Omega}^J$ and $\hat{\Omega}^{J^\c}$ satisfy the following condition. Specification of such estimators
to satisfy the conditions is deferred to Section \ref{sec:omega_estimators}.

\begin{conditionC}
\label{con:omage_pred}
For some positive constant $C_2$, and constants $A_{1}$ and $A_{2}$ depending only on $(s_1, s_{1\xi})$,
the following conditions hold for $\hat{\Omega}^J$ and $\hat{\Omega}^{J^\c}$.
\label{precision_matrix_LM}
$\newline$
\noindent  (i) With probability at least $1-\delta_2(n)$, $\hat{\Omega}^J$ satisfies,
\begin{align*}
    \Vert \hat{\Omega}^{J}_{S_\xi, S_\xi}-\Omega^*_{S_\xi, S_\xi}\Vert_2 &\leq  C_2 A_{1}  \sqrt{\frac{\log(p\vee n)}{n}} , \\
    \Vert \hat{\Omega}^{J}_{S_\xi,S_\xi}-\Omega^*_{S_\xi,S_\xi}\Vert_1 &\leq   C_2 A_{2} \sqrt{\frac{\log(p\vee n)}{n}},
\end{align*}
where $\delta_2(n)$ satisfies $\delta_2(n)\rightarrow 0$ as $n\rightarrow \infty$.

\noindent (ii) With probability at least $1-\delta_3(n)$, $\hat{\Omega}^{J^\c}$ satisfies,
\begin{align*}
   \Vert \hat{\Omega}^{J^\c}_{S_\xi,S_\xi}-\Omega^*_{S_\xi,S_\xi}\Vert_2  &\leq   C_2 A_{1} \sqrt{\frac{\log(p\vee n)}{n}} ,\\
    \Vert \hat{\Omega}^{J^\c}_{S_\xi,S_\xi}-\Omega^*_{S_\xi,S_\xi}\Vert_1 &\leq   C_2 A_{2} \sqrt{\frac{\log(p\vee n)}{n}} ,
\end{align*}
where $\delta_3(n)$ satisfies $\delta_3(n)\rightarrow 0$ as $n\rightarrow \infty$.
\end{conditionC}

As mentioned earlier, we consider one of the following conditions on $X_0$, Assumption \ref{ass:sub_gaussian_entries} or \ref{ass:sub_gaussian_LM}.
If  Assumption \ref{ass:sub_gaussian_LM} holds,
then the entries of $X_0$ are sub-gaussian with norm bounded by $\sigma_1$, and $\xi^\T\Omega^* X_0$ is sub-gaussian with norm bounded by $\sigma_1\Vert \Omega^*\xi\Vert_2$. In other words,
Assumption \ref{ass:sub_gaussian_entries} is implied by (hence weaker than) Assumption \ref{ass:sub_gaussian_LM}.

\begin{assumptionA}
\label{ass:sub_gaussian_entries}
$\newline$
\noindent (i) The entries of $X_0$ are sub-gaussian random variables with sub-gaussian norm bounded by $\sigma_1$.
$\newline$
\noindent (ii) $\xi^\T\Omega^* X_0$ is a sub-gaussian random variable with sub-gaussian norm bounded by $\sigma_1\Vert \Omega^* \xi\Vert_2$.

\end{assumptionA}

\begin{assumptionA}
\label{ass:sub_gaussian_LM}
$X_0$ is sub-gaussian random vector with sub-gaussian norm bounded by $\sigma_1$.
\end{assumptionA}

We also make the following assumptions on $\Omega^*$ and the noise $Y_0-X_0^\T\bar\beta$.
Assumption \ref{sigma_eigenvalue} is also used in \cite{cai2020optimal}, and similar assumptions as Assumption \ref{sub_gaussian_residual} are used in \cite{Ning}
and \cite{cai2020optimal} among others.

\begin{assumptionA}
\label{sigma_eigenvalue}
$\newline$
\noindent (i) $\Omega^*$ satisfies $C_3 \leq \lambda_{\min}(\Omega^*)$.

\noindent (ii) $\Omega^*$ satisfies $\lambda_{\max}(\Omega^*) \leq C_4$.
$\newline$
 \noindent  Here $C_3$ and $C_4$ are some positive constants.
\end{assumptionA}

\begin{assumptionA}
\label{sub_gaussian_residual}
    $Y_0-X_0^\T\bar\beta$ is a sub-gaussian random variable with sub-gaussian norm bounded by $\sigma_2$ and
    $\mathbb{E}[(Y_0-X_0^\T\bar\beta)^2|X_0] \geq C_5$, where $C_5$ is a positive constant.
\end{assumptionA}

To present the theoretical properties of our method,
we use the following expressions. For $\widehat{\xi^\T\beta}$ in (\ref{debiased_est_1}) and $\hat{V}$ in (\ref{eq:definition_V_hat}), let
\begin{subequations}
\begin{equation}
 \frac{\widehat{\xi^{\T}\beta} - \xi^\T\bar\beta}{\sqrt{V}}= \frac{1}{n}\sum_{i=1}^n \frac{(\xi^\T \Omega^* X_i)(Y_i-X_i^\T\bar\beta)}{\sqrt{V}} + \Delta_1  n^{-\frac{1}{2}} , \label{eq:transform_lm_1}
\end{equation}
\begin{equation}
\frac{\hat{V}}{V} = 1+\Delta_2,\label{eq:transform_lm_2}
\end{equation}
\end{subequations}
where $V=\mathbb{E}[(\xi^\T \Omega^* X_i)^2(Y_i-X_i^\T\bar\beta)^2]$.
The theoretical properties of $\widehat{\xi^\T\beta}$ in (\ref{debiased_est_1}) and $\hat{V}$ in (\ref{eq:definition_V_hat}) are established by characterizing $\Delta_1$ and $\Delta_2$ in the theorem below.

\begin{thm}
\label{thm_1}
\noindent Suppose that Conditions \ref{beta_hat_condition} and \ref{precision_matrix_LM} are satisfied. Then
the following results hold for $\Delta_1$ in (\ref{eq:transform_lm_1}) and $\Delta_2$ in (\ref{eq:transform_lm_2}).

\noindent (i) If further Assumptions \ref{ass:sub_gaussian_entries}, \ref{sigma_eigenvalue}, and \ref{sub_gaussian_residual} are satisfied, $A_2\frac{\Vert \xi\Vert_1}{\Vert \xi\Vert_2} \tau_n \leq 1$ and $s_0 \tau_n \leq \min\{\frac{c(1)\rho\sqrt{n} }{3},\frac{c(\frac{2}{3})n^{\frac{1}{6}}}{4},1\}$, then we have with probability at least $1-\delta_1(n)-\delta_2(n)-\delta_3(n)-\frac{4}{n}-\frac{14}{p \vee n}- 2\exp\{-c(\frac{1}{2})n^{\frac{1}{3}})\}-4\exp\{-\frac{1}{2}\sqrt{n}\}$,
\begin{align*}
  |\Delta_1|& \leq   C (A_2 \frac{\Vert \xi\Vert_1}{\Vert \xi\Vert_2} +A_1 \sqrt{s_0}  +s_0)\tau_n, \\
    |\Delta_2| &\leq  C (n^{-\frac{1}{3}}  +  \sqrt{s_0} \tau_n  +A_2 \frac{\Vert \xi\Vert_1}{\Vert \xi\Vert_2} \tau_n ).
\end{align*}

\noindent (ii) If further Assumptions \ref{ass:sub_gaussian_LM}, \ref{sigma_eigenvalue}, and \ref{sub_gaussian_residual} are satisfied,
$A_1 \tau_n \leq 1$ and $s_0\tau_n \leq \min\{c(1)\rho \sqrt{n}, \frac{c(1)\sqrt{n}}{2}, \frac{c(\frac{2}{3})n^{\frac{1}{6}}}{3},1\}$,
then we have with probability at least $1-\delta_1(n)-\delta_2(n)-\delta_3(n)-\frac{8}{n}-\frac{16}{p \vee n}- 2\exp\{-c(\frac{1}{2})n^{\frac{1}{3}})\}-4\exp\{-\frac{1}{2}\sqrt{n}\}$,
\begin{align*}
  |\Delta_1|& \leq   C ( A_1 \sqrt{s_0} +s_0)\tau_n, \\
    |\Delta_2| &\leq  C (n^{-\frac{1}{3}}  +  \sqrt{s_0} \tau_n
    +A_1  \tau_n).
\end{align*}
Here $\tau_n=\frac{\log(p\vee n)}{\sqrt{n}}$, $c(\frac{1}{2})$, $c(\frac{2}{3})$ and $c(1)$ are constants from Lemma \ref{lem:alpha_sub_exponential_concentration} and $C$ is a constant depending on $(C_1, C_2, C_3, C_4, C_5, \sigma_1,\sigma_2,\rho)$ only.
\end{thm}

\subsection{Estimators of $\Omega^*$ and theoretical properties}

\label{sec:omega_estimators}

We present two specific estimators of $\Omega^*$ such that Condition \ref{con:omage_pred} is satisfied.
The first method is directly CLIME \citep{clime} with cross-fitting, as shown in Algorithm \ref{alg:omega_est}.

\begin{algorithm}[t]
    \caption{CLIME with cross-fitting}
      \vspace{-0.8em}
    \begin{flushleft}
    \hspace*{\algorithmicindent} \textbf{Input: }$X_J\in \mathbb{R}^{|J|\times p}$, $X_{J^\c}\in \mathbb{R}^{|J^\c|\times p}$ and $\rho_n >0$. \\
 \hspace*{\algorithmicindent} \textbf{Output: }$\hat{\Omega}^J\in\mathbb{R}^{p\times p}$ and $\hat{\Omega}^{J^\c}\in \mathbb{R}^{p\times p}$.
  \end{flushleft}
    \vspace{-1.2em}
    \begin{algorithmic}[1]
     % Input:
     % Output:
    \State  Compute $\hat{\Sigma}_{J}=\frac{1}{|J|}\sum_{i\in J}X_iX_i^\T$ and $\hat{\Sigma}_{J^\c}=\frac{1}{|J^\c|}\sum_{i\in J^\c}X_iX_i^\T$.

    \State  Obtain $\hat{G}^J\in \mathbb{R}^{p\times p}$ by solving
        \begin{align*}
            \hat{G}^{J} = &\argmin_{\text{$G\in\mathbb{R}^{p\times p}$}}  \sum_{j=1}^p \Vert G_{\cdot j}\Vert_{1}\\
            &\text{s.t. } \Vert \hat{\Sigma}_J G - I_p\Vert_{\max} \leq \rho_n.
        \end{align*}

   \State Compute $\hat{\Omega}^J\in\mathbb{R}^{p\times p}$ such that $\hat{\Omega}^J_{ij}=\hat{G}^J_{ij}\textbf{1}_{[|\hat{G}^J_{ij}| \leq |\hat{G}^J_{ji}|]} + \hat{G}^J_{ji}\textbf{1}_{[|\hat{G}^J_{ij}| > |\hat{G}^J_{ji}|]}$.

    \State Repeat steps 2-3 with $\hat{\Sigma}_{J}$ replaced by $\hat{\Sigma}_{J^\c}$ to compute $\hat{\Omega}^{J^\c}$.
    \end{algorithmic}
    \label{alg:omega_est}
\end{algorithm}

The theoretical properties of the estimators in Algorithm \ref{alg:omega_est} are shown in Lemma \ref{lem:error_rate_omega_clime}, based on \cite{clime}. The following assumption is required on $\Omega^*$.
Note that Assumption \ref{sigma_eigenvalue}(ii) with $C_4=L$ can be derived from Assumption \ref{ass:sparse_omega} due to the symmetry of $\Omega^*$ and the fact that $ \Vert A \Vert_2 \le \Vert A \Vert_1$ for a symmetric matrix $A$. See Theorem 5.6.9 in \cite{horn2012matrix}.

\begin{assumptionA}
    \label{ass:sparse_omega}
     $\Omega^*$ satisfies
    $\Vert \Omega^* \Vert_1 = \max_{1\leq j\leq p}\Vert \Omega^*_{\cdot j}\Vert_1\leq L$.
\end{assumptionA}

\begin{restatable}{lem}{errorrateomegaclime}
\label{lem:error_rate_omega_clime}
Suppose that Assumptions \ref{ass:sub_gaussian_entries}(i) and \ref{ass:sparse_omega} are satisfied
and $\frac{\log(p\vee n)}{\sqrt{n}}\leq\frac{c(1)\rho }{3}\sqrt{n}$. Then
for Algorithm \ref{alg:omega_est} with $\rho_n=\bar{A}_0\bar{\lambda}_n$ and any $\bar{A}_0\geq 1$,
we have with probability at least $1-\frac{4}{ p \vee n }$,
\begin{align}
   \Vert \hat{\Omega}^{J}_{\cdot j}-\Omega^*_{\cdot j}\Vert_2^2  &\leq  C s_{1j} \frac{\log(p\vee n)}{n}\quad\text{for $j=1,\ldots, p$ ,} \notag \\
    \Vert \hat{\Omega}^{J}_{\cdot j}-\Omega^*_{\cdot j}\Vert_1 &\leq   C s_{1j}\sqrt{\frac{\log(p\vee n)}{n}}\quad\text{for $j=1,\ldots, p$ ,}   \notag
\end{align}
where $\bar{\lambda}_n=\sqrt{\frac{3}{c(1)\rho}}K(1)(1+(d(1) \log(2))^{-1})   \sigma_1^2\sqrt{\frac{\log(p\vee n)}{n}}$, $c(1)$ is constants from Lemma \ref{lem:alpha_sub_exponential_concentration}, $K(1)$ and $d(1)$ are constants from Lemma \ref{lem:rv_centering} and $C$ is a constant depending on $(\bar{A}_0, \sigma_1, \rho, L)$ only. Furthermore, due to the symmetry of $\hat{\Omega}^J$, we have
\begin{align}
    \Vert \hat{\Omega}^J_{S_\xi,S_\xi} -\Omega^*_{S_\xi,S_\xi} \Vert_2 &\leq  C \min\{\sqrt{\Vert s_{1\xi}\Vert_1}, \Vert s_{1\xi}\Vert_\infty  \} \sqrt{\frac{\log(p\vee n)}{n}} , \label{inq:omega_operatpr_norm}\\
     \Vert \hat{\Omega}^J_{S_\xi,S_\xi} -\Omega^*_{S_\xi,S_\xi} \Vert_1&\leq C \Vert s_{1\xi}\Vert_\infty  \sqrt{\frac{\log(p\vee n)}{n}}. \notag
\end{align}
A similar result holds for $\hat{\Omega}^{J^c}$ in Algorithm \ref{alg:omega_est}.
\end{restatable}

\begin{rem}
\label{rem:omega_hat_sq_two_scenarios}

The two entries of $\min\{\cdot, \cdot\}$ on the right-hand side of (\ref{inq:omega_operatpr_norm}) are obtained from two different arguments.
First, the term $\Vert \hat{\Omega}^J_{S_\xi,S_\xi} -\Omega^*_{S_\xi,S_\xi} \Vert_2$ is upper bounded by $\Vert \hat{\Omega}^J_{S_\xi,S_\xi} -\Omega^*_{S_\xi,S_\xi}\Vert_1$ due to the symmetry of $\hat{\Omega}^J_{S_\xi,S_\xi}$ and $\Omega^*_{S_\xi,S_\xi}$. See Theorem 5.6.9 in \cite{horn2012matrix}. Then $\Vert \hat{\Omega}^J_{S_\xi,S_\xi} -\Omega^*_{S_\xi,S_\xi}\Vert_1$ is upper bounded by $\max_{j\in S_\xi}\Vert \hat{\Omega}^{J}_{\cdot j}-\Omega^*_{\cdot j} \Vert_1$.
Second, $\Vert \hat{\Omega}^J_{S_\xi,S_\xi} -\Omega^*_{S_\xi,S_\xi} \Vert_2$ is also upper bounded by $\Vert \hat{\Omega}^J_{S_\xi,S_\xi} -\Omega^*_{S_\xi,S_\xi}\Vert_\F$ which is upper bounded by $\sqrt{\sum_{j\in S_\xi}\Vert \hat{\Omega}^{J}_{\cdot j}-\Omega^*_{\cdot j}\Vert_2^2}$.

\end{rem}

Combining Theorem \ref{thm_1} and Lemma \ref{lem:error_rate_omega_clime} leads to
the following theorem for the theoretical properties of $\widehat{\xi^\T \beta}$  in (\ref{debiased_est_1}) and
$\hat{V}$ in (\ref{eq:definition_V_hat}), using $\hat{\Omega}^J$ and $\hat{\Omega}^{J^\c}$ from Algorithm \ref{alg:omega_est}.
Proposition \ref{prop:lm} can be deduced as a variation of Theorem \ref{thm_clime}, when applied to the setup with two independent samples. From the discussion of the first scenario in Section \ref{sec:high-level-LM},
results (i) can also be established if the CLIME estimator of $\Omega^*$ is used without cross-fitting.
For convenience, both results (i) and (ii) in Theorem \ref{thm_clime} are presented for the cross-fitted CLIME estimator.

\begin{thm}
\label{thm_clime}
\noindent Suppose that  Assumption \ref{sigma_eigenvalue},  \ref{sub_gaussian_residual}, \ref{ass:sparse_omega}, and Condition \ref{beta_hat_condition} are satisfied,
and the estimators $\hat{\Omega}^J$ and $\hat{\Omega}^{J^\c}$ in (\ref{debiased_est_1}) and (\ref{eq:definition_V_hat}) are taken from Algorithm \ref{alg:omega_est}
with $\rho_n=A_0\lambda_n$ and a sufficiently large constant $A_0$.
Then the following results hold for $\widehat{\xi^\T \beta}$ in (\ref{debiased_est_1}) and $\hat{V}$ in (\ref{eq:definition_V_hat}).\\
\noindent (i)
If further Assumption \ref{ass:sub_gaussian_entries} is satisfied,  then (\ref{result_LM}) holds provided $\max\{s_0, \Vert s_{1\xi}\Vert_\infty \frac{\Vert \xi\Vert_1}{\Vert \xi\Vert_2}, \min\{\sqrt{\Vert s_{1\xi}\Vert_1},\Vert s_{1\xi}\Vert_\infty\} \sqrt{ s_0}\}\tau_n=o(1)$.\\
\noindent (ii) If further Assumption \ref{ass:sub_gaussian_LM} is satisfied, then (\ref{result_LM}) holds provided $\max\{s_0, \min\{\sqrt{\Vert s_{1\xi}\Vert_1},\Vert s_{1\xi}\Vert_\infty \} \sqrt{ s_0} \}\tau_n=o(1)$. \\
Here $\tau_n=\frac{\log(p\vee n)}{\sqrt{n}}$ and $\lambda_n=\sqrt{\frac{\log(p\vee n)}{n}}$.
\end{thm}

We discuss the implications of Theorem \ref{thm_clime} and comparison with existing sparsity-based results in the following two remarks.

\begin{rem}[On sparsity conditions]
\label{rem:on_sparisity_condition}
As mentioned in the Introduction, existing results from \cite{zhang2014confidence} and \cite{van_de_Geer_2014}
deal with only inference about a single coefficient $\bar\beta_j$ (i.e., $\xi$ is set to $e_j$), while allowing model misspecification.
Nevertheless,
application of Theorem \ref{thm_clime} (i) to $\xi=e_j$ (hence $\Vert s_{1\xi}\Vert_\infty=\Vert s_{1\xi}\Vert_1=s_{1j}$)
shows that, similarly to the existing results,
valid inference is obtained about $\bar\beta_j$ under the same sparsity condition, namely $(s_0 \vee s_{1j} )\frac{\log(p\vee n)}{\sqrt{n}} =o(1)$.
Application of Theorem \ref{thm_clime} (ii) to $\xi=e_j$
shows valid inference about $\bar\beta_j$ under a weaker sparsity condition, $(s_0 \vee \sqrt{s_{1j}s_0} )\frac{\log(p\vee n)}{\sqrt{n}} =o(1)$,
provided the additional assumption that $X_0$ is a sub-gaussian vector.
By combining such entry-wise results for $j=1,\ldots,p$, valid inference about
\textit{each} $\bar\beta_j$ can be obtained under $(s_0 \vee s_{1} )\frac{\log(p\vee n)}{\sqrt{n}} =o(1)$ from existing results and Theorem \ref{thm_clime} (i),
and under $(s_0 \vee \sqrt{s_{1}s_0} )\frac{\log(p\vee n)}{\sqrt{n}} =o(1)$ from Theorem \ref{thm_clime} (ii) with a sub-gaussian covariate vector $X_0$,
where $s_1 = \max_{j=1,\ldots,p} \{s_{1j}\}$.
However, the significance of Theorem \ref{thm_clime} (ii) is that it also indicates valid inference about $\xi^\T\bar\beta$ for a general, non-sparse loading $\xi$
under the sparsity condition $(s_0 \vee s_1\sqrt{s_0})\frac{\log(p\vee n)}{\sqrt{n}} =o(1)$ with a sub-gaussian vector $X_0$.
\end{rem}

\begin{rem}[On CLIME vs nodewise regression]
\label{rem:on_clime_nodewise}
    We discuss why CLIME is chosen for estimating $\Omega^*$ in our method.
    In \cite{zhang2014confidence} and \cite{van_de_Geer_2014}, for inference about single coefficient $\bar\beta_j$,
    Lasso nodewise regression (i.e., Lasso linear regression of $j$th covariate on remaining covariates) is used to estimate $j$th column of $\Omega^*$.
    Let $\hat{\Omega}_{\mytext{node}} \in \mathbb{R}^{p\times p}$
    be the estimator of $\Omega^*$ obtained column by column through Lasso nodewise regression for $j=1,\ldots,p$.
    Similar column-wise error bounds as in Lemma \ref{lem:error_rate_omega_clime} can be shown for $\hat{\Omega}_{\mytext{node}}$ (without Assumption \ref{ass:sparse_omega}). However, $\hat{\Omega}_{\mytext{node}}$ is in general not symmetric and, to our knowledge, may not attain a similar \textit{matrix} $\Vert\cdot\Vert_2$-norm error bound to \eqref{inq:omega_operatpr_norm} for CLIME in Lemma \ref{lem:error_rate_omega_clime}. A direct symmetrization of $\hat{\Omega}_{\mytext{node}}$, for example, as in CLIME does not seem to overcome the issue. As a consequence,
    if $\hat{\Omega}_{\mytext{node}}$ is used to define the debiased estimator of $\xi^\T\bar\beta$ and estimator of $V$, then we would only obtain weaker results than Theorem \ref{thm_clime}, with $\min\{\sqrt{\Vert s_{1\xi}\Vert_1},\Vert s_{1\xi}\Vert_\infty\}$ replaced by $\sqrt{\Vert s_{1\xi}\Vert_1}$,
    because only the second argument in Remark~\ref{rem:omega_hat_sq_two_scenarios} would be applicable.
    In particular, without a desired result similar to Theorem \ref{thm_clime} (ii),
    when $\frac{\Vert\xi\Vert_1}{\Vert\xi\Vert_2}$ is large (namely $\xi$ is dense), the debiased method using $\hat{\Omega}_{\mytext{node}}$ may not enable valid inference about $\xi^\T \bar\beta$ even with  Assumption \ref{ass:sub_gaussian_LM} that $X_0$ is a sub-gaussian vector.
\end{rem}

We point out that any estimator of $\Omega^*$ satisfying Condition \ref{precision_matrix_LM} can be applied in our method to obtain Theorem \ref{thm_1}.
In the following, we briefly introduce a new two-stage algorithm to estimate $\Omega^*$, as shown in Algorithm \ref{alg:omega_est_two_stage} with cross-fitting.
The first-stage estimator $\hat{G}^{J}$ can be verified to be a variation of Lasso nodewise estimator  $\hat{\Omega}_{\mytext{node}}$ in Remark \ref{rem:on_clime_nodewise} but applied to $X_J$ (see Supplement Section \ref{sec:first_stage_lasso_nodel}).
Hence the two-stage algorithm serves as a \textit{symmetric} correction of asymmetric variant Lasso nodewise estimator  $\hat{G}^J$
to achieve comparable both column-wise and matrix-norm error bounds as in Lemma \ref{lem:error_rate_omega_clime} for CLIME and, when used in the debiased estimator for $\xi^\T\bar\beta$, to obtain
valid inference under comparable sparsity conditions as in Theorem \ref{thm_clime}.

\begin{algorithm}[t]
    \caption{Two-stage estimator with cross-fitting}
    \vspace{-0.8em}
     \begin{flushleft}
    \hspace*{\algorithmicindent} \textbf{Input: }$X_J\in \mathbb{R}^{|J|\times p}$, $X_{J^\c}\in \mathbb{R}^{|J^\c|\times p}$ and $\rho_n >0$. \\
 \hspace*{\algorithmicindent} \textbf{Output: }$\hat{\Omega}^J\in\mathbb{R}^{p\times p}$ and $\hat{\Omega}^{J^\c}\in \mathbb{R}^{p\times p}$.
 \end{flushleft}
 \vspace{-1.2em}
    \begin{algorithmic}[1]
     % Input:
     % Output:
    \State Compute $\hat{\Sigma}_{J}=\frac{1}{|J|}\sum_{i\in J}X_iX_i^\T$ and $\hat{\Sigma}_{J^\c}=\frac{1}{|J^\c|}\sum_{i\in J^\c}X_iX_i^\T$.

    \State (\textbf{Stage 1}) Obtain $\hat{G}^J\in \mathbb{R}^{p\times p}$ by solving
         \begin{align*}
            \hat{G}^{J} = \argmin_{G\in\mathbb{R}^{p\times p}} \frac{1}{2} G^\T\hat{\Sigma}_J G-\Tr(G) +\rho_n \Vert G\Vert_1.
        \end{align*}

    \State (\textbf{Stage 2}) Obtain $\hat{\Omega}^J\in \mathbb{R}^{p\times p}$ by solving
        \begin{align*}
            \hat{\Omega}^{J} = &\argmin_{\text{$G\in\mathbb{R}^{p\times p}$ and $G$ is symmetric}}  \Vert \hat{G}^J-G\Vert_{1}\\
            &\text{s.t. } \Vert \hat{\Sigma}_J G - I_p\Vert_{\max} \leq \rho_n.
        \end{align*}

    \State Repeat steps 2-3 with $\hat{\Sigma}_{J}$ replaced by $\hat{\Sigma}_{J^\c}$ to compute $\hat{\Omega}^{J^\c}$.
    \end{algorithmic}
    \label{alg:omega_est_two_stage}
\end{algorithm} \vspace{.1in}

The theoretical properties of the estimators in Algorithm \ref{alg:omega_est_two_stage} are shown in Lemma \ref{lem:error_rate_omega_two_stage} below. Assumption \ref{ass:entries_omega_x_0} is required on $X_0$ and $\Omega^*$.

\begin{assumptionA}
\label{ass:entries_omega_x_0}
 The entries of $\Omega^* X_0$ are sub-gaussian random variables with sub-gaussian norm bounded by $\sigma_1$.
\end{assumptionA}

\begin{restatable}{lem}{errorrateomegatwostage}
\label{lem:error_rate_omega_two_stage}
Suppose that Assumptions \ref{ass:sub_gaussian_entries}(i), \ref{sigma_eigenvalue}(ii) and \ref{ass:entries_omega_x_0} are satisfied and
$\frac{\log(p\vee n)}{\sqrt{n}} \leq\frac{c(1)\rho}{3} \sqrt{n}$,
$\delta = C_4 s_1 \bar{\lambda}_n(1+\eta_0)^2 \leq 1$ for some constants
$\eta_0>1$.
Then for Algorithm \ref{alg:omega_est_two_stage} with $\rho_n=\bar{A}_0\bar{\lambda}_n$ and $\bar{A}_0\geq\frac{\eta_0+1}{\eta_0-1}$, we have with probability at least $1-\frac{4}{ p \vee n }$,
\begin{align*}
   \Vert \hat{\Omega}^{J}_{\cdot j}-\Omega^*_{\cdot j}\Vert_2^2 &\leq  C s_{1} \frac{\log(p\vee n)}{n} \quad \text{for $j=1,\ldots, p$,} \\
    \Vert \hat{\Omega}^{J}_{\cdot j}-\Omega_{\cdot j}\Vert_1 &\leq   C s_{1}\sqrt{\frac{\log(p\vee n)}{n}}\quad \text{for $j=1,\ldots, p$,}
\end{align*}
where $\bar\lambda_n =\sqrt{\frac{3}{c(1)\rho}}K(1)(1+(d(1) \log(2))^{-1})   \sigma_1^2\sqrt{\frac{\log(p\vee n)}{n}} $, $c(1)$ is a constant from Lemma \ref{lem:alpha_sub_exponential_concentration}, $K(1)$ and $d(1)$ are constants from Lemma \ref{lem:rv_centering} and $C$ is a constant depending on $(\bar{A}_0,C_4,  \sigma_1, \delta, \rho)$ only. Furthermore, due to the symmetry of $\hat{\Omega}^J$, we have
\begin{align*}
    \Vert \hat{\Omega}^J_{S_\xi,S_\xi} -\Omega^*_{S_\xi, S_\xi} \Vert_2 &\leq  C \min\{\sqrt{\|s_{1\xi}\|_0 s_1}  , s_1  \} \sqrt{\frac{\log(p\vee n)}{n}} ,\\
     \Vert \hat{\Omega}^J_{S_\xi,S_\xi} -\Omega^*_{S_\xi,S_\xi} \Vert_1&\leq C s_1 \sqrt{\frac{\log(p\vee n)}{n}} .
\end{align*}
A similar result holds for $\hat{\Omega}^{J^c}$ in Algorithm \ref{alg:omega_est_two_stage}.
\end{restatable}

Compared to Lemma \ref{lem:error_rate_omega_clime} about the CLIME algorithm, Lemma \ref{lem:error_rate_omega_two_stage} about
the two-stage algorithm does not require Assumption \ref{ass:sparse_omega} to control $\Vert \Omega^*\Vert_{1}$.
On the other hand, the column-wise $\|\cdot\|_2$-norm and $\|\cdot\|_1$-norm error bounds,
and consequently the matrix $\|\cdot\|_2$-norm and $\|\cdot\|_1$-norm error bounds,
in Lemma~\ref{lem:error_rate_omega_two_stage}
are slightly weaker than those in Lemma \ref{lem:error_rate_omega_clime}, because $s_{1j} \le s_1$ for $j=1,\ldots,p$
by the definition $s_1 = \max_{j=1,\ldots, p} \{s_{1j}\}$.

Combining Theorem \ref{thm_1} and Lemma \ref{lem:error_rate_omega_two_stage} leads to the following theorem for the theoretical properties of $\widehat{\xi^\T \beta}$ in (\ref{debiased_est_1}) and
$\hat{V}$ in (\ref{eq:definition_V_hat}), using $\hat{\Omega}^J$ and $\hat{\Omega}^{J^\c}$ from Algorithm \ref{alg:omega_est_two_stage}.
The sparsity conditions in Theorem~\ref{thm_two_stage} are similar to those in Theorem \ref{thm_clime},
with $\Vert s_{1\xi}\Vert_\infty$ and $\Vert s_{1\xi}\Vert_1$ replaced by $s_1$ and $|S_\xi|s_1$ respectively.

\begin{thm}
\label{thm_two_stage}
\noindent Suppose that Assumption \ref{sigma_eigenvalue},  \ref{sub_gaussian_residual}, \ref{ass:entries_omega_x_0} and Condition \ref{beta_hat_condition} are satisfied,
and the estimators $\hat{\Omega}^J$ and $\hat{\Omega}^{J^\c}$ in (\ref{debiased_est_1}) and (\ref{eq:definition_V_hat}) are taken from Algorithm \ref{alg:omega_est_two_stage} with $\rho_n=A_0\lambda_n$ and a sufficiently large constant $A_0$. Then the following results hold for $\widehat{\xi^\T \beta}$ in (\ref{debiased_est_1}) and $\hat{V}$ in (\ref{eq:definition_V_hat}).\\
\noindent (i) If further Assumption \ref{ass:sub_gaussian_entries} is satisfied, then (\ref{result_LM}) holds provided $\max\{s_0, s_1 \frac{\Vert \xi\Vert_1}{\Vert \xi\Vert_2}, \min\{ \sqrt{\|s_{1\xi}\|_0 s_{1}},s_1\} \sqrt{ s_0} \}\tau_n=o(1)$.\\
\noindent (ii) If further Assumption \ref{ass:sub_gaussian_LM} is satisfied, then (\ref{result_LM}) holds provided $\max\{s_0,\min\{ \sqrt{\|s_{1\xi}\|_0 s_{1}},s_1\} \sqrt{ s_0} \}\tau_n=o(1)$. \\
Here $\tau_n=\frac{\log(p\vee n)}{\sqrt{n}}$ and $\lambda_n=\sqrt{\frac{\log(p\vee n)}{n}}$.
\end{thm}

The first stage in Algorithm \ref{alg:omega_est_two_stage} can be computed using proximal gradient descent and the second stage in Algorithm \ref{alg:omega_est_two_stage} can be computed using linear programming. However, as the dimension $p$ increases, Algorithm \ref{alg:omega_est_two_stage} is computationally more costly than Algorithm \ref{alg:omega_est}. The high computational cost may inhibit the application of Algorithm \ref{alg:omega_est_two_stage} when dimension $p$ is large.

\section{Debiased prediction and inference in misspecified GLM settings}
\label{sec:high_dimen_glm}

\subsection{High-level ideas}

\label{sec:high_level_ideas_M}

We discuss the main ideas behind
our method and theory as presented in the Introduction in the GLM setting. Suppose that $\Vert \hat{\beta}-\bar\beta\Vert_1=O_p(s_0\sqrt{\frac{\log(p \vee n)}{n}})$, $(\hat{\beta}-\bar\beta)^\T\hat{\Sigma} (\hat{\beta}-\bar\beta)=O_p(s_0\frac{\log(p \vee n)}{n})$,
and suitable regularity conditions are satisfied as in Proposition \ref{pro:glm_setting}.
We focus on the justification of asymptotic expansion (\ref{result_GLM_1}). Consider the following decomposition of the difference between $\widehat{\xi^\T\beta}$ in (\ref{eq:debiased_glm}) and $\xi^\T\bar\beta$:
\begin{align}
   & \quad \widehat{\xi^\T\beta} - \xi^\T\bar\beta \notag \\
    &= -v^{*\T}\frac{1}{n}\sum_{i=1}^n X_i\phi_1(Y_i,X_i^\T\bar\beta) -(\hat{H} v^*-\xi)^\T(\hat{\beta}-\bar\beta) - (\tilde{v}-v^{*})^\T\frac{1}{n}\sum_{i=1}^n X_i\phi_1(Y_i,X_i^\T\bar\beta)\notag \\&-(\tilde{v}-v^{*})^\T\frac{1}{n}\sum_{i=1}^n X_i(\phi_1(Y_i,X_i^\T\hat{\beta})-\phi_1(Y_i,X_i^\T\bar\beta))
   + R , \label{debiased_est_diff_glm}
\end{align}
where $\hat{H}=\frac{1}{n}\sum_{i=1}^n \phi_2(Y_i,X_i^\T\bar\beta)X_iX_i^\T$ and $R=v^{*\T}\frac{1}{n}\sum_{i=1}^nX_i(\phi_1(Y_i,X_i^{\T}\bar\beta)-\phi_1(Y_i,X_i^{\T}\hat{\beta}))+v^{*\T}\hat{H}(\hat{\beta}-\bar\beta)$ can be calculated by the fundamental theorem of calculus: $R=-\frac{1}{n}\sum_{i=1}^n (v^{*\T} X_i)\int_0^1 (1-t)\phi_3(Y_i,X_i^\T\bar\beta+t X_i^\T(\hat{\beta}-\bar\beta))dt (X_i^\T(\hat{\beta}-\beta))^2$. The first term on the right-hand side of (\ref{debiased_est_diff_glm}) gives the leading term in (\ref{result_GLM_1}), by the definition $v^*= M^* \xi$. By standard arguments,
\begin{align*}
   -\frac{\sqrt{n}}{\sqrt{V}} v^{* \T}\frac{1}{n}\sum_{i=1}^n X_i\phi_1(Y_i,X_i^\T\bar\beta)  \rightarrow_d \mathcal{N}(0, 1),
\end{align*}
and $V\asymp \Vert \xi\Vert^2_2$. To obtain asymptotic expansion (\ref{result_GLM_1}), it suffices to control the remaining terms on the right-hand side of (\ref{debiased_est_diff_glm}) as $o_p(\frac{\Vert \xi\Vert_2}{\sqrt{n}})$. We distinguish two different scenarios, mainly depending on different sub-gaussian assumptions on $X_0$, Assumption \ref{ass:sub_gaussian_entries_GLM} or \ref{ass:sub_gaussian_GLM}.

\textbf{First scenario.}\;
Suppose that Assumption \ref{ass:sub_gaussian_entries_GLM} holds, i.e., (i) the entries of $X_0$ are sub-gaussian random variables and
(ii) $X_0^\T M^*\xi$ is a sub-gaussian random variable with sub-gaussian norm being $O(\Vert M^*\xi \Vert_2)$, which is $O(\Vert \xi\Vert_2)$ provided Assumption \ref{ass:H_eigen_value}(ii).
The independence of $\tilde{M}$ with $\{(X_i, Y_i)\}_{i=1,\ldots, n}$ is not required.

The term $(\hat{H} v^* -\xi)^\T(\bar\beta-\hat{\beta})$ can be controlled by the H\"{o}lder inequality: $|(\hat{H} v^* -\xi)^\T(\bar\beta-\hat{\beta})|\leq \Vert \hat{H} v^* -\xi\Vert_\infty \Vert \bar\beta-\hat{\beta}\Vert_1$. The $j$th entry of $\hat{H} u^* -\xi$ can be expressed as $\frac{1}{n}\sum\limits_{i=1}^n(\phi_2(Y_i, X_i^\T\bar\beta)X_{ij}(X_i^\T M^* \xi)-\mathbb{E}[\phi_2(Y_i, X_i^\T\bar\beta)X_{ij}(X_i^\T M^* \xi)])$.
Under Assumption \ref{ass:sub_gaussian_entries_GLM}, \ref{ass:H_eigen_value}(ii) and \ref{bounded_higher_differentiate}, it can be shown that
 $\Vert \hat{H} v^* -\xi\Vert_\infty=O_p(\sqrt{\frac{\log(p \vee n)}{n}}\Vert \xi\Vert_2)$ (Lemma \ref{lem:bound_of_random_quan} (v)).
Together with the assumption on $\Vert \hat{\beta}- \bar\beta\Vert_1$, we have
\begin{align}
    (\hat{H} v^* -\xi)^\T(\bar\beta-\hat{\beta})=O_p(s_0\frac{\log(p\vee n)}{n}\Vert\xi\Vert_2) . \label{eq:sigma_v_beta}
\end{align}

The term  $(\tilde{v}-v^*)^\T\frac{1}{n} \sum\limits_{i=1}^n X_i\phi_1(Y_i, X_i^\T\bar\beta)$ can also be controlled by the H\"{o}lder inequality:
\begin{align*}
 \left| (\tilde{v}-v^*)^\T\frac{1}{n} \sum\limits_{i=1}^n X_i\phi_1(Y_i,X_i^\T\bar\beta) \right| \leq \Vert \tilde{v}-v^*\Vert_1\Vert \frac{1}{n} \sum\limits_{i=1}^n X_i\phi_1(Y_i,X_i^\T\bar\beta)\Vert_\infty.
\end{align*}
By the definition of $\bar\beta$ in (\ref{eq:target_beta_GLM}), we have $\mathbb{E}[X_{0j}\phi_1(Y_0,X_0^\T\bar\beta)] = 0$ for $j=1,\ldots, p$.
Under Assumption \ref{ass:sub_gaussian_entries_GLM} (i) and \ref{ass:residulal_glm_lower}, it can be shown that $\Vert \frac{1}{n} \sum\limits_{i=1}^n X_i\phi_1(Y_i, X_i^\T\bar\beta)\Vert_\infty=O_p(\sqrt{\frac{\log(p\vee n)}{n}})$ (Lemma \ref{lem:bound_of_random_quan} (iv)). Hence
\begin{align}
    (\tilde{v}-v^*)^\T\frac{1}{n} \sum\limits_{i=1}^n X_i\phi_1(Y_i,X_i^\T \bar\beta) = O_p(\Vert \tilde{v}-v^*\Vert_1 \sqrt{\frac{\log(p\vee n)}{n}}) . \label{eq:v_Y_Xbeta_GLM}
\end{align}

For the term $(\tilde{v}-v^*)^\T\hat{H}(\bar\beta-\hat{\beta})$, it suffices to control $(\tilde{v}-v^*)^\T\hat{\Sigma}(\bar\beta-\hat{\beta})$
provided $\phi_2(Y_i,X_i^\T \bar\beta)$ is bounded for $i=1,\ldots,n$. Then by the Cauchy--Schwartz inequality,
\begin{align}
    (\tilde{v}-v^*)^\T\hat{H}(\hat{\beta}-\bar\beta)&= O(|(\tilde{v}-v^*)^\T\hat{\Sigma}(\hat{\beta}-\bar\beta)|) \notag \\
    &= O(\sqrt{(\tilde{v}-v^*)^\T\hat{\Sigma} (\tilde{v}-v^*)} \sqrt{(\hat{\beta}-\bar\beta)^\T\hat{\Sigma}(\hat{\beta}-\bar\beta)} ).\label{inq:cross_term_v_beta}
\end{align}

Under Assumption \ref{ass:sub_gaussian_entries_GLM} (i), it can be shown that $\Vert \hat{\Sigma}-\Sigma\Vert_{\max}=O_p(\sqrt{\frac{\log(p\vee n)}{n}})$ (Lemma \ref{lem:bound_of_random_quan} (i)).
Then $(\tilde{v}-v^*)^\T\hat{\Sigma} (\tilde{v}-v^*)=(\tilde{v}-v^*)^\T(\hat{\Sigma}-\Sigma) (\tilde{v}-v^*)+(\tilde{v}-v^*)^\T\Sigma (\tilde{v}-v^*) \leq  \Vert\tilde{v}-v^*\Vert_1^2 \Vert\hat{\Sigma}-\Sigma\Vert_{\max}+(\tilde{v}-v^*)^\T\Sigma (\tilde{v}-v^*)=O_p(\Vert\tilde{v}-v^*\Vert_1^2 \sqrt{\frac{\log(p\vee n)}{n}} + \Vert\tilde{v}-v^*\Vert_2^2 )$. Together with the assumption on $(\hat{\beta}-\bar\beta)^\T\hat{\Sigma} (\hat{\beta}-\bar\beta)$, we have
\begin{align}
(\tilde{v}-v^*)^\T\hat{H}(\hat{\beta}-\bar\beta)&=
    O(|(\tilde{v}-v^*)^\T\hat{\Sigma}(\hat{\beta}-\bar\beta)|)\notag \\
    &=O_p(\sqrt{\Vert\tilde{v}-v^*\Vert_1^2 \sqrt{\frac{\log(p\vee n)}{n}} + \Vert\tilde{v}-v^*\Vert_2^2} \sqrt{s_0\frac{\log(p\vee n) }{n}}) . \label{eq:cross_term_glm}
\end{align}

For the term $R$, provided that $\phi_3(a, b)$ is bounded for $(a,b)\in\mathbb{R}^2$, we have
\begin{align}
    R &=\frac{1}{n}\sum_{i=1}^n (v^{*\T} X_i)\int_0^1 (1-t)\phi_3(Y_i,X_i^\T\beta+t X_i^\T(\hat{\beta}-\bar\beta))dt (X_i^\T(\hat{\beta}-\beta))^2\notag\\
&=O(\frac{1}{n}\sum_{i=1}^n |(v^{*\T} X_i)(X_i^\T(\hat{\beta}-\bar\beta))^2|) \notag  \\
&=O(\max_{i=1,\ldots, n}\{|v^{*\T}X_i|\}\frac{1}{n}\sum_{i=1}^n (X_i^\T(\hat{\beta}-\beta))^2)  .  \label{inq:R_glm_intermediate}
\end{align}
By Assumption \ref{ass:sub_gaussian_entries_GLM} (ii) and \ref{ass:H_eigen_value}(ii), it can be shown that $\max_{i=1,\ldots, n}\{|v^{*\T}X_i|\}=O_p(\sqrt{\log(n)}\Vert\xi\Vert_2)$ (Lemma \ref{lem:bound_of_random_quan}(vi)). Together with (\ref{inq:R_glm_intermediate}) and the assumption on $(\hat{\beta}-\bar\beta)\hat{\Sigma}(\hat{\beta}-\bar\beta)$, we have
\begin{align}
    R = O_p(s_0\frac{\log(p\vee n)\sqrt{\log(n)}}{n}\Vert \xi\Vert_2). \label{inq:R_glm}
\end{align}

From (\ref{eq:sigma_v_beta}), (\ref{eq:v_Y_Xbeta_GLM}), (\ref{eq:cross_term_glm}) and (\ref{inq:R_glm}),
the non-leading terms on the right-hand side of (\ref{debiased_est_diff_glm}) can be controlled through $\Vert \tilde{v}-v^*\Vert_{1}$ and $\Vert \tilde{v}-v^*\Vert_{2}$,
which in turn can be controlled by
$\Vert \tilde{M}_{S_\xi,S_\xi}-M^*_{S_\xi,S_\xi}\Vert_1\Vert\xi\Vert_1$ and $\Vert \tilde{M}_{S_\xi,S_\xi}-M^*_{S_\xi,S_\xi}\Vert_2\Vert\xi\Vert_2$ respectively.
Therefore, if $\frac{\Vert\xi\Vert_1}{\Vert\xi\Vert_2}$ is bounded,
then  asymptotic expansion (\ref{result_GLM_1}) can be obtained from (\ref{debiased_est_diff_glm}) by controlling
$\Vert \tilde{M}_{S_\xi,S_\xi}-M^*_{S_\xi,S_\xi}\Vert_1$ and $\Vert \tilde{M}_{S_\xi,S_\xi}-M^*_{S_\xi,S_\xi}\Vert_2$.

\vspace{.05in}
\textbf{Second scenario.}\;
With possibly large $\frac{\Vert\xi\Vert_1}{\Vert \xi\Vert_2}$, we develop a different approach
to control the terms $(\tilde{v}-v^*)^\T\frac{1}{n} \sum\limits_{i=1}^n X_i\phi_1(Y_i,X_i^\T \bar\beta)$ and $(\tilde{v}-v^*)^\T\hat{H}(\bar\beta-\hat{\beta})$.
Suppose that
Assumption \ref{ass:sub_gaussian_GLM} holds (i.e., $X_0$ is a sub-gaussian vector) and $\tilde{M}$ is independent of $\{(X_i, Y_i)\}_{i=1,\ldots, n}$.

To control the term $(\tilde{v}-v^*)^\T\frac{1}{n} \sum\limits_{i=1}^n X_i\phi_1(Y_i,X_i^\T\bar\beta)$, we note that
conditionally on $\tilde{v}$,
$\mathbb{E}[(\tilde{v}-v^*)^\T X_{i}\phi_1(Y_i,X_i^\T \bar\beta)] = 0$ by the definition of $\bar\beta$ and the independence of $\tilde{M}$ with $(X_i,Y_i)$,
and $(\tilde{v}-v^*)^\T X_{i}$ is a sub-gaussian variable with sub-gaussian norm being $O(\Vert \tilde{v}-v^* \Vert_2)$ by Assumption \ref{ass:sub_gaussian_GLM}.
Then conditionally on $\tilde{v}$,
 $(\tilde{v}-v^*)^\T X_{i}\phi_1(Y_i,X_i^\T\bar\beta)$ is a centered sub-exponential random variable with sub-exponential norm being $O(\Vert \tilde{v}-v^*\Vert_2)$
 provided that $\phi_1(Y_i,X_i^\T \bar\beta)$ is sub-gaussian for $i=1,\ldots,n$.
By sub-exponential concentration properties, we have (Lemma \ref{lem:bound_of_random_quan} (x))
\begin{align}
    (\tilde{v}-v^*)^\T\frac{1}{n} X_i\phi_1(Y_i,X_i^\T \bar\beta) = O_p(\Vert \tilde{v}-v^*\Vert_2 \sqrt{\frac{\log(p\vee n)}{n}}). \label{inq:v_y_xi_beta_vector_ass}
\end{align}
Compared to (\ref{eq:v_Y_Xbeta_GLM}), the right-hand side of (\ref{inq:v_y_xi_beta_vector_ass}) depends on $\Vert \tilde{v}-v^*\Vert_2$, but not $\Vert \tilde{v}-v^*\Vert_1$.

To control the term $(\tilde{v}-v^*)^\T\hat{H}(\hat{\beta}-\bar\beta)$,
it suffices, from (\ref{inq:cross_term_v_beta}), to control $(\tilde{v}-v^*)^\T\hat{\Sigma} (\tilde{v}-v^*)$, which can be expressed as $\frac{1}{n}\sum_{i=1}^n (( \tilde{v}-v^*)^\T X_i)^2 $.
By Assumption \ref{ass:sub_gaussian_GLM} and the independence of $\tilde{M}$ with $X_i$, we have that $((\tilde{v}-v^*)^\T X_0)^2$ is a sub-exponential random variable with sub-exponential norm being $O(\Vert \tilde{v}-v^*\Vert^2_2)$, and hence $( \tilde{v}-v^*)^\T\hat{\Sigma} (\tilde{v}-v^*) =O_p(\Vert \tilde{v}-v^*\Vert_2^2)$ (Lemma \ref{lem:bound_of_random_quan} (xi)).
Together with (\ref{inq:cross_term_v_beta}) and the assumption on $(\hat{\beta}-\bar\beta)^\T\hat{\Sigma} (\hat{\beta}-\bar\beta)$, we have
\begin{align}
    (\tilde{v}-v^*)^\T\hat{H}(\hat{\beta}-\bar\beta)=O_p(\Vert\tilde{v}-v^*\Vert_2\sqrt{s_0\frac{\log(p\vee n) }{n}} ) . \label{eq:cross_v_sigma_beta}
\end{align}
Compared to (\ref{eq:cross_term_glm}), the right-hand side of (\ref{eq:cross_v_sigma_beta}) depends on $\Vert \tilde{v}-v^*\Vert_2$, but not $\Vert \tilde{v}-v^*\Vert_1$.

From (\ref{eq:sigma_v_beta}), (\ref{inq:R_glm}), (\ref{inq:v_y_xi_beta_vector_ass}) and (\ref{eq:cross_v_sigma_beta}),
the non-leading terms on the right-hand side of (\ref{debiased_est_diff_glm}) can be controlled through $\Vert \tilde{v}-v^*\Vert_{2}$,
independently of $\Vert \tilde{v}-v^*\Vert_{1}$. Hence asymptotic expansion (\ref{result_GLM_1}) can be obtained from (\ref{debiased_est_diff_glm}) by controlling
$\Vert \tilde{M}_{S_\xi,S_\xi}-M^*_{S_\xi,S_\xi}\Vert_2$ regardless of the magnitude of $\frac{\Vert \xi\Vert_1}{\Vert\xi\Vert_2}$.

\vspace{.05in}
\textbf{Summary.}\; To complete the justification of asymptotic expansion (\ref{result_GLM_1}), the key is to find a suitable estimator $\tilde{M}$ such that
the \textit{matrix norms} $\Vert \tilde{M}_{S_\xi,S_\xi}-M^*_{S_\xi,S_\xi}\Vert_1$ and $\Vert \tilde{M}_{S_\xi,S_\xi}-M^*_{S_\xi,S_\xi}\Vert_2$ are
$O_p( \sqrt{\frac{\log(p\vee n)}{n}})$, up to some sparsity pre-factors.
We take $\tilde{M}$ to be the extended CLIME estimator in Algorithm \ref{alg:est_M_prod} and show that
$\Vert \tilde{M}_{S_\xi,S_\xi}-M^*_{S_\xi,S_\xi}\Vert_1=O_p(\Vert s_{1\xi}\Vert_\infty\sqrt{s_0\frac{\log(p\vee n)}{n}})$ and $\Vert \tilde{M}_{S_\xi,S_\xi}-M^*_{S_\xi,S_\xi}\Vert_2=O_p(\min\{\sqrt{\Vert s_{1\xi}\Vert_1}, \Vert s_{1\xi}\Vert_\infty\}\sqrt{s_0\frac{\log(p\vee n)}{n}})$ (Lemma \ref{lem:error_rate_M}).
Then we obtain the sparsity conditions as stated in Proposition \ref{pro:glm_setting} (i) and (ii),
corresponding to respectively the first and second scenarios discussed above.

\subsection{Debiased prediction with cross-fitting in GLM settings}

\label{sec:sparse_M}

We give our debiased estimator of $\xi^\T\bar\beta$ with cross-fitting in the misspecified GLM setting. We consider the cross-fitting scheme
as described in Section \ref{sec:sparse_omega}.

Let $\hat{\beta}\in \mathbb{R}^p$ be
an initial regularized sparse estimator of $\bar\beta$, for example, a Lasso estimator computed from $\{(X_i,Y_i)\}_{i=1,\ldots,n}$.
Our debiased estimator of $\xi^\T\bar\beta$ in GLM setting is defined as
\begin{align}
    \widehat{\xi^{\T}\beta} = \xi^\T \hat{\beta} - \xi^\T\hat{M}^{J^\c} \frac{1}{n}\sum_{i\in J}X_i\phi_1(Y_i, X_i^\T \hat{\beta}) - \xi^\T\hat{M}^{J} \frac{1}{n}\sum_{i\in J^\c}X_i\phi_1(Y_i,X_i^\T \hat{\beta}) .  \label{debiased_est_sparse_M}
\end{align}
where $\hat{M}^J$ and $\hat{M}^{J^\c}$ are estimators of $M^*$ based on $(X_J, Y_J)$ and $(X_{J^\c},Y_{J^\c})$ respectively, which are specified later in Section \ref{sec:M_est_glm}.
By construction, $\hat{M}^{J^\c}$ is independent of $(X_J,Y_J)$ and $\hat{M}^J$ is independent of $(X_{J^\c}, Y_{J^\c})$,
which is important for theoretical analysis in the second scenario as discussed in Section~\ref{sec:high_level_ideas_M}.
To construct a Wald confidence interval for $\xi^\T\bar\beta$, our estimator of the asymptotic variance for $\widehat{\xi^{\T}\beta}$ is defined as
\begin{align}
    \hat{\V} = \frac{1}{n} \sum_{i\in J} (\xi^\T\hat{M}^{J^\c}X_i)^2\phi_1^2(Y_i,X_i^\T\hat{\beta}) + \frac{1}{n}\sum_{i\in J^\c}(\xi^\T\hat{M}^{J}X_i)^2\phi_1^2(Y_i,X_i^\T\hat{\beta}). \label{eq:V_estimator}
\end{align}

\subsection{Theoretical properties of debiased prediction in GLM settings}
\label{sec:theory_GLM}

We give the theoretical properties of our debiased estimator of $\xi^\T\bar\beta$ and the variance estimator $\hat V$ introduced in Section \ref{sec:sparse_M},
assuming certain convergence rates on $\hat\beta$, $\hat{M}^J$ and $\hat{M}^{J^\c}$.
Specifically, we assume that $\hat{\beta}$ satisfies Condition \ref{beta_hat_condition} as stated in Section \ref{sec:theorey_LM}.

We also assume that the estimators $\hat{M}^J$ and $\hat{M}^{J^\c}$ satisfy the following condition. Specification of such estimators
to satisfy the conditions is deferred to Section \ref{sec:M_est_glm}.

\begin{conditionC}
\label{con:M_precisoin}
For some positive constant $C_2$, and constants $A_{1}$ and $A_{2}$ depending only on $(s_0, s_1, s_{1\xi})$,
the following conditions hold for $\hat{M}^J$ and $\hat{M}^{J^\c}$.
\label{precision_matrix_GLM}
$\newline$
\noindent  (i) With probability at least $1-\delta_2(n)$, $\hat{M}^J$ satisfies,
\begin{align*}
   \Vert \hat{M}^{J}_{S_\xi,S_\xi}-M^*_{S_\xi,S_\xi}\Vert_2 &\leq C_2 A_{1} \sqrt{\frac{\log(p\vee n)}{n}} \text{,} \\
    \Vert \hat{M}^{J}_{S_\xi,S_\xi}-M^*_{S_\xi,S_\xi}\Vert_1 &\leq   C_2 A_{2}\sqrt{\frac{\log(p\vee n)}{n}} \text{,}
\end{align*}
where $\delta_2(n)$ satisfies $\delta_2(n)\rightarrow 0$ as $n\rightarrow \infty$.

\noindent (ii) With probability at least $1-\delta_3(n)$, $\hat{M}^{J^\c}$ satisfies,
\begin{align*}
   \Vert \hat{M}^{J^\c}_{S_\xi,S_\xi}-M^*_{S_\xi,S_\xi} \Vert_2 &\leq   C_2 A_{1} \sqrt{\frac{\log(p\vee n)}{n}}\text{,} \\
    \Vert \hat{M}^{J^\c}_{S_\xi,S_\xi}-M^*_{S_\xi,S_\xi}\Vert_1 &\leq   C_2 A_{2}\sqrt{\frac{\log(p\vee n)}{n}}\text{,}
\end{align*}
where $\delta_3(n)$ satisfies $\delta_3(n)\rightarrow 0$ as $n\rightarrow \infty$.

\end{conditionC}

We consider one of the following conditions on $X_0$, Assumption \ref{ass:sub_gaussian_entries_GLM} or \ref{ass:sub_gaussian_GLM}.
If Assumption \ref{ass:sub_gaussian_GLM} holds,
then the entries of $X_0$ are sub-gaussian with norm bounded by $\sigma_1$, and $\xi^\T M^* X_0$ is sub-gaussian with norm bounded by $\sigma_1\Vert M^*\xi\Vert_2$. In other words, Assumption \ref{ass:sub_gaussian_entries_GLM} is implied by (hence weaker than) Assumption \ref{ass:sub_gaussian_GLM}.

\begin{assumptionA}
\label{ass:sub_gaussian_entries_GLM}
$\newline$
\noindent (i) The entries of $X_0$ are sub-gaussian random variables with sub-gaussian norm bounded by $\sigma_1$. \\
\noindent (ii) $\xi^\T M^* X_0$ is a sub-gaussian random variable with sub-gaussian norm bounded by $\sigma_1\Vert M^* \xi\Vert_2$.

\end{assumptionA}

\begin{assumptionA}
\label{ass:sub_gaussian_GLM}
$X_0$ is a sub-gaussian random vector with sub-gaussian vector norm bounded by $\sigma_1$.
\end{assumptionA}

We also make the following assumptions on $\Omega^*$, $M^*$ and the derivatives of $\phi (Y_0, X_0^\T\bar\beta)$.

\begin{assumptionA}
\label{ass:H_eigen_value}
$\newline$
  \noindent (i) $\Omega^*$ satisfies $C_3\leq \lambda_{\min}(\Omega^*)$ and $M^*$ satisfies $C_3\leq \lambda_{\min}(M^*)$ and $C_3\leq \lambda_{\min}(M^*\Sigma M^* )$.
  $\newline$
  \noindent (ii) $M^*$ satisfies $\lambda_{\max}(M^*) \leq C_4$.
$\newline$
 \noindent Here $C_3$ and $C_4$ are some positive constants.
\end{assumptionA}
\begin{assumptionA}
\label{ass:residulal_glm_lower}
  $\phi_1(Y_0,X_0^\T \bar\beta)$ is a sub-gaussian random variable with sub-gaussian norm $\sigma_2$ and
   \begin{align*}
       \mathbb{E}[\phi^2_1(Y_0,X_0^\T \bar\beta)|X_0]\geq C_5,
   \end{align*}
where $C_5$ is some positive constant.
\end{assumptionA}

\begin{assumptionA}
\label{bounded_higher_differentiate}
   $\phi_2(Y_0,X_0^\T\bar\beta)\leq C_6$ and $|\phi_3(Y_0,X_0^\T \bar\beta)|\leq C_7$ where $C_6$ and $C_7$ are some positive constants.
\end{assumptionA}

In the LM setting where $\phi(a,b)  =\frac{1}{2}(a-b)^2$, Assumption \ref{ass:sub_gaussian_entries_GLM}, \ref{ass:sub_gaussian_GLM}, \ref{ass:H_eigen_value}, \ref{ass:residulal_glm_lower} reduce to Assumption \ref{ass:sub_gaussian_entries}, \ref{ass:sub_gaussian_LM}, \ref{sigma_eigenvalue}, \ref{sub_gaussian_residual} respectively, whereas Assumption \ref{bounded_higher_differentiate} automatically holds with $C_6 = 1$ and $C_7=0$.

To present the theoretical properties of our method,
we use the following expressions. For $\widehat{\xi^\T\beta}$ in (\ref{debiased_est_sparse_M}) and $\hat{V}$ in (\ref{eq:V_estimator}), let
\begin{subequations}
\begin{equation}
 \frac{\widehat{\xi^{\T}\beta} - \xi^\T\bar\beta}{\sqrt{V}} = -\frac{1}{n}\sum_{i=1}^n \frac{(\xi^\T M^* X_i)\phi_1(Y_i,X_i^\T\bar\beta)}{\sqrt{V}}+ \Delta_1  n^{-\frac{1}{2}} , \label{eq:transform_glm_1}
\end{equation}
\begin{equation}
\frac{\hat{V}}{V} = 1+\Delta_2,\label{eq:transform_glm_2}
\end{equation}
\label{eq:transform_glm}
\end{subequations}
where $V=\mathbb{E}[(\xi^\T M^* X_0)^2\phi_1^2(Y_0,X_0^\T\bar\beta)]$.
The theoretical properties of $\widehat{\xi^\T\beta}$ in (\ref{debiased_est_sparse_M}) and $\hat{V}$ in (\ref{eq:V_estimator}) are established by characterizing $\Delta_1$ and $\Delta_2$ in the theorem below.

\begin{thm}
\label{thm_2}

Suppose that Conditions \ref{beta_hat_condition} and \ref{precision_matrix_GLM} are satisfied. Then
the following results hold for $\Delta_1$ in (\ref{eq:transform_glm_1}) and $\Delta_2$ in (\ref{eq:transform_glm_2}).

\noindent (i) If further Assumption \ref{ass:sub_gaussian_entries_GLM}, \ref{ass:H_eigen_value}, \ref{ass:residulal_glm_lower} and  \ref{bounded_higher_differentiate} are satisfied, $A_2 \frac{\Vert\xi\Vert_1}{\Vert\xi\Vert_2}\tau_n\leq1$ and $s_0\tau_n \leq \min\{\frac{c(1)\rho\sqrt{n} }{3},\frac{c(\frac{2}{3})n^{\frac{1}{6}}}{4},1\}$, then we have with probability at least $1-\delta_1(n)-\delta_2(n)-\delta_3(n)-\frac{2}{n}-\frac{12}{p \vee n}- 2\exp\{-c(\frac{1}{2})n^{\frac{1}{3}})\}-4\exp\{-\frac{1}{2}\sqrt{n}\}$,
\begin{align*}
    |\Delta_1|& \leq   C (A_2 \frac{\Vert \xi\Vert_1}{\Vert \xi\Vert_2}+A_1 \sqrt{ s_0}+\sqrt{\log( n)}s_0 )\tau_n
     ,\\
    |\Delta_2| &\leq  C (n^{-\frac{1}{3}}  +  \sqrt{s_0} \tau_n
    +A_2\frac{\Vert \xi\Vert_1}{\Vert \xi\Vert_2}\tau_n).
\end{align*}

\noindent (ii) If further Assumption \ref{ass:sub_gaussian_GLM}, \ref{ass:H_eigen_value}, \ref{ass:residulal_glm_lower} and  \ref{bounded_higher_differentiate} are satisfied, $A_1 \tau_n\leq1$ and $s_0\tau_n \leq \min\{c(1)\rho \sqrt{n}, \frac{c(1)\sqrt{n}}{2}, \frac{c(\frac{2}{3})n^{\frac{1}{6}}}{3},1\}$, then we have with probability at least $1-\delta_1(n)-\delta_2(n)-\delta_3(n)-\frac{2}{n}-\frac{20}{p \vee n}- 2\exp\{-c(\frac{1}{2})n^{\frac{1}{3}})\}-4\exp\{-\frac{1}{2}\sqrt{n}\}$,
\begin{align*}
    |\Delta_1|& \leq   C (
 A_1 \sqrt{ s_0}  + \sqrt{\log( n)}s_0 )\tau_n,\\
    |\Delta_2| &\leq  C (n^{-\frac{1}{3}}  +  \sqrt{s_0} \tau_n + A_1 \tau_n).
\end{align*}
Here $\tau_n=\frac{\log(p\vee n)}{\sqrt{n}}$, $c(\frac{1}{2})$, $c(\frac{2}{3})$ and $c(1)$ are constants from Lemma \ref{lem:alpha_sub_exponential_concentration} and $C$ is a constant depending on $(C_1, C_2, C_3, C_4, C_5, C_6,C_7, \sigma_1,\sigma_2,\rho)$ only.
\end{thm}

\begin{rem} \label{rem:x_max_bound_benefit}
Compared to Theorem \ref{thm_1} in LM setting, there is an additional log factor, $\sqrt{\log( n)}$, in the last term of the bound on $\Delta_1$,
$\sqrt{\log( n)} s_0 \tau_n $, in Theorem \ref{thm_2}.
This can be attributed to the extra term $R$ in the discussion of Section \ref{sec:high_level_ideas_M}.
     If, instead of Assumption \ref{ass:sub_gaussian_entries}, we assume $\Vert M^*X_0\Vert_\infty \leq K_1$ for some $K_1>0$, then the last term of the bound on $\Delta_1$ can be shown to be $K_1 s_0 \tau_n\frac{\Vert \xi\Vert_1}{\Vert\xi\Vert_2}$, which becomes $K_1 s_0 \tau_n$ when $\xi=e_j$ and hence $\xi^\T\bar\beta=\beta_j$ for $j=1,\ldots, p$.
    Such a result is used in existing works for inference about single coefficients \citep{van_de_Geer_2014, Ning}.
    See Remark \ref{rem:on_sparsity_conditoin_M} for related discussion.
\end{rem}

\subsection{Estimators of $M^*$ and theoretical properties}

\label{sec:M_est_glm}

We present two specific estimators of $M^*$ such that Condition \ref{con:M_precisoin} is satisfied.
The first method is a direct extension of CLIME estimator with cross-fitting, as shown in Algorithm \ref{alg:est_M}.

\begin{algorithm}[t]
    \caption{Extended CLIME estimator with cross-fitting}
         \vspace{-0.8em}
     \begin{flushleft}
    \hspace*{\algorithmicindent} \textbf{Input: }$\hat{\beta}^\prime_J\in \mathbb{R}^p$ estimated by $(X_J, Y_J)$, $\hat{\beta}^\prime_{J^\c}\in \mathbb{R}^p$ estimated by $(X_{J^\c}, Y_{J^\c})$,
    $X_J\in \mathbb{R}^{|J|\times p}$, $Y_J\in \mathbb{R}^{|J|}$, $X_{J^\c}\in \mathbb{R}^{|J^\c|\times p}$, $Y_{J^\c} \in \mathbb{R}^{|J^\c|}$ and $\rho_n >0$. \\
 \hspace*{\algorithmicindent} \textbf{Output: }$\hat{M}^J$ and $\hat{M}^{J^\c}$.
 \end{flushleft}
 \vspace{-1.2em}
    \begin{algorithmic}[1]
     % Input:
     % Output:
     \State Compute $\hat{H}_{J}=\frac{1}{|J|}\sum\limits_{i\in J}\phi_2(Y_i,X_i^\T \hat{\beta}^\prime_J)X_iX_i^\T$ and $\hat{H}_{J^\c}=\frac{1}{|J^\c|}\sum\limits_{i\in J^\c}\phi_2(Y_i,X_i^\T \hat{\beta}^\prime_{J^\c})X_iX_i^\T$.
    \State Obtain $\hat{G}^J\in \mathbb{R}^{p\times p}$ by solving
        \begin{align*}
            \hat{G}^{J} = &\argmin_{\text{$G\in\mathbb{R}^{p\times p}$}}  \sum_{j=1}^p \Vert G_{\cdot j}\Vert_{1}\\
            &\text{s.t. } \Vert \hat{H}_J G - I_p\Vert_{\max} \leq \rho_n.
        \end{align*}

    \State Compute $\hat{M}^{J}\in\mathbb{R}^{p\times p}$ such that $(\hat{M}^{J})_{ij}=\hat{G}^J_{ij}\textbf{1}_{|\hat{G}^J_{ij}| \leq |\hat{G}^J_{ji}|]} + \hat{G}^J_{ji}\textbf{1}_{[|\hat{G}^J_{ij}| > |\hat{G}^J_{ji}|]}$.
    \State Repeat steps 2-3 with $\hat{H}_{J}$ replaced by $\hat{H}_{J^\c}$ to compute $\hat{M}^{J^\c}\in \mathbb{R}^{p\times p}$.
    \end{algorithmic}
    \label{alg:est_M}
\end{algorithm}

The theoretical properties of the estimators in Algorithm \ref{alg:est_M} are shown in Lemma \ref{lem:error_rate_M}.
The following assumption is required on $M^*$.
Note that Assumption \ref{ass:H_eigen_value}(ii) with $C_4=L$ can be derived from Assumption \ref{ass:sparse_M} due to symmetry of $M^*$ and the fact that $ \Vert A \Vert_2 \le \Vert A \Vert_1$ for a symmetric matrix $A$. See Theorem 5.6.9 in \cite{horn2012matrix}.

\begin{assumptionA}
    \label{ass:sparse_M}
     $M^*$ satisfies
    $\max_{1\leq j\leq q}\Vert M^*_{\cdot j}\Vert_1\leq L$.
\end{assumptionA}

For convenience, we state the following condition about $\hat{\beta}^\prime_J$ and $\hat{\beta}^\prime_{J^\c}$ in Algorithm \ref{alg:est_M},
which is the instantiation of Condition \ref{beta_hat_condition} to the subsample $J$ and $J^c$.

\begin{conditionB}
\label{con:beta_hat_J}
$\newline$
\noindent (i) With probability at least $1-\delta_4(n)$, the estimator $\hat{\beta}^\prime_J$ satisfies
\begin{align*}
    \Vert \hat{\beta}^\prime_J-\bar\beta\Vert_1 &\leq C_8 s_0\sqrt{\frac{\log(p\vee n)}{n}}  ,\\
   \frac{1}{|J|} \sum_{i\in J}^n  (X_i^\T \hat{\beta}^\prime_{J} -X_i^\T \bar\beta)^2 &\leq C_8  s_0\frac{\log(p \vee n)}{n},
\end{align*}
where $\delta_4(n)$ satisfies $\delta_4(n)\rightarrow 0$ as $n\rightarrow \infty$ and $C_8$ is a positive constant.

\noindent (ii) With probability at least $1-\delta_5(n)$, the estimator $\hat{\beta}^\prime_{J^\c}$ satisfies
\begin{align*}
    \Vert \hat{\beta}^\prime_{J^\c}-\bar\beta\Vert_1 &\leq C_8 s_0 \sqrt{\frac{\log(p \vee n)}{n}},\\
    \frac{1}{|J^\c|}\sum_{i\in J^\c}^n(X_i^\T \hat{\beta}^\prime_{J^\c} -X_i^\T \bar\beta)^2 &\leq C_8 s_0\frac{\log(p \vee n)}{n},
\end{align*}
where $\delta_5(n)$ satisfies $\delta_5(n)\rightarrow 0$ as $n\rightarrow \infty$ and $C_8$ is a positive constant.
\end{conditionB}

\begin{restatable}{lem}{lemerrorrateM}
\label{lem:error_rate_M}
Suppose that Assumption \ref{ass:sub_gaussian_entries_GLM}(i), \ref{bounded_higher_differentiate} and \ref{ass:sparse_M} and Condition \ref{con:beta_hat_J} are satisfied, and $\frac{\log(p\vee n)}{\sqrt{n}}\leq\min\{\frac{c(1)\rho \sqrt{n}}{3}, \frac{c(\frac{1}{2})}{3}\}$. Then for $\hat{M}^J$ in Algorithm \ref{alg:est_M} with $\rho_n=\sqrt{s_0}\bar{A}_0\bar{\lambda}_n$ and any $\bar{A}_0\geq 1$, we have with probability at least $1-\delta_4(n)-\frac{8}{ p \vee n}$,
\begin{align*}
   \Vert \hat{M}^{J}_{\cdot j}-M^*_{\cdot j}\Vert_2^2 &\leq  C s_{1j}s_0 \frac{\log(p\vee n)}{n}
\quad\text{for $j=1,\ldots, p$},  \\
    \Vert \hat{M}^{J}_{\cdot j}-M^*_{\cdot j}\Vert_1 &\leq   C s_{1j}\sqrt{s_0\frac{\log(p\vee n)}{n}}\quad\text{for $j=1,\ldots, p$},
\end{align*}
where $\bar{\lambda}_n=(\sqrt{\frac{3}{c(1)\rho}}K(1)(1+(d(1) \log(2))^{-1})  C_6\sigma_1^2 + C_7C_8\sqrt{(K(\frac{1}{2})(1+(d(\frac{1}{2}) \log(2))^{-2})+\frac{1}{d(\frac{1}{2})})\sigma_1^4})*\sqrt{\frac{\log(p\vee n)}{n}}$, $c(\frac{1}{2})$ and $c(1)$ are constants from Lemma \ref{lem:alpha_sub_exponential_concentration},  $K(\frac{1}{2}),K(1), d(\frac{1}{2})$ and $d(1)$ are constants from Lemma \ref{lem:rv_centering} and $C$ is a constant depending on $(\bar{A}_0,C_6, C_7, C_8, \rho, \sigma_1,L)$ only. Furthermore, due to the symmetry of $\hat{M}^J$, we have
\begin{align}
    \Vert \hat{M}^J_{S_\xi,S_\xi} -M^*_{S_\xi,S_\xi} \Vert_2 &\leq  C \min\{\sqrt{\Vert s_{1\xi}\Vert_1}, \Vert s_{1\xi}\Vert_\infty  \} \sqrt{s_0} \sqrt{\frac{\log(p\vee n)}{n}} ,  \label{inq:M_error_bounds}\\
     \Vert \hat{M}^J_{S_\xi,S_\xi} -M^*_{S_\xi,S_\xi} \Vert_1&\leq C \Vert s_{1\xi}\Vert_\infty  \sqrt{s_0}\sqrt{\frac{\log(p\vee n)}{n}}.   \notag
\end{align}
A similar result holds for $\hat{M}^{J^c}$ in Algorithm \ref{alg:est_M}.
\end{restatable}

\begin{rem}
\label{rem:est_error}
Compared to Lemma \ref{lem:error_rate_omega_clime} in the LM setting,
the column-wise and matrix $\Vert\cdot\Vert_2$-norm and $\Vert\cdot\Vert_1$-norm error bounds in Lemma \ref{lem:error_rate_M}
involve extra factors related to $s_0$, which can be attributed to the estimation error of $\hat\beta^\prime_J$ and $\hat\beta^\prime_{J^\c}$ in the weights $\phi_2(Y_i,X_i^\T\hat{\beta}^\prime_J)$ and $\phi_2(Y_i,X_i^\T\hat{\beta}^\prime_{J^\c})$
for $\hat H_J$ and $\hat H_{J^\c}$ in Algorithm \ref{alg:est_M}.
See Remark \ref{rem:CLIME_weighted_nodewise_regression_M} for further discussion.
\end{rem}

Combining Theorem \ref{thm_2} and Lemma \ref{lem:error_rate_M} leads to
the following theorem for the theoretical properties of $\widehat{\xi^\T \beta}$  in (\ref{debiased_est_sparse_M}) and
$\hat{V}$ in (\ref{eq:V_estimator}), using $\hat{M}^J$ and $\hat{M}^{J^\c}$ from Algorithm \ref{alg:est_M}.
Proposition \ref{pro:glm_setting} can be deduced as a variation of Theorem \ref{thm_clime_M}, when applied to the setup with two independent samples.

\begin{thm}
\label{thm_clime_M}
\noindent Suppose that Assumption \ref{ass:H_eigen_value}, \ref{ass:residulal_glm_lower}, \ref{bounded_higher_differentiate}, \ref{ass:sparse_M} and Condition \ref{beta_hat_condition}, \ref{con:beta_hat_J} are satisfied, and the estimators $\hat{M}^J$ and $\hat{M}^{J^\c}$ in (\ref{debiased_est_sparse_M}) and (\ref{eq:V_estimator}) are taken from Algorithm \ref{alg:est_M} with $\rho_n = \sqrt{s_0}A_0\lambda_n$ and a sufficiently large constant $A_0$. Then the following results hold for $\widehat{\xi^\T \beta}$ in (\ref{debiased_est_sparse_M}) and $\hat{V}$ in (\ref{eq:V_estimator}).\\
\noindent (i)
If further Assumption \ref{ass:sub_gaussian_entries_GLM} is satisfied, then (\ref{result_GLM}) holds provided $\max\{s_0\sqrt{\log( n)}, \Vert s_{1\xi}\Vert_\infty \sqrt{s_0}\frac{\Vert \xi\Vert_1}{\Vert \xi\Vert_2} ,\min\{ \sqrt{ \Vert s_{1\xi}\Vert_1 }  ,  \Vert s_{1\xi}\Vert_\infty \}s_0\}\tau_n=o(1)$. \\
\noindent
(ii)
If further Assumption \ref{ass:sub_gaussian_GLM} is satisfied, then (\ref{result_GLM}) holds provided $\max\{s_0 \sqrt{\log( n)},\min\{ \sqrt{ \Vert s_{1\xi}\Vert_1}  ,  \Vert s_{1\xi}\Vert_\infty \}s_0 \}\tau_n = o(1)$.\\
Here we denote  $\tau_n=\frac{\log(p \vee
 n)}{\sqrt{n}}$ and $\lambda_n=\sqrt{\frac{\log(p \vee n)}{n}}$.
\end{thm}

\begin{rem}[On sparsity conditions]
\label{rem:on_sparsity_conditoin_M}
As mentioned in the Introduction, existing results from \cite{van_de_Geer_2014} and \cite{Ning}
deal with only inference about a single coefficient $\bar\beta_j$ (i.e., $\xi$ is set to $e_j$), while allowing model misspecification.
Nevertheless, it is instructive to compare the sparsity conditions involved in those results and entry-wise application of our results.
By combining the entry-wise results for $j=1,\ldots,p$ as in Remark \ref{rem:on_sparisity_condition}, valid inference about
\textit{each} $\bar\beta_j$ can be obtained under $\max \{s_1,s_0\} \tau_n=o(1)$ from existing results,
under $\max \{s_0\sqrt{\log(n)},s_1\sqrt{s_0}, \sqrt{s_1}s_0\} \tau_n=o(1)$ from Theorem \ref{thm_clime_M} (i),
and under $\max\{s_0 \sqrt{\log(n)}, \sqrt{s_1}s_0 \}$ from Theorem \ref{thm_clime_M} (ii) with a sub-gaussian covariate vector $X_0$.
However, the significance of Theorem \ref{thm_clime_M} (ii) is that it also indicates valid inference about $\xi^\T\bar\beta$ for a general, non-sparse loading $\xi$
under the sparsity condition $\max\{s_0 \sqrt{\log(n)},s_1s_0 \}$ with a sub-gaussian vector $X_0$.
The extra log factor, $\sqrt{\log(n)}$, in the latter three conditions is discussed in Remark \ref{rem:x_max_bound_benefit}.
The differences in the dependency on $(s_0,s_1)$ are explained in the next remark.
\end{rem}

\begin{rem}[On extended CLIME vs weighted nodewise regression]
\label{rem:CLIME_weighted_nodewise_regression_M}
The entry-wise results for inference about $\bar\beta_j$ from \cite{van_de_Geer_2014} and \cite{Ning} as mentioned in Remark \ref{rem:on_sparsity_conditoin_M} are obtained under boundedness assumptions
similar to $\Vert X_0\Vert_\infty \leq K_1$ and $\Vert M^*X_0\Vert_\infty \leq K_1$, while estimating $M^*_{\cdot j}$ ($j$th column of $M^*$) by weighted Lasso nodewise regression
(i.e., Lasso linear regression of $j$th covariate on remaining covariates, weighted by $\phi_2 (Y_i,X_i^\T\hat\beta)$, $i=1,\ldots,n$).
It can be shown that the resulting estimator, denoted as $\hat{M}_{\mytext{node}, \cdot j}$, satisfies
\begin{align}
   \Vert \hat{M}_{\mytext{node}, \cdot j}-M^*_{\cdot j}\Vert_2^2 &\leq C (s_{1j} + s_0) \frac{\log(p\vee n)}{n} , \notag  \\ %\label{inq:M_hat_Lsq}
    \Vert \hat{M}_{\mytext{node}, \cdot j}-M_{\cdot j}\Vert_1 &\leq  C (s_{1j}+s_0)\sqrt{\frac{\log(p\vee n)}{n}}, \notag %\label{inq:M_hat_L1}
\end{align}
Compared with these error bounds, those in Lemma \ref{lem:error_rate_M} from extended CLIME involve pre-factors depending on multiplication (instead of summation) of $s_0$ and $s_{1j}$,
which as mentioned in Remark \ref{rem:est_error} can be attributed to the estimation error in the weights like $\phi_2 (Y_i,X_i^\T\hat\beta)$.
In other words, it seems that weighted Lasso nodewise regression leads to improved error bounds than weighted Dantzig-like regression in extended CLIME.
However, similarly as discussed in Remark \ref{rem:on_clime_nodewise} in LM setting,
the estimator $\hat M_{\mytext{node}}$, defined column by column through weighted Lasso node regression, is not symmetric and may not achieve
a similar matrix $\Vert\cdot\Vert_2$-norm error bound to (\ref{inq:M_error_bounds}) for extended CLIME in Lemma \ref{lem:error_rate_M}.
   Without such a matrix $\Vert\cdot\Vert_2$-norm error bound, the debiased estimator of $\xi^\T\bar\beta$ using $\hat M_{\mytext{node}}$
    may not satisfy Theorem \ref{thm_clime_M} (ii) for a general loading $\xi$.
\end{rem}

In the following, we briefly introduce an extension of the two-stage algorithm to estimate $M^*$ as shown in Algorithm \ref{alg:est_M_two_stage} below.
Similarly as in the LM setting, the first-stage estimator $\hat{G}^{J}$ can be verified to be a variation of weighted Lasso nodewise estimator  $\hat{M}_{\mytext{node}}$ in Remark \ref{rem:CLIME_weighted_nodewise_regression_M} but applied to $X_J$ (see Supplement Section \ref{sec:first_stage_lasso_nodel_glm}).
Hence the two-stage algorithm serves as a \textit{symmetric} correction of asymmetric variant Lasso nodewise estimator  $\hat{G}^J$
to achieve comparable both column-wise and matrix-norm error bounds as in Lemma \ref{lem:error_rate_M} for CLIME and, when used in the debiased estimator for $\xi^\T\bar\beta$, to obtain
valid inference under comparable sparsity conditions as in Theorem \ref{thm_clime_M}.

\begin{algorithm}[t]
    \caption{Extended Two-stage estimator}
    \vspace{-0.8em}
     \begin{flushleft}
    \hspace*{\algorithmicindent} \textbf{Input: }$\hat{\beta}_J\in \mathbb{R}^p$ estimated by $(X_J, Y_J)$, $\hat{\beta}_{J^\c}\in \mathbb{R}^p$ estimated by $(X_{J^\c}, Y_{J^\c})$,
    $X_J\in \mathbb{R}^{|J|\times p}$, $Y_J\in\mathbb{R}^{|J|}$, $X_{J^\c}\in \mathbb{R}^{|J^\c|\times p}$, $Y_{J^\c}\in \mathbb{R}^{|J^\c|}$ and $\rho_n >0$. \\
 \hspace*{\algorithmicindent} \textbf{Output: }$\hat{M}^J\in\mathbb{R}^{p\times p}$ and $\hat{M}^{J^\c}\in \mathbb{R}^{p\times p}$.
 \end{flushleft}
 \vspace{-1.2em}
    \begin{algorithmic}[1]
\State Compute $\hat{H}_{J}=\frac{1}{|J|}\sum_{i\in J}\phi_2(Y_i,X_i^\T \hat{\beta}^\prime_J)X_iX_i^\T$ and $\hat{H}^\prime_{J^\c}=\frac{1}{|J^\c|}\sum_{i\in J^\c}\phi_2(Y_i,X_i^\T \hat{\beta}^\prime_{J^\c})X_iX_i^\T$.

    \State \textbf{Stage 1: } Compute $\hat{G}^J\in \mathbb{R}^{p\times p}$ by solving
        \begin{align*}
            \hat{G}^{J} = \argmin_{G\in\mathbb{R}^{p\times p}} \frac{1}{2} G^\T\hat{H}_J G-\Tr(G) +\rho_n \Vert G\Vert_1.
        \end{align*}

    \State \textbf{Stage 2: }Compute $\hat{M}^J\in \mathbb{R}^{p\times p}$ by solving
        \begin{align*}
            \hat{M}^{J} = &\argmin_{\text{$G\in\mathbb{R}^{p\times p}$ and $G$ is symmetric }} \Vert \hat{G}^J-G\Vert_{1} \\
            &\text{s.t. } \Vert \hat{H}_J G - I_p\Vert_{\max} \leq  \rho_n.
        \end{align*}

    \State Repeat steps 2-3 with $\hat{H}_{J}$ replaced by $\hat{H}_{J^\c}$ to compute $\hat{M}^{J^\c}\in \mathbb{R}^{p\times p}$.
    \end{algorithmic}
    \label{alg:est_M_two_stage}
\end{algorithm}

The theoretical properties of the estimators in Algorithm \ref{alg:est_M_two_stage} are shown in the Lemma \ref{lem:error_rate_M_two_stage} below. Assumption \ref{ass:entries_M_x_0} is required on $X_0$ and $M^*$.

\begin{assumptionA}
\label{ass:entries_M_x_0}
 The entries of $M^* X_0$ are sub-gaussian random variables with sub-gaussian norm bounded by $\sigma_1$.
\end{assumptionA}

\begin{restatable}{lem}{errorrateMtwostage}
\label{lem:error_rate_M_two_stage} Suppose that Assumption \ref{ass:sub_gaussian_entries_GLM}(i), \ref{ass:H_eigen_value}(ii), \ref{bounded_higher_differentiate}, and \ref{ass:entries_M_x_0} and Condition \ref{con:beta_hat_J} are satisfied, and $\frac{\log(p\vee n)}{\sqrt{n}}\leq\min\{ \frac{c(1)\rho \sqrt{n}}{3}, \frac{c(\frac{1}{2})}{3}\}$, $\delta=C_4 s_1\sqrt{s_0} \bar{\lambda}_n(1+\eta_0)^2 \leq 1$. Then for Algorithm \ref{alg:est_M_two_stage} with $\rho_n=\sqrt{s_0}\bar{A}_0\bar{\lambda}_n$ and $\bar{A}_0\geq\frac{\eta_0+1}{\eta_0-1}$, we have with probability at least $1-\delta_4(n)-\frac{8}{ p \vee n }$,
\begin{align*}
   \Vert \hat{M}^{J}_{\cdot j}-M_{\cdot j}\Vert_2^2 &\leq C s_1 s_0 \frac{\log(p\vee n)}{n} \quad \text{for $j=1,\ldots, p$,} \\
    \Vert \hat{M}^{J}_{\cdot j}-M_{\cdot j}\Vert_1 &\leq  C s_1\sqrt{s_0\frac{\log(p\vee n)}{n}}\quad \text{for $j=1,\ldots, p$,}
\end{align*}
where  $\bar{\lambda}_n=C_7 C_8 \sqrt{(K(\frac{1}{2})(1+(d(\frac{1}{2}) \log(2))^{-2}) +\frac{1}{d(\frac{1}{2})}) \sigma_1^4}+\sqrt{\frac{3}{c(1)\rho}}K(1)(1+(d(1) \log(2))^{-1})C_6\sigma_1^2\sqrt{\frac{\log(p\vee n)}{n}}$, $c(\frac{1}{2})$ and $c(1)$ are constants from Lemma \ref{lem:alpha_sub_exponential_concentration}, $K(\frac{1}{2})$ and $d(\frac{1}{2})$ are constants from Lemma \ref{lem:rv_centering} and $C$ is a constant depending on $(C_4, C_6, C_7, C_8,\bar{A}_0,\delta,\sigma_1,\rho)$.
only. Furthermore, due to the symmetry of $\hat{M}^J$, we have
\begin{align*}
    \Vert \hat{M}^J_{S_\xi,S_\xi} -M^*_{S_\xi,S_\xi} \Vert_2 &\leq  C \min\{\sqrt{\Vert s_{1\xi}\Vert_0s_1}, s_1  \} \sqrt{s_0} \sqrt{\frac{\log(p\vee n)}{n}} ,\\
     \Vert \hat{M}^J_{S_\xi,S_\xi} -M^*_{S_\xi,S_\xi} \Vert_1&\leq C s_1  \sqrt{s_0}\sqrt{\frac{\log(p\vee n)}{n}}.
\end{align*}
A similar result holds for $\hat{M}^{J^c}$ in Algorithm \ref{alg:est_M_two_stage}.
\end{restatable}

Compared to Lemma \ref{lem:error_rate_M} about the extended CLIME algorithm, Lemma \ref{lem:error_rate_M_two_stage} about the extended two-stage algorithm does not require Assumption \ref{ass:sparse_M} to control $\Vert M^*\Vert_{1}$.
On the other hand, the column-wise $\|\cdot\|_2$-norm and $\|\cdot\|_1$-norm error bounds,
and consequently the matrix $\|\cdot\|_2$-norm and $\|\cdot\|_1$-norm error bounds, in Lemma \ref{lem:error_rate_M_two_stage} are slightly weaker than
those in Lemma \ref{lem:error_rate_M} because $s_{1j} \le s_1$ for $j=1,\ldots,p$
by the definition $s_1 = \max_{j=1,\ldots, p} \{s_{1j}\}$.

Combining Theorem \ref{thm_2} and Lemma \ref{lem:error_rate_M_two_stage} leads to the following theorem for the theoretical properties of $\widehat{\xi^\T \beta}$ in (\ref{debiased_est_sparse_M}) and
$\hat{V}$ in (\ref{eq:V_estimator}), using $\hat{M}^J$ and $\hat{M}^{J^\c}$ from Algorithm \ref{alg:est_M_two_stage}.

\begin{thm}
\label{thm_two_stage_M}
Suppose that Assumption \ref{ass:H_eigen_value}, \ref{ass:residulal_glm_lower}, \ref{bounded_higher_differentiate}, \ref{ass:entries_M_x_0} and Condition \ref{beta_hat_condition}, \ref{con:beta_hat_J} are satisfied, and the estimators $\hat{M}^J$ and $\hat{M}^{J^\c}$ in (\ref{debiased_est_sparse_M}) and (\ref{eq:V_estimator}) are taken from Algorithm \ref{alg:est_M_two_stage}  with $\rho_n=\sqrt{s_0}A_0\lambda_n$ and a sufficiently large constant $A_0$. Then the following results hold for $\widehat{\xi^\T \beta}$ in (\ref{debiased_est_sparse_M}) and $\hat{V}$ in (\ref{eq:V_estimator}).\\
\noindent (i) If further Assumption \ref{ass:sub_gaussian_entries_GLM} is satisfied, then (\ref{result_GLM}) holds provided $\max\{s_0\sqrt{\log( n)}, s_1\sqrt{s_0}\frac{\Vert \xi\Vert_1}{\Vert \xi\Vert_2}, \min\{ \sqrt{ \|s_{1\xi}\|_0 s_1}  ,  s_1 \}s_0 \}\tau_n=o(1)$.

\noindent (ii) If further Assumption  \ref{ass:sub_gaussian_GLM} is satisfied, then (\ref{result_GLM}) holds provided $\max\{s_0 \sqrt{\log( n)},\min\{ \sqrt{ |S_\xi |s_1}  ,  s_1 \}s_0 \}\tau_n = o(1)$.\\
Here we denote  $\tau_n=\frac{\log(p\vee n)}{\sqrt{n}}$ and $\lambda_n=\sqrt{\frac{\log(p\vee n)}{n}}$.
\end{thm}

Similarly to Algorithm \ref{alg:omega_est} in the LM setting,
as the dimension $p$ increases, Algorithm \ref{alg:est_M_two_stage} is computationally more costly than Algorithm \ref{alg:est_M}.
The high computational cost may inhibit the application of Algorithm \ref{alg:est_M_two_stage} when dimension $p$ is large.

\section{Proof outline of main results}

\label{sec:outline_of_proof}

We give proof outlines for our main results, Theorem \ref{thm_1} and \ref{thm_2}.
We give an initial result (Theorem \ref{thm:initial_thm}) for the debiased estimator of $\xi^\T\bar\beta$ in the GLM setting,
depending on certain error terms related to estimators of $M^*$, which are controlled in Lemmas \ref{lem:cross_term_control_1} later, subject to Condition \ref{con:M_precisoin} on estimation errors for $M^*$. Theorem \ref{thm_1} in the LM setting can be proved by
incorporating Condition \ref{con:omage_pred} and applying the initial theoretical result with $\phi(a,b) = \frac{1}{2}(a-b)^2$.
Theorem \ref{thm_2} in the GLM setting can be proved by incorporating Condition \ref{con:M_precisoin} and applying the initial theoretical result.
In the following discussion, $C$ is a constant that may differ from line to line, and w.h.p. indicates that an event holds with probability converging to $1$ as the sample size $n$ tends to infinity.

\label{sec:theoretical_results_oracle_M}

\subsection{Initial theoretical result in GLM setting}

We give an initial result (Theorem \ref{thm:initial_thm}) for the debiased estimator of $\xi^\T\bar\beta$ in the GLM setting.
We use the same notation as in Section \ref{sec:high_dimen_glm}. In particular, recall the definitions of $\widehat{\xi^\T\beta}$ and $\hat{V}$ in (\ref{debiased_est_sparse_M}) and (\ref{eq:V_estimator}) respectively
and the expressions in (\ref{eq:transform_glm}) where $V=\mathbb{E}[(\xi^\T M^* X_0)^2\phi_1^2(Y_0,X_0^\T\bar\beta)]$.

We give the theoretical properties of $\widehat{\xi^{\T}\beta}$ in (\ref{debiased_est_sparse_M}) and $\hat{V}$ in (\ref{eq:V_estimator}) by characterizing the bounds of the terms $\Delta_1$ in (\ref{eq:transform_glm_1}) and $\Delta_2$ in (\ref{eq:transform_glm_2}) as below.

\begin{restatable}{thm}{initialthm}
\label{thm:initial_thm}
Suppose that Condition \ref{beta_hat_condition} is satisfied. Then the following results hold for $\Delta_1$ in (\ref{eq:transform_glm_1}) and $\Delta_2$ in (\ref{eq:transform_glm_2}).

\noindent (i) If further Assumptions \ref{ass:sub_gaussian_entries_GLM}, \ref{ass:H_eigen_value}, \ref{ass:residulal_glm_lower} and \ref{bounded_higher_differentiate} are satisfied and $\frac{\log(p\vee n)}{\sqrt{n} }  \leq\min\{\frac{c(1)\rho\sqrt{n} }{2},\frac{c(\frac{2}{3})n^{\frac{1}{6}}}{4}\}$, then we have with probability at least $1-\delta_1(n)-\frac{4}{n}-\frac{8}{p \vee n}- 2 \exp\{-c(\frac{1}{2})n^{\frac{1}{3}})\}-4\exp\{-\frac{1}{2}\sqrt{n}\}$,
\begin{align*}
   |\Delta_1|& \leq   C ( s_0\frac{\log(p\vee n)}{\sqrt{n}}+\sqrt{\log(p\vee n)}\frac{L_1 + L_1^\prime}{\Vert\xi\Vert_2}+  \frac{\sqrt{Q_1}+\sqrt{Q_1^\prime}}{\Vert\xi\Vert_2} \sqrt{s_0\log(p\vee n)}  \\
   &+C_7  s_0 \frac{\log(p\vee n)\sqrt{\log(n)}}{\sqrt{n} } ),
\end{align*}
\begin{align*}
     |\Delta_2| &\leq  C (n^{-\frac{1}{3}}  +  \sqrt{s_0} \sqrt{\frac{\log(n) \log(p\vee n)}{n}}+  s_0 \frac{\log(p\vee n)\log(n)}{n}\\
    &+\frac{\sqrt{\log(n)}(\sqrt{Q_1}+\sqrt{Q_1^\prime})}{\Vert \xi\Vert_2}+ s^2_0 \frac{\log(p\vee n)\sqrt{\log(n)}}{n} \frac{L_1 + L_1^\prime}{\Vert \xi\Vert_2} +  \frac{(Q_1 + Q_1^\prime) \log(n)}{\Vert \xi\Vert_2^2} \\
    &+ s^2_0  \frac{\log(p\vee n) }{n} \frac{L^2_1 + L_1^{\prime 2} }{\Vert \xi\Vert_2^2}).
\end{align*}

\noindent (ii) If further assume Assumptions \ref{ass:sub_gaussian_GLM}, \ref{ass:H_eigen_value}, \ref{ass:residulal_glm_lower} and \ref{bounded_higher_differentiate} are satisfied and $\frac{\log(p\vee n)}{\sqrt{n} }\leq\min\{c(1)\rho \sqrt{n}, \frac{c(1)\sqrt{n}}{2}, \frac{c(\frac{2}{3})n^{\frac{1}{6}}\rho^{\frac{2}{3}}}{3}\}$, then we have with probability at least $1-\delta_1(n)-\frac{8}{n}-\frac{12}{p \vee n}- 2exp\{-c(\frac{1}{2})n^{\frac{1}{3}}\}-4exp\{-\frac{1}{2}\sqrt{n}\}$,
\begin{align*}
   |\Delta_1|& \leq   C (s_0\frac{\log(p\vee n)}{\sqrt{n}}+\sqrt{\log(p\vee n)}\frac{L_2 + L_2^\prime}{\Vert\xi\Vert_2} +  \frac{\sqrt{Q_1}+\sqrt{Q_1^\prime}}{\Vert\xi\Vert_2} \sqrt{s_0\log(p\vee n)}  \\
   &+C_7 s_0 \frac{\log(p\vee n)\sqrt{\log(n)}}{\sqrt{n}} ),
\end{align*}
\begin{align*}
    |\Delta_2| &\leq  C (n^{-\frac{1}{3}}  +  \sqrt{s_0} \sqrt{\frac{\log(n) \log(p\vee n)}{n}}+  s_0 \frac{\log(p\vee n)\log(n)}{n} \\
&+\frac{\sqrt{\log(n)}(\sqrt{Q_1}+\sqrt{Q_1^\prime})}{\Vert\xi\Vert_2}+ s^2_0 \frac{\log(p\vee n)\sqrt{\log(n)}}{n}\frac{L_2 + L_2^\prime}{\Vert\xi\Vert_2} +  \frac{(Q_1+Q_1^\prime) \log(n)}{\Vert \xi\Vert_2^2} \\
    &+ s^2_0
 \frac{\log(p\vee n) \sqrt{\log(n)} }{n}\frac{L^2_2 + L_2^{\prime 2} }{\Vert\xi\Vert_2^2} ).
\end{align*}
Here $c(\frac{1}{2})$, $c(\frac{2}{3})$ and $c(1)$ are constants from Lemma \ref{lem:alpha_sub_exponential_concentration}, $C$ is a constant depending only on $(C_1, C_3, C_4, C_5, C_6, \sigma_1,\sigma_2,\rho)$ and
$Q_1=\xi^\T (\hat{M}^{J^\c}- M^*)\sum_{i\in J}\frac{1}{n}X_iX_i^\T(\hat{M}^{J^\c}- M^*)\xi$,
$Q_1^\prime = \xi^\T (\hat{M}^{J}- M^*)\sum_{i\in J^\c}\frac{1}{n}X_iX_i^\T(\hat{M}^{J}- M^*)\xi$, $L_1 = \Vert (\hat{M}^{J}-M^*)\xi\Vert_1$, $L^\prime_1 = \Vert (\hat{M}^{J^\c}-M^*)\xi\Vert_1$, $L_2 = \Vert (\hat{M}^{J}-M^*)\xi\Vert_2$, $L^\prime_2 = \Vert (\hat{M}^{J^\c}-M^*)\xi\Vert_2$.
\end{restatable}

\subsubsection{Proof outline of Theorem \ref{thm:initial_thm}}

\label{sec:outline_thm_7}

We give a proof outline of Theorem \ref{thm:initial_thm}, first the bound of $\Delta_1$ and then that of $\Delta_2$.

\noindent \textbf{(a) Bound of $\Delta_1$.} By the definition (\ref{debiased_est_sparse_M}) and algebraic calculation, we have
\begin{align*}
    & \quad\widehat{\xi^\T\beta} - \xi^\T\bar\beta \notag \\
    &= -\xi^\T M^*\frac{1}{n}\sum_{i=1}^n X_i\phi_1(Y_i, X_i^\T\bar\beta)-\underbrace{(\frac{1}{n}\sum_{i=1}^n \phi_2(Y_i, X_i^\T\bar\beta) X_iX_i^\T M^*\xi-\xi)^\T(\bar\beta-\hat{\beta})}_{I} \notag \\
    &-\underbrace{\xi^\T (\hat{M}^{J^\c}- M^*)\frac{1}{n}\sum_{i\in J}X_i\phi_1(Y_i, X_i^\T\bar\beta)}_{II} -\underbrace{\xi^\T (\hat{M}^{J^\c}- M^*)\frac{1}{n}\sum_{i\in J}X_i(\phi_1(Y_i, X_i^\T\hat{\beta}) -  \phi_1(Y_i, X_i^\T\bar\beta) )}_{III}  \notag \\
   &-\underbrace{\xi^\T (\hat{M}^{J}- M^*)\frac{1}{n}\sum_{i\in J^\c }X_i\phi_1(Y_i, X_i^\T\bar\beta)}_{IV}  -\underbrace{\xi^\T (\hat{M}^{J}- M^*)\frac{1}{n}\sum_{i\in J^\c }X_i(\phi_1(Y_i, X_i^\T\hat{\beta}) -  \phi_1(Y_i, X_i^\T\bar\beta) )}_{V} \notag  \\&+ R,
\end{align*}
where
$R = -\sum_{i=1}^n \frac{1}{n} (\xi^\T M^* X_i)\int_0^1 (1-t)\phi_3(Y_i,X_i^\T\hat{\beta}+t X_i^\T(\hat{\beta}-\bar\beta))dt  (X_i^\T(\hat{\beta}-\bar\beta))^2$.

We discuss how to control the following terms, $I$, $II$, $III$ and $R$. The term $IV$ can be controlled similarly as the term $II$,
and the term $V$ can be controlled similarly as the term $III$.

\noindent \textcircled{1} For the term $I$, we have by H\"{o}lder's inequality,
\begin{align}
   |I| &= |(\frac{1}{n}\sum_{i=1}^n \phi_2(Y_i, X_i^\T\bar\beta) X_iX_i^\T M^*\xi-\xi)^\T(\bar\beta-\hat{\beta}) | \notag \\&\leq \Vert (\frac{1}{n}\sum_{i=1}^n \phi_2(Y_i, X_i^\T\bar\beta) X_iX_i^\T M^* - I_p)\xi \Vert_\infty \Vert \bar\beta-\hat{\beta}\Vert_1 . \label{inq:V_initial}
\end{align}
The $j$th entry of $ (\frac{1}{n}\sum_{i=1}^n \phi_2(Y_i, X_i^\T\bar\beta) X_iX_i^\T M^* - I_p)\xi $ can be expressed as $\frac{1}{n}\sum\limits_{i=1}^n(\phi_2(Y_i, X_i^\T\bar\beta)X_{ij}(X_i^\T M^* \xi)-\mathbb{E}[\phi_2(Y_i, X_i^\T\bar\beta)X_{ij}(X_i^\T M^* \xi)])$. Under Assumption \ref{ass:sub_gaussian_entries_GLM}, \ref{ass:H_eigen_value}(ii) and \ref{ass:residulal_glm_lower}, we have that from Lemma \ref{lem:product-two-sub-gaussian-rv-norm}, $ \phi_2(Y_i, X_i^\T\bar\beta)X_{ij}(X_i^\T M^* \xi)$ for $i=1,\ldots,n$, $j=1,\ldots,p$ are sub-exponential random variables with bounded sub-exponential norm being $O(\Vert \xi\Vert_2)$, and hence from Lemma S12, $(\phi_2(Y_i, X_i^\T\bar\beta)X_{ij}(X_i^\T M^* \xi)-\mathbb{E}[\phi_2(Y_i, X_i^\T\bar\beta)X_{ij}(X_i^\T M^* \xi)])$ for $i=1,\ldots,n$, $j=1,\ldots,p$ are centered sub-exponential random variables with sub-exponential norm being $O(\Vert\xi\Vert_2)$. By concentration properties of centered sub-exponential random variables we have w.h.p. that for $j=1,\ldots, p$,
\begin{align*}
  | \frac{1}{n} \sum\limits_{i=1}^n (\phi_2(Y_i, X_i^\T\bar\beta)X_{ij}(X_i^\T M^* \xi)-\mathbb{E}[\phi_2(Y_i, X_i^\T\bar\beta)X_{ij}(X_i^\T M^* \xi)])| \leq C \sqrt{\frac{\log(p\vee n)}{n}}\Vert \xi\Vert_2.
\end{align*}
Hence w.h.p.
\begin{align}
    \Vert (\frac{1}{n}\sum_{i=1}^n \phi_2(Y_i, X_i^\T\bar\beta) X_iX_i^\T M^* - I_p)\xi\Vert_\infty \leq C \sqrt{\frac{\log(p\vee n)}{n}} \Vert \xi\Vert_2. \quad\text{(Lemma \ref{lem:bound_of_random_quan} (v))}\label{inq:max_norm}
\end{align}
From (\ref{inq:V_initial}), (\ref{inq:max_norm}) and Condition \ref{beta_hat_condition}, we have w.h.p.
\begin{align}
    |I|&\leq C s_0\frac{\log(p\vee n)}{n}  \Vert \xi\Vert_2. \label{inq:V_term_outline}
\end{align}

\noindent \textcircled{2} For term $II$, we have
\begin{align*}
    |II|=|\xi^\T (\hat{M}^{J^\c}- M^*)\frac{1}{n}\sum_{i\in J}X_i\phi_1(Y_i, X_i^\T\bar\beta)|\leq |\frac{1}{|J|}\sum_{i\in J} \xi^\T (\hat{M}^{J^\c}- M^*) X_i \phi_1(Y_i,X_i^\T\bar\beta)|.
\end{align*}
To control the term $II$, it is suffice to control $|\frac{1}{|J|}\sum_{i\in J}^n\xi^\T (\hat{M}^{J^\c}- M^*) X_i \phi_1(Y_i,X_i^\T\bar\beta)|$.

First, in the scenario where Assumption \ref{ass:sub_gaussian_entries_GLM} is satisfied, we have by H\"{o}lder's inequality,
\begin{align*}
    |\frac{1}{|J|}\sum_{i\in J}\xi^\T (\hat{M}^{J^\c}- M^*) X_i \phi_1(Y_i,X_i^\T\bar\beta)|&\leq \Vert (\hat{M}^{J^\c}- M^*)\xi\Vert_1\Vert \frac{1}{|J|} \sum\limits_{i\in J} X_i\phi_1(Y_i,X_i^\T\bar\beta)\Vert_\infty.
\end{align*}

By the definition of $\bar\beta$ in (\ref{eq:target_beta_GLM}), we have $E [X_{0j}\phi_1(Y_0,X_0^\T\bar\beta)] = 0$ for $j=1,\ldots, p$.
Under Assumption \ref{ass:sub_gaussian_entries_GLM} (i) and \ref{ass:residulal_glm_lower}, we have that $X_{0j}\phi_1(Y_0,X_0^\T\bar\beta)$ for $j=1,\ldots, p$ are centered sub-exponential random variables with bounded expoential norm. By the concentration properties of centered sub-exponential random variables, we have w.h.p.
\begin{align*}
    \Vert \frac{1}{|J|} \sum\limits_{i\in J} X_i\phi_1(Y_i, X_i^\T\bar\beta)\Vert_\infty\leq C \sqrt{\frac{\log(p\vee n)}{n}} .\quad\text{(Lemma \ref{lem:bound_of_random_quan} (iv))}
\end{align*}
Hence we have w.h.p.
\begin{align}
    |II|\leq C \sqrt{\frac{\log(p\vee n)}{n}} \Vert  (\hat{M}^{J^\c}- M^*)\xi\Vert_1 .\label{inq:bound_of_term_I_outline_sce1}
\end{align}

Second, we discuss the scenario where Assumption \ref{ass:sub_gaussian_GLM} is satisfied. Conditionally on $\hat{M}^{J^\c}$, we have $\mathbb{E}[\xi^\T (\hat{M}^{J^\c}-M^*) X_{i}\phi_1(Y_i,X_i^\T\bar\beta)] = 0$ for $i\in J$ by the definition of $\bar\beta$. Then conditionally on $\hat{M}^{J^\c}$, under Assumption \ref{ass:sub_gaussian_GLM} and \ref{ass:residulal_glm_lower}, we have $\xi^\T(\hat{M}-M^*) X_{i}\phi_1(Y_i,X_i^\T\bar\beta)$ for $i\in J$ are centered sub-exponential random variables with sub-exponential norm being $O(\Vert (\hat{M}-M^*)\xi\Vert_2)$. By sub-exponential concentration properties, we have w.h.p.
\begin{align}
  |II| = | \xi^\T(\hat{M}^{J^\c}-M^*)\frac{1}{n}\sum_{i=1}^n X_i\phi_1(Y_i,X_i^\T \bar\beta)| \leq C \Vert (\hat{M}-M^*)\xi\Vert_2 \sqrt{\frac{\log(p\vee n)}{n}}. \quad\text{(Lemma \ref{lem:bound_of_random_quan} (x))}\label{inq:bound_of_term_I_outline_sce2}
\end{align}

\noindent \textcircled{3} For the term $III$, by the Cauchy--Schwartz inequality and under Assumption \ref{bounded_higher_differentiate} and Condition \ref{beta_hat_condition}, we have w.h.p.
\begin{align}
   |III|&=|\xi^\T (\hat{M}^{J^\c}- M^*)\frac{1}{n}\sum_{i\in J}X_i(\phi_1(Y_i, X_i^\T\hat{\beta}) -  \phi_1(Y_i, X_i^\T\bar\beta) )| \notag \\
   &\leq   \sqrt{\xi^\T (\hat{M}^{J^\c}- M^*)\sum_{i\in J}\frac{1}{n}X_iX_i^\T(\hat{M}^{J^\c}- M^*)\xi}
    \sqrt{\frac{1}{n}\sum_{i\in J}(\phi_1(Y_i, X_i^\T\hat{\beta}) -  \phi_1(Y_i, X_i^\T\bar\beta) )^2} \notag  \\
    &\leq \sqrt{\xi^\T (\hat{M}^{J^\c}- M^*)\sum_{i\in J}\frac{1}{n}X_iX_i^\T(\hat{M}^{J^\c}- M^*)\xi}
    \sqrt{\frac{1}{n}\sum_{i=1}^n(\phi_1(Y_i, X_i^\T\hat{\beta}) -  \phi_1(Y_i, X_i^\T\bar\beta) )^2} \notag \\
    &\leq  C  \sqrt{\xi^\T (\hat{M}^{J^\c}- M^*)\sum_{i\in J}\frac{1}{n}X_iX_i^\T(\hat{M}^{J^\c}- M^*)\xi}
    \sqrt{(\hat{\beta}-\bar\beta)^\T\sum_{i=1}^n\frac{1}{n}X_iX_i^\T(\hat{\beta}-\bar\beta)} \notag \\
    &\leq  C\sqrt{Q_1} \sqrt{s_0\frac{\log(p\vee n)}{n}}      .\label{inq:bound_of_term_II_outline}
\end{align}
We control $Q_1$ in the scenario where Assumption \ref{ass:sub_gaussian_entries_GLM} is satisfied and respectively in the scenario where Assumption \ref{ass:sub_gaussian_GLM} is satisfied. This is deferred to Lemma \ref{lem:cross_term_control_1} (i) and (ii).

\noindent \textcircled{4} For the term $R$, under Assumption \ref{bounded_higher_differentiate}, we have
\begin{align*}
   |R| &= |\sum_{i=1}^n \frac{1}{n} (\xi^\T M^* X_i)\int_0^1 (1-t)\phi_3(Y_i,X_i^\T\bar\beta+t X_i^\T(\hat{\beta}-\bar\beta))dt * (X_i^\T(\hat{\beta}-\bar\beta))^2| \notag \\
   &\leq C_7 \frac{1}{n}\sum_{i=1}^n   |\xi^\T M^* X_i| (X_i^\T(\hat{\beta}-\bar\beta))^2      \\
   &\leq C_7 \max_{1\leq i\leq n}\{|\xi^\T M^* X_i|\} \frac{1}{n} \sum_{i=1}^n (X_i^\T(\hat{\beta}-\bar\beta))^2.
\end{align*}
Under Assumption \ref{ass:sub_gaussian_entries_GLM} (ii) and \ref{ass:H_eigen_value}(ii),
$\xi^\T M^*X_0$ is a sub-gaussian random variable with the sub-gaussian norm being $O(\Vert\xi\Vert_2)$. Then we have w.h.p.
\begin{align}
    \max_{i=1,\ldots, n}\{|\xi^\T M^*X_i|\}\leq C \sqrt{\log( n)}\Vert\xi\Vert_2. \quad\text{(Lemma \ref{lem:bound_of_random_quan} (vi))
}  \label{inq:max_M_xi}
\end{align}
Hence we have w.h.p. that
\begin{align}
    |R| \leq C C_7 \frac{s_0\log(p\vee n)\sqrt{\log(n)}}{n} \Vert \xi\Vert_2. \label{inq:R_term_outline}
\end{align}

In summary, the bound of $\Delta_1$ in Theorem \ref{thm:initial_thm} (i) can be derived from (\ref{inq:V_term_outline}), (\ref{inq:bound_of_term_I_outline_sce1}), (\ref{inq:bound_of_term_II_outline}) and (\ref{inq:R_term_outline}). The bound of $\Delta_2$ in Theorem \ref{thm:initial_thm} (ii) can be derived from (\ref{inq:V_term_outline}), (\ref{inq:bound_of_term_I_outline_sce2}), (\ref{inq:bound_of_term_II_outline}) and (\ref{inq:R_term_outline}).

$\newline$
\noindent \textbf{(b) Bound of $\Delta_2$.} With (\ref{eq:V_estimator}), we have
\begin{align*}
    & |\hat{V} - V| \leq \underbrace{|\frac{1}{n} \sum_{i=1}^n (\xi^\T M^* X_i)^2\phi^2_1(Y_i,X_i^\T\bar\beta)- \mathbb{E}[(\xi^\T M^* X_1)^2\phi^2_1(Y_1,X_1^\T\bar\beta)] |}_{I} \\
    &+\underbrace{|\frac{2}{n}\sum_{i=1}^n (\xi^\T M^* X_i)^2\phi_1(Y_i,X_i^\T\bar\beta)(\phi_1(Y_i,X_i^\T\bar\beta)-\phi_1(Y_i,X_i^\T\hat{\beta}))|}_{II} \\ &+ \underbrace{|\frac{1}{n}\sum_{i=1}^n (\xi^\T M^* X_i)^2(\phi_1(Y_i,X_i^\T\bar\beta)-\phi_1(Y_i,X_i^\T\hat{\beta}))^2 )|}_{III} \\
    &+\underbrace{\frac{4}{n}\sum_{i\in J} |(\xi^\T M^* X_i)( \xi^\T(\hat{M}^{J^\c}  -M^*)X_i)|\phi^2_1(Y_i,X_i^\T\bar\beta)}_{IV} \\
    &+ \underbrace{\frac{4}{n}\sum_{i\in J} |(\xi^\T M^* X_i)( \xi^\T(\hat{M}^{J^\c}-M^*)X_i)|(\phi_1(Y_i,X_i^\T\bar\beta)-\phi_1(Y_i,X_i^\T\hat{\beta}))^2}_{V} \notag \\
    &+\underbrace{\frac{2}{n}\sum_{i\in J} (\xi^\T(\hat{M}^{J^\c}-M^*)X_i )^2\phi^2_1(Y_i,X_i^\T\bar\beta)}_{VI} \\
    &+ \underbrace{\frac{2}{n}\sum_{i\in J} (\xi^\T(\hat{M}^{J^\c}-M^*)X_i )^2(\phi_1(Y_i,X_i^\T\bar\beta)-\phi_1(Y_i,X_i^\T\hat{\beta}))^2}_{VII} \\
    &+VIII+IX+X+XI,
\end{align*}
where $VIII=\frac{4}{n}\sum_{i\in J^\c} |(\xi^\T M^* X_i)( \xi^\T(\hat{M}^{J}-M^*)X_i)|\phi^2_1(Y_i,X_i^\T\bar\beta)$, \\
$IX=\frac{4}{n}\sum_{i\in J^\c} |(\xi^\T M^* X_i)( \xi^\T(\hat{M}^{J}-M^*)X_i)|(\phi_1(Y_i,X_i^\T\bar\beta)-\phi_1(Y_i,X_i^\T\hat{\beta}))^2$, \\
$X=\frac{2}{n}\sum_{i\in J^\c} (\xi^\T(\hat{M}^{J}-M^*)X_i )^2\phi^2_1(Y_i,X_i^\T\bar\beta)$ \\
and $XI=\frac{2}{n}\sum_{i\in J^\c} (\xi^\T(\hat{M}^{J}-M^*)X_i )^2(\phi_1(Y_i,X_i^\T\bar\beta)-\phi_1(Y_i,X_i^\T\hat{\beta}))^2$.

We discuss how to control the following terms, $I$, $II$, $III$, $IV$, $V$, $VI$ and $VII$.
The remaining four terms $VIII$, $IX$, $X$ and $XI$ can be controlled similarly as $IV$, $V$, $VI$ and $VII$ respectively.
\\
\noindent \textcircled{1} We first discuss the term $I$. From Assumption \ref{ass:sub_gaussian_entries_GLM} (ii) and \ref{ass:residulal_glm_lower}, with Lemma \ref{lem:product-two-sub-gaussian-rv-norm} and \ref{lem:product-two-sub-exponential-rv-norm}, we have $(\xi^\T M^* X_i)^2\phi^2_1(Y_i,X_i^\T\bar\beta)$ for $i=1,\ldots, p$ are $\frac{1}{2}$-sub-exponential random variables with $\frac{1}{2}$-sub-exponential norm being $O(\Vert \xi\Vert_2^2)$. Hence $(\xi^\T M^* X_i)^2\phi^2_1(Y_i,X_i^\T\bar\beta)-\mathbb{E}[(\xi^\T M^* X_i)^2\phi^2_1(Y_i,X_i^\T\bar\beta)]$ for $i=1,\ldots, p$ are centered $\frac{1}{2}$-sub-exponential random variables with $\frac{1}{2}$-sub-exponential norm being $O(\Vert \xi\Vert_2^2)$. By concentration properties of centered $\frac{1}{2}$-sub-exponential random variables, we have w.h.p. that
\begin{align}
    |I|&= |\frac{1}{n} \sum_{i=1}^n (\xi^\T M^* X_i)^2\phi^2_1(Y_i,X_i^\T\bar\beta)- \mathbb{E}[(\xi^\T M^* X_1)^2\phi^2_1(Y_1,X_1^\T\bar\beta)] | \notag \\
    &\leq C n^{-\frac{1}{3}} \Vert\xi\Vert_2^2 .\quad\text{(Lemma \ref{lem:bound_of_random_quan} (vii))}   \label{inq:term_I_V_outline}
\end{align}

\noindent \textcircled{2} For the term $II$, under Assumption \ref{bounded_higher_differentiate}, we have
\begin{align*}
 |II| &= |\frac{2}{n}\sum_{i=1}^n (\xi^\T M^* X_i)^2\phi_1(Y_i,X_i^\T\bar\beta)(\phi_1(Y_i,X_i^\T\bar\beta)-\phi_1(Y_i,X_i^\T\hat{\beta}))|  \\
 &\leq C |\frac{1}{n}\sum_{i=1}^n (\xi^\T M^* X_i)^2\phi_1(Y_i,X_i^\T\bar\beta) (X_i^\T\bar\beta-X_i^\T\hat{\beta})| \\
 &\leq C \max_{1\leq i\leq p}\{|\phi_1(Y_i,X_i^\T\bar\beta)|\}\sqrt{\frac{1}{n}\sum_{i=1}^n (\xi^\T M^* X_i)^4}\sqrt{\frac{1}{n}\sum_{i=1}^n (X_i^\T\bar\beta-X_i^\T\hat{\beta})^2}.
\end{align*}
Under Assumption \ref{ass:sub_gaussian_entries_GLM}(ii) and \ref{ass:H_eigen_value}(ii), by Lemma \ref{lem:product-two-sub-gaussian-rv-norm} and \ref{lem:product-two-sub-exponential-rv-norm}, we have $(\xi^\T M^* X_i)^4$ are $\frac{1}{2}$-sub-exponential random variables with $\frac{1}{2}$-sub-exponential norm being $O(\Vert \xi\Vert_2^4)$. From the concentration property of $\frac{1}{2}$-sub-exponential random variables, we have w.h.p. that
\begin{align*}
  \frac{1}{n}\sum_{i=1}^n (\xi^\T M^* X_i)^4 \leq  C \Vert \xi\Vert_2^4. \quad\text{(Lemma \ref{lem:bound_of_random_quan} (vii))}
\end{align*}
Provided Assumption \ref{ass:residulal_glm_lower}, we have w.h.p. that
\begin{align}
    \max_{1\leq i\leq p} \{ |\phi_1(Y_i, X_i^\T \bar\beta)| \} \leq  C\sqrt{\log(n)} .\quad\text{(Lemma \ref{lem:bound_of_random_quan} (ii))} \label{inq:max_phi_1}
\end{align}
Combining the above and Condition \ref{beta_hat_condition}, we have, w.h.p. that
\begin{align}
    |II| &\leq C   \frac{\sqrt{s_0\log(p\vee n)}\sqrt{\log(n)}}{\sqrt{n}} \Vert \xi\Vert_2^2   .\label{inq:term_II_V_outline}
\end{align}

\noindent \textcircled{3} For the term $III$, under Assumption \ref{bounded_higher_differentiate}, we have
\begin{align*}
    |III|&=|\frac{1}{n}\sum_{i=1}^n (\xi^\T M^* X_i)^2(\phi_1(Y_i,X_i^\T\bar\beta)-\phi_1(Y_i,X_i^\T\hat{\beta}))^2 )| \\
    &\leq  C (\max_{1\leq i\leq n}  |\xi^\T M^* X_i|)^2 \frac{1}{n}\sum_{i=1}^n(X_i^\T\bar\beta-X_i^\T\hat{\beta})^2.
\end{align*}
Provided Assumption \ref{ass:sub_gaussian_entries_GLM}(ii) and \ref{ass:H_eigen_value}(ii), we have w.h.p. that
\begin{align}
    \max_{1\leq i\leq n}  |\xi^\T M^* X_i|\leq C\sqrt{\log(n)}\Vert \xi\Vert_2 .\quad\text{(see Lemma \ref{lem:bound_of_random_quan} (vi))}  \label{inq:max_xi_M_Xi}
\end{align}
 Combining the above and Condition \ref{beta_hat_condition}, we have w.h.p. that
\begin{align}
    |III| &\leq  C\frac{s_0\log(p\vee n)\log(n)}{n}\Vert \xi\Vert_2^2   .\label{inq:term_III_V_outline}
\end{align}

\noindent \textcircled{4} For the term $IV$, we have
\begin{align*}
    |IV|&=|\frac{4}{n}\sum_{i\in J} |(\xi^\T M^* X_i)( \xi^\T(\tilde{M^*}_{J^\c}-M^*)X_i)|\phi^2_1(Y_i,X_i^\T\bar\beta)|\\
    &\leq C \max_{1\leq i\leq n}|\xi^\T M^* X_i|\sqrt{\frac{1}{|J|}\sum_{i\in J} (\xi^\T (\hat{M}^{J^\c}-M^*)X_i)^2 }\sqrt{\frac{1}{|J|}\sum_{i\in J}\phi_1^4(Y_i,X_i^\T\bar\beta)} \\
    &\leq C \max_{1\leq i\leq n}|\xi^\T M^* X_i|\sqrt{\frac{1}{|J|}\sum_{i\in J} (\xi^\T (\hat{M}^{J^\c}-M^*)X_i)^2 }\sqrt{\frac{1}{n}\sum_{i=1 }^n \phi_1^4(Y_i,X_i^\T\bar\beta)}.
\end{align*}
Provide Assumption \ref{ass:residulal_glm_lower}, by Lemma \ref{lem:product-two-sub-gaussian-rv-norm} and \ref{lem:product-two-sub-exponential-rv-norm}, we have $\phi^4_1(Y_i, X_i^\T \bar\beta)$ are $\frac{1}{2}$-sub-exponential random variables, from the concentration properties of $\frac{1}{2}$-sub-exponential random variables, we have w.h.p.
\begin{align*}
    \frac{1}{n}\sum_{i=1 }^n \phi_1^4(Y_i, X_i^\T \bar\beta) \leq C. \quad\text{(Lemma \ref{lem:bound_of_random_quan} (iii))}
\end{align*}
Combining the above and (\ref{inq:max_M_xi}), we have w.h.p.
\begin{align}
    |IV|&\leq C \sqrt{Q_1}\sqrt{\log(n)}.\label{inq:term_IV_V_outline}
\end{align}

\noindent \textcircled{5} For the term $V$, we have
\begin{align*}
    |V|&=|\frac{4}{n}\sum_{i\in J} |(\xi^\T M^* X_i)( \xi^\T(\hat{M}^{J^\c}-M^*)X_i)|(\phi_1(Y_i,X_i^\T\bar\beta)-\phi_1(Y_i,X_i^\T\hat{\beta}))^2|\\
   &\leq C \frac{1}{n}\sum_{i\in J} |(\xi^\T M^* X_i)( \xi^\T(\hat{M}^{J^\c}-M^*)X_i)|(X_i^\T\bar\beta-X_i^\T\hat{\beta})^2 \notag \\ &\leq  C \Vert \bar\beta-\hat{\beta} \Vert^2_1 \max_{1\leq l \leq p, 1\leq m\leq p}\{\frac{1}{n}\sum_{i\in J}| (\xi^\T(\hat{M}^{J^\c}-M^*)X_i)(\xi^\T M^* X_i)X_{il} X_{im}| \}.
\end{align*}

First, we discuss the scenario where Assumption \ref{ass:sub_gaussian_entries_GLM} is satisfied. We have
\begin{align*}
    &|V|  \\
    &\leq C \Vert \bar\beta-\hat{\beta} \Vert^2_1 \max_{1\leq l \leq p, 1\leq m\leq p}\{\frac{1}{n}\sum_{i\in J}| (\xi^\T(\hat{M}^{J^\c}-M^*)X_i)(\xi^\T M^* X_i)X_{il} X_{im}| \} \\
    &\leq   C \Vert \bar\beta-\hat{\beta} \Vert^2_1 \Vert (\hat{M}^{J^\c}-M^*)\xi\Vert_1  \max_{i=1,\ldots,n} \{|\xi^\T M^* X_i|\}  \max_{1\leq k \leq p,1\leq l \leq p, 1\leq m\leq p}\{\frac{1}{n}\sum_{i\in J}| X_{ik}X_{il} X_{im}| \} \\
    &\leq   C \Vert \bar\beta-\hat{\beta} \Vert^2_1 \Vert (\hat{M}^{J^\c}-M^*)\xi\Vert_1  \max_{i=1,\ldots,n} \{|\xi^\T M^* X_i|\}  \max_{1\leq k \leq p,1\leq l \leq p, 1\leq m\leq p}\{\frac{1}{n}\sum_{i=1}^n| X_{ik}X_{il} X_{im}| \} .
\end{align*}
Under Assumption \ref{ass:sub_gaussian_entries_GLM} (i), by Lemma \ref{lem:product-two-sub-gaussian-rv-norm} and \ref{lem:product-three-rv-norm}, we have $|X_{ik} X_{il} X_{im}|$ for $i=1,\ldots,n$ and $1\leq k,l,m\leq p$ are $\frac{2}{3}$-sub-exponential random variables. By concentration properties of $\frac{2}{3}$-sub-exponential random variables, we have w.h.p.
\begin{align*}
    \max_{1\leq k\leq p,1\leq l\leq p,1\leq m\leq p,1\leq q\leq p} \{\frac{1}{n} \sum_{i=1}^n |X_{ik} X_{il} X_{im}| \}\leq C .\quad\text{(Lemma \ref{lem:bound_of_random_quan} (ii))}
\end{align*}
Combining the (\ref{inq:max_xi_M_Xi}) and Condition \ref{beta_hat_condition}, hence we have w.h.p. that
\begin{align}
    |V| \leq  C s^2_0 \frac{\log(p\vee n)\sqrt{\log(n)}}{n} \Vert (\hat{M}^{J^\c}-M^*)\xi\Vert_1 \Vert \xi\Vert_2   .\label{inq:term_V_V_outline_sce_1}
\end{align}

Second, in the scenario where Assumption \ref{ass:sub_gaussian_GLM} is satisfied, we have
\begin{align}
   |V|   &\leq   C \Vert \bar\beta-\hat{\beta} \Vert^2_1 \max_{1\leq l \leq p, 1\leq m\leq p}\{\frac{1}{n}\sum_{i\in J}| (\xi^\T(\hat{M}^{J^\c}-M^*)X_i)(\xi^\T M^* X_i)X_{il} X_{im}| \}        \notag \\
   &\leq  C  \Vert \bar\beta-\hat{\beta} \Vert^2_1 \max_{i \in J}\{|\xi^\T(\hat{M}^{J^\c}-M^*)X_i|\}\max_{1\leq l \leq p, 1\leq m\leq p}\{\frac{1}{n}\sum_{i\in J}| (\xi^\T M^* X_i)X_{il} X_{im}| \}   .
\end{align}
With Assumption \ref{ass:sub_gaussian_GLM}, conditionally on $\hat{M}^{J^\c}$, we have $\xi^\T(\hat{M}^{J^\c}-M^*)X_i$ for $i\in J$ are sub-gaussian random variable with sub-gaussian norm being $O(\Vert (\hat{M}^{J^\c}-M^*)\xi\Vert_2)$. Hence we have w.h.p.
\begin{align}
    \max_{i \in J}\{|\xi^\T(\hat{M}^{J^\c}-M^*)X_i|\} \leq C \sqrt{\log(n)}\Vert (\hat{M}^{J^\c}-M^*)\xi\Vert_2.  \quad\text{(Lemma \ref{lem:bound_of_random_quan} (xiii))}\label{inq:max_M_J_M_xi}
\end{align}
With Assumption \ref{ass:sub_gaussian_GLM} and \ref{ass:H_eigen_value}(ii), by Lemma \ref{lem:product-two-sub-gaussian-rv-norm} and \ref{lem:product-three-rv-norm},we have $|(\xi^\T M^* X_i)X_{il} X_{im}|$ for $i\in J$ are $\frac{2}{3}$-sub-exponential random variables with $\frac{2}{3}$-sub-exponential norm being $O(\Vert \xi\Vert_2)$. From the concentration property of $\frac{2}{3}$-sub-exponential random variables, we have w.h.p.
\begin{align*}
   \max_{1\leq l\leq p, 1\leq m\leq p} \frac{1}{n}\sum_{i=1}^n |(\xi^\T M^* X_i) X_{il}X_{im}  | \leq C\Vert\xi\Vert_2. \quad\text{(Lemma \ref{lem:bound_of_random_quan} (vi))}
\end{align*}
Combining the above and Condition \ref{beta_hat_condition}, we have w.h.p. that
\begin{align}
    |V|   &\leq  C s^2_0 \frac{\log(p\vee n)\sqrt{\log(n)}}{n} \Vert (\hat{M}^{J^\c}-M^*)\xi\Vert_2 \Vert \xi\Vert_2   . \label{inq:term_V_V_outline_sce_2}
\end{align}

 \noindent \textcircled{6} For the term $VI$, we have
\begin{align*}
 |VI| &=  |\frac{2}{n}\sum_{i\in J} (\xi^\T(\hat{M}^{J^\c}-M^*)X_i )^2\phi^2_1(Y_i,X_i^\T\bar\beta) |\\
 &\leq  C(\max_{1\leq i\leq p}\{|\phi_1(Y_i,X_i^\T \bar\beta)|\})^2 \frac{1}{n}\sum_{i\in J}(\xi^\T(\hat{M}^{J^\c}-M^*)X_i )^2.
\end{align*}
Under Assumption \ref{ass:residulal_glm_lower}, from (\ref{inq:max_phi_1}), we have w.h.p.
\begin{align}
    |VI|&\leq Q_1\log( n)     .\label{inq:term_VI_V_outline}
\end{align}

\noindent \textcircled{7} For the term $VII$, under Assumption \ref{bounded_higher_differentiate}, we have
\begin{align*}
    |VII| &= |\frac{2}{n}\sum_{i\in J} (\xi^\T(\hat{M}^{J^\c}-M^*)X_i )^2(\phi_1(Y_i,X_i^\T\bar\beta)-\phi_1(Y_i,X_i^\T\hat{\beta}))^2| \notag\\
  &\leq C \frac{2}{n}\sum_{i\in J} (\xi^\T(\hat{M}^{J^\c}-M^*)X_i )^2(X_i^\T\bar\beta-X_i^\T\hat{\beta})^2   \notag \\
  &\leq   C \Vert\bar\beta-\hat{\beta}\Vert_1^2  \max_{1\leq l\leq p,1\leq m\leq p} \{ \frac{1}{n}\sum_{i=1}^p(\xi^\T(\hat{M}^{J^\c}-M^*)X_i )^2 |X_{il}X_{im}|\} .
\end{align*}

In the scenario with Assumption \ref{ass:sub_gaussian_entries_GLM} (i), we have
\begin{align}
  |VII| &\leq  C \Vert\bar\beta-\hat{\beta}\Vert_1^2  \max_{1\leq l\leq p,1\leq m\leq p} \{ \frac{1}{n}\sum_{i=1}^p(\xi^\T(\hat{M}^{J^\c}-M^*)X_i )^2 |X_{il}X_{im}|\}  \notag\\
  &\leq C \Vert\bar\beta-\hat{\beta}\Vert_1^2   \Vert (\hat{M}^{J^\c}-M^*)\xi \Vert^2_1\max_{1\leq k\leq p,1\leq l\leq p,1\leq m\leq p,1\leq q\leq p} \{\frac{1}{n}\sum_{i=1}^p| X_{ik} X_{il} X_{im}X_{iq}| \}. \notag
\end{align}
Under Assumption \ref{ass:sub_gaussian_entries_GLM} (i), with Lemma \ref{lem:product-two-sub-gaussian-rv-norm} and \ref{lem:product-two-sub-exponential-rv-norm}, we have $| X_{ik} X_{il} X_{im}X_{iq}|$ are $\frac{1}{2}$-sub-exponential random variables. With the concentration property of $\frac{1}{2}$-sub-exponential random variables, we have w.h.p. that
\begin{align*}
    \max_{1\leq k\leq p,1\leq l\leq p,1\leq m\leq p,1\leq q\leq p} \{\frac{1}{n}\sum_{i=1}^p| X_{ik} X_{il} X_{im}X_{iq}| \} \leq C .  \quad\text{(Lemma \ref{lem:bound_of_random_quan} (ii))}
\end{align*}
Combining Condition \ref{beta_hat_condition}, we have
\begin{align}
    |VII| \leq Cs_0^2\frac{\log(p\vee n)}{n} \Vert (\hat{M}^{J^\c}-M^*)\xi \Vert^2_1 .\label{inq:term_VII_V_outline_sce_1}
\end{align}

In the scenario with Assumption \ref{ass:sub_gaussian_GLM}, we have
\begin{align}
  &|VII| \\
  &\leq C \Vert\bar\beta-\hat{\beta}\Vert_1^2  \max_{1\leq l\leq p,1\leq m\leq p} \{ \frac{1}{n}\sum_{i=1}^p(\xi^\T(\hat{M}^{J^\c}-M^*)X_i )^2 |x_{il}x_{im}|\}  \notag \\
  &\leq C \Vert\bar\beta-\hat{\beta}\Vert_1^2 \max_{i=1,\ldots,n}\{|\xi^\T(\hat{M}^{J^\c}-M^*)X_i|\}\max_{1\leq l\leq p, 1\leq m\leq p} \{ \frac{1}{|J|}\sum_{i\in J} | (\xi^\T (\hat{M}^{J^\c}- M^*) X_i)X_{il} X_{im}| \} .  \notag
\end{align}
Under Assumption \ref{ass:sub_gaussian_GLM} and conditionally on $\hat{M}^{J^\c}$, by Lemma \ref{lem:product-two-sub-gaussian-rv-norm} and \ref{lem:product-three-rv-norm}, we have $|(\xi^\T (\hat{M}^{J^\c}-M^*) X_i)X_{il} X_{im}|$ for $i\in J$ are $\frac{2}{3}$-sub-exponential random variables with $\frac{2}{3}$-sub-exponential norm being $O(\Vert (\hat{M}^{J^\c}-M^*)\xi\Vert_2)$. From concentration properties of $\frac{2}{3}$-sub-exponential random variables, we have w.h.p. that
\begin{align*}
   \max_{1\leq l\leq p, 1\leq m\leq p} \frac{1}{n}\sum_{i=1}^n |(\xi^\T (\hat{M}^{J^\c} - M^*) X_i) X_{il}X_{im}  | \leq C\Vert (\hat{M}^{J^\c} - M^*)\xi\Vert_2.\quad\text{(Lemma \ref{lem:bound_of_random_quan}(xii))}
\end{align*}
Combining the above, Condition \ref{beta_hat_condition} and (\ref{inq:max_M_J_M_xi}), we have w.h.p.
\begin{align}
    |VII| \leq  C s^2_0 \frac{\log(p\vee n) \sqrt{\log(n)} }{n} \Vert (\hat{M}^{J^\c}-M^*)\xi \Vert^2_2  .\label{inq:term_VII_V_outline_sce_2}
\end{align}

The bound of $\Delta_2$ in Theorem \ref{thm:initial_thm} (i) can be proved by incorporating (\ref{inq:term_I_V_outline}), (\ref{inq:term_II_V_outline}), (\ref{inq:term_III_V_outline}), (\ref{inq:term_IV_V_outline}), (\ref{inq:term_V_V_outline_sce_1}), (\ref{inq:term_VI_V_outline}), (\ref{inq:term_VII_V_outline_sce_1}).
The bound of $\Delta_2$ in Theorem \ref{thm:initial_thm} (ii) can be proved by incorporating (\ref{inq:term_I_V_outline}), (\ref{inq:term_II_V_outline}), (\ref{inq:term_III_V_outline}), (\ref{inq:term_IV_V_outline}), (\ref{inq:term_V_V_outline_sce_2}), (\ref{inq:term_VI_V_outline}), (\ref{inq:term_VII_V_outline_sce_2}).

\subsection{Bounds of $Q_1$ and $Q_1^\prime$  in Theorem \ref{thm:initial_thm}}

We give bounds of $Q_1$ in Theorem \ref{thm:initial_thm}. The bounds of $Q_1^\prime$ are similar.
We consider the GLM setting with general $\phi(a,b)$ below. The LM setting can be treated as a special case of the GLM setting.

\begin{restatable}{lem}{controlqqprime}
\label{lem:cross_term_control_1}
$\newline$
\noindent (i) Suppose that Assumption \ref{ass:sub_gaussian_entries_GLM} (i) and \ref{ass:H_eigen_value} (i) are satisfied
and $\frac{\log(p\vee n)}{\sqrt{n} }\leq \frac{c(1)\rho}{3}\sqrt{n}$. Then we have with probability at least $1-\frac{2}{p\vee n}$,
\begin{align*}
    Q_1 &\leq C  (\Vert (\hat{M}^{J^\c}- M^*)\xi\Vert_2^2 +  \sqrt{\frac{\log(p\vee n)}{n}}\Vert (\hat{M}^{J^\c}- M^*)\xi\Vert_1^2 ),
\end{align*}
where $c(1)$ is a constant from Lemma \ref{lem:rv_centering} and $C$ is a constant depending on $(C_3, \sigma_1, \rho)$ only.

\noindent (ii) Suppose that Assumption  \ref{ass:sub_gaussian_GLM} and \ref{ass:H_eigen_value} (i) are satisfied and $\frac{\log(p\vee n)}{\sqrt{n}}\leq \c(1) \rho\sqrt{n}$, we have with probability at least $1-\frac{2}{p\vee n}$,
\begin{align*}
    Q_1 &\leq C \Vert (\hat{M}^{J^\c}- M^*)\xi\Vert_2^2 ,
\end{align*}
where $c(1)$ is a constant from Lemma \ref{lem:rv_centering} and $C$ is a constant depending on $(C_3, \sigma_1)$ only.
\end{restatable}

\subsection{Proofs of main results}

By the definition of matrix $\Vert \cdot \Vert_1$-norm and matrix $\Vert\cdot \Vert_2$-norm, we have
\begin{equation}
    \begin{aligned}
     L_1 &=  \Vert (\hat{M}^{J}-M^*)\xi\Vert_1 \leq \Vert \hat{M}^{J}_{S_\xi,S_\xi} - M^*_{S_\xi,S_\xi}\Vert_1 \Vert\xi\Vert_1 \\
   L_2&= \Vert (\hat{M}^{J}-M^*)\xi\Vert_2 \leq \Vert \hat{M}^{J}_{S_\xi,S_\xi} - M^*_{S_\xi,S_\xi}\Vert_2 \Vert\xi\Vert_2\\
   L_1^\prime &= \Vert (\hat{M}^{J^\c}-M^*)\xi\Vert_1 \leq \Vert \hat{M}^{J^\c}_{S_\xi,S_\xi} - M^*_{S_\xi,S_\xi}\Vert_1 \Vert\xi\Vert_1\\
  L_2^\prime &=  \Vert (\hat{M}^{J^\c}-M^*)\xi\Vert_2 \leq \Vert \hat{M}^{J^\c}_{S_\xi,S_\xi} - M^*_{S_\xi,S_\xi}\Vert_2 \Vert\xi\Vert_2
    \end{aligned}
    \label{inq:M_hat_M_xi}
\end{equation}

Theorem \ref{thm_1} can be proved by incorporating (\ref{inq:M_hat_M_xi}) with replacing $\hat{M}^J$ and $\hat{M}^{J^\c}$ with $\hat{\Omega}^J$ and $\hat{\Omega}^{J^\c}$ and Condition \ref{precision_matrix_LM} and applying Theorem \ref{thm:initial_thm} and Lemma \ref{lem:cross_term_control_1} with $\phi(a,b)=\frac{1}{2}(a-b)^2$. Theorem \ref{thm_2} can be proved by incorporating (\ref{inq:M_hat_M_xi}) and Condition \ref{precision_matrix_GLM}
and applying Theorem \ref{thm:initial_thm} and Lemma \ref{lem:cross_term_control_1}. The probability control is determined by the intersection of the relevant events included in the results above.

\section{Numerical studies}

\label{sec:numerical_study}

We present numerical studies to support the validity of our confidence intervals and, in the LM setting, to compare our proposed method and the method of \cite{cai2020optimal},
which is called high-dimensional individualized treatment selection (HITS) and implemented by the codes in https://github.com/zijguo/ITE.
In Section \ref{sec:numerical_study_sparse_omega}, we consider the setting of LM. In Section \ref{sec:numerical_study_sparse_M}, we consider the setting of GLM. Throughout, $\hat{\beta}$, $\hat{\beta}_J$ and $\hat{\beta}_{J^\c}$ are computed by the Lasso algorithm in \texttt{glmnet} in R \citep{glmnet1} and the penalty parameters are decided by 10-fold cross-validation. Given high computational cost of Algorithms \ref{alg:omega_est_two_stage} and \ref{alg:est_M_two_stage}, we only apply Algorithm \ref{alg:omega_est} and \ref{alg:est_M}, which are computed by \texttt{flare} in R \citep{flare}.

\subsection{Numerical study in LM Setting}
\label{sec:numerical_study_sparse_omega}

In the numerical study, we let $p=400$ and $n \in \{200, 400, 500, 600, 800\}$. The true regression model for data generation is
\begin{align}
	Y_i = X_i^\T \gamma^* + X^2_{i1}\epsilon_{i1}+ X^2_{i4}\epsilon_{i2},\quad \text{for $i=1,\ldots,n$}, \label{eq:LM_true_model}
\end{align}
where $\gamma^*=(4,2,4,4,-2, 0_{395})^\T \in \mathbb{R}^p$ (with $0_{395}\in \mathbb{R}^{395}$ a vector of all zeros) and
$(\epsilon_{i1}, \epsilon_{i2})_{i=1,\ldots,n}$ are i.i.d.~from $\mathcal{N}(0,1)$. It can be verified that $\bar\beta=\gamma^*=(4,2,4,4,-2, 0_{395})^\T$ is also the target coefficient vector satisfying (\ref{eq:target_beta_LM}),
even though the linear regression model of $Y_i$ on $X_i$ is misspecified.

We let $(X_1,\ldots,X_n)$ be i.i.d.~as $X_0$
where $X_0\sim \mathcal{N}_p(0_p,\Sigma)$ and $\Sigma \in \mathbb{R}^{p\times p}$ satisfying $\Sigma_{ij}=0.3^{|i-j|}$. We consider two choices of loadings. The first choice of loading is generated as $\xi =(-0.5, -0.25, 0.25, 0.5, 0.25,\tilde\xi)^\T$ and the second choice of loading is generated as $\xi =(0.25,0.25,0.25,0.25,0.25,\tilde\xi)^\T$, where $\tilde\xi\in\mathbb{R}^{p-5}$ is randomly generated by
$\tilde\xi/q \sim \mathcal{N}_{p-5}(0_{p-5}, I_{p-5})$ and then fixed in the entire simulations, and $q=0.025$ or $q=0.01$.
For $\xi$ used in the results below,  we have
$\xi^\T\bar\beta=0$, and $\|\xi\|_1=9.61$ and $\|\xi_2\|_2= 0.96$ with $q=0.025$ or $\|\xi\|_1=4.87$ and $\|\xi_2\|_2= 0.85$ with $q=0.01$ for the first choice of loading,
and $\xi^\T\bar\beta = 3$ and $\|\xi\|_1=8.88$ and $\|\xi_2\|_2= 0.73$ with $q=0.025$ or $\|\xi\|_1=4.67$ and $\|\xi_2\|_2= 0.60$ with $q=0.01$ for the second choice of loading.

We compute the debiased estimator $\widehat{\xi^\T\beta}$ in (\ref{debiased_est_1}) and estimated variance $\hat{V}$ in (\ref{eq:definition_V_hat}). We construct $95\%$ confidence intervals using (\ref{CI}). We let $|J|=|J^\c|=\frac{n}{2}$ and $\rho_n=\sqrt{\frac{\log(p)}{n}}$ in Algorithm \ref{alg:omega_est}.
As mentioned earlier, we also apply the method HITS from \cite{cai2020optimal}.
We report the simulation results based on 200 replications in Table \ref{table:sparse_omega_LM}.

From Table \ref{table:sparse_omega_LM}, we see that both our method and HITS achieve significant bias reduction as $n$ increases. Moreover, when $n$ is large enough, our method yields coverage proportions closer to the nominal level $95\%$ than HITS.
This comparison numerically demonstrates potential advantage of our proposed method over HITS in the presence of model misspecification.

\begin{table}[!ht]
\scriptsize		
\caption{Summary of $95\%$ prediction confidence intervals (CIs) in LM setting}
\centering
\begin{tabular}{|c|c|c|cccc|cccc|cc|}
  \hline
&& &\multicolumn{4}{c|}{Proposed method}&\multicolumn{4}{c|}{HITS}& \multicolumn{2}{c|}{Original Lasso} \\
$\xi^\T\bar\beta$& $q$ &$n$ & Bias & Sd & Cov prop & CI length & Bias & Sd & Cov prop & CI length & Bias & Sd \\
  \hline
\multirow{8}{*}{0}&\multirow{4}{*}{0.025}& 200& 0.025 & 0.244 & 0.890 & 0.831 & 0.018 & 0.242 & 0.700 & 0.542 & 0.048 & 0.228 \\
  &&400& 0.006 & 0.178 & 0.930 & 0.638 & 0.002 & 0.181 & 0.785 & 0.434 & 0.044 & 0.167 \\
 &&500 & 0.004 & 0.145 & 0.960 & 0.579 & 0.001 & 0.145 & 0.815 & 0.390 & 0.038 & 0.148 \\
  && 600& 0.037 & 0.132 & 0.935 & 0.525 & 0.037 & 0.134 & 0.800 & 0.353 & 0.007 & 0.123 \\
 &&800& 0.008 & 0.122 & 0.970 & 0.480 & 0.007 & 0.118 & 0.815 & 0.308 & 0.028 & 0.111 \\
  \cline{2-13}
&\multirow{4}{*}{0.01}&200 & 0.011 & 0.236 & 0.910 & 0.760 & 0.012 & 0.236 & 0.665 & 0.460 & 0.049 & 0.234 \\
 &&400 & 0.005 & 0.169 & 0.935 & 0.616 & 0.007 & 0.166 & 0.765 & 0.378 & 0.054 & 0.163 \\
  && 500& 0.010 & 0.130 & 0.960 & 0.575 & 0.011 & 0.138 & 0.820 & 0.341 & 0.032 & 0.130 \\
 &&600 & 0.010 & 0.120 & 0.975 & 0.510 & 0.012 & 0.117 & 0.810 & 0.308 & 0.027 & 0.120 \\
  &&800 & 0.005 & 0.117 & 0.940 & 0.440 & 0.004 & 0.116 & 0.750 & 0.263 & 0.026 & 0.113 \\
   \hline
\multirow{8}{*}{3}&\multirow{4}{*}{0.025}& 200& 0.058 & 0.141 & 0.965 & 0.585 & 0.014 & 0.133 & 0.840 & 0.372 & 0.143 & 0.124 \\
  &&400 & 0.047 & 0.100 & 0.965 & 0.445 & 0.024 & 0.095 & 0.915 & 0.307 & 0.098 & 0.084 \\
  && 500& 0.029 & 0.091 & 0.955 & 0.401 & 0.012 & 0.087 & 0.885 & 0.275 & 0.093 & 0.075 \\
 &&600 & 0.028 & 0.082 & 0.975 & 0.362 & 0.013 & 0.077 & 0.915 & 0.253 & 0.082 & 0.068 \\
  &&800 & 0.016 & 0.076 & 0.950 & 0.307 & 0.010 & 0.072 & 0.870 & 0.218 & 0.071 & 0.062 \\   \cline{2-13}
&\multirow{4}{*}{0.01}&200 & 0.048 & 0.120 & 0.980 & 0.556 & 0.014 & 0.115 & 0.805 & 0.285 & 0.150 & 0.119 \\
  &&400 & 0.042 & 0.088 & 0.955 & 0.400 & 0.009 & 0.087 & 0.805 & 0.223 & 0.104 & 0.085 \\
  && 500& 0.018 & 0.079 & 0.965 & 0.356 & 0.005 & 0.077 & 0.835 & 0.203 & 0.104 & 0.076 \\
 &&600 & 0.023 & 0.069 & 0.980 & 0.317 & 0.004 & 0.066 & 0.810 & 0.183 & 0.084 & 0.063 \\
  &&800 & 0.020 & 0.065 & 0.945 & 0.267 & 0.008 & 0.061 & 0.815 & 0.158 & 0.071 & 0.061 \\
   \hline
\end{tabular}
\label{table:sparse_omega_LM}
\end{table}

\subsection{Numerical study in GLM setting}

\label{sec:numerical_study_sparse_M}

In the GLM numerical study, we let $p=300$ and $n \in \{150, 300, 500, 600, 800\}$. The true regression model for data generation is
\begin{align}
 P(Y_i=1|X_i) = \frac{1}{2}\frac{1}{1+\exp\{-X_i^\T\gamma^*+2\}} + \frac{1}{2}\frac{1}{1+\exp\{-X_i^\T\gamma^*-2\}}, \quad \text{for $i=1,\ldots,n$},  \label{eq:GLM_true_model_misspecified}
\end{align}
where $\gamma^*=(0.5,0.5,0.5,0.5,0.075,0_{295})^\T \in \mathbb{R}^p$. We let $(X_1,\ldots,X_n)$ be i.i.d.~as $X_0$ and $X_0=\Sigma^{\frac{1}{2}} Z_0$ where the components of $Z_0$ are i.i.d.~Rademacher variables and $\Sigma \in \mathbb{R}^{p\times p}$ satisfying
\begin{equation}
    \begin{cases}
    \Sigma_{ij}= 5 \times 0.1^{|i-j|},\quad\text{if $i\leq 5$, $j\leq 5$ or $i>5, j>5$,} \\
    \Sigma_{ij}= 0,\quad\text{Otherwise.}
    \end{cases}  \notag
\end{equation}
In such a setting, the logistic regression model of $Y_i$ on $X_i$ is misspecified.
It can be derived that the target coefficients vector is
$\bar\beta=(0.2974, 0.2968, 0.2968, 0.2974, 0.0445, 0_{295})^\T$. See Supplement Section \ref{sec:calculation_of_beta_star}.

We consider two choices of loadings. The first choice of loading is generated as $\xi=(-1.15, 1, -1, 1, 1,\tilde\xi)$ and the second choice
of loading is generated as $\xi=(-1, 1, -1, 1, 3,\tilde\xi)$, where $\tilde\xi\in\mathbb{R}^{p-5}$ is randomly generated by $\tilde\xi/q \sim \mathcal{N}_{p-5}(0_{p-5}, I_{p-5})$ and then fixed in the entire simulations, and $q=0.1$ or $q=0.05$.
For $\xi$ used in the results below,  we have
$\xi^\T\bar\beta=-0.00016$, and $\|\xi\|_1=27.20$ and $\|\xi_2\|_2= 2.81$ with $q=0.1$ or $\|\xi\|_1=16.60$ and $\|\xi_2\|_2= 2.45$ with $q=0.05$ for the first choice of loading,
and $\xi^\T\bar\beta =  0.1335$, and $\|\xi\|_1=30.13$ and $\|\xi_2\|_2= 3.99$ with $q=0.1$ or $\|\xi\|_1=18.55$ and $\|\xi \|_2= 3.7$ with $q=0.05$ for the second choice of loading.

We compute the debiased estimator $\widehat{\xi^\T\beta}$ in (\ref{debiased_est_sparse_M}) and estimated variance $\hat{V}$ in (\ref{eq:V_estimator}). We construct $95\%$ confidence intervals using (\ref{CI}). We let $|J|=|J^\c|=\frac{n}{2}$ and $\rho_n=\sqrt{\frac{\log(p)}{n}}$ in Algorithm \ref{alg:est_M}. We report the simulation results based on 200 replications in Table \ref{table:sparse_M}.

From Table \ref{table:sparse_M}, when $n$ is large enough, our method achieves significant bias reduction
and yield coverage proportions reasonably close to $95\%$ even for dense loadings $\xi$ in the presence of model misspecification.

\begin{table}[!ht]
\scriptsize		
\caption{Summary of $95\%$ prediction confidence intervals (CIs) in GLM setting}
\centering
\begin{tabular}{|c|c|c|cccc|cc|}
\hline
&&&\multicolumn{4}{c|}{Proposed method}&\multicolumn{2}{c|}{Original Lasso} \\
$\xi^\T\bar\beta$ &$q$& $n$& Bias & Sd & Cov prop & CI length & Bias & sd \\
  \hline
\multirow{8}{*}{-0.00016}&\multirow{4}{*}{0.1}&150& 0.013 & 0.288 & 0.940 & 0.961 & 0.011 & 0.192 \\
&&300& 0.005 & 0.185 & 0.945 & 0.708 & 0.033 & 0.129 \\
&& 500& 0.016 & 0.139 & 0.975 & 0.573 & 0.031 & 0.104 \\
&& 600& 0.010 & 0.131 & 0.965 & 0.534 & 0.026 & 0.090 \\
  &&800& 0.001 & 0.105 & 0.965 & 0.464 & 0.022 & 0.072 \\
\cline{2-9}
&\multirow{4}{*}{0.05}& 150& 0.016 & 0.229 & 0.940 & 0.848 & 0.012 & 0.191 \\
 && 300& 0.017 & 0.156 & 0.940 & 0.618 & 0.033 & 0.127 \\
  && 500& 0.008 & 0.140 & 0.935 & 0.506 & 0.031 & 0.103 \\
  && 600& 0.005 & 0.122 & 0.935 & 0.469 & 0.026 & 0.089 \\
  && 800& 0.003 & 0.101 & 0.965 & 0.406 & 0.022 & 0.072 \\
  \hline
\multirow{8}{*}{0.1335}&\multirow{4}{*}{0.1}& 150& 0.010 & 0.380 & 0.930 & 1.363 & 0.095 & 0.206 \\
&&300 & 0.022 & 0.255 & 0.950 & 0.997 & 0.107 & 0.134 \\
 && 500& 0.004 & 0.200 & 0.945 & 0.803 & 0.093 & 0.116 \\
 && 600& 0.007 & 0.194 & 0.945 & 0.742 & 0.080 & 0.104 \\
  &&800 & 0.000 & 0.163 & 0.960 & 0.641 & 0.068 & 0.097 \\
\cline{2-9}
&\multirow{4}{*}{0.05}&150 & 0.030 & 0.464 & 0.940 & 1.290 & 0.095 & 0.205 \\
 &&300 & 0.008 & 0.254 & 0.945 & 0.938 & 0.108 & 0.133 \\
  &&500 & 0.007 & 0.189 & 0.965 & 0.755 & 0.093 & 0.116 \\
 &&600 & 0.000 & 0.171 & 0.960 & 0.697 & 0.080 & 0.103 \\
 &&800 & 0.007 & 0.149 & 0.970 & 0.602 & 0.068 & 0.097 \\
  \hline
\end{tabular}
\label{table:sparse_M}
\end{table}

\newpage
\bibliographystyle{apacite}
\bibliography{references}

\begin{thebibliography}{}

\bibitem [\protect \citeauthoryear {%
Bickel%
, Ritov%
\BCBL {}\ \BBA {} Tsybakov%
}{%
Bickel%
\ \protect \BOthers {.}}{%
{\protect \APACyear {2009}}%
}]{%
lasso_and_dantzig}
\APACinsertmetastar {%
lasso_and_dantzig}%
\begin{APACrefauthors}%
Bickel, P\BPBI J.%
, Ritov, Y.%
\BCBL {}\ \BBA {} Tsybakov, A\BPBI B.%
\end{APACrefauthors}%
\unskip\
\newblock
\APACrefYearMonthDay{2009}{}{}.
\newblock
{\BBOQ}\APACrefatitle {{Simultaneous analysis of Lasso and Dantzig selector}}
  {{Simultaneous analysis of Lasso and Dantzig selector}}.{\BBCQ}
\newblock
\APACjournalVolNumPages{The Annals of Statistics}{37}{}{1705 -- 1732}.
\PrintBackRefs{\CurrentBib}

\bibitem [\protect \citeauthoryear {%
Cai%
, Cai%
\BCBL {}\ \BBA {} Guo%
}{%
Cai%
\ \protect \BOthers {.}}{%
{\protect \APACyear {2021}}%
}]{%
cai2020optimal}
\APACinsertmetastar {%
cai2020optimal}%
\begin{APACrefauthors}%
Cai, T.%
, Cai, T\BPBI T.%
\BCBL {}\ \BBA {} Guo, Z.%
\end{APACrefauthors}%
\unskip\
\newblock
\APACrefYearMonthDay{2021}{}{}.
\newblock
{\BBOQ}\APACrefatitle {Optimal statistical inference for individualized
  treatment effects in high-dimensional models} {Optimal statistical inference
  for individualized treatment effects in high-dimensional models}.{\BBCQ}
\newblock
\APACjournalVolNumPages{Journal of the Royal Statistical Society Series B:
  Statistical Methodology}{83}{}{669-719}.
\PrintBackRefs{\CurrentBib}

\bibitem [\protect \citeauthoryear {%
Cai%
, Liu%
\BCBL {}\ \BBA {} Luo%
}{%
Cai%
\ \protect \BOthers {.}}{%
{\protect \APACyear {2011}}%
}]{%
clime}
\APACinsertmetastar {%
clime}%
\begin{APACrefauthors}%
Cai, T.%
, Liu, W.%
\BCBL {}\ \BBA {} Luo, X.%
\end{APACrefauthors}%
\unskip\
\newblock
\APACrefYearMonthDay{2011}{}{}.
\newblock
{\BBOQ}\APACrefatitle {A constrained $\ell_1$ minimization approach to sparse
  precision matrix estimation} {A constrained $\ell_1$ minimization approach to
  sparse precision matrix estimation}.{\BBCQ}
\newblock
\APACjournalVolNumPages{Journal of the American Statistical
  Association}{106}{}{594--607}.
\PrintBackRefs{\CurrentBib}

\bibitem [\protect \citeauthoryear {%
Friedman%
, Tibshirani%
\BCBL {}\ \BBA {} Hastie%
}{%
Friedman%
\ \protect \BOthers {.}}{%
{\protect \APACyear {2010}}%
}]{%
glmnet1}
\APACinsertmetastar {%
glmnet1}%
\begin{APACrefauthors}%
Friedman, J.%
, Tibshirani, R.%
\BCBL {}\ \BBA {} Hastie, T.%
\end{APACrefauthors}%
\unskip\
\newblock
\APACrefYearMonthDay{2010}{}{}.
\newblock
{\BBOQ}\APACrefatitle {Regularization paths for generalized linear models via
  coordinate descent} {Regularization paths for generalized linear models via
  coordinate descent}.{\BBCQ}
\newblock
\APACjournalVolNumPages{Journal of Statistical Software}{33}{}{1--22}.
\PrintBackRefs{\CurrentBib}

\bibitem [\protect \citeauthoryear {%
Guo%
, Rakshit%
, Herman%
\BCBL {}\ \BBA {} Chen%
}{%
Guo%
\ \protect \BOthers {.}}{%
{\protect \APACyear {2021}}%
}]{%
guo2021inference}
\APACinsertmetastar {%
guo2021inference}%
\begin{APACrefauthors}%
Guo, Z.%
, Rakshit, P.%
, Herman, D\BPBI S.%
\BCBL {}\ \BBA {} Chen, J.%
\end{APACrefauthors}%
\unskip\
\newblock
\APACrefYearMonthDay{2021}{}{}.
\newblock
{\BBOQ}\APACrefatitle {Inference for the case probability in high-dimensional
  logistic regression} {Inference for the case probability in high-dimensional
  logistic regression}.{\BBCQ}
\newblock
\APACjournalVolNumPages{Journal of Machine Learning Research}{22}{}{1--54}.
\PrintBackRefs{\CurrentBib}

\bibitem [\protect \citeauthoryear {%
Hastie%
, Tibshirani%
\BCBL {}\ \BBA {} Wainwright%
}{%
Hastie%
\ \protect \BOthers {.}}{%
{\protect \APACyear {2015}}%
}]{%
hastie_book}
\APACinsertmetastar {%
hastie_book}%
\begin{APACrefauthors}%
Hastie, T.%
, Tibshirani, R.%
\BCBL {}\ \BBA {} Wainwright, M.%
\end{APACrefauthors}%
\unskip\
\newblock
\APACrefYear{2015}.
\newblock
\APACrefbtitle {Statistical {L}earning with {S}parsity: {T}he Lasso and
  {G}eneralizations} {Statistical {L}earning with {S}parsity: {T}he lasso and
  {G}eneralizations}.
\newblock
\APACaddressPublisher{}{Chapman \& Hall/CRC}.
\PrintBackRefs{\CurrentBib}

\bibitem [\protect \citeauthoryear {%
Horn%
\ \BBA {} Johnson%
}{%
Horn%
\ \BBA {} Johnson%
}{%
{\protect \APACyear {2012}}%
}]{%
horn2012matrix}
\APACinsertmetastar {%
horn2012matrix}%
\begin{APACrefauthors}%
Horn, R\BPBI A.%
\BCBT {}\ \BBA {} Johnson, C\BPBI R.%
\end{APACrefauthors}%
\unskip\
\newblock
\APACrefYear{2012}.
\newblock
\APACrefbtitle {{M}atrix {A}nalysis} {{M}atrix {A}nalysis}.
\newblock
\APACaddressPublisher{}{Cambridge University Press}.
\PrintBackRefs{\CurrentBib}

\bibitem [\protect \citeauthoryear {%
Javanmard%
\ \BBA {} Montanari%
}{%
Javanmard%
\ \BBA {} Montanari%
}{%
{\protect \APACyear {2014}}%
}]{%
JMLR:v15:javanmard14a}
\APACinsertmetastar {%
JMLR:v15:javanmard14a}%
\begin{APACrefauthors}%
Javanmard, A.%
\BCBT {}\ \BBA {} Montanari, A.%
\end{APACrefauthors}%
\unskip\
\newblock
\APACrefYearMonthDay{2014}{}{}.
\newblock
{\BBOQ}\APACrefatitle {Confidence intervals and hypothesis testing for
  high-dimensional regression} {Confidence intervals and hypothesis testing for
  high-dimensional regression}.{\BBCQ}
\newblock
\APACjournalVolNumPages{Journal of Machine Learning
  Research}{15}{}{2869--2909}.
\PrintBackRefs{\CurrentBib}

\bibitem [\protect \citeauthoryear {%
Li%
, Zhao%
, Wang%
, Yuan%
\BCBL {}\ \BBA {} Liu%
}{%
Li%
\ \protect \BOthers {.}}{%
{\protect \APACyear {2022}}%
}]{%
flare}
\APACinsertmetastar {%
flare}%
\begin{APACrefauthors}%
Li, X.%
, Zhao, T.%
, Wang, L.%
, Yuan, X.%
\BCBL {}\ \BBA {} Liu, H.%
\end{APACrefauthors}%
\unskip\
\newblock
\APACrefYearMonthDay{2022}{}{}.
\newblock
\APACrefbtitle {flare: Family of Lasso regression.} {flare: Family of lasso
  regression.}
\newblock
\begin{APACrefURL} \url{https://CRAN.R-project.org/package=flare}
  \end{APACrefURL}
\newblock
\APACrefnote{R package version 1.7.0.1}
\PrintBackRefs{\CurrentBib}

\bibitem [\protect \citeauthoryear {%
Ning%
\ \BBA {} Liu%
}{%
Ning%
\ \BBA {} Liu%
}{%
{\protect \APACyear {2017}}%
}]{%
Ning}
\APACinsertmetastar {%
Ning}%
\begin{APACrefauthors}%
Ning, Y.%
\BCBT {}\ \BBA {} Liu, H.%
\end{APACrefauthors}%
\unskip\
\newblock
\APACrefYearMonthDay{2017}{}{}.
\newblock
{\BBOQ}\APACrefatitle {{A general theory of hypothesis tests and confidence
  regions for sparse high-dimensional models}} {{A general theory of hypothesis
  tests and confidence regions for sparse high-dimensional models}}.{\BBCQ}
\newblock
\APACjournalVolNumPages{The Annals of Statistics}{45}{}{158 -- 195}.
\PrintBackRefs{\CurrentBib}

\bibitem [\protect \citeauthoryear {%
van~de Geer%
, Bühlmann%
, Ritov%
\BCBL {}\ \BBA {} Dezeure%
}{%
van~de Geer%
\ \protect \BOthers {.}}{%
{\protect \APACyear {2014}}%
}]{%
van_de_Geer_2014}
\APACinsertmetastar {%
van_de_Geer_2014}%
\begin{APACrefauthors}%
van~de Geer, S.%
, Bühlmann, P.%
, Ritov, Y.%
\BCBL {}\ \BBA {} Dezeure, R.%
\end{APACrefauthors}%
\unskip\
\newblock
\APACrefYearMonthDay{2014}{}{}.
\newblock
{\BBOQ}\APACrefatitle {On asymptotically optimal confidence regions and tests
  for high-dimensional models} {On asymptotically optimal confidence regions
  and tests for high-dimensional models}.{\BBCQ}
\newblock
\APACjournalVolNumPages{The Annals of Statistics}{42}{}{1166 -- 1202}.
\PrintBackRefs{\CurrentBib}

\bibitem [\protect \citeauthoryear {%
Zhang%
\ \BBA {} Zhang%
}{%
Zhang%
\ \BBA {} Zhang%
}{%
{\protect \APACyear {2014}}%
}]{%
zhang2014confidence}
\APACinsertmetastar {%
zhang2014confidence}%
\begin{APACrefauthors}%
Zhang, C\BHBI H.%
\BCBT {}\ \BBA {} Zhang, S\BPBI S.%
\end{APACrefauthors}%
\unskip\
\newblock
\APACrefYearMonthDay{2014}{}{}.
\newblock
{\BBOQ}\APACrefatitle {Confidence intervals for low-dimensional parameters in
  high-dimensional linear models} {Confidence intervals for low-dimensional
  parameters in high-dimensional linear models}.{\BBCQ}
\newblock
\APACjournalVolNumPages{Journal of the Royal Statistical Society Series B:
  Statistical Methodology}{76}{}{217--242}.
\PrintBackRefs{\CurrentBib}

\bibitem [\protect \citeauthoryear {%
Zhao%
, Zhou%
\BCBL {}\ \BBA {} Liu%
}{%
Zhao%
\ \protect \BOthers {.}}{%
{\protect \APACyear {2023}}%
}]{%
yliu}
\APACinsertmetastar {%
yliu}%
\begin{APACrefauthors}%
Zhao, J.%
, Zhou, Y.%
\BCBL {}\ \BBA {} Liu, Y.%
\end{APACrefauthors}%
\unskip\
\newblock
\APACrefYearMonthDay{2023}{}{}.
\newblock
{\BBOQ}\APACrefatitle {Estimation of linear functionals in high-dimensional
  linear models: from sparsity to non-sparsity} {Estimation of linear
  functionals in high-dimensional linear models: from sparsity to
  non-sparsity}.{\BBCQ}
\newblock
\APACjournalVolNumPages{Journal of the American Statistical
  Association}{119}{}{1-27}.
\PrintBackRefs{\CurrentBib}

\bibitem [\protect \citeauthoryear {%
Zhu%
\ \BBA {} Bradic%
}{%
Zhu%
\ \BBA {} Bradic%
}{%
{\protect \APACyear {2018}}%
}]{%
zhu2018linear}
\APACinsertmetastar {%
zhu2018linear}%
\begin{APACrefauthors}%
Zhu, Y.%
\BCBT {}\ \BBA {} Bradic, J.%
\end{APACrefauthors}%
\unskip\
\newblock
\APACrefYearMonthDay{2018}{}{}.
\newblock
{\BBOQ}\APACrefatitle {Linear hypothesis testing in dense high-dimensional
  linear models} {Linear hypothesis testing in dense high-dimensional linear
  models}.{\BBCQ}
\newblock
\APACjournalVolNumPages{Journal of the American Statistical
  Association}{113}{}{1583--1600}.
\PrintBackRefs{\CurrentBib}

\end{thebibliography}


\begin{thebibliography}{}

\bibitem [\protect \citeauthoryear {%
G{\"o}tze%
, Sambale%
\BCBL {}\ \BBA {} Sinulis%
}{%
G{\"o}tze%
\ \protect \BOthers {.}}{%
{\protect \APACyear {2021}}%
}]{%
alpha_sub_exponential}
\APACinsertmetastar {%
alpha_sub_exponential}%
\begin{APACrefauthors}%
G{\"o}tze, F.%
, Sambale, H.%
\BCBL {}\ \BBA {} Sinulis, A.%
\end{APACrefauthors}%
\unskip\
\newblock
\APACrefYearMonthDay{2021}{}{}.
\newblock
{\BBOQ}\APACrefatitle {{Concentration inequalities for polynomials in
  $\alpha$-sub-exponential random variables}} {{Concentration inequalities for
  polynomials in $\alpha$-sub-exponential random variables}}.{\BBCQ}
\newblock
\APACjournalVolNumPages{Electronic Journal of Probability}{26}{}{1 -- 22}.
\PrintBackRefs{\CurrentBib}

\bibitem [\protect \citeauthoryear {%
van~de Geer%
, Bühlmann%
, Ritov%
\BCBL {}\ \BBA {} Dezeure%
}{%
van~de Geer%
\ \protect \BOthers {.}}{%
{\protect \APACyear {2014}}%
}]{%
van_de_Geer_2014}
\APACinsertmetastar {%
van_de_Geer_2014}%
\begin{APACrefauthors}%
van~de Geer, S.%
, Bühlmann, P.%
, Ritov, Y.%
\BCBL {}\ \BBA {} Dezeure, R.%
\end{APACrefauthors}%
\unskip\
\newblock
\APACrefYearMonthDay{2014}{}{}.
\newblock
{\BBOQ}\APACrefatitle {On asymptotically optimal confidence regions and tests
  for high-dimensional models} {On asymptotically optimal confidence regions
  and tests for high-dimensional models}.{\BBCQ}
\newblock
\APACjournalVolNumPages{The Annals of Statistics}{42}{}{1166 -- 1202}.
\PrintBackRefs{\CurrentBib}

\bibitem [\protect \citeauthoryear {%
Vershynin%
}{%
Vershynin%
}{%
{\protect \APACyear {2018}}%
}]{%
Vershynin_2018}
\APACinsertmetastar {%
Vershynin_2018}%
\begin{APACrefauthors}%
Vershynin, R.%
\end{APACrefauthors}%
\unskip\
\newblock
\APACrefYear{2018}.
\newblock
\APACrefbtitle {{H}igh-Dimensional {P}robability: An {I}ntroduction with
  {A}pplications in {D}ata {S}cience} {{H}igh-dimensional {P}robability: An
  {I}ntroduction with {A}pplications in {D}ata {S}cience}.
\newblock
\APACaddressPublisher{}{Cambridge University Press}.
\PrintBackRefs{\CurrentBib}

\end{thebibliography}

\clearpage
\setcounter{page}{1}

\setcounter{section}{0}
\setcounter{equation}{0}

\setcounter{figure}{0}
\setcounter{table}{0}

\renewcommand{\theequation}{S\arabic{equation}}
\renewcommand{\thesection}{\Roman{section}}

\renewcommand\thefigure{S\arabic{figure}}
\renewcommand\thetable{S\arabic{table}}

\setcounter{lem}{0}
\renewcommand{\thelem}{S\arabic{lem}}

\newtheorem{lemmain}{Lemma}
\renewcommand{\thelemmain}{\arabic{lemmain}}

\begin{center}
{\Large Supplementary Material for}

{\Large ``Debiased prediction inference with non-sparse loadings in misspecified high-dimensional regression models''}

\vspace{.1in} {\large Libin Liang\footnotemark[1] and Zhiqiang Tan\footnotemark[1]}

\vspace{.1in}
\today
\end{center}

\footnotetext[1]{Department of Statistics, Rutgers University. Address: 110 Frelinghuysen Road,
Piscataway, NJ 08854. E-mails: ll866@stat.rutgers.edu, ztan@stat.rutgers.edu.}

\tableofcontents

\addtocontents{toc}{\protect\setcounter{tocdepth}{2}}

\vspace*{1in}

\setcounter{section}{0}
\section{Definitions}

\subsection{Definition of sub-gaussian random variable and vector}

A random variable $z\in \mathbb{R}$ is sub-gaussian if it has a finite sub-gaussian norm defined as
\begin{align*}
    \Vert z\Vert_{\psi_2} = \inf\{t>0:\mathbb{E}\exp(z^2/t^2)\leq 2\}.
\end{align*}
A random vector $Z\in \mathbb{R}^p$ is sub-gaussian if it has a finite sub-gaussian vector norm defined as
\begin{align*}
    \Vert Z\Vert_{\psi_2} = \sup_{s\neq 0}\Vert \frac{ \langle s, Z\rangle}{\Vert s\Vert} \Vert_{\psi_2}.
\end{align*}

\subsection{Definition of sub-exponential random variable}

A random variable $z\in \mathbb{R}$ is sub-exponential if it has a finite sub-exponential norm defined as
\begin{align*}
    \Vert z\Vert_{\psi_1} = \inf\{t>0:\mathbb{E}\exp(|z|/t)\leq 2\}.
\end{align*}

\subsection{Definition of $\alpha$-sub-exponential random variable}

A random variable $z\in \mathbb{R}$ is $\alpha$-sub-exponential where $0\leq\alpha \leq 1$ if it has a finite $\alpha$-sub-exponential norm defined as
\begin{align*}
    \Vert z\Vert_{\psi_\alpha} = \inf\{t>0:\mathbb{E}\exp(|z|^\alpha/t^\alpha)\leq 2\}.
\end{align*}

\section{Proofs of theoretical results}

We mainly give proofs of the theoretical results of the debiased estimator and the estimator of the asymptotic variance(as outlined in main paper Section \ref{sec:outline_of_proof}) in the high-dimensional GLM setting. The theoretical results of LM setting can be regarded as a special case of GLM setting.

\subsection{Initial theoretical result in GLM setting}
\label{sec:theoretical_tilde_M}

We give the proof of the initial theoretical result (Theorem \ref{thm:initial_thm}), depending on certain estimators of $M^*$.

\initialthm*

\begin{prf}
We first give the bound of $\Delta_1$ in (\ref{eq:transform_glm_1}) and then we give the bound of $\Delta_2$ in (\ref{eq:transform_glm_2}).
\noindent \textbf{(a)}. Bound of $\Delta_1$ in (\ref{eq:transform_glm_1}).
\begin{align*}
    &\widehat{\xi^\T\beta} - \xi^\T\bar\beta \\
    &= -\xi^\T M^*\frac{1}{n}\sum_{i=1}^n X_i\phi_1(Y_i, X_i^\T\bar\beta) \notag \\
    &-\xi^\T M^* \frac{1}{n}\sum_{i=1}^n X_i (\phi_1(Y_i, X_i^\T\hat{\beta}) - \phi_1(Y_i,X_i^\T \bar\beta)) + \xi^\T(\hat{\beta}-\bar\beta) \notag \\
    &-\xi^\T (\hat{M}^{J^\c}- M^*)\frac{1}{n}\sum_{i\in J}X_i\phi_1(Y_i, X_i^\T\bar\beta) -\xi^\T (\hat{M}^{J^\c}- M^*)\frac{1}{n}\sum_{i\in J}X_i(\phi_1(Y_i, X_i^\T\hat{\beta}) -  \phi_1(Y_i, X_i^\T\bar\beta) )  \notag \\
   &-\xi^\T (\hat{M}^{J}- M^*)\frac{1}{n}\sum_{i\in J^\c }X_i\phi_1(Y_i, X_i^\T\bar\beta)  -\xi^\T (\hat{M}^{J}- M^*)\frac{1}{n}\sum_{i\in J^\c }X_i(\phi_1(Y_i, X_i^\T\hat{\beta}) -  \phi_1(Y_i, X_i^\T\bar\beta) ).  \notag
\end{align*}
Hence we have
\begin{align*}
    &\widehat{\xi^\T\beta} - \xi^\T\bar\beta \\
    &= -\xi^\T M^*\frac{1}{n}\sum_{i=1}^n X_i\phi_1(Y_i, X_i^\T\bar\beta) -\underbrace{(\frac{1}{n}\sum_{i=1}^n \phi_2(Y_i, X_i^\T\bar\beta) X_iX_i^\T M^*\xi-\xi)^\T(\bar\beta-\hat{\beta})}_{I}\notag \\
    &-\underbrace{\xi^\T (\hat{M}^{J^\c}- M^*)\frac{1}{n}\sum_{i\in J}X_i\phi_1(Y_i, X_i^\T\bar\beta)}_{II} -\underbrace{\xi^\T (\hat{M}^{J^\c}- M^*)\frac{1}{n}\sum_{i\in J^\c }X_i(\phi_1(Y_i, X_i^\T\hat{\beta}) -  \phi_1(Y_i, X_i^\T\bar\beta) )}_{III}  \notag \\
   &-\underbrace{\xi^\T (\hat{M}^{J}- M^*)\frac{1}{n}\sum_{i\in J^\c }X_i\phi_1(Y_i, X_i^\T\bar\beta)}_{IV}  -\underbrace{\xi^\T (\hat{M}^{J}- M^*)\frac{1}{n}\sum_{i\in J^\c }X_i(\phi_1(Y_i, X_i^\T\hat{\beta}) -  \phi_1(Y_i, X_i^\T\bar\beta) )}_{V} \notag  \\
    &+ R.
\end{align*}
where
$R = -\sum_{i=1}^n \frac{1}{n} (\xi^\T M^* X_i)\int_0^1 (1-t)\phi_3(Y_i,X_i^\T\hat{\beta}+t X_i^\T(\hat{\beta}-\bar\beta))dt  (X_i^\T(\hat{\beta}-\bar\beta))^2$.

We discuss how to control the following terms, $I$, $II$, $III$ and $R$. The term $IV$ can be controlled similarly as the term $II$, the term $V$ can be controlled similarly as term $III$.

\noindent \textcircled{1} Term $I$: We have
\begin{align}
   |I| &= |(\frac{1}{n}\sum_{i=1}^n \phi_2(Y_i, X_i^\T\bar\beta) X_iX_i^\T M^*\xi-\xi)^\T(\bar\beta-\hat{\beta}) | \notag\\
   &\leq \Vert (\frac{1}{n}\sum_{i=1}^n \phi_2(Y_i, X_i^\T\bar\beta) X_iX_i^\T M^* - I_p)\xi \Vert_\infty \Vert \bar\beta-\hat{\beta}\Vert_1 . \notag %\label{inq:lem_1_3}
\end{align}
Under the event $\mathcal{A}_{5}$ (Lemma \ref{lem:bound_of_random_quan}(v)), we have
\begin{align}
    \Vert (\frac{1}{n}\sum_{i=1}^n \phi_2(Y_i, X_i^\T\bar\beta) X_iX_i^\T M^* - I)\xi \Vert_\infty \leq  \sqrt{\frac{2}{c(1)}}K(1)(1+(d(1) \log(2))^{-1}) C_6\sigma^2_1 \sqrt{\frac{\log(p\vee n)}{n}}\Vert \xi\Vert_2. \notag %\label{inq:lem_1_4}
\end{align}
Then given Condition \ref{beta_hat_condition}, we have
\begin{align}
   |I|= |(\frac{1}{n}\sum_{i=1}^n \phi_2(Y_i, X_i^\T\bar\beta) X_iX_i^\T M^*\xi-\xi)^\T(\bar\beta-\hat{\beta})| \leq   Cs_0\frac{\log(p\vee n)}{n}\Vert  \xi\Vert_2  , \label{inq:part_5}
\end{align}
where $C$ is a constant depending on $(C_1, C_6, \sigma_1)$ only.

\noindent \textcircled{2} Term $II$: We have
\begin{align*}
    II=\xi^\T (\hat{M}^{J^\c}- M^*)\frac{1}{n}\sum_{i\in J}X_i\phi_1(Y_i, X_i^\T\bar\beta)=\frac{|J|}{n} \frac{1}{|J|}\sum_{i\in J} \xi^\T (\hat{M}^{J^\c}- M^*) X_i \phi_1(Y_i,X_i^\T\bar\beta).
\end{align*}
Under the event $\mathcal{A}_4$ (Lemma \ref{lem:bound_of_random_quan}(iv)), we have
\begin{align}
   |II| = |\xi^\T (\hat{M}^{J^\c}- M^*)\frac{1}{n}\sum_{i\in J}X_i\phi_1(Y_i, X_i^\T\bar\beta)| &\leq \Vert (\hat{M}^{J^\c}- M^*)\xi\Vert_1\Vert \frac{1}{n}\sum_{i\in J}X_i\phi_1(Y_i, X_i^\T\bar\beta)\Vert_\infty   \notag \\
   &\leq \rho \Vert (\hat{M}^{J^\c}- M^*)\xi\Vert_1\Vert \frac{1}{|J|}\sum_{i\in J}X_i\phi_1(Y_i, X_i^\T\bar\beta)\Vert_\infty   \notag \\
   &\leq C  \Vert (\hat{M}^{J^\c}- M^*)\xi\Vert_1 \sqrt{\frac{\log(p\vee n)}{n}}, \label{inq:lem_1_part_1}
\end{align}
where $C$ is a constant depending on $(\sigma_1,\sigma_2, \rho)$ only.

\noindent If under the event $\mathcal{B}_{1}$ (Lemma \ref{lem:bound_of_random_quan}(x)), we have
\begin{align}
   |II| = |\xi^\T (\hat{M}^{J^\c}- M^*)\frac{1}{n}\sum_{i\in J}X_i\phi_1(Y_i, X_i^\T\bar\beta)|
    \leq C \Vert (\hat{M}^{J^\c}-M) \xi\Vert_2 \sqrt{\frac{\log(p\vee n)}{n}}, \label{inq:lem_1_part_1_2}
\end{align}
where $C$ is a constant depending on $(\sigma_1,\sigma_2, \rho)$.

\noindent \textcircled{3} Term $III$: Given Assumption \ref{bounded_higher_differentiate}, we have
\begin{align}
   | III| \\
   &=|\xi^\T (\hat{M}^{J^\c}- M^*)\frac{1}{n}\sum_{i\in J}(\phi_1(Y_i, X_i^\T\hat{\beta}) -  \phi_1(Y_i, X_i^\T\bar\beta) )| \notag \\
   &\leq   \sqrt{\xi^\T (\hat{M}^{J^\c}- M^*)\sum_{i\in J}\frac{1}{n}X_iX_i^\T(\hat{M}^{J^\c}- M^*)\xi}
    \sqrt{\frac{1}{n}\sum_{i\in J}(\phi_1(Y_i, X_i^\T\hat{\beta}) -  \phi_1(Y_i, X_i^\T\bar\beta) )^2} \notag  \\
    &\leq   \sqrt{\xi^\T (\hat{M}^{J^\c}- M^*)\sum_{i\in J}\frac{1}{n}X_iX_i^\T(\hat{M}^{J^\c}- M^*)\xi}
    \sqrt{\frac{1}{n}\sum_{i=1}^n (\phi_1(Y_i, X_i^\T\hat{\beta}) -  \phi_1(Y_i, X_i^\T\bar\beta) )^2} \notag  \\
    &\leq  C_6  \sqrt{\xi^\T (\hat{M}^{J^\c}- M^*)\sum_{i\in J}\frac{1}{n}X_iX_i^\T(\hat{M}^{J^\c}- M^*)\xi}
    \sqrt{(\hat{\beta}-\bar\beta)^\T\sum_{i=1}^n\frac{1}{n}X_iX_i^\T(\hat{\beta}-\bar\beta)} .\notag
\end{align}
Under Condition \ref{beta_hat_condition}, we have
\begin{align}
    |III| \leq C\sqrt{Q_1} \sqrt{s_0\frac{\log(p\vee n)}{n}} ,\label{inq:part_2}
\end{align}
where $C$ is a constant depending on $(C_1,C_6)$.

\noindent \textcircled{4} Term $R$: Given Assumption \ref{bounded_higher_differentiate}, we have
\begin{align*}
   |R| &= |\sum_{i=1}^n \frac{1}{n} (\xi^\T M^* X_i)\int_0^1 (1-t)\phi_3(Y_i,X_i^\T\bar\beta+t X_i^\T(\hat{\beta}-\bar\beta))dt * (X_i^\T(\hat{\beta}-\bar\beta))^2| \notag \\
   &\leq C_7 \frac{1}{n}\sum_{i=1}^n   |\xi^\T M^* X_i| (X_i^\T(\hat{\beta}-\bar\beta))^2      \\
   &\leq C_7 \max_{1\leq i\leq n}\{|\xi^\T M^* X_i|\} \frac{1}{n} \sum_{i=1}^n (X_i^\T(\hat{\beta}-\bar\beta))^2.
\end{align*}
Under the event $\mathcal{A}_{62}$ (Lemma \ref{lem:bound_of_random_quan}(vi)) and given Condition \ref{beta_hat_condition}, we have
\begin{align}
    |R|&=|\sum_{i=1}^n \frac{1}{n} (\xi^\T M^* X_i)\int_0^1 (1-t)\phi_3(Y_i,X_i^\T\bar\beta+t X_i^\T(\hat{\beta}-\bar\beta))dt * (X_i^\T(\hat{\beta}-\bar\beta))^2| \notag \\
    &\leq  C_7 C  s_0 \frac{\log(p\vee n)\sqrt{\log(n)}}{n} \Vert \xi\Vert_2,  \label{inq:part_R}
\end{align}
where $C$ is a constant depending on $(C_1, C_4, \sigma_1)$.

\noindent From Assumption \ref{ass:H_eigen_value}(i) and \ref{ass:residulal_glm_lower}, we have
\begin{align}
    V &= \mathbb{E}[(\xi^\T M^* X_i)^2\phi_1^2(Y_i, X_i^\T\bar\beta)]  \notag\\
    &\geq C_3C_5^2 \Vert \xi\Vert_2^2. \label{inq:lower_bound_V}
\end{align}
$\newline$
\noindent From (\ref{inq:part_5}),(\ref{inq:lem_1_part_1}), (\ref{inq:part_2}), (\ref{inq:part_R}) and (\ref{inq:lower_bound_V}), suppose that Assumption \ref{ass:sub_gaussian_entries_GLM}, \ref{ass:H_eigen_value},  \ref{ass:residulal_glm_lower} and \ref{bounded_higher_differentiate} and Condition \ref{beta_hat_condition} is satisfied and $\frac{\log(p\vee n)}{\sqrt{n} }\leq \frac{c(1)\rho \sqrt{n}}{2}$ and $\log(n)>1$, we have with probability at least $1-\delta_1(n)-\frac{6}{p \vee n}-\frac{2}{n}$,
\begin{align*}
   |\Delta_1|& \leq   C (\sqrt{\log(p\vee n)}\frac{\Vert (\hat{M}^{J^\c}- M^*)\xi\Vert_1 + \Vert (\hat{M}^{J}- M^*)\xi\Vert_1}{\Vert\xi\Vert_2} +  \frac{\sqrt{Q_1}+\sqrt{Q_1^\prime}}{\Vert\xi\Vert_2} \sqrt{s_0\log(p \vee n)} \\ &+s_0\frac{\log(p\vee n)}{\sqrt{n}} +C_7   s_0 \frac{\log(p\vee n)\sqrt{\log(n)} }{\sqrt{n} } ),
\end{align*}
where $C$ is a constant depending on $(C_1, C_2, C_3, C_4, C_5, C_6, \sigma_1,\sigma_2,\rho)$ only. This gives the bound of $\Delta_1$ in Theorem \ref{thm:initial_thm} (i).

$\newline$
\noindent From (\ref{inq:part_5}),(\ref{inq:lem_1_part_1_2}), (\ref{inq:part_2}), (\ref{inq:part_R}) and (\ref{inq:lower_bound_V}), suppose that Assumption \ref{ass:sub_gaussian_GLM}, \ref{ass:H_eigen_value} \ref{ass:residulal_glm_lower} and \ref{bounded_higher_differentiate} and Condition \ref{beta_hat_condition} is satisfied and $\frac{\log(p\vee n)}{\sqrt{n} }\leq\min\{c(1)\rho \sqrt{n}, \frac{c(1)\sqrt{n}}{2}\}$ and $\log(n)>1$, we have with probability at least $1-\delta_1(n)-\frac{6}{p \vee n}-\frac{2}{n}$,
\begin{align*}
   |\Delta_1|& \leq   C (\sqrt{\log(p\vee n)}\frac{\Vert (\hat{M}^{J^\c}- M^*)\xi\Vert_2 + \Vert (\hat{M}^{J}- M^*)\xi\Vert_2}{\Vert\xi\Vert_2}+  \frac{\sqrt{Q_1}+\sqrt{Q_1^\prime}}{\Vert\xi\Vert_2} \sqrt{s_0\log(p\vee n)} \\&+s_0\frac{\log(p\vee n)}{\sqrt{n}}
   +C_7  s_0 \frac{\log(p\vee n)\sqrt{\log(n)}}{\sqrt{n} } ),
\end{align*}
where $C$ is a constant depending on $(C_1, C_2,C_3, C_4, C_5, C_6, \sigma_1,\sigma_2,\rho)$. This gives the bound of $\Delta_1$ in Theorem \ref{thm:initial_thm} (ii).

$\newline$
\noindent \textbf{(b)}. Bound of $\Delta_2$ in (\ref{eq:transform_glm_2}).

With (\ref{eq:V_estimator}), we have
\begin{align*}
    &\hat{V} - V \\
    &= \frac{1}{n} \sum_{i=1}^n ((\xi^\T M^* X_i)^2\phi^2_1(Y_i,X_i^\T\bar\beta)- \mathbb{E}[(\xi^\T M^* X_i)^2\phi^2_1(Y_i,X_i^\T\bar\beta)]) \\
    &+\frac{1}{n}\sum_{i=1}^n (\xi^\T M^* X_i)^2(2\phi_1(Y_i,X_i^\T\bar\beta)(\phi_1(Y_i,X_i^\T\bar\beta)-\phi_1(Y_i,X_i^\T\hat{\beta})) + (\phi_1(Y_i,X_i^\T\bar\beta)-\phi_1(Y_i,X_i^\T\hat{\beta}))^2 ) \\
    &+\frac{2}{n}\sum_{i\in J} (\xi^\T M^* X_i)( \xi^\T(\hat{M}^{J^\c}-M^*)X_i)\phi^2_1(Y_i,X_i^\T\hat{\beta})+\frac{1}{n}\sum_{i\in J} (\xi^\T(\hat{M}^{J^\c}-M^*)X_i )^2\phi^2_1(Y_i,X_i^\T\hat{\beta}) \\
    &+\frac{2}{n}\sum_{i\in {J^C}} (\xi^\T M^* X_i)( \xi^\T (\hat{M}^J - M^*)X_i)(\phi^2_1(Y_i,X_i^\T\hat{\beta})
  +\frac{1}{n}\sum_{i\in {J^C}} (\xi^\T(\hat{M}^J - M^*)X_i)^2\phi^2_1(Y_i,X_i^\T\hat{\beta}).
\end{align*}
Hence we have
\begin{align*}
    & |\tilde{V} - V| \\
     &\leq \underbrace{|\frac{1}{n} \sum_{i=1}^n ((\xi^\T M^* X_i)^2\phi^2_1(Y_i,X_i^\T\bar\beta)- \mathbb{E}[(\xi^\T M^* X_i)^2\phi^2_1(Y_i,X_i^\T\bar\beta)]) |}_{I} \\
    &+\underbrace{|\frac{2}{n}\sum_{i=1}^n (\xi^\T M^* X_i)^2\phi_1(Y_i,X_i^\T\bar\beta)(\phi_1(Y_i,X_i^\T\bar\beta)-\phi_1(Y_i,X_i^\T\hat{\beta}))|}_{II} \\
    &+ \underbrace{|\frac{1}{n}\sum_{i=1}^n (\xi^\T M^* X_i)^2(\phi_1(Y_i,X_i^\T\bar\beta)-\phi_1(Y_i,X_i^\T\hat{\beta}))^2 )|}_{III} \\
    &+\underbrace{\frac{4}{n}\sum_{i\in J} |(\xi^\T M^* X_i)( \xi^\T(\hat{M}^{J^\c}-M^*)X_i)|\phi^2_1(Y_i,X_i^\T\bar\beta)}_{IV} \\
    &+ \underbrace{\frac{4}{n}\sum_{i\in J} |(\xi^\T M^* X_i)( \xi^\T(\hat{M}^{J^\c}-M^*)X_i)(\phi_1(Y_i,X_i^\T\bar\beta)-\phi_1(Y_i,X_i^\T\hat{\beta}))^2}_{V} \notag \\
    &+\underbrace{\frac{2}{n}\sum_{i\in J} (\xi^\T(\hat{M}^{J^\c}-M^*)X_i )^2\phi^2_1(Y_i,X_i^\T\bar\beta)}_{VI} \\
    &+ \underbrace{\frac{2}{n}\sum_{i\in J} (\xi^\T(\hat{M}^{J^\c}-M^*)X_i )^2(\phi_1(Y_i,X_i^\T\bar\beta)-\phi_1(Y_i,X_i^\T\hat{\beta}))^2}_{VII} \\
    &+VIII+IX+X+XI,
\end{align*}
where $VIII=\frac{4}{n}\sum_{i\in J^\c} |(\xi^\T M^* X_i)( \xi^\T(\hat{M}^{J}-M^*)X_i)|\phi^2_1(Y_i,X_i^\T\bar\beta)$, \\
$IX=\frac{4}{n}\sum_{i\in J^\c} |(\xi^\T M^* X_i)( \xi^\T(\hat{M}^{J}-M^*)X_i)|(\phi_1(Y_i,X_i^\T\bar\beta)-\phi_1(Y_i,X_i^\T\hat{\beta}))^2$, \\
$X=\frac{2}{n}\sum_{i\in J^\c} (\xi^\T(\hat{M}^{J}-M^*)X_i )^2\phi^2_1(Y_i,X_i^\T\bar\beta)$ \\
and $XI=\frac{2}{n}\sum_{i\in J^\c} (\xi^\T(\hat{M}^{J}-M^*)X_i )^2(\phi_1(Y_i,X_i^\T\bar\beta)-\phi_1(Y_i,X_i^\T\hat{\beta}))^2$.

We discuss how to control the following terms, $I$, $II$, $III$, $IV$, $V$, $VI$ and $VII$. The term $VIII$ can be controlled similarly as the term $IV$. The term $IX$ can be controlled similarly as the term $V$. The term $X$ can be controlled similarly as term $VI$ and the term $XI$ can be controlled similarly as the term $VII$.

\noindent \textcircled{1} Term I: Under event $\mathcal{A}_{71}$ (Lemma \ref{lem:bound_of_random_quan}(vii)), we have
\begin{align}
  |I|= | \frac{1}{n}\sum_{i=1}^n(\xi^\T M^* X_i)^2\phi^2_1(Y_i,X_i^\T\bar\beta) -\mathbb{E}[\Vert (\xi^\T M^* X_i)^2\phi^2_1(Y_i,X_i^\T\bar\beta) ]|\leq C n^{-\frac{1}{3}} \Vert \xi\Vert^2_2 , \label{inq:bound_var_part_I}
\end{align}
where $C$ is a constant depending on $(C_4,\sigma_1,\sigma_2)$ only.

\noindent \textcircled{2} Term II: Under Assumption \ref{bounded_higher_differentiate}, we have
\begin{align}
   |II| &=|\frac{2}{n}\sum_{i=1}^n (\xi^\T M^* X_i)^2\phi_1(Y_i,X_i^\T\bar\beta)(\phi_1(Y_i,X_i^\T\bar\beta)-\phi_1(Y_i,X_i^\T\hat{\beta}))| \notag \\
    &\leq 2 \sqrt{\frac{1}{n}\sum_{i=1}^n (\xi^\T M^* X_i)^4}\sqrt{\frac{1}{n}\sum_{i=1}^n \phi^2_1(Y_i,X_i^\T\bar\beta)(\phi_1(Y_i,X_i^\T\bar\beta)-\phi_1(Y_i,X_i^\T\hat{\beta}))^2}| \notag \\
    &\leq 2 \max_{1\leq i\leq n}\{|\phi_1(Y_i,X_i^\T\bar\beta)|\}\sqrt{\frac{1}{n}\sum_{i=1}^n (\xi^\T M^* X_i)^4}\sqrt{\frac{1}{n}\sum_{i=1}^n (\phi_1(Y_i,X_i^\T\bar\beta)-\phi_1(Y_i,X_i^\T\hat{\beta}))^2}| \notag \\
    &\leq  2 C_6 \max_{1\leq i\leq n}\{|\phi_1(Y_i,X_i^\T\bar\beta)|\}\sqrt{\frac{1}{n}\sum_{i=1}^n (\xi^\T M^* X_i)^4}\sqrt{\frac{1}{n}\sum_{i=1}^n (X_i^\T\bar\beta-X_i^\T\hat{\beta})^2} . \notag %\label{inq:bounded_part_II}
\end{align}
Under event $\mathcal{A}_{24}\cap\mathcal{A}_{72}$ (Lemma \ref{lem:bound_of_random_quan}(ii), (vii)) and with Condition \ref{beta_hat_condition}, we have
\begin{align}
    |II|&=|\frac{2}{n}\sum_{i=1}^n (\xi^\T M^* X_i)^2\phi_1(Y_i,X_i^\T\bar\beta)(\phi_1(Y_i,X_i^\T\bar\beta)-\phi_1(Y_i,X_i^\T\hat{\beta}))| \notag \\
    &\leq C \sqrt{s_0} \sqrt{\frac{\log(n) \log(p\vee n)}{n}}\Vert  \xi\Vert^2_2  ,\label{inq:bound_var_part_II}
\end{align}
where $C$ is a constant depending on $(C_1,C_4,C_6,\sigma_1,\sigma_2)$ only.

\noindent \textcircled{3} Term III: Under Assumption \ref{bounded_higher_differentiate}, we have
\begin{align}
   |III| &= |\frac{1}{n}\sum_{i=1}^n (\xi^\T M^* X_i)^2(\phi_1(Y_i,X_i^\T\bar\beta)-\phi_1(Y_i,X_i^\T\hat{\beta}))^2 )| \notag \\
   &\leq  C^2_6 (\max_{1\leq i\leq p}  |\xi^\T M^* X_i|)^2 \frac{1}{n}\sum_{i=1}^n(X_i^\T\bar\beta-X_i^\T\hat{\beta})^2 . \notag %\label{inq:bound_part_III}
\end{align}
Given Condition \ref{beta_hat_condition} and under event $\mathcal{A}_{62}$ (Lemma \ref{lem:bound_of_random_quan}(vi)), we have
\begin{align}
  |III| = |\frac{1}{n}\sum_{i=1}^n (\xi^\T M^* X_i)^2(\phi_1(Y_i,X_i^\T\bar\beta)-\phi_1(Y_i,X_i^\T\hat{\beta}))^2 )| \leq C  s_0 \frac{\log(p\vee n)\log(n)}{n}\Vert \xi\Vert^2_2,  \label{inq:bound_var_part_III}
\end{align}
where $C$ is a constant depending on $(C_1, C_4, C_6,\sigma_1)$ only.

\noindent \textcircled{4} Term IV: We have
\begin{align}
    |IV|&\leq \frac{4}{n}\sum_{i\in J} |(\xi^\T M^* X_i)( \xi^\T(\hat{M}^{J^\c}-M^*)X_i)|\phi^2_1(Y_i,X_i^\T\bar\beta) \notag \\
    &\leq  4\sqrt{\frac{1}{|J|}\sum_{i\in J} (\xi^\T M^* X_i)^2 (\xi^\T (\hat{M}^{J^\c}-M^*)X_i)^2 }\sqrt{\frac{1}{|J|}\sum_{i\in J}\phi_1^4(Y_i,X_i^\T\bar\beta)} \notag \\
    &\leq  4 \max_{1\leq i\leq n}|\xi^\T M^* X_i|\sqrt{\frac{1}{|J|}\sum_{i\in J} (\xi^\T (\hat{M}^{J^\c}-M^*)X_i)^2 }\sqrt{\frac{1}{|J|}\sum_{i\in J}\phi_1^4(Y_i,X_i^\T\bar\beta)}\notag \\
    &\leq \frac{4}{\sqrt{\rho}} \max_{1\leq i\leq n}|\xi^\T M^* X_i|\sqrt{\frac{1}{|J|}\sum_{i\in J} (\xi^\T (\hat{M}^{J^\c}-M^*)X_i)^2 }\sqrt{\frac{1}{n}\sum_{i=1 }^n \phi_1^4(Y_i,X_i^\T\bar\beta)}. \notag %\label{inq:bounded_part_IV}
\end{align}Under event $\mathcal{A}_{3}\cap \mathcal{A}_{62}$ (Lemma \ref{lem:bound_of_random_quan}(iii), (vi)), we have
\begin{align}
    |IV| &\leq \frac{4}{n}\sum_{i\in J} |(\xi^\T M^* X_i)( \xi^\T(\hat{M}^{J^\c}-M^*)X_i)|\phi^2_1(Y_i,X_i^\T\bar\beta) \notag \\
     &\leq
   C \sqrt{\log(n)} \sqrt{Q_1} ,  \label{inq:bound_var_part_IV}
\end{align}
where $C$ is a constant depending on $(C_4, \sigma_2,\rho)$ only.

\noindent \textcircled{5} Term V: Given Assumption \ref{bounded_higher_differentiate}, we have
\begin{align}
   |V| &\leq \frac{4}{n}\sum_{i\in J} |(\xi^\T M^* X_i)( \xi^\T(\hat{M}^{J^\c}-M^*)X_i)|(\phi_1(Y_i,X_i^\T\bar\beta)-\phi_1(Y_i,X_i^\T\hat{\beta}))^2 \notag\\
   &\leq C_6 \frac{4}{n}\sum_{i\in J} |(\xi^\T M^* X_i)( \xi^\T(\hat{M}^{J^\c}-M^*)X_i)|(X_i^\T\bar\beta-X_i^\T\hat{\beta})^2 \notag \\
   &\leq 4 C_6 \Vert \bar\beta-\hat{\beta} \Vert^2_1 \max_{1\leq l \leq p, 1\leq m\leq p}\{\frac{1}{n}\sum_{i\in J}| (\xi^\T(\hat{M}^{J^\c}-M^*)X_i)(\xi^\T M^* X_i)X_{il} X_{im}| \}.
    \label{inq:bounded_part_V}
\end{align}

\noindent

\noindent Under event $\mathcal{A}_{22}\cap\mathcal{A}_{62}$ (Lemma \ref{lem:bound_of_random_quan}(ii), (vi)) and given Condition \ref{beta_hat_condition}, from (\ref{inq:bounded_part_V}), we have
\begin{align}
   &|V|\notag\\
   &\leq  4 C_6 \Vert \bar\beta-\hat{\beta} \Vert^2_1 \max_{1\leq l \leq p, 1\leq m\leq p}\{\frac{1}{n}\sum_{i\in J}| (\xi^\T(\hat{M}^{J^\c}-M^*)X_i)(\xi^\T M^* X_i)X_{il} X_{im}| \}\notag \\
   &\leq  4 C_6 \Vert \bar\beta-\hat{\beta} \Vert^2_1 \max_{i=1,\ldots,n} \{|\xi^\T M^* X_i|\}  \max_{1\leq l \leq p, 1\leq m\leq p}\{\frac{1}{n}\sum_{i\in J}| (\xi^\T(\hat{M}^{J^\c}-M^*)X_i)X_{il} X_{im}| \}   \notag \\
   &\leq   4 C_6 \Vert \bar\beta-\hat{\beta} \Vert^2_1 \Vert (\hat{M}^{J^\c}-M^*)\xi\Vert_1  \max_{i=1,\ldots,n} \{|\xi^\T M^* X_i|\}  \max_{1\leq k \leq p,1\leq l \leq p, 1\leq m\leq p}\{\frac{1}{n}\sum_{i\in J}| X_{ik}X_{il} X_{im}| \}      \notag \\
    &\leq  C s^2_0 \frac{\log(p\vee n)\sqrt{\log( n)}}{n} \Vert (\hat{M}^{J^\c}-M^*)\xi\Vert_1 \Vert \xi\Vert_2,\label{inq:bound_var_V_1}
\end{align}
where $C$ is a constant depending on $(C_1,C_4,C_6,\sigma_1)$.

$\newline$
\noindent Under event $\mathcal{A}_{61}\cap\mathcal{B}_{4}$ (Lemma \ref{lem:bound_of_random_quan}(vi), (xiii)), from (\ref{inq:bounded_part_V}), we have
\begin{align}
   |V| &\leq \frac{4}{n}\sum_{i\in J} |(\xi^\T M^* X_i)( \xi^\T(\hat{M}^{J^\c}-M^*)X_i)|(\phi_1(Y_i,X_i^\T\bar\beta)-\phi_1(Y_i,X_i^\T\hat{\beta}))^2 \notag \\
   &\leq    4 C_6 \Vert \bar\beta-\hat{\beta} \Vert^2_1 \max_{1\leq l \leq p, 1\leq m\leq p}\{\frac{1}{n}\sum_{i\in J}| (\xi^\T(\hat{M}^{J^\c}-M^*)X_i)(\xi^\T M^* X_i)X_{il} X_{im}| \}        \notag \\
   &\leq  4 C_6 \Vert \bar\beta-\hat{\beta} \Vert^2_1 \max_{i \in J}\{|\xi^\T(\hat{M}^{J^\c}-M^*)X_i|\}\max_{1\leq l \leq p, 1\leq m\leq p}\{\frac{1}{n}\sum_{i\in J}| (\xi^\T M^* X_i)X_{il} X_{im}| \}    \notag \\
    &\leq  C s^2_0 \frac{\log(p\vee n)\sqrt{\log(n)}}{n}\Vert (\hat{M}^{J^\c}-M^*)\xi\Vert_2 \Vert \xi\Vert_2,\label{inq:bound_var_V_2}
\end{align}
where $C$ is a constant depending on $(C_1,C_4,C_6,\sigma_1)$.

\noindent \textcircled{6} Term VI: We have
\begin{align}
   |VI| &= |\frac{2}{n}\sum_{i\in J} (\xi^\T(\hat{M}^{J^\c}-M^*)X_i )^2\phi^2_1(Y_i,X_i^\T\bar\beta)| \notag \\
   &\leq  2(\max_{1\leq i\leq p}\{|\phi_1(Y_i,X_i^\T \bar\beta)|\})^2 \frac{1}{n}\sum_{i\in J}(\xi^\T(\hat{M}^{J^\c}-M^*)X_i )^2. \notag
\end{align}
Under event $\mathcal{A}_{24}$ (Lemma \ref{lem:bound_of_random_quan}(ii)), we have
\begin{align}
  |VI| =  |\frac{2}{n}\sum_{i\in J} (\xi^\T(\hat{M}^{J^\c}-M^*)X_i )^2\phi^2_1(Y_i,X_i^\T\bar\beta)| \leq  C Q_1 \log(n), \label{inq:bound_var_part_VI}
\end{align}
where $C$ is a constant depending on $\sigma_2$.

\noindent \textcircled{7} Term VII: We have
\begin{align}
  |VII| &= \frac{2}{n}\sum_{i\in J} (\xi^\T(\hat{M}^{J^\c}-M^*)X_i )^2(\phi_1(Y_i,X_i^\T\bar\beta)-\phi_1(Y_i,X_i^\T\hat{\beta}))^2 \notag\\
  &\leq C_6\frac{2}{n}\sum_{i\in J} (\xi^\T(\hat{M}^{J^\c}-M^*)X_i )^2(X_i^\T\bar\beta-X_i^\T\hat{\beta})^2   \notag \\
  &\leq   2 C_6 \Vert\bar\beta-\hat{\beta}\Vert_1^2  \max_{1\leq l\leq p,1\leq m\leq p} \{ \frac{1}{n}\sum_{i\in J} (\xi^\T(\hat{M}^{J^\c}-M^*)X_i )^2 |X_{il}X_{im}|\} . \notag
\end{align}
Under event $\mathcal{A}_{22}$ (Lemma \ref{lem:bound_of_random_quan}(ii)), under Condition \ref{beta_hat_condition}, we have
\begin{align}
  |VII| &\leq  2 C_6 \Vert\bar\beta-\hat{\beta}\Vert_1^2  \max_{1\leq l\leq p,1\leq m\leq p} \{ \frac{1}{n}\sum_{i\in J} (\xi^\T(\hat{M}^{J^\c}-M^*)X_i )^2 |x_{il}x_{im}|\}  \notag\\
  &\leq  2 C_6 \Vert\bar\beta-\hat{\beta}\Vert_1^2   \Vert (\hat{M}^{J^\c}-M^*)\xi \Vert^2_1*\max_{1\leq k\leq p,1\leq l\leq p,1\leq m\leq p,1\leq q\leq p} \{\frac{1}{n}\sum_{i=1}^n |X_{ik} X_{il} X_{im}X_{iq}| \} \notag \\
  &\leq C s^2_0 \frac{\log(p\vee n)  }{n} \Vert (\hat{M}^{J^\c}-M^*)\xi \Vert^2_1, \label{inq:bound_var_part_VII_1}
\end{align}
where $C$ is a constant depending on $(C_1, C_6,\sigma_1)$.

$\newline$
\noindent Under event $\mathcal{B}_{3}\cap\mathcal{B}_{4}$ (Lemma \ref{lem:bound_of_random_quan}(xii)(xiii)), we have
\begin{align}
  &|VII| \notag \\&\leq  2 C_6 \Vert\bar\beta-\hat{\beta}\Vert_1^2  \max_{1\leq l\leq p,1\leq m\leq p} \{ \frac{1}{n}\sum_{i\in J} (\xi^\T(\hat{M}^{J^\c}-M^*)X_i )^2 |X_{il}X_{im}|\},  \notag \\
  &\leq 2 C_6 \Vert\bar\beta-\hat{\beta}\Vert_1^2 \max_{i\in J}\{|\xi^\T(\hat{M}^{J^\c}-M^*)X_i|\} \notag \\
  &\times \max_{1\leq l\leq p, 1\leq m\leq p} \{ \frac{1}{|J|}\sum_{i\in J} | (\xi^\T (\hat{M}^{J^\c}- M^*) X_i)X_{il} X_{im}| \}   \notag \\
  &\leq  C s^2_0 \frac{\log(p\vee n) \sqrt{\log(n)} }{n} \Vert (\hat{M}^{J^\c}-M^*)\xi \Vert^2_2, \label{inq:bound_var_part_VII_2}
\end{align}
where $C$ is a constant depending on $(C_1,C_6,\sigma_1)$.

$\newline$
\noindent From (\ref{inq:lower_bound_V}), (\ref{inq:bound_var_part_I}), (\ref{inq:bound_var_part_II}), (\ref{inq:bound_var_part_III}), (\ref{inq:bound_var_part_IV}), (\ref{inq:bound_var_V_1}), (\ref{inq:bound_var_part_VI}) and (\ref{inq:bound_var_part_VII_1}), suppose that Assumption \ref{ass:sub_gaussian_entries_GLM}, \ref{ass:H_eigen_value}, \ref{ass:residulal_glm_lower} and \ref{bounded_higher_differentiate} and Condition \ref{beta_hat_condition} are satisfied and $\frac{\log(p\vee n)}{\sqrt{n}}\leq\frac{c(\frac{2}{3})n^{\frac{1}{6}}}{4}$ and $\log(n)>1$, we have with probability at least $1-\delta_1(n)-\frac{4}{n}-\frac{2}{p\vee n} - 2exp\{-c(\frac{1}{2})n^{\frac{1}{3}})\}-4exp\{-\frac{1}{2}\sqrt{n}\}$
\begin{align*}
    |\Delta_2| &\leq  C (n^{-\frac{1}{3}}  +  \sqrt{s_0} \sqrt{\frac{\log(n) \log(p\vee n)}{n}}+  s_0 \frac{\log(p\vee n)\log(n)}{n}\\
    &+\frac{\sqrt{\log(n)}(\sqrt{Q_1}+\sqrt{Q_1^\prime})}{\Vert \xi\Vert_2}+ s^2_0 \frac{\log(p\vee n)\sqrt{\log(n)}}{n} \frac{\Vert (\hat{M}^{J^\c}-M^*)\xi\Vert_1+\Vert (\hat{M}^{J}-M^*)\xi\Vert_1}{\Vert \xi\Vert_2} \\&+  \frac{(Q_1 + Q_1^\prime) \log(n)}{\Vert \xi\Vert_2^2}
    + s^2_0  \frac{\log(p\vee n) }{n} \frac{\Vert (\hat{M}^{J^\c}-M^*)\xi\Vert^2_1+\Vert (\hat{M}^{J}-M^*)\xi\Vert^2_1}{\Vert \xi\Vert_2^2}),
\end{align*}
where $C$ is a constant depending on $(C_1,C_3, C_4,C_5,C_6,\sigma_1,\sigma_2,\rho)$. This give the bound of $\Delta_2$ in Theorem \ref{thm:initial_thm} (i).

$\newline$
\noindent From (\ref{inq:lower_bound_V}), (\ref{inq:bound_var_part_I}), (\ref{inq:bound_var_part_II}), (\ref{inq:bound_var_part_III}), (\ref{inq:bound_var_part_IV}), (\ref{inq:bound_var_V_2}), (\ref{inq:bound_var_part_VI}) and (\ref{inq:bound_var_part_VII_2}), suppose that Assumption  \ref{ass:sub_gaussian_GLM}, \ref{ass:H_eigen_value}, \ref{ass:residulal_glm_lower} and \ref{bounded_higher_differentiate} and  Condition \ref{beta_hat_condition} are satisfied and $\frac{\log(p\vee n)}{\sqrt{n}}\leq\frac{c(\frac{2}{3})n^{\frac{1}{6}}\rho^{\frac{2}{3}}}{3}$ and $\log(n)>1$, we have with probability at least $1-\delta_1(n)-\frac{8}{n}-\frac{6}{p\vee n} - 2exp\{-c(\frac{1}{2})n^{\frac{1}{3}}\}-4exp\{-\frac{1}{2}\sqrt{n}\}$
\begin{align*}
    |\Delta_2| &\leq  C (n^{-\frac{1}{3}}  +  \sqrt{s_0} \sqrt{\frac{\log(n) \log(p\vee n)}{n}}+  s_0 \frac{\log(p\vee n)\log(n)}{n} \\
&+\frac{\sqrt{\log(n)}(\sqrt{Q_1}+\sqrt{Q_1^\prime})}{\Vert\xi\Vert_2}+ s^2_0 \frac{\log(p\vee n)\sqrt{\log(n)}}{n}\frac{\Vert (\hat{M}^{J^\c}-M^*)\xi\Vert_2+\Vert (\hat{M}^{J}-M^*)\xi\Vert_2}{\Vert\xi\Vert_2} \\&+  \frac{(Q_1+Q_1^\prime) \log(n)}{\Vert \xi\Vert_2^2} + s^2_0
 \frac{\log(p\vee n) \sqrt{\log(n)} }{n}\frac{\Vert (\hat{M}^{J^\c}-M^*)\xi\Vert^2_2+\Vert (\hat{M}^{J}-M^*)\xi\Vert^2_2}{\Vert\xi\Vert_2^2} ),
\end{align*}
where $C$ is a constant depending on $(C_1,C_3,C_4,C_5,C_6,\sigma_1,\sigma_2,\rho)$. This give the bound of $\Delta_2$ in Theorem \ref{thm:initial_thm} (ii).
\end{prf}

\subsection{Control of $Q_1$ and $Q_1^\prime$  in Theorem \ref{thm:initial_thm}}

We provide the proof of the bounds of $Q_1$ in Theorem \ref{thm:initial_thm}. The bounds of $Q_1^\prime$ can be proved similarly.

\controlqqprime*
\begin{prf}
For $Q_1$, from Assumption \ref{ass:H_eigen_value}(i), we have
\begin{align}
    Q_1
    &= \frac{|J|}{n} \xi^\T (\hat{M}^{J^\c}- M^*)\Sigma(\hat{M}^{J^\c}- M^*)\xi +  \frac{|J|}{n}\xi^\T (\hat{M}^{J^\c}- M^*)(\sum_{i\in J}\frac{1}{|J|}X_iX_i^\T -\Sigma)(\hat{M}^{J^\c}- M^*)\xi \notag \\
    &\leq  \frac{|J|}{n} \frac{1}{C_3}\xi^\T (\hat{M}^{J^\c}- M^*)(\hat{M}^{J^\c}- M^*)\xi  + \Vert  \sum_{i\in J}\frac{1}{|J|}X_iX_i^\T - \Sigma\Vert_{\max} \Vert (\hat{M}^{J^\c}- M)_{S_{\xi},S_{\xi}}\Vert_1^2 \Vert\xi\Vert_1^2  \notag \\
    &\leq  \frac{|J|}{n} \frac{1}{C_3} \Vert (\hat{M}^{J^\c}- M^*)\xi\Vert_2^2 + \Vert  \sum_{i\in J}\frac{1}{|J|}X_iX_i^\T - \Sigma\Vert_{\max} \Vert (\hat{M}^{J^\c}- M)\xi\Vert_1^2 .\label{inq:cross_term_diff_4}
\end{align}

From (\ref{inq:cross_term_diff_4}), under the event $\mathcal{A}_{14}$ (Lemma \ref{lem:bound_of_random_quan}(i)), we have
\begin{align}
  Q_1 &=\xi^\T (\hat{M}^{J^\c}- M^*)\sum_{i\in J}\frac{1}{n}X_iX_i^\T(\hat{M}^{J^\c}- M^*)\xi \notag \\
   &\leq  C (\Vert (\hat{M}^{J^\c}- M^*)\xi\Vert_2^2 +  \sqrt{\frac{\log(p\vee n)}{n}}\Vert (\hat{M}^{J^\c}- M^*)\xi\Vert_1^2 )    ,\label{inq:corss_term_bound(1)}
\end{align}
where $C$ is a constant depending on $(C_3, \sigma_1, \rho)$ only.

\noindent If we further assume Assumption \ref{ass:sub_gaussian_GLM}, under the event $\mathcal{B}_{2}$ (Lemma \ref{lem:bound_of_random_quan}(xi)), we have
\begin{align}
  Q_1 &\leq \xi^\T (\hat{M}^{J^\c}- M^*)\hat{\Sigma}_J(\hat{M}^{J^\c}- M^*)\xi  \notag \\&\leq C \Vert (\hat{M}^{J^\c}- M^*)\xi\Vert_2^2   ,\label{inq:cross_term_diff_2_1}
\end{align}
where $C$ is a constant depending on $\sigma_1$ only.

\end{prf}

\subsection{Error bounds of estimators in Algorithms \ref{alg:omega_est} and \ref{alg:est_M}}

We give errors bounds for the estimators in Algorithms \ref{alg:omega_est} and \ref{alg:est_M}. We first give the following lemma which is important in the following analysis. For any symmetric $W\in \mathbb{R}^{p \times p}$, $\rho_n>0$, we consider the  estimator $\hat{\Gamma}\in\mathbb{R}^{p\times p}$ as below,
\begin{align}
    \hat{G} = &\argmin_{\text{$G\in\mathbb{R}^{p\times p}$}}  \sum_{j=1}^p \Vert G_{\cdot j}\Vert_{1} \notag \\
    &\text{s.t. } \Vert W G - I_p\Vert_{\max} \leq \rho_n, \notag \\
   \hat{\Gamma}=&(\hat{\Gamma}_{ij})\quad \text{where $\hat{\Gamma}_{ij}=\hat{G}_{ij}\textbf{1}_{[|\hat{G}_{ij}| \leq |\hat{G}_{ji}|]} + \hat{G}_{ji}\textbf{1}_{[|\hat{G}_{ij}| > |\hat{G}_{ji}|]}$ for $1\leq i,j\leq p$. }  \label{def:est_gamma_clime}
\end{align}

\begin{lem}
\label{lem:clime_est_proof}
For any symmetric matrix $W\in \mathbb{R}^{p \times p}$, non-singular $\Gamma \in \mathbb{R}^{p\times p}$, $L>0$ and $\rho_n>0$, if
\begin{align}
    \Vert W \Gamma - I \Vert_{\max} &\leq \rho_n ,  \label{cstr:M_max}\\
    \Vert W - \Gamma^{-1}\Vert_{\max}  &\leq  \rho_n, \label{cstr:M_inv}\\
    \Vert \Gamma\Vert_{1}&\leq L, \label{cstr:M_L1}
\end{align}
then we have for $\hat{\Gamma}$ in (\ref{def:est_gamma_clime}),
\begin{align*}
    \Vert \hat{\Gamma}_{\cdot j}-\Gamma_{\cdot j}\Vert_2^2 &\leq C \Vert \Gamma_{\cdot j}\Vert_0 \rho^2_n\quad\text{for $j=1,\ldots, p$},\\
    \Vert \hat{\Gamma}_{\cdot j}-\Gamma_{\cdot j}\Vert_1 &\leq  C \Vert \Gamma_{\cdot j}\Vert_0\rho_n\quad\text{for $j=1,\ldots, p$},
\end{align*}
where $C$ is a constant depending on $L$.
\end{lem}
\begin{prf}
For the definition of $\hat{G}$ and (\ref{cstr:M_max}), we have
\begin{align*}
   \Vert \hat{G}_{\cdot j}\Vert_{1}&\leq \Vert \Gamma_{\cdot j}\Vert_{1} \quad \text{for $j=1,\ldots, p$,}  \\
    \Vert \hat{\Gamma}_{\cdot j}\Vert_{1}&\leq \Vert \Gamma_{\cdot j}\Vert_{1} \quad \text{for $j=1,\ldots, p$,} \\
    \Vert W(\hat{G}-M) \Vert_{\max} &\leq  \Vert W \hat{G} - I\Vert_{\max} +  \Vert  W \Gamma -I_p  \Vert_{\max} \\
    &\leq 2\rho_n.
\end{align*}
We have
\begin{align*}
    \Vert \Gamma^{-1} (\hat{G}-\Gamma)\Vert_{\max} &\leq \Vert W (\hat{G}-\Gamma)\Vert_{\max} + \Vert (W - \Gamma^{-1}) (\hat{G}-\Gamma)\Vert_{\max} \\
    &\leq  \Vert W (\hat{G}-\Gamma)\Vert_{\max} + \Vert W - \Gamma^{-1}\Vert_{\max}\Vert\hat{G}-\Gamma\Vert_{1} \\
    &\leq (2+2L)\rho_n .
\end{align*}
Hence we have
\begin{align*}
    \Vert \hat{G}-\Gamma\Vert_{\max} &\leq  \Vert \Gamma \Vert_{1} \Vert \Gamma^{-1} (\hat{G}-\Gamma)\Vert_{\max} \\
    &\leq  (2L+2L^2)\rho_n ,
\end{align*}
and
\begin{align}
    \Vert \hat{\Gamma}-\Gamma\Vert_{\max} &\leq  (2L+2L^2)\rho_n .   \label{inq:matrix_max_norm}
\end{align}
From the triangle inequality, we have for $j=1,\ldots, p$,
\begin{align}
    \sum_{l\notin \supp(\Gamma_{\cdot j})} |\hat{\Gamma}_{lj}-\Gamma_{lj}| &= \sum_{l\notin \supp(\Gamma_{\cdot j})} |\hat{\Gamma}_{lj}| \notag \\
    &\leq  \sum_{l \in \supp(\Gamma_{\cdot j})}  |\Gamma_{lj}| - \sum_{l \in \supp(\Gamma_{\cdot j})}  |\hat{\Gamma}_{lj}| \notag \\
    &\leq   \sum_{l \in \supp(\Gamma_{\cdot j})} |\hat{\Gamma}_{lj}-\Gamma_{lj}|.\notag
\end{align}
Hence we have
\begin{align}
    \Vert \hat{\Gamma}_{\cdot j} - \Gamma_{\cdot j}\Vert_1 &\leq  2\sum_{l \in \supp(\Gamma_{\cdot j})} |\hat{\Gamma}_{lj}-\Gamma_{lj}| \notag \\
    &\leq  2\Vert \hat{\Gamma} -\Gamma\Vert_{\max} \Vert \Gamma_{\cdot j}\Vert_0  \notag  \\
    &\leq   (4L + 4L^2) \Vert \Gamma_{\cdot j}\Vert_0 \rho_n, \quad\text{(From \ref{inq:matrix_max_norm})}   \label{inq:L_1_norm}
\end{align}
and
\begin{align}
     \Vert \hat{\Gamma}_{\cdot j} - \Gamma_{\cdot j}\Vert_2^2  &\leq  \Vert \hat{\Gamma}_{\cdot j} - \Gamma_{\cdot j}\Vert_1 \Vert \hat{\Gamma}_{\cdot j} - \Gamma_{\cdot j}\Vert_\infty \notag \\
    &\leq  \Vert \hat{\Gamma}_{\cdot j} - \Gamma_{\cdot j}\Vert_1 \Vert \hat{\Gamma} - \Gamma\Vert_{\max} \notag \\
    &\leq (4L + 4L^2)(2L+2L^2) \Vert \Gamma_{\cdot j}\Vert_0 \rho^2_n. \quad\text{(From \ref{inq:matrix_max_norm}) and (\ref{inq:L_1_norm}))}   \notag
\end{align}

\end{prf}

\subsubsection{Error bound of the estimator in Algorithm \ref{alg:omega_est}}

 \errorrateomegaclime*

\begin{prf}
With $\phi(a,b)=\frac{1}{2}(a-b)^2$, under event $\mathcal{A}_{11}\cap\mathcal{A}_{8}$ (Lemma \ref{lem:bound_of_random_quan}(i), (viii)) and Assumption \ref{ass:sparse_omega}, we have
\begin{align*}
    \Vert \hat{\Sigma}_J \Omega^* - I\Vert_{\max} &\leq  \rho_n, \\
    \Vert \hat{\Sigma}_J - \Sigma\Vert_{\max} &\leq  \rho_n, \\
    \Vert \Omega^*\Vert_{1} &\leq L.
\end{align*}
From Lemma \ref{lem:clime_est_proof}, we have for $j=1,\ldots, p$,
\begin{align}
    \Vert \hat{\Omega}^J_{\cdot j}-\Omega^*_{\cdot j}\Vert_2^2 &\leq C s_{1j} \frac{\log(p\vee n)}{n}, \label{inq:clime_omega_1}\\
    \Vert \hat{\Omega}^J_{\cdot j}-\Omega^*_{\cdot j}\Vert_1 &\leq  C s_{1j}\sqrt{\frac{\log(p\vee n)}{n}}, \label{inq:clime_omega_2}
\end{align}
where $C$ is a constant depending on $(\bar{A}_0, \sigma_1,\rho,L)$ only.

From \eqref{inq:clime_omega_2}, we have
\begin{align*}
    \Vert \hat{\Omega}^{J}_{S_\xi,S_\xi} -\Omega^*_{S_\xi,S_\xi} \Vert_1
     &\leq \max_{j\in S_{\xi}} \Vert \hat{\Omega}^{J}_{\cdot j} - \Omega^*_{\cdot_j}\Vert_1 \\
     &\leq C \Vert s_{1\xi} \Vert_\infty  \sqrt{\frac{\log(p\vee n)}{n}} . \\
\end{align*}
From the symmetry of $\hat{\Omega}^{J}$ and $\Omega^*$, we have
\begin{align}
    \Vert \hat{\Omega}^{J}_{S_\xi,S_\xi} -\Omega^*_{S_\xi,S_\xi} \Vert_2
     &\leq \Vert \hat{\Omega}^{J}_{S_\xi,S_\xi} -\Omega^*_{S_\xi,S_\xi} \Vert_1 \notag \\
     &\leq C \Vert s_{1\xi} \Vert_\infty \sqrt{\frac{\log(p\vee n)}{n}}. \label{inq:clime_omega_3}
\end{align}

From (\ref{inq:clime_omega_1}), we also have
\begin{align}
    \Vert \hat{\Omega}^{J}_{S_\xi,S_\xi}-\Omega^*_{S_\xi,S_\xi}\Vert_2 &\leq  \Vert \hat{\Omega}^{J^\c}_{S_\xi,S_\xi}-\Omega^*_{S_\xi,S_\xi} \Vert_2 \notag \\
    &\leq  \Vert \hat{\Omega}^{J}_{S_\xi,S_\xi}-\Omega^*_{S_\xi,S_\xi} \Vert_{\F} \notag \\
    &\leq \sqrt{\Tr( (\hat{\Omega}^{J}_{S_\xi,S_\xi}- \Omega^*_{S_\xi,S_\xi})(\hat{\Omega}^{J}_{S_\xi,S_\xi}- \Omega^*_{S_\xi,S_\xi}))}\notag\\
    &\leq \sqrt{\sum_{j\in S_{\xi}} (\hat{\Omega}^{J}_{\cdot j} - \Omega^*_{\cdot_j})^\T (\hat{\Omega}^{J}_{\cdot j} - \Omega^*_{\cdot_j})}  \notag\\
    &\leq \sqrt{C} \sqrt{\Vert s_{1\xi}\Vert_1 } \sqrt{\frac{\log(p\vee n)}{n}}. \label{inq:clime_omega_4}
\end{align}

From (\ref{inq:clime_omega_3}) and (\ref{inq:clime_omega_4}), we have
\begin{align}
    \Vert (\hat{\Omega}^{J}-\Omega^*)_{S_\xi,S_\xi} \Vert_2 \leq (C + \sqrt{C})\min\{\sqrt{\Vert s_{1\xi}\Vert_1}, \Vert s_{1\xi}\Vert_\infty\}\sqrt{\frac{\log(p\vee n)}{n}}.
\end{align}

\end{prf}

\subsubsection{Error bound of the estimator in Algorithm \ref{alg:est_M}}

\lemerrorrateM*

\begin{prf}
Under the event $\mathcal{A}_{8}\cap\mathcal{A}_{9}$ (Lemma \ref{lem:bound_of_random_quan}(viii)(ix)) and given Condition \ref{con:beta_hat_J}, we have for $j=1,\ldots, p$,
\begin{align}
    &\Vert \hat{H}_J M_{\cdot j}^* -e_j\Vert_\infty \notag \\
    &\leq \Vert \frac{1}{|J|}\sum_{i\in J} \phi_2(Y_i,X_i^\T\bar\beta)X_iX_i^\T M_{\cdot j}^*-e_j\Vert_\infty + \Vert \frac{1}{|J|}\sum_{i\in J} (\phi_2(Y_i,X_i^\T\bar\beta)-\phi_2(Y_i,X_i^\T \hat{\beta}_J))X_iX_i^\T M_{\cdot j}^*\Vert_\infty  \notag\\
    &\leq  \Vert \frac{1}{|J|}\sum_{i\in J} \phi_2(Y_i,X_i^\T\bar\beta)X_iX_i^\T M_{\cdot j}^*-e_j\Vert_\infty \notag \\
    &+ C_7  \sqrt{\frac{1}{|J|}\sum_{i\in J}(X_i^\T(\bar\beta-\hat{\beta}_J))^2}\sqrt{\max_{1\leq l\leq p}\{\frac{1}{|J|}\sum_{i\in J} x_{il}^2(M_{\cdot j}^{*\T}X_i)^2\}} \quad\text{(From Assumption \ref{bounded_higher_differentiate})}\notag \\
    &\leq (\sqrt{\frac{3}{c(1)\rho}}K(1)(1+(d(1) \log(2))^{-1})  C_6\sigma_1^2 + C_7C_8\sqrt{(K(\frac{1}{2})(1+(d(\frac{1}{2}) \log(2))^{-2})+\frac{1}{d(\frac{1}{2})})\sigma_1^4}) \notag \\
    &\times \sqrt{s_0\frac{\log(p\vee n)}{n}}. \label{inq:bound_max_entries}
\end{align}
Hence we have
\begin{align*}
    &\Vert \hat{H}_J M^* -I_p\Vert_{\max} \\ &\leq (\sqrt{\frac{3}{c(1)\rho}}K(1)(1+(d(1) \log(2))^{-1})  C_6\sigma_1^2 + C_7C_8\sqrt{(K(\frac{1}{2})(1+(d(\frac{1}{2}) \log(2))^{-2})+\frac{1}{d(\frac{1}{2})})\sigma_1^4})\\
    &\times \sqrt{s_0\frac{\log(p\vee n)}{n}}.
\end{align*}
From Lemma \ref{lem:H_diff}, under event $\mathcal{A}_{11}\cap \mathcal{A}_{21}$(Lemma \ref{lem:bound_of_random_quan}(i), (ii)) and Assumption \ref{bounded_higher_differentiate} and Condition \ref{con:beta_hat_J} are satisfied, we have
\begin{align*}
    &\Vert \hat{H}_J -H\Vert_{\max}\\
    &\leq (\sqrt{\frac{3}{c(1)\rho}}K(1)(1+(d(1) \log(2))^{-1})  C_6\sigma_1^2 + C_7C_8\sqrt{(K(\frac{1}{2})(1+(d(\frac{1}{2}) \log(2))^{-2})+\frac{1}{d(\frac{1}{2})})\sigma_1^4})\\
    &\times \sqrt{s_0\frac{\log(p\vee n)}{n}}.
\end{align*}
Further with $\Vert M^*\Vert_{1}\leq L$, from Lemma \ref{lem:clime_est_proof}, for $j=1,\ldots p$, we have
\begin{align}
   \Vert \hat{M}^{J}_{\cdot j}-M^*_{\cdot j}\Vert_2^2 &\leq  C s_{1j}s_0 \frac{\log(p\vee n)}{n},  \label{inq:clime_M_1} \\
    \Vert \hat{M}^{J}_{\cdot j}-M^*_{\cdot j}\Vert_1 &\leq   C s_{1j}\sqrt{s_0\frac{\log(p\vee n)}{n}}, \label{inq:clime_M_2}
\end{align}
where $C$ is a constant depending on $(\bar{A}_0, C_6, C_7, C_8, \rho, \sigma_1,L)$.

From \eqref{inq:clime_M_2}, we have
\begin{align*}
    \Vert \hat{M}^{J}_{S_\xi,S_\xi}-M^*_{S_\xi,S_\xi} \Vert_1
     &\leq \max_{j\in S_{\xi}} \Vert \hat{M}^{J}_{\cdot j} - M^*_{\cdot_j}\Vert_1 \\
     &\leq C \Vert s_{1\xi}\Vert_\infty \sqrt{s_0}\sqrt{\frac{\log(p\vee n)}{n}} . \\
\end{align*}

From the symmetry of $\hat{M}^{J}$ and $M^*$, we have
\begin{align}
    \Vert \hat{M}^{J}_{S_\xi,S_\xi}-M^*_{S_\xi,S_\xi} \Vert_2
     &\leq \Vert \hat{M}^{J}_{S_\xi,S_\xi}-M^*_{S_\xi,S_\xi} \Vert_1 \notag \\
     &\leq C \Vert s_{1\xi}\Vert_\infty \sqrt{s_0}\sqrt{\frac{\log(p\vee n)}{n}}  . \label{inq:M_J_3}
\end{align}

From (\ref{inq:clime_M_1}), we also have
\begin{align}
    \Vert \hat{M}^{J}_{S_\xi,S_\xi}-M^*_{S_\xi,S_\xi}\Vert_2 &\leq  \Vert \hat{M}^{J^\c}_{S_\xi,S_\xi}-M^*_{S_\xi,S_\xi} \Vert_2 \notag \\
    &\leq  \Vert \hat{M}^{J}_{S_\xi,S_\xi}-M^*_{S_\xi,S_\xi} \Vert_{\F} \notag \\
    &\leq \sqrt{\Tr( (\hat{M}^{J}_{S_\xi,S_\xi}- M^*_{S_\xi,S_\xi})(\hat{M}^{J}_{S_\xi,S_\xi}- M^*_{S_\xi,S_\xi}))}\notag\\
    &\leq \sqrt{\sum_{j\in S_{\xi}} (\hat{M}^{J}_{\cdot j} - M^*_{\cdot_j})^\T (\hat{M}^{J}_{\cdot j} - M^*_{\cdot_j})}  \notag\\
    &\leq \sqrt{C} \sqrt{\Vert s_{1\xi}\Vert_1 } \sqrt{s_0} \sqrt{\frac{\log(p\vee n)}{n}}. \label{inq:M_J_4}
\end{align}

From (\ref{inq:M_J_3}) and (\ref{inq:M_J_4}), we have
\begin{align}
    \Vert \hat{M}^{J}_{S_\xi,S_\xi}-M^*_{S_\xi,S_\xi} \Vert_2 \leq (C + \sqrt{C})\min\{\sqrt{\Vert s_{1\xi}\Vert_1}, \Vert s_{1\xi}\Vert_\infty\}\sqrt{s_0}\sqrt{\frac{\log(p\vee n)}{n}}.
\end{align}
\end{prf}

\subsection{Error bounds of estimators in Algorithms \ref{alg:omega_est_two_stage} and \ref{alg:est_M_two_stage}}

We give error bounds for the estimators in Algorithms \ref{alg:omega_est_two_stage} and \ref{alg:est_M_two_stage}. We first give the following conditions that are important in our analysis.

\begin{conditionE}[Compatibility Condition] For any symmetric matrix $W\in \mathbb{R}^{p \times p}$, $S\subset \{1,\ldots, p\}$, some constants $\nu_0>0$ and $\eta_0>1$, we have
\begin{align*}
    \nu_0^2 (\sum_{j\in S} |b_j|)^2 \leq |S|(b^\T W b),
\end{align*}
for any vector $b\in \mathbb{R}^p$ satisfying the cone condition below,
\begin{align*}
    \sum_{j \notin S} |b_j|\leq \eta_0 \sum_{j\in S}|b_j|.
\end{align*}
And we denote $W$, $S$, $\nu_0$ and $\eta_0$ satisfy the $CompCondition(W, S, \nu_0,\eta_0)$.
\end{conditionE}

\begin{conditionE}[Bounded infinity norm]
\label{con:bounded_infinity_norm}
    For any symmetric $W\in \mathbb{R}^{p \times p}$, $\Gamma^*\in \mathbb{R}^{p\times p}$ and $\bar{\lambda}_n>0$, we have
    \begin{align*}
        \Vert W \Gamma^* -I_p\Vert_{\max}  \leq \bar{\lambda}_n .
    \end{align*}
   And we denote $W$, $\Gamma^*$ and $\bar{\lambda}_n$ satisfy the $InfinityNorm(W, \Gamma^*,\bar{\lambda}_n)$.
\end{conditionE}

For any symmetric $W\in \mathbb{R}^{p \times p}$ , $\bar{\lambda}_n>0$ and $\bar{A}_0>1$, we consider the following estimator,
\begin{align}
   \hat{\Gamma} &= \argmin_{G\in\mathbb{R}^{p\times p}} \frac{1}{2}G^\T W G -\Tr(G) +\bar{A}_0 \bar{\lambda}_n \Vert G\Vert_1. \label{def:est_gamma}
\end{align}

\begin{lem}
\label{lem:estimating_v}
For symmetric $W\in \mathbb{R}^{p \times p}$, $\nu_0>0$, $\eta_0>1$, $\Gamma^* \in \mathbb{R}^{p\times p}$, $\bar{\lambda}_n>0$ and let $S_j=\{1\leq l\leq p:\Gamma^*_{lj}\neq 0\}$ for $j=1,\ldots,p$, given $CompCondition(W, S_j, \nu_0,\eta_0)$ for $j=1,\ldots, p$ are satisfied and $InfinityNorm(W,  \Gamma^*,\bar{\lambda}_n)$ is satisfied and $\bar{A}_0\geq \frac{\eta_0+1}{\eta_0-1}$, we have for $\hat{\Gamma}$ in (\ref{def:est_gamma}),
\begin{align*}
    \Vert \hat{\Gamma}_{\cdot j}-{\Gamma^*}_{\cdot j}\Vert_1 & \leq C|S_j|\bar{\lambda}_n\quad\text{for $j=1,\ldots p$,}
\end{align*}
where $C>0$ is a constant depending only on $(\bar{A}_0, \nu_0)$ only.
\end{lem}
\begin{prf}
Define $\Obj(G)=\frac{1}{2} \Tr(G^\T W G) - \Tr(G)+\lambda\Vert G\Vert_1$. Because $\hat{\Gamma}$ is a minimizer of $\Obj(G)$, we have $\hat{\Gamma}$ satisfies the stationary condition below
\begin{align}
    0\in W \hat{\Gamma} - I +\bar{A}_0\bar{\lambda}_n \frac{\Vert \hat{\Gamma}\Vert_{L_1}}{\partial \hat{\Gamma}} .\label{con:stationary_con}
\end{align}

\noindent For any $\Gamma \in \mathbb{R}^{p\times p}$, we have
\begin{align*}
    \Tr((\hat{\Gamma}-\Gamma)^\T(W\hat{\Gamma}-I_p))
    &= \sum_{i,j} (\hat{\Gamma}_{ij} -\Gamma_{ij}) (W\hat{\Gamma}-I_p)_{ij} \\
    &=\sum_{i,j} (\hat{\Gamma}_{ij} -\Gamma_{ij}) (-\bar{A}_0\bar{\lambda}_n\frac{\partial \Vert \hat{\Gamma}\Vert_{L_1}}{\partial \hat{\Gamma}})_{ij}.\quad(\text{from (\ref{con:stationary_con})})
\end{align*}

\noindent $\newline$
By direct calculation, we have
\begin{align*}
   & \text{If $\hat{\Gamma}_{ij}$ > 0, } (\hat{\Gamma}_{ij} -\Gamma_{ij}) (-\bar{A}_0\bar{\lambda}_n\frac{\partial \Vert \hat{\Gamma}\Vert_{L_1}}{\partial \hat{\Gamma}})_{ij} = -\bar{A}_0\bar{\lambda}_n (\hat{\Gamma}_{ij} -\Gamma_{ij}) \leq \bar{A}_0\bar{\lambda}_n (|\Gamma_{ij}|-\hat{\Gamma}_{ij})= \bar{A}_0\bar{\lambda}_n (|\Gamma_{ij}|-|\hat{\Gamma}_{ij}|).\\
   & \text{If $\hat{\Gamma}_{ij} \leq 0$, }  (\hat{\Gamma}_{ij} -\Gamma_{ij}) (-\bar{A}_0\bar{\lambda}_n\frac{\partial \Vert \hat{\Gamma}\Vert_{L_1}}{\partial \hat{\Gamma}})_{ij}=\bar{A}_0\bar{\lambda}_n(\hat{\Gamma}_{ij} -\Gamma_{ij})\leq \bar{A}_0\bar{\lambda}_n(|\Gamma_{ij}|+\hat{\Gamma}_{ij}) =\bar{A}_0\bar{\lambda}_n(|\Gamma_{ij}|-|\hat{\Gamma}_{ij}|)  . \\
   &\text{If $\hat{\Gamma}_{ij}=0$, }  (\hat{\Gamma}_{ij} -\Gamma_{ij}) (-\bar{A}_0\bar{\lambda}_n\frac{\partial \Vert \hat{\Gamma}\Vert_{L_1}}{\partial \hat{\Gamma}})_{ij} \leq \bar{A}_0\bar{\lambda}_n |\Gamma_{ij}| =\bar{A}_0\bar{\lambda}_n(|\Gamma_{ij}|-|\hat{\Gamma}_{ij}|).
\end{align*}

\noindent In conclusion, for any $\Gamma\in \mathbb{R}^{p\times p}$, we have
\begin{align}
     \Tr((\hat{\Gamma}-\Gamma)^\T (W\hat{\Gamma}-I_p))   \leq  \bar{A}_0\bar{\lambda}_n (\Vert \Gamma\Vert_{L_1} - \Vert \hat{\Gamma}\Vert_{L_1}) . \label{inq:basic_inq}
\end{align}

From the definition of $\hat{\Gamma}$, for $1\leq j \leq p$, we have
\begin{align}
   &\quad  \Tr( \hat{\Gamma}_{\cdot j}^\T W \hat{\Gamma}_{\cdot j} ) - \hat{\Gamma}_{jj} + \bar{A}_0\bar{\lambda}_n \Vert \hat{\Gamma}_{\cdot j}\Vert_1 \leq  \notag \Tr( \Gamma^{*\T}_{\cdot j} W \Gamma^*_{\cdot j} ) - \Gamma^*_{jj} + \bar{A}_0\bar{\lambda}_n \Vert \Gamma^*_{\cdot j}\Vert_1\\
   &\Rightarrow  (\hat{\Gamma}_{\cdot j}-\Gamma_{\cdot j}^*)^\T W(\hat{\Gamma}_{\cdot j}-\Gamma_{\cdot j}^*) + (\hat{\Gamma}_{\cdot j}-\Gamma_{\cdot j}^* )^\T(W \Gamma^{*}_{\cdot j}-e_j) \leq \bar{A}_0\bar{\lambda}_n(\Vert \Gamma_{\cdot j}^*\Vert_{1} - \Vert \hat{\Gamma}_{\cdot j}\Vert_{1}) .\label{inq:basic_inq_transformed}
\end{align}

Now we focus on the term $(\hat{\Gamma}_{\cdot j}-\Gamma_{\cdot j}^* )^\T(W \Gamma^{*}_{\cdot j}-e_j)$, we have
\begin{align}
    &\quad (\hat{\Gamma}_{\cdot j}-\Gamma^*_{\cdot j})^\T(W\Gamma^{*}_{\cdot j}-e_j) \notag \\
    &\leq \Vert \hat{\Gamma}_{\cdot j}-\Gamma^*_{\cdot j}\Vert_1 \Vert W\Gamma^{*}_{\cdot j}-e_j\Vert_{\infty}  \notag\\
    &\leq \Vert \hat{\Gamma}_{\cdot j}-\Gamma^*_{\cdot j}\Vert_1 \Vert W\Gamma^{*}-I_p\Vert_{\max} \notag\\
    &\leq \bar{\lambda}_n \Vert \hat{\Gamma}_{\cdot j}-\Gamma^*_{\cdot j}\Vert_1. \quad(\text{from Condition \ref{con:bounded_infinity_norm}}) \label{inq:mediate_term_basic_inq}
\end{align}
From (\ref{inq:basic_inq_transformed}) and (\ref{inq:mediate_term_basic_inq}), we have
\begin{align*}
    (\hat{\Gamma}_{\cdot j}-\Gamma_{\cdot j}^*)^\T W(\hat{\Gamma}_{\cdot j}-\Gamma_{\cdot j}^*) \leq \bar{A}_0\bar{\lambda}_n(\Vert \Gamma_j^*\Vert_1 - \Vert \hat{\Gamma}_{\cdot j}\Vert_{1})+\bar{\lambda}_n \Vert \hat{\Gamma}_{\cdot j}-\Gamma^*_{\cdot j}\Vert_1.
\end{align*}
Hence we have
\begin{align*}
    (\hat{\Gamma}_{\cdot j}-\Gamma_{\cdot j}^*)^\T W(\hat{\Gamma}_{\cdot j}-\Gamma_{\cdot j}^*) +  (\bar{A}_0-1) \bar{\lambda}_n  \Vert \hat{\Gamma}_{\cdot j}-\Gamma^*_{\cdot j}\Vert_1 + \bar{A}_0\bar{\lambda}_n \Vert \hat{\Gamma}_{\cdot j}\Vert_1 \leq \bar{A}_0\bar{\lambda}_n \Vert \Gamma^*_{\cdot j}\Vert_1 + \bar{A}_0 \bar{\lambda}_n  \Vert \hat{\Gamma}_{\cdot  j}-\Gamma^*_{\cdot j}\Vert_1 .
\end{align*}
After applying the equation and inequality below,
\begin{align*}
    &|\hat{\Gamma}_{l j}| = |\hat{\Gamma}_{l j}-\Gamma^*_{l j}| \quad \text{if $l\notin S_j$},\\
    &|\hat{\Gamma}_{l j}| \geq |\Gamma^*_{l j}| - |\hat{\Gamma}_{l j}-\Gamma^*_{l j}| \quad \text{if $l\in S_j$} ,
\end{align*}
we have
\begin{align*}
     &\bar{A}_0\bar{\lambda}_n \Vert \Gamma^*_{\cdot j}\Vert_1 + \bar{A}_0 \bar{\lambda}_n\Vert \hat{\Gamma}_{\cdot j}-\Gamma^*_{\cdot j}\Vert_1 \\
     &\geq  (\hat{\Gamma}_{\cdot j}-\Gamma_{\cdot j}^*)^\T W(\hat{\Gamma}_{\cdot j}-\Gamma_{\cdot j}^*)+  (\bar{A}_0-1) \bar{\lambda}_n  \Vert \hat{\Gamma}_{\cdot j}-\Gamma^*_{\cdot j}\Vert_1 + \bar{A}_0\bar{\lambda}_n \Vert \hat{\Gamma}_{\cdot j}\Vert_1 \\
     &\geq   (\hat{\Gamma}_{\cdot j}-\Gamma_{\cdot j}^*)^\T W(\hat{\Gamma}_{\cdot j}-\Gamma_{\cdot j}^*)+(\bar{A}_0-1) \bar{\lambda}_n  \Vert \hat{\Gamma}_{\cdot j}-\Gamma^*_{\cdot j}\Vert_1 \\
     &\quad\quad\quad\quad + \bar{A}_0\bar{\lambda}_n (\sum_{l\notin S_j}(|\hat{\Gamma}_{l j}-\Gamma^*_{l j}|) + \sum_{l\in S_j}(|\Gamma^*_{l j}|-|\hat{\Gamma}_{l j}-\Gamma^*_{l j}|)).
\end{align*}
That is, we have
\begin{align}
    2\bar{A}_0\bar{\lambda}_n  \sum_{l\in S_j} |\hat{\Gamma}_{l j}-\Gamma_{l j}^*|  \geq   (\hat{\Gamma}_{\cdot j}-\Gamma_{\cdot j}^*)^\T W(\hat{\Gamma}_{\cdot j}-\Gamma_{\cdot j}^*) + (\bar{A}_0-1)\bar{\lambda}_n \Vert \hat{\Gamma}_{\cdot j}-\Gamma^*_{\cdot j}\Vert_1 . \label{inq:bound_S_S_C}
\end{align}
From (\ref{inq:bound_S_S_C}), we have
\begin{align*}
   & 2\bar{A}_0\bar{\lambda}_n  \sum_{l\in S_j} |\hat{\Gamma}_{lj}-\Gamma^*_{lj}| \geq  (\bar{A}_0-1)\bar{\lambda}_n \Vert \hat{\Gamma}_{\cdot j}-\Gamma^*_{\cdot j}\Vert_1. \notag
\end{align*}
Hence we have,
\begin{align*}
    \sum_{l\notin S_j} |\hat{\Gamma}_{lj}-\Gamma_{lj}^*|\leq \frac{\bar{A}_0+1}{\bar{A}_0-1} \sum_{l\in S_j} |\hat{\Gamma}_{lj}-\Gamma_{lj}^*| .
\end{align*}
From $\bar{A}_0\geq \frac{\eta_0+1}{\eta_0-1}$, we have $\frac{\bar{A}_0+1}{\bar{A}_0-1}\leq \eta_0$, so we have $\hat{\Gamma}_{\cdot j}-\Gamma_{\cdot j}$ satisfies the cone condition. From condition $CompCondition(W, S_j, \nu_0,\eta_0)$, we have
\begin{align}
    \nu^2_0 (\sum_{l\in S_j}|\hat{\Gamma}_{l j}-\Gamma^*_{l j}|)^2\leq 2|S_j| (\hat{\Gamma}_{\cdot j}-\Gamma_{\cdot j}^*)^\T W(\hat{\Gamma}_{\cdot j}-\Gamma_{\cdot j}^*).  \label{inq:comp_con_inq_1}
\end{align}
Denote $\Delta =  (\hat{\Gamma}_{\cdot j}-\Gamma_{\cdot j}^*)^\T W(\hat{\Gamma}_{\cdot j}-\Gamma_{\cdot j}^*) + (\bar{A}_0-1)\bar{\lambda}_n \Vert \hat{\Gamma}_{\cdot j}-\Gamma^*_{\cdot j}\Vert_1$, we have
\begin{align*}
   \Delta &= (\hat{\Gamma}_{\cdot j}-\Gamma_{\cdot j}^*)^\T W(\hat{\Gamma}_{\cdot j}-\Gamma_{\cdot j}^*) + (\bar{A}_0-1)\bar{\lambda}_n \Vert \hat{\Gamma}_{\cdot j}-\Gamma^*_{\cdot j}\Vert_1 \\&\leq   2\bar{A}_0\bar{\lambda}_n  \sum_{l\in S_j} |\hat{\Gamma}_{lj}-\Gamma_{lj}^*|\quad\text{(From (\ref{inq:bound_S_S_C}))}\\
    &\leq  2\sqrt{2}\bar{A}_0\bar{\lambda}_n  \nu_0^{-1}|S_j|^{\frac{1}{2}}(\hat{\Gamma}_{\cdot j}-\Gamma_{\cdot j}^*)^\T W(\hat{\Gamma}_{\cdot j}-\Gamma_{\cdot j}^*)\\
    &\leq   2\sqrt{2}\bar{A}_0\bar{\lambda}_n  \nu_0^{-1}|S_j|^{\frac{1}{2}} \Delta^{\frac{1}{2}}. \quad\text{(From (\ref{inq:comp_con_inq_1}))}
\end{align*}
Hence we have
\begin{align*}
     (\hat{\Gamma}_{\cdot j}-\Gamma_{\cdot j}^*)^\T W(\hat{\Gamma}_{\cdot j}-\Gamma_{\cdot j}^*) + (\bar{A}_0-1)\bar{\lambda}_n \Vert \hat{\Gamma}_{\cdot j}-\Gamma^*_{\cdot j}\Vert_1 \leq  8\bar{A}_0^2 \nu_0^{-2}\bar{\lambda}_n^2 |S_j|,
\end{align*}
and
\begin{align*}
    \Vert \hat{\Gamma}_{\cdot j}-\Gamma^*_{\cdot j}\Vert_1 &\leq   \frac{8\bar{A}_0^2 \nu_0^{-2}}{\bar{A}_0-1} |S_j|\bar{\lambda}_n .
\end{align*}
\end{prf}

\subsubsection{Error bound of the estimators in Algorithms \ref{alg:omega_est_two_stage} }

\begin{lem}
\label{lem:empirical_compatibility_condition}

Under event $\mathcal{A}_{14}$(Lemma \ref{lem:bound_of_random_quan}(i)), for any $\eta_0>1$ and any $S\subset\{1,\ldots p\}$, let $W=\frac{1}{|J|}\sum_{i\in J}X_iX_i^\T$ and $\delta=C_4|S|\bar{\lambda}_n (1+\eta_0)^2 \leq1$, then condition $CompCondition(W, S, \sqrt{\frac{1}{C_4}(1-\delta)}, \eta_0)$ holds, where $\bar{\lambda}_n=\sqrt{\frac{3}{c(1)\rho}}K(1)(1+(d(1) \log(2))^{-1})\sigma_1^2\sqrt{\frac{\log(p\vee n)}{n}}$.
\end{lem}
\begin{prf} Under the event $\mathcal{A}_{14}$(Lemma \ref{lem:bound_of_random_quan}(i)), we have
\begin{align}
    \Vert W -\Sigma\Vert_{\max} \leq \sqrt{\frac{3}{c(1)}}\frac{K(1)(1+(d(1) \log(2))^{-1})}{\sqrt{\rho}}\sigma_1^2\sqrt{\frac{\log(p\vee n)}{n}}.  \notag %\label{inq:bounded_sigma_hat_entries}
\end{align}
For any $S\subset \{1,\ldots,p\}$ and $\eta_0 > 1$, let $b\in\mathbb{R}^p$ satisfies
\begin{align*}
    \sum_{j\notin S} |b_j|\leq \eta_0 (\sum_{j\in S}|b_j|).
\end{align*}
We have,
\begin{align*}
    \frac{1}{C_4}  (\sum_{j\in S}|b_j|)^2 &\leq   \frac{1}{C_4} |S| (\sum_{j\in S}b_j^2) \\
    &\leq  \frac{1}{C_4} |S|  \Vert b\Vert_2^2 \\
    &\leq   |S| b^\T \Sigma b \quad\text{(from Assumption \ref{sigma_eigenvalue}(ii))} \\
    &=  |S| b^\T W b + |S| b^\T (\Sigma-W)b.
\end{align*}
For $b^\T (\Sigma-W)b$, we have
\begin{align*}
    b^\T (\Sigma-W)b &\leq  \Vert b\Vert_1^2 \Vert W-\Sigma\Vert_{\max} \\
    &\leq   \Vert b\Vert_1^2 \bar{\lambda}_n.
\end{align*}

\noindent Hence we have
\begin{align*}
     \frac{1}{C_4}  (\sum_{j\in S}|b_j|)^2&\leq  |S| b^\T W b + |S|\bar{\lambda}_n \Vert b\Vert_1^2 \\
     &\leq   |S| b^\T W b + |S|\bar{\lambda}_n (1+\eta_0)^2(\sum_{j\in S}|b_j|)^2.
\end{align*}
That is, we have
\begin{align*}
    \frac{1}{C_4}(1-C_4|S|\bar{\lambda}_n (1+\eta_0)^2) (\sum_{j\in S}|b_j|)^2\leq  |S| b^\T W b.
\end{align*}
Hence
\begin{align*}
    \frac{1}{C_4}(1-\delta) (\sum_{j\in S}|b_j|)^2\leq  |S| b^\T W b.
\end{align*}
\end{prf}

\begin{lem}
\label{lem:inf_norm_con_u}
Under event $\mathcal{A}_{8}$(Lemma \ref{lem:bound_of_random_quan}(viiii)) with $\phi(a,b)=\frac{1}{2}(a-b)^2$, let $W=\frac{1}{|J|}\sum_{i\in J}X_iX_i^\T$, then condition $InfinityNorm(W, \Omega^*,\bar{\lambda}_n)$ holds, where $\bar{\lambda}_n=\sqrt{\frac{3}{c(1)\rho}}K(1)(1+(d(1) \log(2))^{-1})   \sigma_1^2\sqrt{\frac{\log(p\vee n)}{n}}$.
\end{lem}
\begin{prf}
Let $\phi(a,b)=\frac{1}{2}(a-b)^2$, the Assumption \ref{bounded_higher_differentiate} holds with $C_6=1$ and $C_7=0$. Under the event $\mathcal{A}_{8}$ (Lemma \ref{lem:bound_of_random_quan}(viiii)) with $\phi(a,b)=\frac{1}{2}(a-b)^2$, we have
\begin{align*}
    & \Vert W \Omega^* -I_p \Vert_\infty \leq  \sqrt{\frac{3}{c(1)\rho}}K(1)(1+(d(1) \log(2))^{-1})   \sigma_1^2\sqrt{\frac{\log(p\vee n)}{n}} .
\end{align*}
\end{prf}

\errorrateomegatwostage*

\begin{prf}

\noindent From Lemma \ref{lem:estimating_v}, \ref{lem:empirical_compatibility_condition} and \ref{lem:inf_norm_con_u}, under event $\mathcal{A}_{14}\cap \mathcal{A}_8$ (Lemma \ref{lem:bound_of_random_quan}(i), (viii)), we have for $j=1,\ldots, p$,
\begin{align*}
    \Vert \hat{G}_{\cdot j}-{\Omega^*}_{\cdot j}\Vert_{1} & \leq C^\prime s_1\bar{\lambda}_n,
\end{align*}
where $C^\prime$ is a constant depending on $(\bar{A}_0,C_4,\delta)$ only.

\noindent Hence we have
\begin{align*}
    \Vert \hat{G}^J-{\Omega^*}\Vert_{1} & \leq C^\prime s_1\bar{\lambda}_n.
\end{align*}

\noindent Under event $\mathcal{A}_{8}$(Lemma \ref{lem:bound_of_random_quan}(viii)) with $\phi(a,b)=\frac{1}{2}(a-b)^2$, we have
\begin{align*}
   \Vert \hat{\Sigma}_J\Omega^* - I_p \Vert \leq \bar{A}_0\bar{\lambda}_n.
\end{align*}
\noindent From the definition of $\hat{\Omega}^J$ in Algorithm \ref{alg:omega_est_two_stage} for $j=1,\ldots p$, we have
\begin{align*}
    \Vert \hat{\Omega}^{J}_{\cdot j}  - \hat{G}^J_{\cdot j}\Vert_{1}\leq \Vert \hat{G}^J_{\cdot j}  - \Omega^*_{\cdot j}\Vert_{1}.
\end{align*}
Hence we have,
\begin{align}
    \Vert \hat{\Omega}^{J}_{\cdot j}  - \Omega^*_{\cdot j}\Vert_{1}&\leq \Vert \hat{\Omega}^{J}_{\cdot j}  - \hat{G}^J_{\cdot j}\Vert_{1}+\Vert \hat{G}^J_{\cdot j}  - \Omega^*_{\cdot j}\Vert_{1} \notag \\
    &\leq     2C^\prime s_1\bar{\lambda}_n\notag \\
    &\leq C^{\prime\prime} s_1\sqrt{\frac{\log(p\vee n)}{n}},\label{inq:omega_J_1}
\end{align}
where $C^{\prime\prime}$ is a constant depending on $(\bar{A}_0, C_4, \sigma_1,\delta,\rho)$ only.

For $j=1,\ldots, p$, we have
\begin{align*}
     (\hat{\Omega}^{J}_{\cdot j}-\Omega^*_{\cdot j})^\T \hat{\Sigma}_{J}  (\hat{\Omega}^{J}_{\cdot j}-\Omega^*_{\cdot j}) &\leq \Vert \hat{\Omega}^{J}_{\cdot j}-\Omega^*_{\cdot j}\Vert_1 \Vert \hat{\Sigma}_{J}  (\hat{\Omega}^{J}_{\cdot j}-\Omega^*_{\cdot j})\Vert_\infty \\
     &\leq \Vert \hat{\Omega}^{J}_{\cdot j}-\Omega^*_{\cdot j}\Vert_1 (\Vert \hat{\Sigma}_{J}  \hat{\Omega}^{J}-I_p\Vert_{\max} + \Vert \hat{\Sigma}_{J}  \Omega^*-I_p\Vert_{\max}) \\
     &\leq \Vert \hat{\Omega}^{J}  - \Omega^*\Vert_{1}(\Vert \hat{\Sigma}_{J}  \hat{\Omega}^{J}-I_p\Vert_{\max} + \Vert \hat{\Sigma}_{J}  \Omega^*-I_p\Vert_{\max}) \\
     &\leq  2C^\prime \bar{A}_0s_1\bar{\lambda}^2_n.
\end{align*}

\noindent Hence for $j=1,\ldots,p$, we have
\begin{align}
   \Vert \hat{\Omega}^{J}_{\cdot j}-\Omega^*_{\cdot j}\Vert_2^2 &\leq  C_4(\hat{\Omega}^{J}_{\cdot j}-\Omega^*_{\cdot j})^\T \Sigma  (\hat{\Omega}^{J}_{\cdot j}-\Omega^*_{\cdot j}) \quad\text{(from Assumption \ref{sigma_eigenvalue} (ii))}  \notag \\
  &\leq C_4 (\hat{\Omega}^{J}_{\cdot j}-\Omega^*_{\cdot j})^\T \hat{\Sigma}_J  (\hat{\Omega}^{J}_{\cdot j}-\Omega^*_{\cdot j}) +  C_4(\hat{\Omega}^{J}_{\cdot j}-\Omega^*_{\cdot j})^\T (\Sigma- \hat{\Sigma}_J ) (\hat{\Omega}^{J}_{\cdot j}-\Omega^*_{\cdot j}) \notag \\
  &\leq C_4(\hat{\Omega}^{J}_{\cdot j}-\Omega^*_{\cdot j})^\T \hat{\Sigma}_J  (\hat{\Omega}^{J}_{\cdot j}-\Omega^*_{\cdot j}) +  C_4 \Vert \hat{\Omega}^{J}_{\cdot j}-\Omega^*_{\cdot j}\Vert_1^2 \Vert\Sigma- \hat{\Sigma}_J \Vert_{\max}  \notag \\
  &\leq  2C_4 C^\prime \bar{A}_0s_1\bar{\lambda}^2_n + 4C_4C^{\prime 2}s^2_1\bar{\lambda}^3_n \notag \\
  &\leq   C^{\prime\prime\prime} s_1\frac{\log(p\vee n)}{n} , \label{inq:omega_J_2} \quad\text{(From $\delta = C_4 s_1 \bar{\lambda}_n(1+\eta_0)^2 \leq 1$)}
\end{align}
where $C^{\prime\prime\prime}$ is a constant depending on $(\bar{A}_0, C_4, \delta,\sigma_1,\rho)$ only.

From \eqref{inq:omega_J_1}, we have
\begin{align*}
    \Vert \hat{\Omega}^{J}_{S_\xi,S_\xi}-\Omega^*_{S_\xi,S_\xi} \Vert_1 &\leq \max_{j\in S_\xi} \Vert \hat{\Omega}^J_{\cdot j}-\Omega_{\cdot j}^*\Vert_1 \\
     &\leq C^{\prime\prime} s_1\sqrt{\frac{\log(p\vee n)}{n}}.
\end{align*}

From the symmetry of $\hat{\Omega}^{J}$ and $\Omega^*$, we have
\begin{align}
    \Vert \hat{\Omega}^{J}_{S_\xi,S_\xi}-\Omega^*_{S_\xi,S_\xi} \Vert_2
     &\leq \Vert (\hat{\Omega}^{J}-\Omega^*)_{S_\xi,S_\xi} \Vert_1 \notag \\
     &\leq  C^{\prime\prime} s_1\sqrt{\frac{\log(p\vee n)}{n}}. \label{inq:omega_J_3}
\end{align}

From (\ref{inq:omega_J_2}), we also have
\begin{align}
    \Vert \hat{\Omega}^{J}_{S_\xi,S_\xi}-\Omega^*_{S_\xi,S_\xi}\Vert_2
    &\leq  \Vert \hat{\Omega}^{J}_{S_\xi,S_\xi}-\Omega^*_{S_\xi,S_\xi} \Vert_{\F} \notag \\
    &\leq \sqrt{\Tr( (\hat{\Omega}^{J}_{S_{\xi},S_{\xi}}- \Omega^*_{S_{\xi},S_{\xi}})(\hat{\Omega}^{J}_{S_{\xi},S_{\xi}}- \Omega^*_{S_{\xi},S_{\xi}})}\notag\\
    &\leq \sqrt{\sum_{j\in S_{\xi}} (\hat{\Omega}^{J}_{\cdot j} - \Omega^*_{\cdot_j})^\T (\hat{\Omega}^{J}_{\cdot j} - \Omega^*_{\cdot_j})}  \notag\\
    &\leq \sqrt{C^{\prime\prime\prime}} \sqrt{\Vert s_{1\xi}\Vert_0 s_1 }\sqrt{\frac{\log(p\vee n)}{n}}. \label{inq:omega_J_4}
\end{align}

From (\ref{inq:omega_J_3}) and (\ref{inq:omega_J_4}), we have
\begin{align}
    \Vert \hat{\Omega}^{J}_{S_\xi,S_\xi}-\Omega^*_{S_\xi,S_\xi} \Vert_2 \leq (C^{\prime\prime} + \sqrt{C^{\prime\prime\prime}})\min\{\sqrt{\Vert s_{1\xi}\Vert_0s_1}, s_1\}\sqrt{\frac{\log(p\vee n)}{n}}.
\end{align}

\end{prf}

\subsubsection{Error bound of the estimators in Algorithms \ref{alg:est_M_two_stage} }

\begin{lem}
\label{lem:empirical_compatibility_condition_v}

Under the event $\mathcal{A}_{11}\cap\mathcal{A}_{21}$ (Lemma \ref{lem:bound_of_random_quan} (i), (ii)) and suppose that Condition \ref{con:beta_hat_J} is satisfied, for any $\eta_0>1$ and any $S\subset\{1,\ldots n\}$, let $W = \frac{1}{|J|}\sum_{i\in J}\phi_2(Y_i,X_i^\T\hat{\beta}_J)X_iX_i^\T \in \mathbb{R}^{p \times p}$ and $\delta=C_4|S|\sqrt{s_0}\bar{\lambda}_n  (1+\eta_0)^2 \leq1$, then the condition $CompCondition(W, S, \sqrt{\frac{1}{C_4} (1-\delta)}, \eta_0)$ holds, where $\bar{\lambda}_n=\kappa_1\sqrt{\frac{\log(p\vee n)}{n}}$ and $\kappa_1=C_7 C_8 \sqrt{(K(\frac{1}{2})(1+(d(\frac{1}{2}) \log(2))^{-2}) +\frac{1}{d(\frac{1}{2})}) \sigma_1^4}+\sqrt{\frac{3}{c(1)\rho}}K(1)(1+(d(1) \log(2))^{-1})C_6\sigma_1^2$.

\end{lem}
\begin{prf}
From Lemma \ref{lem:H_diff}, under event $\mathcal{A}_{11}\cap\mathcal{A}_{21}$ (Lemma \ref{lem:bound_of_random_quan}(i), (ii)) and suppose that Condition \ref{con:beta_hat_J} is satisfied, we have
\begin{align*}
    &\Vert W - H\Vert_{\max}  \\
    &\leq ( C_7 C_1 \sqrt{(K(\frac{1}{2})(1+(d(\frac{1}{2}) \log(2))^{-2}) +\frac{1}{d(\frac{1}{2})}) \sigma_1^4}+\sqrt{\frac{3}{c(1)\rho}}K(1)(1+(d(1) \log(2))^{-1})C_6\sigma_1^2)\\
    &\times \sqrt{s_0\frac{\log(p\vee n)}{n}}.
\end{align*}
For any $S\subset \{1,\ldots,p\}$ and $\eta_0 > 1$, let $b\in\mathbb{R}^p$ satisfying
\begin{align*}
    \sum_{j\notin S} |b_j|\leq \eta_0 (\sum_{j\in S}|b_j|).
\end{align*}
We have
\begin{align*}
    \frac{1}{C_4}  (\sum_{j\in S}|b_j|)^2 &\leq   \frac{1}{C_4}  |S| (\sum_{j\in S}b_j^2) \\
    &\leq  \frac{1}{C_4}  |S|  \Vert b\Vert_2^2 \\
    &\leq   |S| b^\T H b \quad\text{(From Assumption \ref{ass:H_eigen_value}(ii))} \\
    &=  |S| b^\T W b + b^\T (H-W)b.
\end{align*}
For $b^\T (H-W)b$, we have
\begin{align*}
    b^\T (H-W)b \leq \Vert b \Vert_1^2 \sqrt{s_0}\bar{\lambda}_n.
\end{align*}

\noindent Hence we have
\begin{align*}
     \frac{1}{C_4}   (\sum_{j\in S}|b_j|)^2&\leq  |S| b^\T W b + |S|\sqrt{s_0}\bar{\lambda}_n \Vert b\Vert_1^2 \\
     &\leq   |S| b^\T W b + |S|\sqrt{s_0}\bar{\lambda}_n (1+\eta_0)^2(\sum_{j\in S}|b_j|)^2.
\end{align*}
That is, we have
\begin{align*}
    \frac{1}{C_4} (1-C_4|S|\sqrt{s_0}\bar{\lambda}_n (1+\eta_0)^2) (\sum_{j\in S}|b_j|)^2\leq  |S| b^\T W b.
\end{align*}
Hence
\begin{align*}
    \frac{1}{C_4} (1-\delta) (\sum_{j\in S}|b_j|)^2\leq  |S| b^\T W b.
\end{align*}
\end{prf}

\begin{lem}
\label{lem:inf_norm_con_GLM}
Under event $\mathcal{A}_{8}\cap \mathcal{A}_{9}$(Lemma \ref{lem:bound_of_random_quan}(viii)(ix)) and suppose that Condition \ref{con:beta_hat_J} is satisfied, let $W=\frac{1}{|J|}\sum_{i\in J}\phi_2(Y_i,X_i^\T\hat{\beta}_J)X_iX_i^\T\in \mathbb{R}^{p\times p}$, then condition $InfinityNorm(W,  M^*,\sqrt{s_0}\bar{\lambda}_n)$ holds, where  $\bar{\lambda}_n=\kappa_1\sqrt{\frac{\log(p\vee n)}{n}}$ and $\kappa_1=C_7 C_8 \sqrt{(K(\frac{1}{2})(1+(d(\frac{1}{2}) \log(2))^{-2}) +\frac{1}{d(\frac{1}{2})}) \sigma_1^4}+\sqrt{\frac{3}{c(1)\rho}}K(1)(1+(d(1) \log(2))^{-1})C_6\sigma_1^2$.
\end{lem}

\begin{prf}
From (\ref{inq:bound_max_entries}), under the event $\mathcal{A}_{8}\cap\mathcal{A}_{9}$(Lemma \ref{lem:bound_of_random_quan}(viii)(ix)) and given Condition \ref{con:beta_hat_J}, we have for $j=1,\ldots, p$,
\begin{align*}
    &\Vert W M_{\cdot j}^* -e_j\Vert_\infty \notag \\
    &\leq (\sqrt{\frac{3}{c(1)\rho}}K(1)(1+(d(1) \log(2))^{-1})  C_6\sigma_1^2 + C_7C_8\sqrt{(K(\frac{1}{2})(1+(d(\frac{1}{2}) \log(2))^{-2})+\frac{1}{d(\frac{1}{2})})\sigma_1^4})\\
    &\times\sqrt{s_0\frac{\log(p\vee n)}{n}}.
\end{align*}
\end{prf}

\errorrateMtwostage*

\begin{prf}
From Lemma \ref{lem:estimating_v}, \ref{lem:empirical_compatibility_condition_v} and \ref{lem:inf_norm_con_GLM}, under event $\mathcal{A}_{11}\cap \mathcal{A}_{21}\cap \mathcal{A}_{8}\cap \mathcal{A}_{9}$(Lemma \ref{lem:bound_of_random_quan} (i), (ii), (viii), (ix)) and suppose that Condition \ref{con:beta_hat_J} is satisfied, we have for $j=1,\ldots, p$,
\begin{align*}
    \Vert \hat{G}^J_{\cdot j}-{M^*}_{\cdot j}\Vert_{1} & \leq C^\prime s_1\sqrt{s_0}\bar{\lambda}_n,
\end{align*}
where $C^\prime $ is a constant depending on $(\bar{A}_0,C_4,\delta)$.

\noindent Hence we have
\begin{align*}
    \Vert \hat{G}^J-{M^*}\Vert_{1} & \leq C^\prime s_1\sqrt{s_0}\bar{\lambda}_n.
\end{align*}

\noindent  Under the event $\mathcal{A}_{8}\cap \mathcal{A}_{9}$ (Lemma \ref{lem:bound_of_random_quan} (viii)(ix)) and suppose that Condition \ref{con:beta_hat_J} is satisfied, from (\ref{inq:bound_max_entries}), we have
\begin{align*}
   \Vert \hat{H}_JM^* - I_p \Vert \leq \bar{A}_0\sqrt{s_0}\bar{\lambda}_n.
\end{align*}
\noindent From the definition of $\hat{M}^J$ in Algorithm \ref{alg:est_M}, for $j=1,\ldots, p$, we have
\begin{align*}
    \Vert \hat{M}^{J}_{\cdot j}  - \hat{G}^J_{\cdot j}\Vert_{1}\leq \Vert \hat{G}^J_{\cdot j}  - M^*_{\cdot j}\Vert_{1}.
\end{align*}
Hence we have,
\begin{align}
    \Vert \hat{M}^{J}_{\cdot j}  - M^*_{\cdot j}\Vert_{1}&\leq \Vert \hat{M}^{J}_{\cdot j}  - \hat{G}^J_{\cdot j}\Vert_{1}+\Vert \hat{G}^J_{\cdot j}  - M^*_{\cdot j}\Vert_{1} \notag \\
    &\leq     2C^\prime s_1\sqrt{s_0}\bar{\lambda}_n\notag\\
    &\leq C^{\prime\prime} s_1\sqrt{s_0}\sqrt{\frac{\log(p\vee n)}{n}}\label{inq:M_two_stage_J_1}
\end{align}
where $C^{\prime\prime}$ is a constant depending on $(\bar{A}_0, C_4,C_6,C_7,C_8,\sigma_1, \delta, \rho)$ only.

For $j=1,\ldots, p$, we have
\begin{align}
     (\hat{M}^{J}_{\cdot j}-M^*_{\cdot j})^\T \hat{H}_{J}  (\hat{M}^{J}_{\cdot j}-M^*_{\cdot j}) &\leq \Vert \hat{M}^{J}_{\cdot j}-M^*_{\cdot j}\Vert_1 \Vert \hat{H}_{J}  (\hat{M}^{J}_{\cdot j}-M^*_{\cdot j})\Vert_\infty \notag \\
     &\leq \Vert \hat{M}^{J}_{\cdot j}-M^*_{\cdot j}\Vert_1 (\Vert \hat{H}_{J}  \hat{M}^{J}-I_p\Vert_{\max} + \Vert \hat{H}_{J}  M^*-I_p\Vert_{\max}) \notag \\
     &\leq \Vert \hat{H}^{J}  - M^*\Vert_{1}(\Vert \hat{H}_{J}  \hat{M}^{J}-I_p\Vert_{\max} + \Vert \hat{H}_{J}  M^*-I_p\Vert_{\max}) \notag  \\
     &\leq  2C^\prime \bar{A}_0s_1s_0\bar{\lambda}^2_n. \notag
\end{align}

\noindent Hence we have
\begin{align}
  \Vert \hat{M}^{J}_{\cdot j}-M^*_{\cdot j}\Vert_2^2 &\leq  C_4(\hat{M}^{J}_{\cdot j}-M^*_{\cdot j})^\T H  (\hat{M}^{J}_{\cdot j}-M^*_{\cdot j}) \quad\text{(From Assumption \ref{ass:H_eigen_value}(ii))} \notag \\
  &\leq C_4 (\hat{M}^{J}_{\cdot j}-M^*_{\cdot j})^\T \hat{H}_J  (\hat{M}^{J}_{\cdot j}-M^*_{\cdot j}) +  C_4(\hat{M}^{J}_{\cdot j}-M^*_{\cdot j})^\T (H- \hat{H}_J ) (\hat{M}^{J}_{\cdot j}-M^*_{\cdot j}) \notag \\
  &\leq C_4(\hat{M}^{J}_{\cdot j}-M^*_{\cdot j})^\T \hat{H}_J  (\hat{M}^{J}_{\cdot j}-M^*_{\cdot j}) + C_4 \Vert \hat{M}^{J}_{\cdot j}-M^*_{\cdot j}\Vert_1^2 \Vert H- \hat{H}_J \Vert_{\max} \notag\\
  &\leq  2C^\prime C_4 \bar{A}_0s_1s_0\bar{\lambda}^2_n + 4C^{\prime 2}C_4s^2_1s_0\sqrt{s_0}\bar{\lambda}^3_n \notag\\
  &\leq   2C^\prime C_4 \bar{A}_0s_1s_0\bar{\lambda}^2_n + 4C^{\prime 2}s_1s_0\bar{\lambda}^2_n,\quad\text{(From $\delta=C_4 s_1\sqrt{s_0} \bar{\lambda}_n(1+\eta_0)^2 \leq 1$)} \notag \\
  &\leq C^{\prime\prime\prime} s_1s_0\frac{\log(p\vee n)}{n} , \label{inq:M_two_stage_J_2}
\end{align}
where $C^{\prime\prime\prime}$ is a constant depending on $(\bar{A}_0,C_4,C_6, C_7,C_8,\sigma_1, \delta,\rho)$ only.

From \eqref{inq:M_two_stage_J_1}, we have
\begin{align*}
    \Vert \hat{M}^{J}_{S_\xi,S_\xi}-M^*_{S_\xi,S_\xi} \Vert_1
     &\leq C^{\prime\prime} s_1\sqrt{s_0}\sqrt{\frac{\log(p\vee n)}{n}}.
\end{align*}

From the symmetry of $\hat{M}^{J}$ and $M^*$, we have
\begin{align}
    \Vert \hat{M}^{J}_{S_\xi,S_\xi}-M^*_{S_\xi,S_\xi} \Vert_2
     &\leq \Vert \hat{M}^{J}_{S_\xi,S_\xi}-M^*_{S_\xi,S_\xi} \Vert_1 \notag \\
     &\leq  C^{\prime\prime} s_1\sqrt{s_0}\sqrt{\frac{\log(p\vee n)}{n}}. \label{inq:M_two_stage_J_3}
\end{align}

From (\ref{inq:M_two_stage_J_2}), we also have
\begin{align}
    \Vert \hat{M}^{J}_{S_\xi,S_\xi}-M^*_{S_\xi,S_\xi}\Vert_2 &\leq  \Vert \hat{M}^{J^\c}_{S_\xi,S_\xi}-M^*_{S_\xi,S_\xi} \Vert_2 \notag \\
    &\leq  \Vert \hat{M}^{J}_{S_\xi,S_\xi}-M^*_{S_\xi,S_\xi} \Vert_{\F} \notag \\
    &\leq \sqrt{\Tr( (\hat{M}^{J}_{S_{\xi},S_{\xi}}- M^*_{S_{\xi},S_{\xi}})(\hat{M}^{J}_{S_{\xi},S_{\xi}}- M^*_{S_{\xi},S_{\xi}})}\notag\\
    &\leq \sqrt{\sum_{j\in S_{\xi}} (\hat{M}^{J}_{\cdot j} - M^*_{\cdot_j})^\T (\hat{M}^{J}_{\cdot j} - M^*_{\cdot_j})}  \notag\\
    &\leq \sqrt{C^{\prime\prime\prime}} \sqrt{\Vert s_{1\xi}\Vert_0 s_1 } \sqrt{s_0}\sqrt{\frac{\log(p\vee n)}{n}}. \label{inq:M_two_stage_J_4}
\end{align}
where $C^{\prime\prime\prime}$ is a constant depending on $(\bar{A}_0,C_4, C_6, C_7,C_8,\sigma_1, \delta,\rho)$ only.

From (\ref{inq:M_two_stage_J_3}) and (\ref{inq:M_two_stage_J_4}), we have
\begin{align}
    \Vert \hat{M}^{J}_{S_\xi,S_\xi}-M^*_{S_\xi,S_\xi} \Vert_2 \leq (C^{\prime\prime} + \sqrt{C^{\prime\prime\prime}})\min\{\sqrt{\Vert s_{1\xi}\Vert_0s_1}, s_1  \} \sqrt{s_0} \sqrt{\frac{\log(p\vee n)}{n}}.
\end{align}
\end{prf}

\begin{lem}
\label{lem:H_diff}
$\newline$
Under event $\mathcal{A}_{11}\cap \mathcal{A}_{21}$(Lemma \ref{lem:bound_of_random_quan}(i), (ii)) and suppose that Assumption \ref{bounded_higher_differentiate} and Condition \ref{con:beta_hat_J} are satisfied, we have
\begin{align*}
    \Vert \frac{1}{|J|}\sum_{i\in J}\phi_2(Y_i,X_i^\T\hat{\beta}_J)X_iX_i^\T-H\Vert_{\max}&\leq C \sqrt{s_0\frac{\log(p\vee n)}{n}}.
\end{align*}
Here we denote \\
$C=( C_7 C_8 \sqrt{(K(\frac{1}{2})(1+(d(\frac{1}{2}) \log(2))^{-2}) +\frac{1}{d(\frac{1}{2})}) \sigma_1^4}+\sqrt{\frac{3}{c(1)\rho}}K(1)(1+(d(1) \log(2))^{-1})C_6\sigma_1^2)$ and $K(\frac{1}{2})$, $K(1)$, $d(\frac{1}{2})$ and $d(1)$ are constants from Lemma \ref{lem:rv_centering}.

\end{lem}

\begin{prf}

\noindent We have for $1\leq l\leq p, 1\leq m \leq p$,
\begin{align}
    &|\frac{1}{|J|}\sum_{i\in J} \phi_2(Y_i,X_i^\T\hat{\beta}_J)X_{il}X_{im} - H_{lm} | \notag \\
    &\leq  |\frac{1}{|J|}\sum_{i\in J}\phi_1(Y_i,X_i^\T\hat{\beta}_J)X_{il}X_{im} - \phi_2(Y_i,X_i^\T\bar\beta)X_{il}X_{im}| + |\frac{1}{|J|}\sum_{i\in J}\phi_1(Y_i,X_i^\T\bar\beta)X_{il}X_{im} - H_{lm}| \notag\\
    &\leq   \frac{1}{|J|}\sum_{i\in J}|\phi_2(Y_i,X_i^\T\hat{\beta}_J)X_{il}X_{im} - \phi_2(Y_i,X_i^\T\bar\beta)X_{il}X_{im}| + |\frac{1}{|J|}\sum_{i\in J}\phi_2(Y_i,X_i^\T\bar\beta)X_{il}X_{im} - H_{lm}| \notag\\
    &\leq  C_7 \frac{1}{|J|}\sum_{i\in J}|(X_i^\T\hat{\beta}_J - X_i^\T\bar\beta)X_{il}X_{im}| + |\frac{1}{|J|}\sum_{i\in J}\phi_2(Y_i,X_i^\T\bar\beta)X_{il}X_{im} - H_{lm}| \quad\text{(From Assumption \ref{bounded_higher_differentiate})}\notag\\
    &\leq C_7 \sqrt{\frac{1}{|J|}\sum_{i\in J}(X_i^\T(\hat{\beta}_J-\bar\beta))^2}\sqrt{\frac{1}{|J|}\sum_{i\in J} X_{il}^2X_{im}^2} + |\frac{1}{|J|}\sum_{i\in J}\phi_2(Y_i,X_i^\T\bar\beta)X_{il}X_{im} - H_{lm}|\notag\\
    &\leq C_7\frac{n}{|J|}\sqrt{\frac{1}{n}\sum_{i=1}^n(X_i^\T(\hat{\beta}_J-\bar\beta))^2}\sqrt{\frac{1}{n}\sum_{i\in J} X_{il}^2X_{im}^2} + |\frac{1}{|J|}\sum_{i\in J}\phi_2(Y_i,X_i^\T\bar\beta)X_{il}X_{im} - H_{lm}|. \label{inq:empirical_H_H_diff}
\end{align}
Under event $\mathcal{A}_{11}\cap\mathcal{A}_{21}$(Lemma \ref{lem:bound_of_random_quan}(i), (ii)) and suppose that Condition \ref{con:beta_hat_J} is satisfied, from (\ref{inq:empirical_H_H_diff}), we have  for $1\leq l\leq p, 1\leq m \leq p$,
\begin{align*}
    &|\frac{1}{|J|}\sum_{i\in J} \phi_2(Y_i,X_i^\T\hat{\beta}_J)X_{il}X_{im} - H_{lm} |  \notag \\
    &\leq ( C_7 C_1\frac{1}{\rho} \sqrt{(K(\frac{1}{2})(1+(d(\frac{1}{2}) \log(2))^{-2}) +\frac{1}{d(\frac{1}{2})}) \sigma_1^4}\\
    &+\sqrt{\frac{3}{c(1)\rho}}K(1)(1+(d(1) \log(2))^{-1})C_6\sigma_1^2)\sqrt{s_0\frac{\log(p\vee n)}{n}} .
\end{align*}
That is, we have
\begin{align*}
    &\Vert \frac{1}{|J|}\sum_{i\in J}\phi_2(Y_i,X_i^\T\hat{\beta}_J)X_iX_i^\T - H\Vert_{\max}  \\
    &\leq ( C_7 C_1 \frac{1}{\rho}\sqrt{(K(\frac{1}{2})(1+(d(\frac{1}{2}) \log(2))^{-2}) +\frac{1}{d(\frac{1}{2})}) \sigma_1^4}\\
    &+\sqrt{\frac{3}{c(1)\rho}}K(1)(1+(d(1) \log(2))^{-1})C_6\sigma_1^2)\sqrt{s_0\frac{\log(p\vee n)}{n}}.
\end{align*}

\end{prf}

\section{Relation between the first-stage estimators in Algorithm \ref{alg:omega_est_two_stage} and \ref{alg:est_M_two_stage} and the Lasso nodewise estimators in \cite{van_de_Geer_2014}}

\subsection{Relation between the first-stage estimator $\hat{G}^J$ in Algorithm \ref{alg:omega_est_two_stage} and Lasso nodewise estimator $\hat{\Omega}_{\mytext{node}}$ in \cite{van_de_Geer_2014}}

\label{sec:first_stage_lasso_nodel}

From the definition of $\hat{G}^J$, for $j=1,\ldots, p$, we have
\begin{align}
    \hat{G}^J_{\cdot j} = \argmin_{\Theta \in\mathbb{R}^{p}} \frac{1}{2} \Theta^\T\hat{\Sigma}_J \Theta-\Theta_j +\rho_n \Vert \Theta \Vert_1.  \label{eq:G_J_computation}
\end{align}

The $j$-th column of Lasso nodewise estimator $\hat{\Omega}_{\mytext{node}}$ applied on $X_J$ is computed as below,
\begin{align}
    \hat{\gamma}_j &= \argmin_{\gamma \in \mathbb{R}^{p-1}} \frac{1}{2 |J| }\Vert X_{J, \cdot j} - X_{J, - j} \gamma\Vert_2^2 + \rho_n\Vert\gamma\Vert_1 , \label{eq:gamma_definition}\\
    \hat{\tau}_j^2 &=  (X_{J, \cdot j} - X_{J, - j} \hat{\gamma}_j)^\T X_{J, \cdot j} /|J| , \notag \\
    \hat{\Omega}_{\mytext{node},\cdot, j} &= \hat{\tau}_j^{-2} \begin{pmatrix} -\hat{\gamma}_{j,1}&\hdots&-\hat{\gamma}_{j,j-1}&1 &\hat{\gamma}_{j,j} &\hdots& -\hat{\gamma}_{j,p-1}\end{pmatrix}^\T  .\notag
\end{align}
We consider a variation of (\ref{eq:gamma_definition}) as below:
\begin{align}
    (\hat{\kappa},\hat{\gamma}_j) &= \argmin_{\kappa \in \mathbb{R},\gamma \in \mathbb{R}^{p-1}} \frac{\kappa^2}{2 |J| }\Vert X_{J, \cdot j} - X_{J, - j} \gamma\Vert_2^2 -\kappa+ |\kappa| \rho_n\Vert\gamma\Vert_1 + |\kappa| \rho_n . \label{eq:new_obj_function}
\end{align}
If we reparameterize $(\kappa, \gamma)$ as
\begin{align*}
    \Theta = \begin{pmatrix}
        -\gamma_1\kappa&\hdots& -\gamma_{j-1}\kappa & \kappa &-\gamma_{j}\kappa&\hdots &-\gamma_{p-1}\kappa
    \end{pmatrix},
\end{align*}
then the objective function in (\ref{eq:new_obj_function}) can be expressed as
\begin{align}
  \frac{\kappa^2}{2 |J| }\Vert X_{J, \cdot j} - X_{J, - j} \gamma\Vert_2^2 -\kappa+ |\kappa| \rho_n\Vert\gamma\Vert_1 + |\kappa| \rho_n=  \frac{1}{2} \Theta^\T\hat{\Sigma}_J \Theta-\Theta_j +\rho_n \Vert \Theta \Vert_1 . \label{eq:new_obj_function_definition}
\end{align}
From (\ref{eq:G_J_computation}) and (\ref{eq:new_obj_function_definition}), $\hat{G}^J_{\cdot j}$ can be viewed as a variation of Lasso nodewise regression in \citesupp{van_de_Geer_2014}.

\subsection{Relation between the first-stage estimator $\hat{G}^J$ in Algorithm \ref{alg:est_M_two_stage} and extended Lasso nodewise estimator $\hat{M}_{\mytext{node}}$ in \cite{van_de_Geer_2014}}

\label{sec:first_stage_lasso_nodel_glm}

From the definition of $\hat{G}^J$, for $j=1,\ldots, p$, we have
\begin{align}
    \hat{G}^J_{\cdot j} = \argmin_{\Theta \in\mathbb{R}^{p}} \frac{1}{2} \Theta^\T\hat{H}_J \Theta-\Theta_j +\rho_n \Vert \Theta \Vert_1.  \label{eq:G_J_computation_GLM}
\end{align}

The $j$-th column of weighted Lasso nodewise estimator $\hat{M}_{\mytext{node}}$ applied on $X_J$ and $\hat{\beta}^\prime_J$ is computed as below,
\begin{align}
    \hat{\gamma}_j &= \argmin_{\gamma \in \mathbb{R}^{p-1}} \frac{1}{2}(\hat{H}_{J; jj} -2\hat{H}_{J; j,-j}\gamma + \gamma^\T \hat{H}_{J; -j,-j}\gamma)+\rho_n \Vert\gamma\Vert_1 , \label{eq:gamma_definition_GLM}\\
    \hat{\tau}_j^2 &=  \hat{H}_{J; jj} - \hat{H}_{J; j,-j} \hat{\gamma}_j ,\notag \\
    \hat{\Omega}_{\mytext{node},\cdot, j} &= \hat{\tau}_j^{-2} \begin{pmatrix} -\hat{\gamma}_{j,1}&\hdots&-\hat{\gamma}_{j,j-1}&1 &\hat{\gamma}_{j,j} &\hdots& -\hat{\gamma}_{j,p-1}\end{pmatrix}^\T . \notag
\end{align}
We consider a variation of (\ref{eq:gamma_definition_GLM}) as
\begin{align}
    (\hat{\kappa},\hat{\gamma}_j) &= \frac{\kappa^2}{2}(\hat{H}_{J; jj} -2\hat{H}_{J; j,-j}\gamma + \gamma^\T \hat{H}_{J; -j,-j}\gamma)-\kappa+|\kappa|\rho_n \Vert\gamma\Vert_1 +|\kappa|\rho_n . \label{eq:new_obj_function_GLM}
\end{align}
If we reparameterize $(\kappa, \gamma)$ as
\begin{align*}
    \Theta = \begin{pmatrix}
        -\gamma_1\kappa&\hdots& -\gamma_{j-1}\kappa & \kappa &-\gamma_{j}\kappa&\hdots &-\gamma_{p-1}\kappa
    \end{pmatrix},
\end{align*}
then the objective function in (\ref{eq:new_obj_function_GLM}) can be expressed as
\begin{align}
  &\frac{\kappa^2}{2}(\hat{H}_{J; jj} -2\hat{H}_{J; j,-j}\gamma + \gamma^\T \hat{H}_{J; -j,-j}\gamma)-\kappa+|\kappa|\rho_n \Vert\gamma\Vert_1 +|\kappa|\rho_n \notag \\
  &=  \frac{1}{2} \Theta^\T\hat{H}_J \Theta-\Theta_j +\rho_n \Vert \Theta \Vert_1 . \label{eq:new_obj_function_definition_GLM}
\end{align}
From (\ref{eq:G_J_computation_GLM}) and (\ref{eq:new_obj_function_definition_GLM}), $\hat{G}^J$ can be viewed as a variation of weighted Lasso nodewise regression in \citesupp{van_de_Geer_2014}.

\section{Calculation of $\bar\beta$ for misspecified logistic regression in Section~\ref{sec:numerical_study_sparse_M}}

\label{sec:calculation_of_beta_star}

In this section, we calculate $\bar\beta$ for the misspecified logistic regression in Section \ref{sec:numerical_study_sparse_M} which is restated as follows.
The true regression model for data generation is
\begin{align*}
 P(Y_i=1|X_i) = \frac{1}{2}\frac{1}{1+\exp\{-X_i^\T\gamma^*+2\}} + \frac{1}{2}\frac{1}{1+\exp\{-X_i^\T\gamma^*-2\}}, \quad \text{for $i=1,\ldots,n$},
\end{align*}
where $\gamma^*=(0.5,0.5,0.5,0.5,0.075,0_{295})^\T \in \mathbb{R}^p$. We let $(X_1,\ldots,X_n)$ be i.i.d.~as $X_0$ and $X_0=\Sigma^{\frac{1}{2}}Z_0$ where the components of $Z_0$ are i.i.d.~Rademacher variables and $\Sigma \in \mathbb{R}^{p\times p}$ satisfying
\begin{equation}
    \begin{cases}
    \Sigma_{ij}=5 \times 0.1^{|i-j|},\quad\text{if $i\leq 5$, $j\leq 5$ or $i>5, j>5$,} \\
    \Sigma_{ij}= 0,\quad\text{otherwise.}
    \end{cases}  \notag
\end{equation}

We consider the logistic regression of $Y_i$ on $X_i$ and $\phi(a,b)=-ab+\log(1+\exp(b))$. From the definition in (\ref{eq:target_beta_GLM}), we have
$\bar\beta = \argmin_{\beta\in\mathbb{R}^{p}} \ell(\beta)$,
where
\begin{align*}
    \ell(\beta) =& \mathbb{E}[(-\frac{1}{2}\frac{1}{1+\exp\{-X_0^\T\gamma^*+2\}} -\frac{1}{2}\frac{1}{1+\exp\{-X_0^\T\gamma^*-2\}} ) X_0^\T\beta\notag \\
    &+ \log(1+\exp\{X_0^\T \beta \})] . %\label{def:beta_all}
\end{align*}

We show that $\bar\beta_j=0$ for $j>5$. Consider $\tilde{\beta}=(\tilde\beta_{1:5}, 0_{295})\in \mathbb{R}^p$ with
\begin{align}
    \tilde\beta_{1:5} = &\argmin_{\beta_{1:5}\in\mathbb{R}^{5}} \ell( \beta_{1:5}, 0_{295}) .\label{def:beta_first}
    %\mathbb{E}[(-\frac{1}{2}\frac{1}{1+\exp\{-X_0^\T\gamma^*+2\}} -\frac{1}{2}\frac{1}{1+\exp\{-X_0^\T\gamma^*-2\}} ) X_{0,1:5}^\T \beta_{1:5}\notag \\
    %&+ \log(1+\exp\{X_{0,1:5}^\T \beta_{1:5} \})]. \notag
\end{align}
From the definition of $\tilde\beta_{1:5}$, we have for $j\leq 5$,
\begin{align}
  &\frac{\partial }{\partial \beta_j}\ell(\beta) |_{\tilde\beta} \notag \\
  &=\mathbb{E}[(\frac{1}{1+\exp\{-X_{0,1:5}^\T \tilde\beta_{1:5} \}}-\frac{1}{2}\frac{1}{1+\exp\{-X_0^\T\gamma^*+2\}} -\frac{1}{2}\frac{1}{1+\exp\{-X_0^\T\gamma^*-2\}}   ) X_{0,j}] =0. \label{eq:partial_diff_beta_part1}
\end{align}
For $j> 5$, we have
\begin{align}
   & \frac{\partial }{\partial {\beta}_j}\ell( \beta) |_{\tilde\beta}  \notag \\
    &= \mathbb{E}[(\frac{1}{1+\exp\{-X_{0,1:5}^\T \tilde\beta_{1:5}\}}-\frac{1}{2}\frac{1}{1+\exp\{-X_{0,1:5}^\T\gamma^*_{1:5}+2\}} -\frac{1}{2}\frac{1}{1+\exp\{-X_{0,1:5}^\T\gamma^*_{1:5}-2\}}  ) X_{0,j}] \notag \\
    &=\mathbb{E}[(\frac{1}{1+\exp\{-X_{1,1:5}^\T \tilde\beta_{1:5} \}}-\frac{1}{2}\frac{1}{1+\exp\{-X_{1,1:5}^\T\gamma^*_{1:5}+2\}}-\frac{1}{2}\frac{1}{1+\exp\{-X_{1,1:5}^\T\gamma^*_{1:5}-2\}} ) ]\mathbb{E}[X_{0,j}]\notag  \\
    &\quad\text{(because $X_{0,j}$ is independent with $X_{0,1:5}$ for $j>5$)}\notag \\
    &=0.  \quad\text{(because $\mathbb{E}[X_{0,j}]=0$)} \label{eq:partial_diff_beta_part2}
\end{align}
From (\ref{eq:partial_diff_beta_part1}) and (\ref{eq:partial_diff_beta_part2}), we have $\bar\beta=\tilde{\beta}$, i.e., $\bar\beta$ can be obtained as $\tilde\beta$.
Hence $\bar\beta_j=0$ for $j> 5$.
From (\ref{def:beta_first}), $\bar\beta_{1:5}=\tilde\beta_{1:5}$ can be numerically computed by gradient descent.

\section{Basic lemmas}

\begin{lem}[Bounds of random quantities] For the events as below, we use $\mathcal{A}_{\cdot\cdot}$ to denote events not depending on the independency of $\hat{M}^{J^\c}$ and we use $\mathcal{B}_{\cdot\cdot}$ to denote events depending on the independency of $\hat{M}^{J^\c}$.
\label{lem:bound_of_random_quan}
$\newline$
\noindent (i)

\noindent \textcircled{1} Denote by $\mathcal{A}_{11}$ the event that
\begin{align*}
    \Vert \frac{1}{|J|}\sum_{i\in J}\phi_2(Y_i, X_i^\T\bar\beta)X_iX_i^\T-H\Vert_{\max} \leq \sqrt{\frac{3}{c(1)}}\frac{K(1)(1+(d(1) \log(2))^{-1})}{\sqrt{\rho}}C_6\sigma_1^2\sqrt{\frac{\log(p\vee n)}{n}}.
\end{align*}
\textcircled{2}
Denote by $\mathcal{A}_{12}$ the event that
\begin{align*}
    \Vert \frac{1}{|J^\c|}\sum_{i\in J^\c}\phi_2(Y_i, X_i^\T\bar\beta)X_iX_i^\T-H\Vert_{\max} \leq \sqrt{\frac{3}{c(1)}}\frac{K(1)(1+(d(1) \log(2))^{-1})}{\sqrt{\rho}}C_6\sigma_1^2\sqrt{\frac{\log(p\vee n)}{n}}.
\end{align*}
\textcircled{3}Denote by $\mathcal{A}_{13}$ the event that
\begin{align}
    \Vert \frac{1}{n}\sum_{i=1}^n\phi_2(Y_i, X_i^\T\bar\beta)X_iX_i^\T-H\Vert_{\max} \leq \sqrt{\frac{3}{c(1)}}K(1)(1+(d(1) \log(2))^{-1})C_6\sigma_1^2\sqrt{\frac{\log(p\vee n)}{n}}. \notag %\label{inq:bounded_H_hat_entries}
\end{align}
\textcircled{4}Denote by $\mathcal{A}_{14}$ the event that
\begin{align*}
    \Vert \frac{1}{|J|}\sum_{i\in J}X_iX_i^\T-\Sigma\Vert_{\max} \leq \sqrt{\frac{3}{c(1)}}\frac{K(1)(1+(d(1) \log(2))^{-1})}{\sqrt{\rho}}\sigma_1^2\sqrt{\frac{\log(p\vee n)}{n}}.
\end{align*}

    \noindent Suppose that Assumption \ref{ass:sub_gaussian_entries_GLM} (i) and \ref{bounded_higher_differentiate} are satisfied and $\frac{\log(p\vee n)}{\sqrt{n} }\leq\frac{c(1)\rho}{3}\sqrt{n}$, we have $\mathbb{P}(\mathcal{A}_{11})\geq 1-\frac{2}{ p \vee n }$ and $\mathbb{P}(\mathcal{A}_{12})\geq 1-\frac{2}{ p \vee n }$. Suppose that Assumption \ref{ass:sub_gaussian_entries_GLM} (i) and \ref{bounded_higher_differentiate} are satisfied and $\frac{\log(p\vee n)}{\sqrt{n}}\leq\frac{c(1)}{3}\sqrt{n}$, we have $\mathbb{P}(\mathcal{A}_{13})\geq 1-\frac{2}{ p \vee n }$. Suppose that Assumption \ref{ass:sub_gaussian_entries_GLM} (i) is satisfied and $\frac{\log(p\vee n)}{\sqrt{n} }\leq\frac{c(1)\rho}{3}\sqrt{n}$, we have $\mathbb{P}(\mathcal{A}_{14})\geq 1-\frac{2}{ p \vee n }$
$\newline$

\noindent (ii)

\noindent
\textcircled{1} Denote by $\mathcal{A}_{21}$ the event that
\begin{align*}
    \max_{1\leq l\leq p,1\leq m\leq p} \{\frac{1}{n}\sum_{i=1}^n X_{il}^2 X_{im}^2 \}  \leq (K(\frac{1}{2})(1+(d(\frac{1}{2}) \log(2))^{-2}) +\frac{1}{d(\frac{1}{2})}) \sigma_1^4 .
\end{align*}
\textcircled{2} Denote by $\mathcal{A}_{22}$ the event that
\begin{align*}
    \max_{1\leq k\leq p,1\leq l\leq p,1\leq m\leq p} \{\frac{1}{n}\sum_{i=1}^n |X_{ik} X_{il} X_{im}| \}  \leq (K(\frac{2}{3})(1+(d(\frac{2}{3}) \log(2))^{-2}) +\frac{1}{d(\frac{2}{3})}) \sigma_1^3  .
\end{align*}
\textcircled{3} Denote by $\mathcal{A}_{23}$ the event that
\begin{align*}
    \max_{1\leq k\leq p,1\leq l\leq p,1\leq m\leq p,1\leq q\leq p} \{\frac{1}{n}\sum_{i=1}^n| X_{ik} X_{il} X_{im}X_{iq}| \}  \leq (K(\frac{1}{2})(1+(d(\frac{1}{2}) \log(2))^{-2}) +\frac{1}{d(\frac{1}{2})}) \sigma_1^4 .
\end{align*}
\textcircled{4} Denote by $\mathcal{A}_{24}$ the event that
\begin{align}
    \max_{1\leq i \leq n} \{|\phi_1(Y_i,X_i^\T\bar\beta)|\} \leq 2K(2)(1+(d(2) \log(2))^{-1/2}) +\frac{1}{d(2)})\sigma_2 \sqrt{\log ( n)}. \notag %\label{inq:bouned_max_phi_1}
\end{align}
Suppose that Assumption \ref{ass:sub_gaussian_entries_GLM}(i) is satisfied and $\frac{\log(p\vee n)}{\sqrt{n}}\leq\frac{c(\frac{1}{2})}{3}$, we have $\mathbb{P}(\mathcal{A}_{21})\geq 1-\frac{2}{ p \vee n }$. Suppose that Assumption \ref{ass:sub_gaussian_entries_GLM}(i) is satisfied and $\frac{\log(p\vee n)}{\sqrt{n}}\leq \frac{c(\frac{2}{3})n^{\frac{1}{6}}}{4}$, we have $\mathbb{P}(\mathcal{A}_{22})\geq 1-\frac{2}{ p \vee n }$. Suppose that Assumption \ref{ass:sub_gaussian_entries_GLM}(i) is satisfied and $\frac{\log(p\vee n)}{\sqrt{n}}\leq\frac{c(\frac{1}{2})}{5}$ , we have $\mathbb{P}(\mathcal{A}_{23})\geq 1-\frac{2}{ p \vee n }$.
Suppose that Assumption \ref{ass:residulal_glm_lower} is satisfied and $\log(n)>1$, we have $\mathbb{P}(\mathcal{A}_{24})\geq 1-\frac{2}{  n }$.

$\newline$
\noindent (iii) Denote by $\mathcal{A}_{3}$ the event that
\begin{align}
    |\frac{1}{n} \sum_{i=1}^n \phi^4_1(Y_i,X_i^\T\bar\beta)| \leq (K(\frac{1}{2})(1+(d(\frac{1}{2}) \log(2))^{-2}) +\frac{1}{d(\frac{1}{2})}) \sigma_2^4.  \notag %\label{inq:phi_1_power_4_concentration}
\end{align}
Suppose that Assumption \ref{ass:residulal_glm_lower} is satisfied, we have $\mathbb{P}(\mathcal{A}_{3})\geq 1-2\exp\{-c(\frac{1}{2}) \sqrt{n}\}$.

\noindent (iv) Denote by $\mathcal{A}_4$ the event that
\begin{align*}
    \Vert \frac{1}{|J|} \sum_{j\in J} X_i \phi_1(Y_i,X_i^{\T}\bar\beta)\Vert_\infty \leq  \sigma_1\sigma_2\sqrt{\frac{2}{c(1)\rho}}\sqrt{\frac{\log(p\vee n)}{n}}.
\end{align*}
Suppose that Assumption \ref{ass:sub_gaussian_entries_GLM} (i) and \ref{ass:residulal_glm_lower} are satisfied and $\frac{\log(p\vee n)}{\sqrt{n}}\leq \frac{c(1)\rho}{2}\sqrt{n}$, $\mathbb{P}(\mathcal{A}_4) \geq 1-\frac{2}{p\vee n}$.

$\newline$
\noindent (v) Denote by $\mathcal{A}_5$ the event that
\begin{align*}
    \Vert (\frac{1}{n}\sum_{i=1}^n \phi_2(Y_i, X_i^\T\bar\beta) X_iX_i^\T M^* - I_p)\xi \Vert_\infty \leq \sqrt{\frac{2}{c(1)}}K(1)(1+(d(1) \log(2))^{-1}) C_6\sigma^2_1 \sqrt{\frac{\log(p\vee n)}{n}}\Vert \xi\Vert_2.
\end{align*}
Suppose that Assumption \ref{ass:sub_gaussian_entries_GLM}, \ref{ass:H_eigen_value}(ii) and \ref{bounded_higher_differentiate} are satisfied and $\frac{\log(p\vee n)}{\sqrt{n}}\leq \frac{c(1)\sqrt{n}}{2}$, we have $\mathbb{P}(\mathcal{A}_5)\geq 1-\frac{2}{ p \vee n }$.

$\newline$
\noindent (vi)

\noindent \textcircled{1} Denote by $\mathcal{A}_{61}$ the event that
\begin{align*}
    \max_{1\leq l \leq p, 1\leq m\leq p} \{ \frac{1}{n}\sum_{i=1}^n | (\xi^\T M^* X_i)X_{il} X_{im}| \} \leq (K(\frac{2}{3})(1+(d(\frac{2}{3}) \log(2))^{-\frac{3}{2}}) + \frac{1}{d(\frac{2}{3})})\sigma_1^3 C_4 \Vert \xi\Vert_2 .
\end{align*}
\textcircled{2} Denote by $\mathcal{A}_{62}$ the event that
\begin{align}
    \max_{1\leq i \leq n} \{|\xi^\T M^* X_i|\} \leq (2 K(2)(1+(d(2) \log(2))^{-1/2}) +\frac{1}{d(2)}) \sigma_1 C_4 \Vert \xi\Vert_2 \sqrt{\log( n)}. \notag %\label{inq:bouned_max_xi0_M_xi}
\end{align}
Suppose that Assumption \ref{ass:sub_gaussian_entries_GLM} and \ref{ass:H_eigen_value}(ii) are satisfied and $\frac{\log(p\vee n)}{\sqrt{n}}\leq\frac{c(\frac{2}{3})n^{\frac{1}{6}}}{3}$ hold, we have $\mathbb{P}(\mathcal{A}_{61})\geq 1-\frac{2}{ p \vee n }$. Suppose that Assumption \ref{ass:sub_gaussian_entries_GLM} (ii) and \ref{ass:H_eigen_value}(ii) are satisfied and $\log(n)>1$, we have $\mathbb{P}(\mathcal{A}_{62})\geq 1-\frac{2}{n}$.

$\newline$
\noindent (vii)

\noindent \textcircled{1} Denote by $\mathcal{A}_{71}$ the event that
\begin{align*}
    &|\frac{1}{n} \sum_{i=1}^n ((\xi^\T M^* X_i)^2\phi_1(Y_i,X_i^\T\bar\beta)^2- \mathbb{E}[(\xi^\T M^* X_i)^2\phi_1(Y_1,X_i^\T\bar\beta)^2]) | \\
    &\leq  K(\frac{1}{2})(1+(d(\frac{1}{2}) \log(2))^{-2}) \sigma^2_1\sigma_2^2 C_4^2\Vert \xi\Vert^2_2 n^{-\frac{1}{3}}.
\end{align*}
\textcircled{2} Denote by $\mathcal{A}_{72}$ the event that
\begin{align*}
    |\frac{1}{n} \sum_{i=1}^n (\xi^\T M^* X_i)^4| \leq (K(\frac{1}{2})(1+(d(\frac{1}{2}) \log(2))^{-2}) +\frac{1}{d(\frac{1}{2})}) \sigma^4_1C_4^4 \Vert  \xi\Vert^4_2.
\end{align*}
Suppose Assumption \ref{ass:sub_gaussian_entries_GLM} (ii), \ref{ass:H_eigen_value}(ii) and \ref{ass:residulal_glm_lower} are satisfied, we have $\mathbb{P}(\mathcal{A}_{71})\geq 1-2\exp\{-c(\frac{1}{2}) n^{\frac{1}{3}}\}$. Suppose Assumption \ref{ass:sub_gaussian_entries_GLM} (ii) and \ref{ass:H_eigen_value}(ii) are satisfied, we have $\mathbb{P}(\mathcal{A}_{72})\geq 1-2\exp\{-c(\frac{1}{2}) \sqrt{n}\}$.

$\newline$
\noindent (viii) Denote by $\mathcal{A}_{8}$ the event that
\begin{align*}
    \Vert \frac{1}{|J|}\sum_{i\in J} \phi_2(Y_i, X_i^\T\bar\beta) X_iX_i^\T M^* - I_p \Vert_\infty \leq \sqrt{\frac{3}{c(1)\rho}}K(1)(1+(d(1) \log(2))^{-1})  C_6 \sigma_1^2\sqrt{\frac{\log(p\vee n)}{n}}.
\end{align*}
Suppose Assumption \ref{ass:sub_gaussian_entries_GLM} (i), \ref{bounded_higher_differentiate} and \ref{ass:entries_M_x_0} are satisfied and $\frac{\log(p\vee n)}{\sqrt{n}}\leq \frac{c(1)\rho \sqrt{n}}{3}$, we have $\mathbb{P}(\mathcal{A}_{8})\geq 1-\frac{2}{ p \vee n }$.

$\newline$
\noindent (ix) Denote by $\mathcal{A}_{9}$ the event that,
\begin{align*}
    \max_{1\leq l\leq p,1\leq j\leq p}\{\frac{1}{|J|} \sum_{i\in J} x^2_{il} (M_{\cdot j}^{*\T}X_i)^2\} \leq (K(\frac{1}{2})(1+(d(\frac{1}{2}) \log(2))^{-2})+\frac{1}{d(\frac{1}{2})})\sigma_1^4.
\end{align*}
Suppose Assumption \ref{ass:sub_gaussian_entries_GLM} (i) and \ref{ass:entries_M_x_0} are satisfied and $\frac{\log(p\vee n)}{\sqrt{n}}\leq \frac{c(\frac{1}{2})\rho^{\frac{1}{2}}}{3}$, we have $\mathbb{P}(\mathcal{A}_{9})=1-\frac{2}{p\vee n}$.

$\newline$
\noindent (x) Denote by $\mathcal{B}_{1}$ the event that
\begin{align}
    |\xi^\T (\hat{M}^{J^\c}- M^*)\frac{1}{|J|}\sum_{i\in J}X_i\phi_1(Y_i, X_i^\T\bar\beta)|\leq\frac{1}{\sqrt{c(1)\rho}} \sqrt{\frac{\log(p\vee n)}{n}}\sigma_1\sigma_2\Vert (\hat{M}^{J^\c}-M) \xi\Vert_2. \notag
\end{align}
Suppose that Assumption \ref{ass:sub_gaussian_GLM} and \ref{ass:residulal_glm_lower} are satisfied and $\frac{\log(p\vee n)}{\sqrt{n}}\leq c(1)\rho\sqrt{n}$, we have $\mathbb{P}(\mathcal{B}_{1})\geq 1-\frac{2}{ p \vee n }$.

$\newline$
\noindent (xi) Denote by $\mathcal{B}_{2}$ the event that
\begin{align}
    &|\xi^\T (\hat{M}^{J^\c}- M^*)\sum_{i\in J}\frac{1}{|J|}X_iX_i^\T(\hat{M}^{J^\c}- M^*)\xi| \\
    &\leq(K(1)(1+(d(1) \log(2))^{-1})+\frac{1}{d(1)})\sigma_1^2
   \Vert  (\hat{M}^{J^\c}- M^*)\xi\Vert_2^2 . \notag
\end{align}
Suppose that Assumption \ref{ass:sub_gaussian_GLM} is satisfied and $\frac{\log(p\vee n)}{\sqrt{n} }\leq c(1)\rho\sqrt{n}$, we have $\mathbb{P}(\mathcal{B}_{2})\geq 1-\frac{2}{ p \vee n }$.

$\newline$
\noindent (xii)

\noindent Denote by $\mathcal{B}_{3}$ the event that
\begin{align*}
    &\max_{1\leq l\leq p, 1\leq m\leq p} \{ \frac{1}{|J|}\sum_{i\in J} | (\xi^\T (\hat{M}^{J^\c}- M^*) X_i)X_{il} X_{im}| \} \\
    &\leq (K(\frac{2}{3})(1+(d(\frac{2}{3}) \log(2))^{-\frac{3}{2}}) + \frac{1}{d(\frac{2}{3})})\sigma_1^3 \Vert (\hat{M}^{J^\c}-M^*)\xi\Vert_2.
\end{align*}
Suppose that Assumption \ref{ass:sub_gaussian_GLM} is satisfied and $\frac{\log(p\vee n)}{\sqrt{n}}\leq\frac{c(\frac{2}{3})n^{\frac{1}{6}}\rho^{\frac{2}{3}}}{3}$, we have $\mathbb{P}(\mathcal{B}_{3})\geq 1-\frac{2}{ p \vee n }$.

$\newline$
\noindent (xiii) Denote by $\mathcal{B}_{4}$ the event that
\begin{align}
    \max_{i \in J}\{|\xi^\T(\hat{M}^{J^\c}-M^*)X_i|\} \leq  2K(2)(1+(d(2) \log(2))^{-1/2}) +\frac{1}{d(2)})\sqrt{\log( n)}\sigma_1 \Vert (\hat{M}^{J^\c}-M^*)\xi\Vert_2 .  \notag
\end{align}
Suppose Assumption \ref{ass:sub_gaussian_GLM} is satisfied and $\log(n)>1$, we have $\mathbb{P}(\mathcal{B}_{4})\geq 1-\frac{2}{n}$.

\end{lem}

\begin{prf}
$\newline$
\noindent (i)

\noindent \textcircled{1} Under Assumption \ref{ass:sub_gaussian_entries_GLM} (i) and \ref{bounded_higher_differentiate}, the entries of $\sqrt{\phi_2(Y_i,X_i^\T\bar\beta)}X_i$ for $i=1,\ldots, n$ are sub-gaussian random variables with sub-gaussian norm bounded by $\sqrt{C_6}\sigma_1$. Given $\frac{\log(p\vee n)}{\rho n}\leq\frac{c(1)}{3}$, from Lemma \ref{lem:max_entries_bounds}, we have with probability at least $1-\frac{2}{p\vee n}$,
\begin{align}
    \Vert \frac{1}{|J|}\phi_2(Y_i, X_i^\T\bar\beta)X_iX_i^\T-H\Vert_{\max} \leq \sqrt{\frac{3}{c(1)}}\frac{K(1)(1+(d(1) \log(2))^{-1})}{\sqrt{\rho}}C_6\sigma_1^2\sqrt{\frac{\log(p\vee n)}{n}}.  \label{inq:entries_empirical}
\end{align}

\noindent \textcircled{2} The bound of $\Vert \frac{1}{|J^\c|}\sum_{i\in J^\c}\phi_2(Y_i, X_i^\T\bar\beta)X_iX_i^\T-H\Vert_{\max}$ can be derived similarly as (\ref{inq:entries_empirical}).

\noindent \textcircled{3} The bound of $ \Vert \frac{1}{n}\phi_2(Y_i, X_i^\T\bar\beta)X_iX_i^\T-H\Vert_{\max}$ can be derived similarly as (\ref{inq:entries_empirical}).

\noindent \textcircled{4} The bound of $ \Vert \frac{1}{|J|}\sum_{i\in J}X_iX_i^\T-\Sigma\Vert_{\max}$ can be derived similarly as (\ref{inq:entries_empirical}).

\noindent (ii)
$\newline$
\noindent \textcircled{1} Under Assumption \ref{ass:sub_gaussian_entries_GLM} (i), from Lemma \ref{lem:product-two-sub-gaussian-rv-norm} and \ref{lem:product-two-sub-exponential-rv-norm}, we have
\begin{align*}
    \Vert X_{il}^2X_{im}^2\Vert_{\psi_{\frac{1}{2}}}\leq \sigma_1^4.
\end{align*}
From Lemma \ref{lem:rv_centering}, we have
\begin{align*}
   & \Vert X_{il}^2X_{im}^2-\mathbb{E}[X_{il}^2X_{im}^2]\Vert_{\psi_{\frac{1}{2}}} \leq  K(\frac{1}{2})(1+(d(\frac{1}{2}) \log(2))^{-2}) \sigma_1^4 ,\\
   & \mathbb{E}[X_{il}^2X_{im}^2]\leq \frac{1}{d(\frac{1}{2})}\sigma_1^4.
\end{align*}
From Lemma \ref{lem:alpha_sub_exponential_concentration}, let $t = K(\frac{1}{2})(1+(d(\frac{1}{2}) \log(2))^{-2})\sigma_1^4$ and $3\log(p\vee n)\leq c(\frac{1}{2})\sqrt{n}$, we have with probability at least $1-\frac{2}{ p \vee n ^3}$,
\begin{align*}
    |\frac{1}{n}\sum_{i=1}^n (X_{il}^2X_{im}^2-\mathbb{E}[X_{il}^2X_{im}^2])| \leq K(\frac{1}{2})(1+(d(\frac{1}{2}) \log(2))^{-2})\sigma_1^4.
\end{align*}
Hence with probability at least $1-\frac{2}{ p \vee n }$,
\begin{align*}
   \max_{1\leq l\leq p,1\leq m \leq p} \{ \frac{1}{n}\sum_{i=1}^n X_{il}^2X_{im}^2 \}\leq (K(\frac{1}{2})(1+(d(\frac{1}{2}) \log(2))^{-2})+\frac{1}{d(\frac{1}{2})})\sigma_1^4.
\end{align*}
$\newline$
\noindent \textcircled{2}  Under Assumption \ref{ass:sub_gaussian_entries_GLM} (i), from Lemma \ref{lem:product-three-rv-norm}, we have
\begin{align*}
    \Vert |X_{ik} X_{il} X_{im} |\Vert_{\psi_{\frac{2}{3}}}\leq \sigma_1^3.
\end{align*}
From Lemma \ref{lem:rv_centering}, we have
\begin{align*}
   & \Vert |X_{ik} X_{il} X_{im}| -\mathbb{E}[|X_{ik} X_{il} X_{im}| ]\Vert_{\psi_{\frac{2}{3}}} \leq  K(\frac{2}{3})(1+(d(\frac{2}{3}) \log(2))^{-2}) \sigma_1^3 ,\\
   & |\mathbb{E}[|X_{ik} X_{il} X_{im}| ]|\leq \frac{1}{d(\frac{2}{3})}\sigma_1^3.
\end{align*}
From Lemma \ref{lem:alpha_sub_exponential_concentration}, let $t = K(\frac{3}{3})(1+(d(\frac{2}{3}) \log(2))^{-2})\sigma_1^3$ and $4\log(p\vee n)\leq c(\frac{2}{3})n^{\frac{2}{3}}$, we have with probability at least $1-\frac{2}{ (p \vee n) ^4}$,
\begin{align*}
   | \frac{1}{n}\sum_{i=1}^n (|X_{ik} X_{il} X_{im}| - \mathbb{E}[|X_{ik} X_{il} X_{im}|] )| \leq K(\frac{2}{3})(1+(d(\frac{2}{3}) \log(2))^{-2})\sigma_1^3.
\end{align*}
Hence with probability at least $1-\frac{2}{ p \vee n }$,
\begin{align*}
    \max_{1\leq k\leq p,1\leq l\leq p,1\leq m\leq p} \{|\frac{1}{n} X_{ik} X_{il} X_{im}| \}  \leq (K(\frac{2}{3})(1+(d(\frac{2}{3}) \log(2))^{-2}) +\frac{1}{d(\frac{2}{3})}) \sigma_1^3 .
\end{align*}

\noindent \textcircled{3} Under Assumption \ref{ass:sub_gaussian_entries_GLM} (i), from Lemma \ref{lem:product-two-sub-gaussian-rv-norm} and \ref{lem:product-two-sub-exponential-rv-norm}, we have
\begin{align*}
    \Vert |X_{ik} X_{il} X_{im} X_{iq}| \Vert_{\psi_{\frac{1}{2}}}\leq \sigma_1^4.
\end{align*}
From Lemma \ref{lem:rv_centering}, we have
\begin{align*}
   & \Vert |X_{ik} X_{il} X_{im} X_{iq}|-\mathbb{E}[|X_{ik} X_{il} X_{im} X_{iq}|]\Vert_{\psi_{\frac{1}{2}}} \leq  K(\frac{1}{2})(1+(d(\frac{1}{2}) \log(2))^{-2}) \sigma_1^4 ,\\
   & |\mathbb{E}[|X_{ik} X_{il} X_{im} X_{iq}|]|\leq \frac{1}{d(\frac{1}{2})}\sigma_1^4.
\end{align*}
From Lemma \ref{lem:alpha_sub_exponential_concentration}, let $t = K(\frac{1}{2})(1+(d(\frac{1}{2}) \log(2))^{-2})\sigma_1^4$ and $5\log(p\vee n)\leq c(\frac{1}{2})\sqrt{n}$, we have with probability at least $1-\frac{2}{ (p \vee n) ^5}$,
\begin{align*}
    |\frac{1}{n}\sum_{i=1}^n (|X_{ik} X_{il} X_{im} X_{iq}| - \mathbb{E}[|X_{ik} X_{il} X_{im} X_{iq}|]) |\leq K(\frac{1}{2})(1+(d(\frac{1}{2}) \log(2))^{-2})\sigma_1^4.
\end{align*}
Hence with probability at least $1-\frac{2}{ p \vee n }$,
\begin{align*}
    \max_{1\leq k\leq p,1\leq l\leq p,1\leq m\leq p,1\leq q\leq p} \{\frac{1}{n} |X_{ik} X_{il} X_{im}X_{iq}| \}  \leq (K(\frac{1}{2})(1+(d(\frac{1}{2}) \log(2))^{-2}) +\frac{1}{d(\frac{1}{2})}) \sigma_1^4 .
\end{align*}

$\newline$
\noindent \textcircled{4} Under Assumption \ref{ass:residulal_glm_lower}, we have $\phi_1(Y_i,X_i^\T\bar\beta)$ are sub-gaussian random variables with sub-gaussian norm bounded by $\sigma_2$. From Lemma \ref{lem:max_sub_gaussian}, we have with probability at least $1-\frac{2}{ n}$,
\begin{align*}
    \max_{1\leq i \leq p} \{|\phi_1(Y_i,X_i^\T\bar\beta)|\} \leq 2K(2)(1+(d(2) \log(2))^{-1/2}) +\frac{1}{d(2)})\sigma_2 \sqrt{\log( n)}.
\end{align*}

\noindent (iii) Under Assumption \ref{ass:residulal_glm_lower}, from Lemma \ref{lem:product-two-sub-gaussian-rv-norm} and \ref{lem:product-two-sub-exponential-rv-norm}, we have $\Vert \phi_1^4(Y_i,X_i^\T\bar\beta)\Vert_{\psi_\frac{1}{2}}\leq \sigma_2^4$. From Lemma \ref{lem:rv_centering}, we have
\begin{align*}
    \Vert \phi_1^4(Y_i,X_i^\T\bar\beta) - \mathbb{E}[\phi_1^4(Y_i,X_i^\T\bar\beta) ] \Vert_{\psi_{\frac{1}{2}}} &\leq K(\frac{1}{2})(1+(d(\frac{1}{2}) \log(2))^{-2}) \sigma_2^4, \\
    \mathbb{E}[\phi_1^4(Y_i,X_i^\T\bar\beta) ] &\leq \frac{1}{d(\frac{1}{2})} \sigma_2^4.
\end{align*}
Then from Lemma \ref{lem:alpha_sub_exponential_concentration}, let $t = K(\frac{1}{2})(1+(d(\frac{1}{2}) \log(2))^{-2})  \sigma_2^4$, we have with probability at least $1-2\exp\{-c(\frac{1}{2}) \sqrt{n}\}$,
\begin{align*}
    |\frac{1}{n}\sum_{i=1}^n\phi_1^4(Y_i,X_i^\T\bar\beta) - \mathbb{E}[\phi_1^4(Y_i,X_i^\T\bar\beta)]| \leq  K(\frac{1}{2})(1+(d(\frac{1}{2}) \log(2))^{-2})  \sigma_2^4.
\end{align*}
Hence,
\begin{align*}
    \frac{1}{n}\sum_{i=1}^n\phi_1^4(Y_i,X_i^\T\bar\beta)  \leq  (K(\frac{1}{2})(1+(d(\frac{1}{2}) \log(2))^{-2})  +\frac{1}{d(\frac{1}{2})}) \sigma_2^4.
\end{align*}

\noindent (iv) Under Assumption \ref{ass:sub_gaussian_entries_GLM} (i) and \ref{ass:residulal_glm_lower}, from Lemma \ref{lem:product-two-sub-gaussian-rv-norm}, we have $X_{ij}\phi_1(Y_i,X_i^\T\bar\beta)$ are sub-expenential random variable with $\Vert X_{ij}\phi_1(Y_i,X_i^\T\bar\beta)\Vert_{\psi_1}\leq \sigma_1\sigma_2$. From the definition of $\bar\beta$, we have $\mathbb{E}[X_{ij}\phi_1(Y_i,X_i^\T\bar\beta)]=0$. From Lemma \ref{lem:alpha_sub_exponential_concentration}, let $t=\sigma_1\sigma_2\sqrt{\frac{2}{c(1)\rho}}\sqrt{\frac{\log(p\vee n)}{n}}$ and $\frac{\log(p\vee n)}{n}\leq \frac{c(1)\rho}{2}$, we have with probability at least $1-\frac{2}{(p\vee n)^2}$,
\begin{align*}
    |\frac{1}{|J|}\sum_{i\in J} X_{ij}\phi_1(Y_i,X_i^\T\bar\beta)|\leq \sigma_1\sigma_2\sqrt{\frac{2}{c(1)\rho}}\sqrt{\frac{\log(p\vee n)}{n}}.
\end{align*}
Hence with probability at least $1-\frac{2}{p\vee n}$,
\begin{align*}
    \Vert \frac{1}{|J|}\sum_{i\in J}X_i \phi_1(Y_i,X_i^\T\bar\beta)\Vert_\infty \leq \sigma_1\sigma_2\sqrt{\frac{2}{c(1)\rho}}\sqrt{\frac{\log(p\vee n)}{n}}.
\end{align*}

\noindent (v) Under Assumption \ref{ass:sub_gaussian_entries_GLM}, \ref{ass:H_eigen_value}(ii) and \ref{bounded_higher_differentiate}, $\phi_2(Y_i,X_i^\T\bar\beta) x_{ij}(X_i^\T M^* \xi)$ for $i=1,\ldots, n$, $j=1,\ldots,p$ are sub-exponential random variables with sub-exponential norm bounded by $C_6\sigma^2_1\Vert \xi\Vert_2$. Moreover, we have $\mathbb{E}[\phi_2(Y_i,X_i^\T\bar\beta) x_{ij}(X_i^\T M^* \xi)]=\xi_{j}$. From Lemma \ref{lem:rv_centering}, we have
\begin{align*}
    &\Vert \phi_2(Y_i,X_i^\T\bar\beta) x_{ij}(X_i^\T M^* \xi) - \xi_{j} \Vert_{\psi_1} \leq K(1)(1+(d(1) \log(2))^{-1}) C_6\sigma^2_1\Vert \xi\Vert_2.
\end{align*}
Hence from Lemma \ref{lem:alpha_sub_exponential_concentration}, let $t = \sqrt{\frac{2}{c(1)} } K(1)(1+(d(1) \log(2))^{-1})C_6\sigma_1^2  \sqrt{\frac{\log(p\vee n)}{n}} \Vert \xi\Vert_2$ and $2(\log(p\vee n))\leq c(1)n$, with probability at least $1-\frac{2}{ (p \vee n) ^2}$,
\begin{align*}
|\frac{1}{n}\sum_{i=1}^n \phi_2(Y_i,X_i^\T\bar\beta) x_{ij}(X_i^\T M^* \xi) - \xi_{j}|  \leq  \sqrt{\frac{2}{c(1)} }K(1)(1+(d(1) \log(2))^{-1}) C_6\sigma^2_1 \sqrt{\frac{\log(p\vee n)}{n}}\Vert \xi\Vert_2.
\end{align*}
Hence, we have with probability at least $1-\frac{2}{ p \vee n }$,
\begin{align*}
    \Vert (\frac{1}{n}\sum_{i=1}^n \phi_2(Y_i, X_i^\T\bar\beta) X_iX_i^\T M^* - I_p)\xi \Vert_\infty \leq \sqrt{\frac{2}{c(1)}}K(1)(1+(d(1) \log(2))^{-1}) C_6\sigma^2_1 \sqrt{\frac{\log(p\vee n)}{n}}\Vert \xi\Vert_2.
\end{align*}

\noindent (vi)

\noindent \textcircled{1} Under Assumption \ref{ass:sub_gaussian_entries_GLM}(ii) and \ref{ass:H_eigen_value}(ii), we have $\xi^\T M^* X_i$ are sub-gaussian random variables with sub-gaussian norm bounded by $\sigma_1C_4\Vert \xi\Vert_2$. Under Assumption \ref{ass:sub_gaussian_entries_GLM}(i), from Lemma \ref{lem:product-three-rv-norm}, we have
\begin{align*}
    \Vert |\xi^\T M^* X_i X_{il}X_{im}|\Vert_{\psi_{\frac{2}{3}}} \leq  \sigma_1^3 C_4 \Vert \xi\Vert_2.
\end{align*}
From Lemma \ref{lem:rv_centering}, we have
\begin{align*}
    &\Vert |\xi^\T M^* X_i  X_{il}X_{im}| - \mathbb{E}[|\xi^\T M X_i X_{il}X_{im}|] \Vert_{\psi_{\frac{2}{3}}} \leq K(\frac{2}{3})(1+(d(\frac{2}{3}) \log(2))^{-3/2}) \Vert |\xi^\T M^* X_i | X_{il}X_{im}\Vert_{\psi_{\frac{2}{3}} },
    \\
    &|\mathbb{E}[ |\xi^\T M^* X_i X_{il}X_{im}|]| \leq \frac{1}{d(\frac{2}{3})}\Vert |\xi^\T M^*  X_i X_{il}X_{im}|\Vert_{\psi_{\frac{2}{3}}}.
\end{align*}
From Lemma \ref{lem:alpha_sub_exponential_concentration}, let $t = K(\frac{2}{3})(1+(d(\frac{2}{3}) \log(2))^{-3/2})\sigma_1^3 C_4 \Vert \xi\Vert_2$ and $3\log(p\vee n)\leq c(\frac{2}{3})n^{\frac{2}{3}}$, with probability at least $1-\frac{2}{ p \vee n ^3}$,
\begin{align*}
|\frac{1}{n}\sum_{i=1}^n | \xi^\T M^* X_i  X_{il}X_{im}| - \mathbb{E}[|\xi^\T M^* X_i X_{il}X_{im}|] |  \leq K(\frac{2}{3})(1+(d(\frac{2}{3}) \log(2))^{-3/2})\sigma_1^3 C_4 \Vert \xi\Vert_2.
\end{align*}
Hence,
\begin{align*}
\frac{1}{n}\sum_{i=1}^n | \xi^\T M^* X_i  X_{il}X_{im}|  &\leq K(\frac{2}{3})(1+(d(\frac{2}{3}) \log(2))^{-3/2})\sigma_1^3C_4 \Vert \xi\Vert_2 + |\mathbb{E}[| |\xi^\T M^* X_i X_{il}X_{im}|]| \\
&\leq   (K(\frac{2}{3})(1+(d(\frac{2}{3}) \log(2))^{-3/2}) + \frac{1}{d(\frac{2}{3})})\sigma_1^3 C_4 \Vert \xi\Vert_2 .
\end{align*}
Hence with probability at least $1-\frac{2}{ p \vee n }$,
\begin{align*}
    \max_{1\leq l\leq p, 1\leq m\leq p} \{ \frac{1}{n}\sum_{i = 1}^n | \xi^\T M^* X_i  X_{il} X_{im}| \} \leq (K(\frac{2}{3})(1+(d(\frac{2}{3}) \log(2))^{-3/2}) + \frac{1}{d(\frac{2}{3})})\sigma_1^3 C_4 \Vert \xi\Vert_2 .
\end{align*}

\noindent \textcircled{2} Under Assumption \ref{ass:sub_gaussian_entries_GLM} (ii) and \ref{ass:H_eigen_value}(ii), we have $\xi^\T M^* X_i$ are sub-gaussian random variables with sub-gaussian norm bounded by $\sigma_1 C_4 \Vert \xi\Vert_2$. From lemma \ref{lem:max_sub_gaussian}, we have with probability at least $1-\frac{2}{n}$,
\begin{align*}
    \max_{1\leq i \leq p} \{|\xi^\T M^* X_i|\} \leq (2K(2)(1+(d(2) \log(2))^{-1/2}) +\frac{1}{d(2)}) \sigma_1 C_4 \Vert \xi\Vert_2 \sqrt{\log( n)}.
\end{align*}

\noindent (vii)

\noindent \textcircled{1} Under Assumption \ref{ass:sub_gaussian_entries_GLM} (ii) and \ref{ass:H_eigen_value}(ii), we have $\xi^\T M^* X_i$ are sub-gaussian random variables with sub-gaussian norm bounded by $\sigma_1 C_4\Vert \xi\Vert_2$. Further with Assumption \ref{ass:residulal_glm_lower}, from Lemma \ref{lem:product-two-sub-gaussian-rv-norm} and \ref{lem:product-two-sub-exponential-rv-norm}, we have \\$\Vert (\xi^\T M^* X_i)^2 \phi^2_1(Y_i,X_i^\T \bar\beta)\Vert_{\psi_\frac{1}{2}}\leq \sigma^2_1\sigma_2^2 C_4^2 \Vert \xi\Vert^2_2$. From Lemma \ref{lem:rv_centering}, we have
\begin{align*}
    \Vert (\xi^\T M^* X_i)^2 \phi^2_1(Y_i,X_i^\T \bar\beta) - \mathbb{E}[(\xi^\T M^* X_i)^2 \phi^2_1(Y_i,X_i^\T \bar\beta)] \Vert_{\psi_{\frac{1}{2}}} \leq K(\frac{1}{2})(1+(d(\frac{1}{2}) \log(2))^{-2}) \sigma^2_1\sigma_2^2 C_4^2\Vert \xi\Vert^2_2.
\end{align*}
From Lemma \ref{lem:alpha_sub_exponential_concentration}, let $t = K(\frac{1}{2})(1+(d(\frac{1}{2}) \log(2))^{-2})  \sigma^2_1\sigma_2^2C_4^2\Vert \xi\Vert^2_2 n^{-\frac{1}{3}}$, we have with probability at least $1-2\exp\{-c(\frac{1}{2}) n^{\frac{1}{3}}\}$,
\begin{align*}
   & |\frac{1}{n}\sum_{i=1}^n(\xi^\T M^* X_i)^2 \phi^2_1(Y_i,X_i^\T \bar\beta) - \mathbb{E}[(\xi^\T M^* X_i)^2 \phi^2_1(Y_i,X_i^\T \bar\beta)]| \\
    &\leq  K(\frac{1}{2})(1+(d(\frac{1}{2}) \log(2))^{-2})  \sigma^2_1\sigma_2^2C_4^2\Vert \xi\Vert^2_2 n^{-\frac{1}{3}}.
\end{align*}

\noindent \textcircled{2}  Under Assumption \ref{ass:sub_gaussian_entries_GLM} (ii) and \ref{ass:H_eigen_value}(ii), we have $\xi^\T M^* X_i$ are sub-gaussian random variables with sub-gaussian norm bounded by $\sigma_1 C_4\Vert \xi\Vert_2$. From Lemma \ref{lem:product-two-sub-gaussian-rv-norm} and \ref{lem:product-two-sub-exponential-rv-norm}, we have $\Vert (\xi^\T M^* X_i)^4 \Vert_{\psi_\frac{1}{2}}\leq \sigma^4_1 C^4_4\Vert \xi\Vert^4_2$. From Lemma \ref{lem:rv_centering}, we have
\begin{align*}
    \Vert (\xi^\T M^* X_i)^4 - \mathbb{E}[(\xi^\T M^* X_i)^4 ] \Vert_{\psi_{\frac{1}{2}}} &\leq K(\frac{1}{2})(1+(d(\frac{1}{2}) \log(2))^{-2})  \sigma^4_1 C^4_4\Vert \xi\Vert^4_2, \\
    \mathbb{E}[(\xi^\T M^* X_i)^4 ] &\leq \frac{1}{d(\frac{1}{2})} \sigma^4_1 C^4_4\Vert \xi\Vert^4_2 .
\end{align*}
Then from Lemma \ref{lem:alpha_sub_exponential_concentration}, let $t = K(\frac{1}{2})(1+(d(\frac{1}{2}) \log(2))^{-2})  \sigma^4_1 C^4_3\Vert \xi\Vert^4_2$, we have with probability at least $1-2\exp\{-c(\frac{1}{2}) \sqrt{n}\}$,
\begin{align*}
    |\frac{1}{n}\sum_{i=1}^n(\xi^\T M^* X_i)^4  - \mathbb{E}[(\xi^\T M^* X_i)^4 ]| \leq  K(\frac{1}{2})(1+(d(\frac{1}{2}) \log(2))^{-2})  \sigma^4_1 C^4_4 \Vert \xi\Vert^4_2.
\end{align*}
Hence,
\begin{align}
    |\frac{1}{n}\sum_{i=1}^n(\xi^\T M^* X_i)^4 | \leq  (K(\frac{1}{2})(1+(d(\frac{1}{2}) \log(2))^{-2})  +\frac{1}{d(\frac{1}{2})}) \sigma^4_1 C^4_4\Vert  \xi\Vert^4_2. \notag  %\label{inq:power_4_xi0_M_xi}
\end{align}

\noindent (viii) Under Assumption \ref{ass:sub_gaussian_entries_GLM} (i), \ref{bounded_higher_differentiate} and \ref{ass:entries_M_x_0}, $\phi_2(Y_i,X_i^\T\bar\beta) X_{il}(X_i^\T M^*_{\cdot j})$ for $i\in J$, $l=1,\ldots,p$ and $j=1,\ldots, p$ are sub-exponential random variables with sub-exponential norm bounded by $C_6\sigma^2_1$. Moreover, we have $\mathbb{E}[\phi_2(Y_i,X_i^\T\bar\beta) X_{il}(X_i^\T M^*_{\cdot j})]=\textbf{1}_{[l=j]}$. From Lemma \ref{lem:rv_centering}, we have
\begin{align*}
    &\Vert  \phi_2(Y_i,X_i^\T\bar\beta) X_{il}(X_i^\T M^*_{\cdot j} ) - \textbf{1}_{[l=j]} \Vert_{\psi_1} \leq K(1)(1+(d(1) \log(2))^{-1}) C_6\sigma^2_1.
\end{align*}
From Lemma \ref{lem:alpha_sub_exponential_concentration}, let $t = \sqrt{\frac{3}{c(1)} } K(1)(1+(d(1) \log(2))^{-1})C_6\sigma_1^2  \sqrt{\frac{\log(p\vee n)}{\rho n}} \Vert \xi\Vert_2$ and $3(\log(p\vee n))\leq c(1)\rho n$, with probability at least $1-\frac{2}{ (p \vee n) ^3}$,
\begin{align*}
|\frac{1}{|J|}\sum_{i\in J} \phi_2(Y_i,X_i^\T\bar\beta) X_{il}(X_i^\T M_{\cdot j} ) - \textbf{1}_{l=m}|  \leq  \sqrt{\frac{3}{c(1)\rho}}K(1)(1+(d(1) \log(2))^{-1}) C_6\sigma^2_1 \sqrt{\frac{\log(p\vee n)}{n}}.
\end{align*}
We have with probability at least $1-\frac{2}{p \vee n}$,
\begin{align*}
    \Vert \frac{1}{|J|}\sum_{i\in J} \phi_2(Y_i, X_i^\T\bar\beta) X_iX_i^\T M^*_{\cdot j} - I_p \Vert_\infty \leq \sqrt{\frac{3}{c(1)\rho}}K(1)(1+(d(1) \log(2))^{-1}) C_6\sigma^2_1 \sqrt{\frac{\log(p\vee n)}{n}}.
\end{align*}

\noindent (ix) Under Assumption \ref{ass:sub_gaussian_entries_GLM}(i) and \ref{ass:entries_M_x_0}, we have $\Vert x^2_{il} (M_{\cdot j}^{*\T}X_i)^2\Vert_{\psi_{\frac{1}{2}}}\leq \sigma_1^4$ for $i\in J$, $l=1,\ldots, p$ and $j=1,\ldots, p$.  From Lemma \ref{lem:rv_centering}, we have
\begin{align*}
    \Vert x^2_{il} (M_{\cdot j}^{*\T}X_i)^2 - \mathbb{E}[x^2_{il} (M_{\cdot j}^{*\T}X_i)^2] \Vert_{\psi_{\frac{1}{2}}} &\leq K(\frac{1}{2})(1+(d(\frac{1}{2}) \log(2))^{-2})  \sigma_1^4,
    \\
    |\mathbb{E}[x^2_{il} (M_{\cdot j}^{*\T}X_i)^2]| &\leq \frac{1}{d(\frac{1}{2})}\sigma_1^4.
\end{align*}
From Lemma \ref{lem:alpha_sub_exponential_concentration}, let $t = K(\frac{1}{2})(1+(d(\frac{1}{2}) \log(2))^{-2})  \sigma^4_1$ and $3\log(p\vee n)\leq c(\frac{1}{2})\rho^{\frac{1}{2}}\sqrt{n}$, with probability at least $1-\frac{2}{ p \vee n ^3}$, we have with probability at least $1-\frac{2}{(p\vee n)^3}$,
\begin{align*}
    |\frac{1}{|J|}\sum_{i\in J} x^2_{il} (M_{\cdot j}^{*\T}X_i)^2  - \mathbb{E}[ x^2_{il} (M_{\cdot j}^{*\T}X_i)^2 ]| \leq  K(\frac{1}{2})(1+(d(\frac{1}{2}) \log(2))^{-2})  \sigma^4_1 .
\end{align*}
Hence with probability at least $1-\frac{2}{p\vee n}$,
\begin{align*}
    \max_{1\leq l\leq p,1\leq j\leq p}\{\frac{1}{|J|} \sum_{i\in J} x^2_{il} (M_{\cdot j}^{*\T}X_i)^2\} \leq (K(\frac{1}{2})(1+(d(\frac{1}{2}) \log(2))^{-2})+\frac{1}{d(\frac{1}{2})})\sigma_1^4.
\end{align*}

\noindent (x) Under Assumption \ref{ass:sub_gaussian_GLM} and \ref{ass:residulal_glm_lower}, we have $\xi^\T (\hat{M}^{J^\c}-M^*)X_i\phi_1(Y_i,X_i^\T \bar\beta)$ for $i\in J$ are centered sub-exponential random variables with sub-exponential norm bounded by $\sigma_1\sigma_2\Vert (\hat{M}^{J^\c}-M) \xi\Vert_2$. From Lemma \ref{lem:alpha_sub_exponential_concentration}, let
$t=\frac{\sqrt{\frac{\log(p\vee n)}{n}}}{\sqrt{c(1)\rho }}\sigma_1\sigma_2\Vert (\hat{M}^{J^\c}-M) \xi\Vert_2$ and $\frac{\log(p\vee n)}{n}\leq c(1)\rho$, we have with probability at least $1-\frac{2}{ p \vee n }$,
\begin{align*}
   |\xi^\T (\hat{M}^{J^\c}-M^*)\frac{1}{n}\sum_{i\in J}X_i\phi_1(Y_i,X_i\T\bar\beta)|
    &=\frac{|J|}{n}|\frac{1}{|J|}\sum_{i\in J}\xi^\T (\hat{M}^{J^\c}-M^*) X_i\phi_1(Y_i,X_i^\T \bar\beta)|\\
    &\leq \frac{1}{\sqrt{c(1)\rho}} \sqrt{\frac{\log(p\vee n)}{n}}\sigma_1\sigma_2\Vert (\hat{M}^{J^\c}-M) \xi\Vert_2.
\end{align*}

$\newline$
$\newline$
\noindent (xi) Under Assumption \ref{ass:sub_gaussian_GLM}, from Lemma \ref{lem:cross_prod_ran_mat} and let $C=\frac{1}{\sqrt{c(1)\rho}}\sqrt{\frac{\log(p\vee n)}{n}}$, we have with probability at least $1-\frac{2}{ p \vee n }$,
\begin{align}
   & |\xi^\T (\hat{M}^{J^\c}- M^*)\sum_{i\in J}\frac{1}{|J|}X_iX_i^\T(\hat{M}^{J^\c}-M^*)\xi| \notag  \\ &\leq (K(1)(1+(d(1) \log(2))^{-1})+\frac{1}{d(1)})\sigma_1^2
   \Vert  (\hat{M}^{J^\c}- M^*)\xi\Vert_2^2.
   \notag
\end{align}

\noindent (xii) Under Assumption \ref{ass:sub_gaussian_GLM}, we have $\xi^\T(\hat{M}^{J^\c}-M^*)X_i$ are sub-gaussian random variables with sub-gaussian norm bounded by $ \sigma_1 \Vert (\hat{M}^{J^\c}-M^*)\xi\Vert_2$. From Lemma \ref{lem:product-three-rv-norm}, we have
\begin{align*}
    \Vert |\xi^\T(\hat{M}^{J^\c}-M^*) X_iX_{il}X_{im}|\Vert_{\psi_{\frac{2}{3}}} \leq  \sigma_1^3 \Vert (\hat{M}^{J^\c}-M^*)\xi\Vert_2.
\end{align*}
From Lemma \ref{lem:rv_centering}, we have
\begin{align*}
    &\Vert |\xi^\T(\hat{M}^{J^\c}-M^*) X_iX_{il}X_{im}| - \mathbb{E}[|\xi^\T(\hat{M}^{J^\c}-M^*) X_iX_{il}X_{im}|] \Vert_{\psi_{\frac{2}{3}}} \\
    &\leq K(\frac{2}{3})(1+(d(\frac{2}{3}) \log(2))^{-\frac{3}{2}})  \Vert |\xi^\T(\hat{M}^{J^\c}-M^*) X_iX_{il}X_{im}|\Vert_{\psi_{\frac{2}{3}} },
    \\
    &|\mathbb{E}[|\xi^\T(\hat{M}^{J^\c}-M^*) X_iX_{il}X_{im}|]| \leq \frac{1}{d(\frac{2}{3})}\Vert |\xi^\T(\hat{M}^{J^\c}-M^*) X_i X_{il}X_{im}|\Vert_{\psi_{\frac{2}{3}}}.
\end{align*}
Further from Lemma \ref{lem:alpha_sub_exponential_concentration}, let $t = K(\frac{2}{3})(1+(d(\frac{2}{3}) \log(2))^{-\frac{3}{2}})\sigma_1^3 \Vert (\hat{M}^{J^\c}-M^*)\xi\Vert_2$ and $3\log(p\vee n)\leq c(\frac{2}{3})n^{\frac{2}{3}}\rho^{\frac{2}{3}}$, with probability at least $1-\frac{2}{ p \vee n ^3}$,
\begin{align*}
&|\frac{1}{|J|}\sum_{i\in J}  |\xi^\T (\hat{M}^{J^\c}- M^*) X_iX_{il} X_{im}|- \mathbb{E}[|\xi^\T(\hat{M}^{J^\c}-M^*) X_iX_{il}X_{im}|]|  \\
&\leq  K(\frac{2}{3})(1+(d(\frac{2}{3}) \log(2))^{-\frac{3}{2}})\sigma_1^3 \Vert (\hat{M}^{J^\c}-M^*)\xi\Vert_2.
\end{align*}
Hence,
\begin{align*}
&\frac{1}{|J|}\sum_{i\in J} | \xi^\T (\hat{M}^{J^\c}- M^*) X_iX_{il} X_{im} |  \\
&\leq  K(\frac{2}{3})(1+(d(\frac{2}{3}) \log(2))^{-\frac{3}{2}})\sigma_1^3 \Vert (\hat{M}^{J^\c}-M^*)\xi\Vert_2 \\
&+ |\mathbb{E}[|\xi^\T(\hat{M}^{J^\c}-M^*)X_iX_{il}X_{im}|]|\\
&\leq  ( K(\frac{2}{3})(1+(d(\frac{2}{3}) \log(2))^{-\frac{3}{2}}) + \frac{1}{d(\frac{2}{3})})\sigma_1^3 \Vert (\hat{M}^{J^\c}-M^*)\xi\Vert_2.
\end{align*}
Considering the union probability bound, we have with probability at least $1-\frac{2}{ p \vee n }$,
\begin{align*}
    &\max_{1\leq l\leq p, 1\leq m\leq p} \{ \frac{1}{|J|}\sum_{i\in J} | \xi^\T (\hat{M}^{J^\c}- M^*) X_iX_{il} X_{im}| \} \\
    &\leq ( K(\frac{2}{3})(1+(d(\frac{2}{3}) \log(2))^{-\frac{3}{2}}) + \frac{1}{d(\frac{2}{3})})\sigma_1^3 \Vert (\hat{M}^{J^\c}-M^*)\xi\Vert_2.
\end{align*}

\noindent (xiii) Under Assumption \ref{ass:sub_gaussian_GLM}, we have $\xi^\T(\hat{M}^{J^\c}-M^*)X_i$ are sub-gaussian random variables with sub-gaussian norm bounded by $\sigma_1 \Vert (\hat{M}^{J^\c}-M^*)\xi\Vert_2$. From lemma \ref{lem:max_sub_gaussian}, we have with probability at least $1-\frac{2}{n}$,
\begin{align*}
     \max_{i\in J}\{|\xi^\T(\hat{M}^{J^\c}-M^*)X_i|\} \leq ( 2K(2)(1+(d(2) \log(2))^{-1/2}) +\frac{1}{d(2)})\sqrt{\log( n)}\sigma_1 \Vert (\hat{M}^{J^\c}-M^*)\xi\Vert_2 .
\end{align*}
\end{prf}

\begin{lem}[Lemma 2.7.7 in \citesupp{Vershynin_2018}]
\label{lem:product-two-sub-gaussian-rv-norm}
Let $X$, $Y$ be sub-gaussian random variables. Then $XY$ is sub-exponential random variable. Moreover,
\begin{align*}
    \Vert XY\Vert_{\psi_1}\leq \Vert X\Vert_{\psi_2}\Vert Y\Vert_{\psi_2}.
\end{align*}
\end{lem}

\begin{lem}
\label{lem:product-two-sub-exponential-rv-norm}
Let $X$, $Y$ be sub-exponential random variables. Then $XY$ is $\frac{1}{2}$-sub-exponential random variable. Moreover,
\begin{align*}
    \Vert XY\Vert_{\psi_{\frac{1}{2}}}\leq \Vert X\Vert_{\psi_1}\Vert Y\Vert_{\psi_1}.
\end{align*}
\end{lem}

\begin{prf}
Without loss of generality, we assume $\Vert X\Vert_{\psi_1}=\Vert Y\Vert_{\psi_1}=1$. We have
\begin{align*}
    \mathbb{E}[\exp\{|XY|^{\frac{1}{2}}\}] &\leq \mathbb{E}[\exp\{\frac{(|X|^{\frac{1}{2}})^{2} }{2} + \frac{(|Y|^{\frac{1}{2}})^{2} }{2} \}]\quad(\text{From Young's inequality}) \\
    &=\mathbb{E}[\exp\{\frac{(|X| }{2} \}\exp\{\frac{|Y| }{2} \}]\\
    &\leq \frac{1}{2} \mathbb{E}[\exp\{|X|\}] +\frac{1}{2} \mathbb{E}[\exp\{|Y|\}]\quad(\text{From Young's inequality})\\
    &\leq 1.
\end{align*}
\end{prf}

\begin{lem}
\label{lem:product-three-rv-norm}
Let $X$, $Y$ and $Z$ be sub-gaussian random variables. Then $XYZ$ is $\frac{2}{3}-$sub-exponential. Moreover,
\begin{align*}
    \Vert XYZ\Vert_{\psi_{\frac{2}{3}}}\leq \Vert X\Vert_{\psi_2}\Vert Y\Vert_{\psi_2}\Vert Z\Vert_{\psi_2}.
\end{align*}
\end{lem}
\begin{prf}
Without loss of generality, we assume $\Vert X\Vert_{\psi_2}=\Vert Y\Vert_{\psi_2}=\Vert Z\Vert_{\psi_2}=1$. Because $X$ and $Y$ are sub-gaussian, from Lemma \ref{lem:product-two-sub-gaussian-rv-norm}, we have $XY$ is sub-exponential and
\begin{align*}
    \Vert XY\Vert_{\psi_1} \leq \Vert X\Vert_{\psi_2} \Vert Y\Vert_{\psi_2} = 1.
\end{align*}
Hence we have
\begin{align*}
    \mathbb{E}[\exp\{|XYZ|^{\frac{2}{3}}\}] &\leq \mathbb{E}[\exp\{\frac{(|XY|^{\frac{2}{3}})^{\frac{3}{2}} }{\frac{3}{2}} + \frac{(|Z|^{\frac{2}{3}})^{\frac{3}{2}*2} }{\frac{3}{2}*2} \}]\quad(\text{From Young's inequality}) \\
    &=\mathbb{E}[\exp\{\frac{(|XY| }{\frac{3}{2}} \}\exp\{\frac{|Z|^2 }{\frac{3}{2}*2} \}]\\
    &\leq \frac{1}{\frac{3}{2}} \mathbb{E}[\exp\{|XY|\}] +\frac{1}{\frac{3}{2}*2} \mathbb{E}[\exp\{|Z|^2\}]\quad(\text{From Young's inequality})\\
    &\leq 1 .
\end{align*}
\end{prf}

\begin{lem}[Lemma A.3 in \citesupp{alpha_sub_exponential}]
\label{lem:rv_centering}
For any $\alpha>0$ and any $\alpha$-sub-exponential random variable $X$, we have
\begin{align*}
    \Vert X - \mathbb{E}[X] \Vert_{\psi_\alpha} &\leq K(\alpha)(1+(d(\alpha) \log(2))^{-1/\alpha})  \Vert X\Vert_{\psi_\alpha},
    \\
    |\mathbb{E}[X]| &\leq \frac{1}{d(\alpha)}\Vert X\Vert_{\psi_\alpha},
\end{align*}
where $d(\alpha) = (2\alpha)^{1/\alpha} \exp\{1/\alpha\}$, $K(\alpha)=2^{1/\alpha}$ if $\alpha\in (0,1)$ or $K(\alpha) = 1$ if $\alpha\geq 1$.
\end{lem}

\begin{lem}
\label{lem:cross_prod_ran_mat}
Suppose that $X_i\in \mathbb{R}^p$ for $i=1,\ldots,n$ are independent sub-gaussian random vectors with sub-gaussian norms bounded by $\sigma_x$ and $\Sigma = Cov(X_i,X_i)$ for $i=1,\ldots, n$. Then for any $u\in\mathbb{R}^p$ and $v\in\mathbb{R}^p$, we have with probability at least $1-2\exp\{-c(1)  \min\{C^2,C\}n\}$,
\begin{align*}
    |w^\T\frac{1}{n}\sum_{i=1}^mX_iX_i^\T v| \leq (K(1)(1+(d(1) \log(2))^{-1})  C +\frac{1}{d(1)})\Vert  w\Vert_2 \Vert  v\Vert_2 \sigma_x^2.
\end{align*}
where $C>0$ is a positive constant, $c(1)$ is a constant from Lemma \ref{lem:alpha_sub_exponential_concentration} and $K(1)$ and $d(1)$ are constants from Lemma \ref{lem:rv_centering}.

\end{lem}
\begin{prf}
Because $X_i\in \mathbb{R}^p$ for $i=1,\ldots,n$ are independent and identically distributed sub-gaussian random vectors with sub-gaussian norms bounded by $\sigma_x$, $w^\T X_i$, $X_i^\T v$ are sub-gaussian random variables with sub-gaussian norm bounded by $\Vert w\Vert_2\sigma_x$ and $\Vert v\Vert_2\sigma_x$ respectively. Then $w^\T X_i X_i^\T v$ for $i=1,\ldots,n$ are independent sub-exponential random variables with sub-exponential norms bounded by $\Vert w\Vert_2 \Vert v\Vert_2 \sigma_x^2$. From Lemma \ref{lem:rv_centering}, we have
$w^\T X_i X_i^\T v - w^\T \Sigma v$ for $i=1,\ldots,n$ are independent centered sub-exponential random variables with sub-exponential norms bounded by $K(1)(1+(d(1) \log(2))^{-1}) \Vert w\Vert_2 \Vert v\Vert_2 \sigma_x^2$ and $|w^\T \Sigma v|\leq \frac{1}{d(1)}\Vert w\Vert_2 \Vert v\Vert_2 \sigma_x^2$. From Lemma \ref{lem:alpha_sub_exponential_concentration}, let $t=K(1)(1+(d(1) \log(2))^{-1}) C\Vert w\Vert_2\Vert\nu\Vert_2\sigma_x^2$, we have with probability at least $1-2\exp\{-c(1)  \min\{C^2,C\}n\}$,
\begin{align*}
    |w^\T\frac{1}{n}\sum_{i=1}^mX_iX_i^\T v-w^\T \Sigma v| \leq K(1)(1+(d(1) \log(2))^{-1})  C \Vert  w\Vert_2 \Vert  v\Vert_2 \sigma_x^2.
\end{align*}
Hence we have
\begin{align*}
    |w^\T\frac{1}{n}\sum_{i=1}^mX_iX_i^\T v| \leq ( K(1)(1+(d(1) \log(2))^{-1})  C + \frac{1}{d(1)})\Vert  w\Vert_2 \Vert  v\Vert_2 \sigma_x^2.
\end{align*}

\end{prf}

\begin{lem}
\label{lem:max_entries_bounds}
Suppose that $X_i\in\mathbb{R}^p$ for $i=1,\ldots,n$ are independent and identically distributed random vectors with the entries of $X_i$ are sub-gaussian random variables with sub-gaussian norm bounded by $\sigma_x$ and $\mathbb{E}[X_iX_i^\T]=\Sigma$ for $i=1,\ldots, n$. For any $D>0$, if $\frac{\log(D)}{\sqrt{n}}\leq\frac{c(1)}{3}\sqrt{n}$, then we have with probability at least $1-\frac{2 p^2}{ D^3 }$,
\begin{align*}
    \max_{1\leq l\leq p;1\leq m\leq p}|\frac{\sum_{i=1}^n X_{il}X_{im}}{n}-\Sigma_{lm}| \leq C\sqrt{\frac{\log(p\vee n)}{n}}.
\end{align*}
where $C=\sqrt{\frac{3}{c(1)}}K(1)(1+(d(1) \log(2))^{-1})\sigma_x^2$ and $c(1)$ is a constant from Lemma \ref{lem:alpha_sub_exponential_concentration} and $K(1)$ and $d(1)$ are constants from Lemma \ref{lem:rv_centering}.

\end{lem}
\begin{prf}
Because the sub-gaussian norm of the entries of $X_i$ are bounded by $\sigma_x$, from Lemma \ref{lem:product-two-sub-gaussian-rv-norm}, we have $X_{il}X_{im}$ are sub-exponential random variables with sub-exponential norm upper bounded by $\sigma_x^2$. From Lemma \ref{lem:rv_centering}, $X_{il}X_{im}-\Sigma_{lm}$ are sub-exponential random variables with sub-exponential norm upper bounded by $K(1)(1+(d(1) \log(2))^{-1}) \sigma_x^2$. From Lemma \ref{lem:alpha_sub_exponential_concentration}, letting $t=C\sqrt{\frac{3}{c(1)} }\sqrt{\frac{\log(D)}{n}}$, $C=K(1)(1+(d(1) \log(2))^{-1})\sigma_x^2$ and $\frac{\log(p\vee n)}{n}\leq\frac{c(1)}{3}$, we have with probability at least $1-\frac{2}{ D ^3}$,
\begin{align*}
    |\frac{\sum_{i=1}^n X_{il}X_{im}}{n} - \Sigma_{lm}| \leq C\sqrt{\frac{\log(p\vee n)}{n}}.
\end{align*}
Hence we have with probability at least $1-\frac{2p^2}{ D^3 }$,
\begin{align*}
    \max_{1\leq l\leq p;1\leq m\leq p}|\frac{\sum_{i=1}^n X_{il}X_{im}}{n} - \Sigma_{lm}| \leq C\sqrt{\frac{\log(p\vee n)}{n}}.
\end{align*}
\end{prf}

\begin{lem}
\label{lem:max_sub_gaussian}
If $X_i\in \mathbb{R}$ for $i=1,\ldots, n$ are independent sub-gaussian random variables with sub-gaussian norms bounded by $\sigma_x$, $D>0$ and $\log(D)>1$, then with probability at least $1-\frac{2n}{D^2}$,
\begin{align*}
    \max_{i=1,\ldots,n}|X_i|\leq ( 2K(2)(1+(d(2) \log(2))^{-1/2}) +\frac{1}{d(2)})\sigma_x \sqrt{\log(D)} .
\end{align*}
where $K(2)$ and $d(2)$ are constants from Lemma \ref{lem:rv_centering}.
\end{lem}

\begin{prf}
From $X_i\in \mathbb{R}$ for $i=1,\ldots, n$ are independent sub-gaussian random variables with sub-gaussian norms bounded by $\sigma_x$, with Lemma \ref{lem:rv_centering}, we have
\begin{align*}
    \Vert X_i - \mathbb{E}[X_i] \Vert_{\psi_2} &\leq K(2)(1+(d(2) \log(2))^{-1/2})  \sigma_x
\quad\text{for $i=1,\ldots, n$,}    \\
    |\mathbb{E}[X_i]| &\leq \frac{1}{d(2)}\sigma_x \quad\text{for $i=1,\ldots, n$.}
\end{align*}

With Lemma \ref{lem:max_sub_gaussian_origin}, let $t=2K(2)(1+(d(2) \log(2))^{-1/2})\sigma_x\sqrt{\log(D)}$, we have with probability at least $1-\frac{2n}{D^2}$,
\begin{align*}
    \max_{i=1,\ldots,n}|X_i - \mathbb{E}[X_i]|\leq 2K(2)(1+(d(2) \log(2))^{-1/2})\sigma_x\sqrt{\log(D)}.
\end{align*}
With $\log(D)>1$, we have
\begin{align*}
     \max_{i=1,\ldots,n}|X_i|\leq (2K(2)(1+(d(2) \log(2))^{-1/2}) +\frac{1}{d(2)})\sigma_x \sqrt{\log(D)} .
\end{align*}
\end{prf}

\begin{lem}[Corollary 1.4 in \citesupp{alpha_sub_exponential}]
\label{lem:alpha_sub_exponential_concentration}
If $X_i\in \mathbb{R}$ for $i=1,\ldots, n$ are independent centered $\alpha$-sub-gaussian random variables with $\alpha-$sub-gaussian norm bounded by $\sigma$ where $0\leq\alpha \leq 1$, then for $t>0$, with probability at least $1-2\exp\{-c(\alpha)\min\{\frac{t^2n}{\sigma^2},\frac{t^{\alpha}n^{\alpha}}{\sigma^{\alpha}}\}\}$,
\begin{align*}
   |\frac{\sum_{i=1}^n X_i}{n}| \leq t,
\end{align*}
where $c(\alpha)$ is a constant only depending on $\alpha$.

\end{lem}

\begin{lem}
\label{lem:max_sub_gaussian_origin}
If $X_i\in \mathbb{R}$ for $i=1,\ldots, n$ are centered independent sub-gaussian random variables with sub-gaussian norms bounded by $\sigma_x$, then for $t>0$, with probability at least $1-2n\exp\{-\frac{t^2}{2\sigma_x^2}\}$,
\begin{align*}
    \max_{i=1,\ldots,n}|X_i|\leq t.
\end{align*}
\end{lem}

\newpage
\bibliographystylesupp{apacite}
\bibliographysupp{references}

\end{document}